\documentclass[12pt]{amsart}
\usepackage{amsmath,amsfonts,amssymb,amscd}
\usepackage{latexsym}

\usepackage[
hyperindex=true,pagebackref=true,bookmarks=true,
colorlinks=true,linkcolor=blue,citecolor=red]
{hyperref}
\def\~{{\rm --}} 

\font\dfont=cmbx10 at 11pt   

\title [Non-semisimple Macdonald polynomials]
{Non-semisimple Macdonald polynomials} 
\author[Ivan Cherednik]{Ivan Cherednik $^\dag$}

\thanks{$^\dag$  \today\ \ \ Partially supported by NSF grant
DMS--0456445}

\address[I. Cherednik]{Department of Mathematics, UNC
Chapel Hill, North Carolina 27599, USA\\
chered@email.unc.edu}

 \def\bysame{{\bf --- }}
 \def\~{{\bf --}}
\newcommand{\comment}[1]{}
\renewcommand{\tilde}{\widetilde}
\renewcommand{\hat}{\widehat}

\renewcommand{\tilde}{\widetilde}
\renewcommand{\hat}{\widehat}

\newcommand{\Z}{{\mathbb Z}}
\newcommand{\Q}{{\mathbb Q}}
\newcommand{\N}{{\mathbb N}}
\newcommand{\C}{{\mathbb C}}
\newcommand{\R}{{\mathbb R}}

\def\HH{\mbox{${\mathcal H}$\kern-5.2pt${\mathcal H}$}}

\newtheorem{theorem}{Theorem}[section]
\newtheorem{maintheorem}[theorem]{Main Theorem}
\newtheorem{proposition}[theorem]{Proposition}
\newtheorem{definition}[theorem]{Definition}
\newtheorem{lemma}[theorem]{Lemma}
\newtheorem{corollary}[theorem]{Corollary}

\newtheorem{theorem }{Theorem}[section]
\newtheorem{maintheorem }[theorem]{Main Theorem}
\newtheorem{proposition }[theorem]{Proposition}
\newtheorem{definition }[theorem]{Definition}
\newtheorem{lemma }[theorem]{Lemma}
\newtheorem{corollary }[theorem]{Corollary}
\newtheorem{notation }[theorem]{Notation}
\newtheorem{remark }[theorem]{Remark}
\newtheorem{example }[theorem]{Example}

\newtheorem{ maintheorem }[theorem]{Main Theorem}
\newtheorem{ theorem}{Theorem}[section]
\newtheorem{ proposition}[theorem]{Proposition}
\newtheorem{ definition}[theorem]{Definition}
\newtheorem{ lemma}[theorem]{Lemma}
\newtheorem{ corollary}[theorem]{Corollary}
\newtheorem{ notation}[theorem]{Notation}
\newtheorem{ remark}[theorem]{Remark}
\newtheorem{ example}[theorem]{Example}

 \newcommand{\rem}{{\bf Comment.\ }}
 \newcommand{\rmk}{{\bf Comment.\ }}

\def\for{\  \hbox{ for } \ }
\def\iif{ \ \hbox{ if } \ }

\def\where{\  \hbox{ where } \ }
\def\and{\  \hbox{ and } \ }
\def\and{\  \hbox{ and } \ }

\def\equal{\stackrel{\,\mathbf{def}}{= \kern-3pt =}}

\def\la{\lambda}
\def\La{\Lambda}
\def\om{\omega}

\def\th{\theta}
\def\al{\alpha}
\def\be{\beta}
\def\ga{\gamma}
\def\ep{\epsilon}

\def\de{\delta}

\def\kapp{\hbox{\bf \ae}}
\def\si{\sigma}

\def\Ga{\Gamma}
\def\ze{\zeta}

\def\vep{\varepsilon}

\def\vth{{\vartheta}}

\newcommand{\bS}{{\mathbf S}}

\def\tal{\tilde{\alpha}}
\def\tbe{\tilde{\beta}}
\def\tde{\tilde{\delta}}

\def\tla{\tilde{\lambda}}
\def\tga{\tilde{\gamma}}
\def\tGa{\tilde{\Gamma}}

\def\tu{\tilde{u}}
\def\tU{\tilde{U}}
\def\tw{\widetilde w}
\def\tW{\widetilde W}

\def\tv{\tilde v}

\def\tz{\tilde z}
\def\tb{\tilde b}

\def\tR{\tilde R}

\def\hw{\widehat{w}}
\def\hW{\widehat{W}}
\def\hu{\hat{u}}

\def\hv{\hat{v}}
\def\hb{\hat{b}}

\def\hE{\widehat{E}}

\def\P{\mathbf{P}}

\def\S{\mathbf{S}}

\def\0{\mathbf{0}}
\def\H{\mathbf{H}}

\def\çF{\mathcal{F}}
\def\o{\mathcal{O}}
\def\t{\mathcal{T}}
\def\r{\mathcal{R}}

\def\p{\mathcal{P}}

\def\h{\mathcal{H}}

\def\y{\mathcal{Y}}
\def\e{\mathcal{E}}
\def\v{\mathcal{V}}

\def\x{\mathcal{X}}
\def\s{\mathcal{S}}
\def\g{\mathcal{G}}

\def\w{\mathcal{W}}
\def\i{\mathcal{I}}

\def\b{\mathcal{B}}

\def\lan{\langle}
\def\llb{(\!(}
\def\ran{\rangle}
\def\rrb{)\!)}
\def\lng{\hbox{\rm{\tiny lng}}}
\def\sht{\hbox{\rm{\tiny sht}}}

\newcommand{\Rad}{\operatorname{Rad}}

\newcommand{\sq}{\phantom{1}\hfill$\qed$}

\newcommand{\eps}{\varepsilon}

\def\HH{\mathfrak{H}}

\def\CC{\mathfrak{C}}
\def\LL{\mathfrak{L}}

\def\HH{\hbox{${\mathcal H}$\kern-5.2pt${\mathcal H}$}}

\font\smm=msbm10 at 12pt 
\def\symbol#1{\hbox{\smm #1}}
\def\lsmash{{\symbol n}}
\def\rsmash{{\symbol o}}
\def\#{\sharp}

\renewcommand{\natural}{\wr}

\begin{document}
\maketitle
\renewcommand{\baselinestretch}{0.94}. 
{\small
\tableofcontents
}          

\vfill\eject
{\bf Basic notations.}

1) $R=\{\al\}\subset \tR=\{[\al,\,\nu_\al j\,], j\in \Z]\}\,:$
root systems in $R^n$ with the
form $(\,,\,)$ and $R^{n+1}$; \ $\nu_\al=(\al,\al)/2$, $\nu_{\sht}=1$.
\smallskip

2) $W\,\subset \tW=\lan s_i,\,0\le i\le n\ran=W\lsmash Q
\,\subset \hW=W\lsmash P$
for the root and weight lattices $Q,P$;\ 
$l(\hw), l_\nu(\hw)\,$ are the length and
partial length;\ $B$ is a lattice between $Q$ and
$P$. 
\smallskip

3) $\Pi=P/Q=\{\pi_r, r\in O\}$, $O\,$ is  the orbit of $\al_0$ 
in the affine Dynkin diagram $\tGa$;\ 
$\al_0 = [-\vth,1]$ for the maximal short root
$\vth$;\ 
$O'=O\setminus \al_0$;\  
$\hW^\flat=W\lsmash B=\tW\rsmash \Pi^\flat $
for the image 
$\Pi^\flat$ of $B$ in $\Pi$. 
\smallskip

4) $\la(\hw)=\tR_+\cap \hw^{-1}(\tR_-),
\hw\in \hW;$\ 
for reduced $\hw=\pi_r s_{i_l}\cdots s_{i_2}s_{i_1}$,\ 
$\la(\hw) = \{\tal^l=\tw^{-1}s_{i_l}(\al_{i_l}),\,
\ldots,\, 
\tal^2=s_{i_1}(\al_{i_2}),\, \tal^1=\al_{i_1} \};$
$\tw=\pi_r^{-1}\,\hw$. 
\smallskip

5) $(wb)\llb z \rrb = w(b+z)$, \ 
$(\,[z,l]\,,\, z'+d) = (z,z')+l$\
for $w\in W,\, b\in P,\, z\in \C^n\,$;\  
$(\,\hw([z,l])\,,\,\hw\llb z' \rrb+d\,) =
([z,l], z'+d) \for \hw\in \hW;$\ \,
$\CC=\{z\in \R^n, (z,\al_i)>0$ as $i>0\},$\ 
$\CC^a = \bigcap_{i=0}^n \LL_{\al_i},$ where 
$\LL_{[\al,\nu_\al j]}=\{z\in \R^n,\
(z,\al)+j>0 \}.$
\smallskip

6) $\rho_\nu = \frac{1}{2}\,\sum_{\nu_{\al}=\nu} \al =
\sum_{\nu_i=\nu}  \om_i$,
where $\al\in R_+$, $\om_i$ are fundamental weights;\ 
$(\rho_\nu)^\vee=\rho_\nu/\nu$,\,
$2\rho_k=\sum_{\al\in R_+}\, k_\al\al;$\ 
$q_\al=q^{\nu_\al},\, t_\al=q_\al^{k_\al}$.
\smallskip

7) $b= \pi_b  u_b \for b\in P,\,u_b\in W$\,, where 
$u_b(b)=b_-\in P_-\,\and \la(\pi_b)\cap R = \emptyset$;\ 
$\pi_r=\pi_{\om_r}$ for $r \in O'$;\ 
$u_b\pi_b = b_-,\ b_+=w_0(b_-)=\varsigma(-b_-);$\ 
$b_\#= b- u_b^{-1}(\rho_k)$.
\smallskip

8) $\v = \Q_{q,t}[X_b]\ =\ \Q_{q,t}[X_b, b\in B]\,:$
the polynomial representation over
$\Q_{q,t} = \Q[q^{\pm 1/m},t^{\pm 1/2}]$
for  $(B,B)=(1/m)\Z;$\  $\Phi_{\hw}$ and  
$\Psi_{\hw}$ are $X$ and $Y$ intertwiners;\ 
$\Psi_i=\tau_+(T_i)+
\frac{t_i^{1/2}-t_i^{-1/2}}{Y_{\al_i}^{-1}-1},\ i\ge 0$,\ 
$P_r=\tau_+(\pi_r),\ r\in O'$ 
for the DAHA-automorphism $\tau_+$.
\smallskip

9) $\widetilde{E}_b, E_b, P_b, \e_b=E_b/E_b(q^{-\rho_k}) \,:$ 
non-semisimple, nonsymmetric, symmetric and
spherical polynomials, $b\in B$;\ 
$\Psi_i^b$ is $\Psi_i$ acting on $E_b$.
\smallskip

10) $\widetilde{V}_b\subset \v_b \subset \v(-b_\#)^\infty\,:$ 
spaces of generalized $Y$\~eigenvectors in $\v$;\ 
$\widetilde{V}_b$ may depend on the reduced
decomposition of $\pi_b (b\in B)$, 
$\v_b=\lim  {\v}^\natural_b$ does not depend on the
decomposition;\ $\widetilde{\Psi}_i$,
$\widetilde{\Psi}_{\hw}$ are the generalized intertwiners
defined in terms of  $\tR^0$.
\smallskip

11) $R^0$ is a root subsystem in $\tR$, 
$\tW^0=\lan s_\al,\al\in R^0\ran;$\ 
in the context of $\v$,
$\tR^0=\{\tal=[\al,\nu_\al j]\in \tR\, \mid \,
q^{\nu_\al j-(\al,\, \rho_k)}=q_\al^{j-(\al^\vee, \rho_k)}=1\};$\ 
$\hW^\flat[\xi]=
\{\,\hw\in \hW^\flat\, \mid\, q^{\hw\llb \xi\rrb}=q^{\xi}\,\};$\  
$\Pi_{\tR}$ is the product of affine exponents. 

\renewcommand{\baselinestretch}{1.0}.
\vfill

\setcounter{section}{-1}
\setcounter{equation}{0}
\section{Introduction}

\comment{
The paper is mainly devoted to the irreducibility
of the polynomial representation of the double
affine Hecke algebra for an arbitrary 
reduced root systems and generic ``central charge" q. 
The technique of intertwiners in the non-semisimple
variant is the main tool. We introduce Macdonald's 
non-semisimple polynomials and 
use them to analyze the reducibility of the polynomial 
representation in terms of the affine exponents,
counterparts of the classical Coxeter exponents.
The focus is on the principal aspects of the technique of 
intertwiners, including related problems in the theory
of reduced decompositions on affine Weyl groups.
}

The paper is mainly devoted to the irreducibility
of the polynomial representation of the Double
affine Hecke algebra, DAHA, for arbitrary 
reduced root systems and generic ``central charge" $q$. 

The technique of intertwiners from \cite{C1} in a
non-semisimple variant is the main tool.
It is important for the
decomposition of the polynomial representation
in terms of irreducible DAHA modules and for its
weight decomposition. We focus on the principal
aspects of the technique of intertwiners and discuss
only basic (and instructional) applications.
Generally, it is more efficient to combine 
the intertwiners with other approaches, to be 
considered in further papers.
including author's                              
(unpublished) construction of the Jantzen\~type 
filtration.                                     
\smallskip

There are several methods that can be used now in the
study of the polynomial representation of DAHA and
its degenerations. Certainly the localization functor
(the KZ\~monodromy) from \cite{GGOR,VV1}
must be mentioned, as well as   
the geometric methods of \cite{VV3} 
and the parabolic induction from recent \cite{BE}.
The technique of intertwiners provides constructive
(relatively elementary)  
tools for managing the irreducibility of the
polynomial representation and its constituents 
for any $q,t$ based on combinatorics of affine Weyl groups.
    
\smallskip
An important general objective of this technique
is finding a counterpart of the classical theory of 
{\em highest vectors} for DAHA (and AHA),
complementary to the geometric method of \cite{KL1}.
It involves difficult combinatorial problems 
and is known only for type $A$ and in some cases of small ranks.
However, the geometric DAHA methods are far from simple too
and explicit theory of DAHA modules is needed in quite a few
applications (see Section \ref{sec:EXPECTEDAP}).     
\medskip

{\bf Main constructions.}
The {\em polynomial representation}, denoted by $\v$ 
in the paper, is well known to be irreducible and semisimple
for generic values of the DAHA-parameters $q$ and $t=q^k.$
It becomes reducible either when $q$ is a root
of unity or for generic $q$ and special $t.$

We perform a complete analysis of the
irreducibility and semisimplicity of $\v$ 
for generic $q.$ Another application
is a construction of the canonical semisimple
submodule in $\v$ generalizing that of type $A$
from \cite{FJMM} (the symmetric variant) 
and \cite{Ka} (the non-symmetric case). 
 
\smallskip

We begin with a description of all {\em singular} 
$t=q^k$ making (by definition) 
the radical of the evaluation pairing 
nonzero. The answer is given in terms of the principal 
values (at $t^{-\rho}$) of the nonsymmetric Macdonald
polynomials \cite{C3,C4,C100} 
for the weights sufficiently large to ensure the
existence of these polynomials. 
The same answer can be obtain using a generalization 
of the method from \cite{O1,O3} in the 
rational case based on the {\em shift-operator}. 
The latter is used to calculate the principal value  of
the $t$\~discriminant $D_Y$ applied to the 
$t$\~discriminant $D_X$ in $\v$, 
where $X,Y$ are the generators of DAHA.   
The approach via the Macdonald polynomials
has no rational counterpart.

The {\em evaluation pairing} is defined as follows: 
$$
\{E,F\}\,=\, E(Y^{-1})(F(X))(t^{-\rho}),\ E,F \in \v\,;
$$
for instance, $D_Y(D_X)(t^{-\rho})=\{D_{X^{-1}},D_X\}$.
In the simply-laced case,
the radical $Rad$ of this pairing is zero if 
and only $\v$ is irreducible. 
This equivalence becomes more subtle in the 
non-simply-laced case, as well as the formula for 
$D_Y(D_X)(t^{-\rho})\,$.
The $q,t$\~theory provides the best (and direct)
method for calculating this formula, including 
managing the rational case through the limiting 
procedure from DAHA to its rational degeneration.
\smallskip

The main objects of this paper are
the {\em chains of the intertwiners\,} and  
{\em non-semisimple Macdonald polynomials\,} 
(defined via such chains).
We define a system of subspaces in $\v$
with the Macdonald polynomials as top elements
in a punctured neighborhood of a singular $t$\,
and then extend this construction to singular $t$.

In contrast to the semisimple case, the related 
combinatorics of the affine root systems and
affine Weyl groups becomes significantly more 
involved. For instance, 
the non-semisimple Macdonald polynomials may depend 
on the choice of
the reduced decomposition of the corresponding elements
in the affine Weyl groups, the {\em relative Bruhat
ordering\,} is needed versus the usual one, and so on. 
\smallskip

There are six Main Theorems in the paper, we will
discuss its contents following these theorems
beginning with those of combinatorial nature. 

\medskip
\subsection{Reduced decompositions}
A significant part of the paper is devoted to the
combinatorics of affine root systems grouped around 
the reduced decompositions. The theory of reduced decompositions
in affine (and non-affine) Weyl groups is far from being 
simple and complete. 

The {\em affine exponents\,} (see below), 
generalizing the classical Coxeter 
exponents, demonstrate that there are many properties
of the reduced decompositions we do not know.

Almost all facts of combinatorial nature
we mention/prove in the paper are really needed here. Some,
like Main Theorem \ref{RANKTWO}, are expected to be used in 
other papers, say, for the
classification of the semisimple representations and
in the theory of Jantzen\~type filtration.
\smallskip

Generally speaking, 
many (if not all) theorems established
in \cite{C101} and previous author's papers 
(including the Macdonald conjectures) are corollaries 
of relatively few facts on affine root systems and 
affine Weyl groups. Extending the list of such basic
facts is very important; any progress here can be
readily translated to the DAHA and AHA theories.
This paper demonstrates it.
\smallskip

We note that an abstract foundation of our approach 
is the notion of a pair of {\em compatible  $r$\~matrices\,} 
associated with an affine root system and its 
root subsystem ($\tR^0\subset \tR$, see below). 
There is a connection with Kauffman's 
axioms of {\em virtual links\,} in topology \cite{GPV}
(when extended to arbitrary affine root systems),
although the relation is direct only for some 
representations of DAHA. 
\smallskip

Let us briefly discuss the main theorems of the
``combinatorial part" of this paper. We begin
with an irreducible reduced root system $R$
and its affine extension 
$\tR=\{[\al,\frac{(\al,\al)}{2}\,j]\}$, where $\al\in R$,
$j\in \Z$; $W$ is the Weyl group of $R$, $\hW$
is the extended affine Weyl group $\hW=W\lsmash P$ 
defined for the weight lattice $P$ of $R$.  
\smallskip

{\em Main Theorem \ref{INTRINLA}}. The key tool we use
in the theory of the extended affine Weyl group $\hW$
is the notion of the $\la$\~set: \ 
$\la(\hw)=\tR_+\cap \hw^{-1}(-\tR_+)$ for $\hw\in \hW$, where
$\tR_+$ is the set of positive roots in $\tR$. It is well-known
that $\hw$ is uniquely determined by  
$\la(\hw)$; many properties of $\hw$ and its
reduced decompositions can be interpreted in terms
of this set. We give an intrinsic description of the 
$\la$\~sets, more generally, the $\la$\~{\em sequences\,},
that are $\la$\~sets with the orderings induced by reduced
decompositions of $\hw$.  
Essentially, the description is given in terms of the 
``triangle triples" 
$\tal,\tal+\tbe,\tbe$.
Say, $\tal,\tbe\in \la(\hw)$
$\Rightarrow$ $\tal+\tbe\in \la(\hw)$ and the latter
root must be between $\tal$ and $\tbe$\, if this set
is treated as a sequence. 

This theorem is not exactly
new (although we cannot give precise references); however, 
we think, it is the most complete one of this
kind. Several adjustments and generalizations 
were needed here since we need to catch the ordering 
of the roots in $\la(\hw)$ and because we are doing the
affine theory. It plays an important technical role in 
the paper; we decided to give its proof. 
It is mainly used as a list of properties of $\la(\hw)$; 
in the opposite direction, it gives that the intersections
of $\la$\~sequences with root subsystems remain 
$\la$\~sequences.  
\smallskip

{\em Main Theorem \ref{RANKTWO}}. 
We need to know when a set of positive 
roots of a rank two subsystem inside a
given {\em sequence} 
$\la(\hw)$ can be made consecutive 
using the Coxeter
transforms in $\la(\hw)$. This problem can be readily 
reduced to considering {\em triangle triples}
$\{\tal,\tal+\tbe,\tbe\}$\, provided 
special conditions (a,b,c) from the theorem. 
Only for affine $A_n, B_2, C_2, G_2$ or when 
$|\tal|\neq|\tbe|$ the answer is always affirmative.
For affine $A_n$, it can be deduced directly 
from the interpretation of \cite{Ch5}) in terms of the lines
on the cylinder. Generally the
{\em admissibility} condition is needed, which is 
formulated in terms of subsystems of $\tR$ of types 
$B_3,C_3,D_4$;
the justification requires using subsystems of types
$B_4,C_4,D_5$. Some technical details of the
proof of this theorem are omitted in the present paper. Namely, 
the obstacles of type $D_4$ are not discussed in full and
the proof is not complete in the case $F_4$ (as $|\tal|=|\tbe|$);
we hope to continue this topic in other works.

\smallskip

{\em Its applications.\,}
Presumably, this theorem  ``explains" 
the difficulties with  
``AHA-DAHA highest vectors",
generalizing Zelevinsky's segments
(known only in the $A$\~case and for some
root systems of small ranks). Generally, we
need the approach from \cite{KL1} and similar 
{\em geometric} constructions. 
The classification of {\em semisimple\,} 
representations of DAHA seems a natural
first step toward the theory of
highest vectors/weights.
  
In the case of affine root system of type $A$,
the classification of such representations
was obtained in \cite{C12} in terms 
of infinite periodic skew Young diagrams (see also \cite{C101}).
Let us also mentioned paper \cite{SuV}, where the result
is the same as in \cite{C12} but the approach is somewhat
different. In \cite{C12}, it was obtained as a corollary 
of the Main Theorem there (for arbitrary root systems) 
based on the technique of intertwiners combined with
Zelevinsky's classification  in
the case of the affine Hecke algebra of type $A$
(adjusted to the semisimple affine case in author's papers).

Theorem \ref{RANKTWO} for the affine $A$-system 
is an important (actually, the key) ingredient of  
the construction from \cite{C12} and its justification.
The classification of the semisimple representations of 
DAHA for arbitrary root systems (not finished at the moment) 
will be addressed in author's future papers.

In this paper, we do not particularly need 
Theorem \ref{RANKTWO} for the theory of $\v$.  
Corollary \ref{CORTALPRIME} 
can be used to justify applying Key Lemma \ref{KEYTRIPLE} 
to the irreducibility of $\v$ 
(Main Theorem \ref{THMIRRV}), but
Lemma \ref{LEMREFLEC} is actually sufficient.
\medskip

\subsection{Affine exponents}
The classical exponents of an irreducible
reduced root system $R$ are given by the formula:
\begin{align*}
&\Pi_R(t)\equal\prod_{\al\in R_+} \frac{1- t^{1+(\al,\rho^\vee)}}
{1-t^{(\al,\rho^\vee)}}=
\prod_{i=1}^n \frac{1-t^{m_i+1}}
{1-t}.
\end{align*}
We will not review their various and important applications in 
algebra, combinatorics, geometry and topology; see, e.g., 
\cite{Bo,Hu} concerning basic (certainly not all) aspects 
of their theory. Algebraically, this formula is about reducing
coinciding terms in the numerator and denominator of the product
in the left-hand side. This viewpoint will be the main in this
paper.

The affine counterpart of $\Pi_R$ is defined in terms
of $q$ and $t_\al=q_{\al}^{k_{\al}}$, 
where $q_\al=q^{\nu_\al}$ for $\nu_\al=(\al,\al)/2$ and 
$k_{\al}$ depends only on  $\nu_\al$: 
\begin{align*}
&\Pi_{\tR}\equal\prod_{\al\in R_+}
\Bigl( (1- q_\al^{k_\al+(\al^\vee,\,\rho+\rho_k)})
\prod_{ j=1}^{(\al^\vee,\,\rho)}
\frac{
(1- q_\al^{j-1+k_\al+(\al^\vee,\,\rho_k)})}
{(1- q_\al^{j-1+(\al^\vee,\,\rho_k)})}\Bigr).
\end{align*}
The normalization here is $\nu_\al=1$ for 
short roots;\, $2\rho_k=\sum_{\al\in R_+}\, k_\al\al.$
Generally, the definition depends on a pair of
affine extensions of $R$ and/or its dual $R^\vee$. 
The definition above is for the pair $\{\tR,\tR\}$ 
where $\tR=\{[\al,\nu_\al j],\, \al\in R,\,j\in \Z\}$.

The {\em root subsystem} 
$\tR^0=\{\,[\al,\nu_\al j],\ q_\al^{j+(\al^\vee,\,\rho_k)}=1\,\}$
(compare with the denominator of the last formula), 
the Weyl group $\tW^0$ of $\tR^0$, and the relative
Bruhat ordering for the pair $\tR^0\subset \tR$
play the key role in our analysis of $\v$. 
The main technical reason is that the roots from 
$\tR^0$ lead to {\em singular intertwiners\,} (see
below); $\tR^0$ appears virtually in all 
statements and formulas of the paper.

We calculate $\Pi_{\tR}$ in terms of the {\em affine
exponents\,} in a way similar to the classical formula; 
$\Pi_{\tR}$ becomes $\Pi_R$ as  $q=0$, $t=t_{\sht}=t_{\lng}$,
(equivalently, $k=k_{\sht}=\nu_{\lng}k_{\lng}$), so our
construction is 
a direct generalization of the classical definition.
\smallskip

In the simply-laced case:  
\begin{align*}
&\Pi_{\tR}=
\prod_{i=1}^n \frac{\prod_{j=0}^{m_i} (1-q^{j} t^{m_i+1})}
{1-t}.
\end{align*}
We give complete lists of the affine exponents for
any (reduced) root systems but their combinatorics is not
considered in full in this paper, as well as 
their various applications. For instance,  
a Langlands-type duality holds, namely, 
Theorem \ref{OTHERAFFINE}:
$$
\Pi_{\widehat{R^{}}}(q,k_{\lng},k_{\sht}) = 
\Pi_{\widehat{R^\vee}}(q,k_{\sht},k_{\lng})
$$
for the affine root system
$\widehat{R}=\{[\al,j],\, \al\in R,j\in \Z\}$
when the pair $\{\widehat{R},\widehat{R^\vee}\}$
is used in the definition of $\Pi_{\widehat{R^{}}}$.

It is directly connected with the DAHA\~Fourier 
transform, to be discussed in further papers.
Generally, there are confirmations 
that, as far as DAHA can be used, Langlands'
correspondence can be associated with the Fourier transform,
which plays the key role in the theory of DAHA.

The functoriality of $\Pi_{\tR}$ and $\Pi_{\widehat{R}}$
with respect to the affine root subsystems of $R$
is not discussed at all. 
The {\em affine exponents\,} for $\Pi_{\widehat{R}}$ 
satisfy a $q\leftrightarrow q'$\~duality, for instance,
rational singular $k$ are in a sense dual to the zeros of the
Poincar\'e polynomial (see (\ref{qqprime})). 
This property is related to Theorem \ref{THMIRRV}, 
but is not discussed systematically. Also, no
interpretation of $\Pi_{\tR}$ and $\Pi_{\widehat{R}}$
as a Poincar\'e series is known.
\smallskip

The main motivation of the affine exponents 
in the paper is that their zeros upon multiplicative
translations by the elements from $q^{\Z_-\,}$ 
constitute  the list of all
{\em singular $t$}, defined as those making the radical 
of the polynomial representation nonzero. For instance, 
$t_{sing}=q^{-j/(m_i+1)}$ for $j>0$ 
in the simply-laced case (the zeros of $\Pi_R$ 
must be excluded).
It is intimately 
related to the following formula.
\medskip

{\bf Calculating {\mathversion{bold}$D_Y(D_X)(t^{-\rho})$}.}
This expression generalizes the one in the
rational setting, which is $D_y(D_x)$ for the 
$y$\~discriminant 
$D_y=\prod_{\al\in R_+} y_\al$ 
in terms of the differential 
{\em rational Dunkl operators\,} \cite{Du},
$$
y_b=\partial_b+\sum_{\al\in R_+}
\frac{k_{\al}(b,\al)}{x_\al}(1-s_\al)
\for b\in P,
$$
applied to the $x$\~discriminant 
$D_x=\prod_{\al\in R_+} x_\al$. 

The latter expression is already  
a constant (depending on the $k$\~parameters),
so the evaluation is not necessary in the rational case.
This formula is due to Opdam \cite{O1} in the crystallographical
case; see \cite{O3, DJO} for explicit
formulas in the non-simply-laced cases and
for $I_2(2m),H_3,H_4$.
A straightforward algebraic verification of this
formula is known (and quite involved) 
only in the $A$\~case (Dunkl, Hanlon).

The methods involved in \cite{O1,O3,DJO} are
the Heckman \~ Opdam theory of Jacobi\~Jack polynomials,
the Macdonald \~ Mehta conjecture (proved by Opdam)
and also the semisimplicity theorem 
for the classical Hecke algebras from \cite{GU}.
Employing the Jack polynomials
is remarkable; it requires differential-trigonometric setting 
(we use the difference-trigonometric
setting in a similar manner).
Using the monodromy method and \cite{GU} 
provides a universal tool, but the least direct.

\comment{
Opdam found the $k$\~zeros of $D_y(D_x)$
(with multiplicities) 
in the case of equal labels $k_{\sht}=\nu_{\lng}k_{\lng}$
 in terms of the classical 
{\em Coxeter exponents\,} of $R$ \,\cite{O1}. 
See \cite{O3, Je1, DJO}
for the non-simply-laced case and non-crystallographic
Coxeter groups. This gives the formula for $D_y(D_x)$
up to a multiplier from $\Q$, that can be 
calculated (numerically, for $H_3,H_4$).  

The methods used by Opdam and de Jeu are
the Macdonald \~ Mehta conjecture (proved by Opdam), 
the Heckman \~ Opdam theory of Jacobi\~Jack polynomials
and also the semisimplicity theorem 
for the classical Hecke algebras from \cite{GU}.
Let us comment on it.

The approach based on the Macdonald \~ Mehta conjecture 
uses its ``analytic" connection with the Bernstein-Sato 
polynomial of $R$. Employing the Jack polynomials
is remarkable; it requires trigonometric setting 
(we need the trigonometric
setting too and also replace the differential operators 
by the difference ones).
Using the monodromy method and \cite{GU} 
provides the most universal tool but the least direct.

Reconstructing $D_y(D_x)$ from its $k$\~zeros requires
knowing their multiplicities, that creates 
technical problems, especially in the non-simply-laced case, 
A straightforward algebraic verification of the 
Opdam \~ de Jeu formula is known (and quite involved) 
only in the $A$\~case (Dunkl, Hanlon).
} 
\medskip

{\em The difference case.}
The calculation of $D_Y(D_X)$ becomes simple
and uniform (for any $k$\~parameters) in the general
$q,t$\~case. See Main Theorem \ref{YOFXEV}.
We check that $D_Y(D_X)$ is proportional to
the symmetric Macdonald  polynomial of weight $\rho$
using the {\em shift operator} from \cite{C2} and 
employ the Macdonald's principal
value conjecture proven in \cite{C12}; the result is
$\Pi_{\tR}$ up to some powers of $q,t$.
Combining this formula with the limiting procedure
from the general DAHA to its rational degeneration gives
another justification of the Opdam formula (in the 
crystallographic case).
\smallskip

Concerning the evaluation (principal value) formula,
its proof is entirely conceptual, a direct corollary of the
DAHA\~duality. We improve its deduction in this paper;
see Proposition \ref{COREBEC}.

Note that the Macdonald polynomials collapse in the rational
limit; the differential rational Dunkl operators are nilpotent
in the polynomial representation and have no eigenvalues
but $0$. In the trigonometric case, the Jack polynomials exist 
but the duality collapse (the evaluation formula holds); 
We think, it explains why the difference theory is
the most relevant to deal with $D_Y(D_X)$ and, correspondingly,
$D_y(D_x).$ 
\smallskip

Last but not the least,
the {\em Coxeter exponents\,} ``naturally" enter  
the $q,t$\~setting, almost directly due to their definition
via $\Pi_R$. It ``explains" why the Coxeter
exponents  appear almost everywhere in the Macdonald \~
Matsumoto theory of $p$\~adic spherical functions;
the limit $q\to 0$ is exactly the passage from DAHA to
the $p$\~adic theory.
Another advantage of the $q,t$\~theory is that the affine
exponents, in contrast to their rational limits,
are generally of multiplicity one. 
The multiplicities of the roots of the Bernstein-Sato 
polynomials are important in their general theory.  
\smallskip

Generally, no simple ways
can be expected for obtaining the Coxeter exponents 
from {\em rational}-differential operators; say,
it was a difficult conjecture (proved by Opdam) 
that the Bernstein-Sato polynomials for the
discriminant of $R$ are given in 
terms of $\{m_i\}$. A metamathematical reason for it is
trigonometric nature of their (main) definitions in 
terms of the $W$\~invariant {\em Laurent} polynomials
and via the product $\Pi_R$.

\comment{
We come to an interesting conclusion:
{\em the classical theory of the radial parts of the
Laplace operators on the symmetric spaces and their
rational degenerations does require difference methods}.

The first demonstration of this fact was \cite{C8}
devoted to the Harish-Chandra transform. The
application of the difference DAHA methods to the 
Bernstein-Sato polynomials of the root systems seems 
equally convincing.} 
\medskip

{\bf The radical.}
The zeros of the {\em rational} $D_y(D_x)$ give the values
of $k$ when the discriminant $D_x$ belongs to the 
{\em radical\,} of the 
evaluation pairing (rational evaluation, at $0$). 
For such $k$, the radical, $Rad_{rat}$, 
is obviously nonzero; therefore the polynomial representation,
$\v_{rat}$, is reducible.
Vice versa, Opdam establishes that the irreducibility
of $\v_{rat}$ occurs exactly
when the (rational) radical $Rad_{rat}$ is nonzero, 
that happens at the zeros of $D_y(D_x)$ up to their
translations by negative integers.

The equivalence $\{\,Rad_{rat}=\{0\}\,\}\, 
\Leftrightarrow\, \{$irreducibility of $\v_{rat}\}$ 
is very simple. Indeed,
it is obvious in the $\Leftarrow$ direction. If $\v'_{rat}$
is a submodule of $\v_{rat}$ then it contains a
$y$\~eigenvector $v'$ (the eigenvalue can be only $0$
for the rational DAHA) and therefore 
$v'-v'(0)1$ belongs to $Rad_{rat}$. 
\smallskip

{\em Main Theorem \ref{THMIRRV}}.
A modified variant of this reasoning can be  
used in the $q,t$\~setting. Generally, $\v$ 
can be reducible when $Rad=\{0\}$ in the non-simply-laced
case; it is of course not true any longer that all 
$Y$\~eigenvalues in $\v$ coincide. The list of exceptional 
cases, when $Rad=\{0\}$ but $\v$ is reducible, is 
given in Theorem \ref{THMIRRV}. Its justification
involves combinatorial case-by-case analysis.
\smallskip

{\em Etingof's theorem}. Recently Etingof \cite{Et} 
obtained the first ``half" of this list
using the rational case and his general reduction theory.
Actually he obtained it for the {\em degenerate\,} DAHA,
the trigonometric limit of the general $q,t$\~DAHA.
The relation between our $q,t$\~theorem 
and his one in the degenerate (trigonometric) case
is as follows. 

The parameter $q,$ the center charge,
is assumed generic in Theorem \ref{THMIRRV}.
However this assumption
is not quite sufficient to connect the $q,t$\~DAHA 
with its trigonometric degeneration. 
One has to impose furthermore that 
$q^a t^b=1\,\Rightarrow\, a+kb=0$ as 
$t=q^k$ for $a,b\in \Q$ in the simply-laced case 
(with two $t$ here for $B,C,F,G$).
Then the reducibility of $\v$ and its (trigonometric)
degeneration will occur exactly at the same $k$. 
This corresponds to the list from (\ref{bfgcondi}) in  
Theorem \ref{THMIRRV}. 

We mention that Etingof obtained 
(\ref{bfgcondi}) practically independently of this paper.
Some exceptional cases when $\v$ is reducible
in spite of $Rad=\{0\}$ were known to him, however he
arrived at their {\em complete} 
list without knowing/using the methods and exact results of 
this paper.
\smallskip

Paper \cite{Et} contains simple and {\em 
conceptual} interpretation of conditions from
(\ref{bfgcondi}). We found this list on the basis of technical 
Lemma \ref{ZFROMM}, which includes a case-by-case
analysis. Actually, this lemma is directly connected with
the {\em zigzag connectivity} from 
Lemma \ref{COMBR01} that seems quite general.
However, the relation of (\ref{bfgcondi}) to
the Borel \~ de Siebenthal algorithm found by Etingof
is remarkable and clarifying. It is
good to have alternative approaches to such problem
that seem of fundamental nature. 

We would like to mention the paper \cite{L} in this context;
considering the ``mixed" products of the normalized
intertwiners and the $T$\~elements ($T_w^c$
from Proposition 8.1 there)
is similar to our approach. However, the main problem
in the theory of polynomial representation is dealing
with {\em non-invertible intertwiners}, where the counterparts of
Lusztig's $T_w^c$ may depend on the particular choice of the reduced
decomposition.
\smallskip

The second ``half" of the exceptional cases, namely    
(\ref{reducasesthm})\,, describes the 
reducibility of $\v$ with zero radical
without imposing the assumption 
$q^a t^b=1,\, a,b\in \Q\, \Rightarrow\, a+kb=0$.
It is, in a sense, {\em dual} to (\ref{bfgcondi}), 
so it is likely that the whole Theorem \ref{THMIRRV} 
can be ``deduced" from the rational case using 
Etingof's approach at greater potential.
\smallskip

Our approach is expected to be applicable to describing 
the cases when $\v/Rad$ is not irreducible (it naturally includes
those from Theorem \ref{THMIRRV}), but we did not establish
Lemma \ref{COMBR01} and related stuff in proper generality
so far. It was proven in \cite{C100} that
$\v/Rad$ is always irreducible for generic $q$
if this quotient is {\em finite-dimensional}.
The technique of intertwiners 
(in the $q,t$\~setting) can be used 
to address this problem in general.
\medskip

\subsection{Non--semisimple polynomials}
A natural approach here is to use the decomposition 
of the polynomial representation $\v$ for generic $q,t$
in terms of the non-symmetric Macdonald polynomials 
$E_b, \, b\in P,\,$ and tend $t$
(or $k$) to a {\em singular\,} values 
$t_{sing}=q^{k_{sing}}$. 
Such decomposition of $\v$ and the related technique
of intertwiners are well known for generic $q,t$; 
see \cite{C1} (\cite{KnS} in the $A$\~case) 
and \cite{Ma3,C101}.
However there are no reasonable formulas for (arbitrary) 
coefficients of $E_b$ and no straight ways to
control directly the limits of these polynomials
coefficient-wise as $t$ becomes singular (unless for
$A_1$ and in some rank 2 examples).

A bypass is in representing $\v$ as a sum of a system of
finite dimensional subspaces $\v_b$ such that their limits 
$t\to t_{sing}$ can be calculated {\em exactly} and
the corresponding ``Gr" is a direct sum of
one-dimensional subspaces, i.e, for $\{\v_b,\, b\in P\}$
constituting a {\em maximal system\,} of vector spaces 
generating  $\v$.

Here the limit is understood in the sense of vector
bundles over a curve (dimension one is important).
{\em All} regular (at $t_{sing}$) linear combinations
of vectors with the coefficients in the field of
rational functions in terms of $(t-t_{sing})$ 
must be considered. 
Such limiting procedure preserves the dimensions
and, given $b$, extends the vector bundle  
$\{\v_b(t)\}$ from generic $t$ to singular $t_{sing}$. 
Generally, there can be many choices of  
{\em one-parametric} limiting procedures, 
but this construction is independent of such choices.
\smallskip

{\em Main Theorem \ref{VINVHAT}.} It contains a 
construction of such system of spaces. Actually,
we simply give a uniform definition of $\v_b$ for all
$t$ (including $t_{sing}$) and check that the dependence
of $t$ is flat. The description of spaces $\v_b$ is 
very explicit; they are given in terms of
$\tR^0$, the subsystem of singular roots for a given 
$t_{sing}$. This construction is entirely combinatorial
in terms of $\tR^0$ ($t$ is arbitrary,
possibly different from $t_{sing}$).

The first application is that 
the Macdonald polynomial $E_b$ for $b\in P$
exists at $t_{sing}$ if ({\em and only if}\, in some sense) 
dim\,$\v_b=1$, which is a pure combinatorial condition 
concerning $b$ and $\tR^0$. 

Given $b\in P$, the definition of the
space $\v_b$ requires a reduced decomposition of 
$\pi_b\in \hW$, the {\em reduction
of $b$ modulo $W$}. To be exact,  $\pi_b$ 
is defined as the element of minimal length in the coset 
$bW\in \hW$ (it is unique). The space $\v_b$ 
does not depend on the particular choice of the reduced 
decomposition of $\pi_b$ and can be calculated 
in terms of the {\em relative Bruhat ordering\,}
defined for $\tR\mod\tR^0$. This space consists of
{\em generalized} $Y$\~eigenvectors for the eigenvalue 
$-b_\#$ (see the definition in Proposition \ref{YONE}); 
generally, not all of them. Its dimension is 
always finite even if $q$ is a root of unity. 

We note that the elements $\{\pi_b\}$ and the
standard Bruhat ordering on $\hW$ govern
the combinatorics of {\em affine Schubert manifolds}, although
$\tR^0$ does not appear in the theory of the affine
Grassmanian and such manifolds. 
\smallskip

{\em Non-semisimple Macdonald
polynomials\,}. They are 
top polynomial $\widetilde{E}_b\in \v_b$.
The construction is explicit;\, $\widetilde{E}_b$ is 
given in terms  of the {\em generalized chain
of intertwiners} (discussed below) corresponding to 
a given reduced decomposition of $\pi_b$. 

The polynomials $\widetilde{E}_b$ are not unique and may depend 
on the choice of this decomposition modulo polynomials 
from ``lower" $\v_c$, although the actual flexibility of
this definition is limited. In many cases they 
are determined uniquely even if dim\,$\v_b>1$.
These polynomials have correct leading terms and 
form a basis of $\v$; \ 
$\widetilde{E}_b=E_b$ if dim\,$\v_b=1$.
\smallskip

{\em Generalized chains\,}. 
The $Y$\~intertwiners corresponding to
simple $\al_i\,(i>0)$ are  
$(T_i+ c(Y_{\al_i}))$ for nonaffine simple roots $\al_i\in R_+$ 
(treated as vectors in $P$) and 
$c(x)=(t_i^{1/2}-t_i^{-1/2})/(x^{-1}-1).$
It is more involved for the affine $\al_0$; namely, 
following \cite{C1}, we need to consider
$\tau_+(T_0)$ instead of $T_0$ for the automorphism
$\tau_+$ of DAHA (formally) corresponding to multiplication
by the Gaussian.  Here $c(Y_{\al_i})$ 
may become infinite in a chain of intertwiners; we 
call the corresponding places {\em singular\,} in the paper,
as well as corresponding simple reflections
and roots from $\la(\pi_b)$.

We define generalized chains of intertwiners
by replacing singular simple intertwiners
$(T_i+ \infty)$ by $T_i$.
See (\ref{hatefinal}), (\ref{hatefinall})
for exact definitions. Choosing $T_i$ here
ensures the proper leading terms of the non-semisimple
Macdonald polynomials. Another motivation is that
the relations between the intertwiners and $T$ 
are of fundamental importance (note a connection with
virtual links).
\smallskip

Two natural problems arise: \\
(a) determining how the non-semisimple
polynomials $\widetilde{E}_b$ depend on the
choice of the reduced decompositions of $\pi_b$, which 
is mainly covered by Theorem \ref{PHIBRUHAT};\\
(b) finding ``large" families of $b\in P$ such that 
dim\,$\v_b=1$ and therefore $\widetilde{E}_b=E_b$, which is
addressed in Proposition \ref{HWALBE} and in the
following theorem.  
\smallskip

{\em Main Theorem \ref{HWALSIMPLE}}. It contains 
efficient tools for solving (b) and exact
calculating the spaces $\v_b$ and the spaces 
$\v(-b_\#)^\infty$. The latter spaces are 
defined as the spaces of 
{\em all generalized eigenvectors} with the 
$Y$\~eigenvalue $-b_\#$ serving $\widetilde{E}_b$;
see the definition in Proposition \ref{YONE}. 
We analyze what happens with $\v_c$ if $\pi_c$ is 
replaced be $\pi_c\hw$ for arbitrary $\hw\in \hW$, 
then do it for $\hw\in \hW$
from the centralizer of the weight $-0_\#=\rho_k$ (the 
$Y$\~eigenvalue of  $1\in \v$) and, 
finally, for $\hw$ from $\tW^0$, the Weyl group of $\tR^0$.
\smallskip
 
An application of Proposition \ref{HWALBE}
and Theorem \ref{HWALSIMPLE} is 
Theorem \ref{FEIGIN} about the canonical
semisimple submodule $\v_{ss\,}\subset\v$. 
It generalizes 
the construction from 
\cite{Ka} in the $A$\~case based, in its turn,
on paper \cite{FJMM}
where the symmetric case was considered. 

We consider the $A$ case in detail because it is 
important to establish the connection with these papers. 
However we do not give the final list of root systems and
parameters $k$ when $\v_{ss\,}$ is nonzero. Obtaining these
conditions and an explicit  
description of $\v_{ss}$ does not look very difficult, 
but involving other methods is more reasonable here. 

To be more precise, we do not touch in the paper the 
{\em wheel condition} from \cite{FJMM,Ka}, an alternative
way (better to say, complementary) to introduce $\v_{ss\,}$
as an ideal in the polynomial
representation. Counterparts of the wheel condition
can be obtained for (many, maybe all) 
root systems within the technique
of this paper; actually, only the evaluation pairing is necessary.
However we prefer to postpone with finalizing the consideration of 
$\v_{ss\,}$ until future papers on the Jantzen\~type filtration.
We only mention here that the Kasatani conjecture \cite{Ka}, 
verified in \cite{En} via the localization functor,
can be almost certainly managed using {\em directly} the 
technique of intertwiners. 

Generally, it is expected that the Jantzen\~type
filtration gives
a natural way to decompose $\v$ with the constituents that
are irreducible in many cases, including the $A$\~case.
We consider only two examples in this paper, the ``highest"
and the ``lowest" constituents, namely, the 
quotient $\v/Rad$ and the submodule $\v_{ss\,}$. Mainly we 
discuss $\v$ subject to $Rad=\{0\}$, but the
same technique can be used for $\v/Rad$.
\medskip

{\bf Further topics.}
The spaces $\v_b$ are not the
smallest ``natural" space of generalized
eigenvectors corresponding to the eigenvalue $-b_\#$ 
that contain $\widetilde{E}_b$.
Following the chain of intertwiners for $\widetilde{E}_b$,
we define the spaces 
$\widetilde{E}_b\in \widetilde{V}_b\subset \v_b,$ 
that form a {\em maximal system\,} of subspaces in $\v$, 
like $\{\v_b\}$, with $\{\widetilde{E}_b\}$ as 
top polynomials. The definition is more direct than that for
$\{\v_b\}$:\ given $b$, the space
$\widetilde{V}_b$ is given {\em exactly} in terms of 
the reduced decomposition of $\pi_b$ that
is used for $\widetilde{E}_b$. However, in contrast to
$\v_b$, this space may depend on a particular choice of 
the reduced decomposition. 

The spaces $\v_b$ and $\widetilde{V}_b$ are $Y$\~modules. 
The space $\widetilde{V}_b$ is
$Y$\~cyclic in many cases; the strongest result in this
direction we have is Proposition \ref{VINDUCEDY} based
on the following theorem.
\smallskip

{\em Main Theorem \ref{THEORCHAIN}}. We prove that
$\widetilde{V}_b$ considered as a $Y$\~module is 
{\em covered} by a certain $Y$\~module defined explicitly 
in terms of the generators and relations with the 
structural constants that are essentially {\em integers\,}. 
It is defined via the Demazure operators and is connected
with Schubert polynomials. 
Due to such integrality, we can switch here from 
DAHA to any degeneration we wish. 
The proof of the theorem is a straightforward
calculation; it gives a natural direct approach to
the $Y$\~cyclicity of $\widetilde{V}_b$. 

In a sense,
$\widetilde{V}_b$ is cyclic ``almost always". To be
more exact, if
the set $\la(\pi_b)$ is sufficiently small, then 
$\tR^0$ contains only pairwise orthogonal roots
and we can use Corollary \ref{HATEECOM}. 
If $\la(\pi_b)$ is sufficiently 
large, then Main Theorem \ref{RANKTWO} can be applied
to collect all singular roots in $\la(\pi_b)$ 
in a connected segment using the Coxeter transforms and 
Proposition \ref{VINDUCEDY} can be employed.
\smallskip

The $Y$\~cyclicity of $\widetilde{V}_b$
(if known) can certainly simplify
using the technique of intertwiners for
$\v$ and other induced modules where it holds.
Concerning general $Y$\~induced modules (free modules
induced $Y$\~eigenvectors), 
we need them only a little in this part of the paper. 
Generalizing $\widetilde{V}_{b}$, we define the spaces 
$\widetilde{V}_{\hw}$; note that using $\hw=\pi_b$ only
is a special feature of $\v$. 
The modules $\widetilde{V}_{\hw}$
are {\em cyclic} for sufficiently large $\la(\hw)$
(not one-dimensional as for $\v$); 
see Proposition \ref{VINDUCEDY}.
\medskip
      
{\em Irreducibility of $\v$.}
The end of the paper is devoted to 
Main Theorem \ref{RADZERO} on irreducibility
of $\v$ subject to the condition $Rad=\{0\}$, which
demonstrates
almost all aspects of the technique of intertwiners.

The condition $Rad=\{0\}$
readily gives that all three spaces  
$\widetilde{V}_b\subset \v_b\subset \v(-b_\#)^\infty$
contain a {\em unique} Macdonald polynomial $E_{b^\circ}$, 
where $b^\circ\in P$ is
defined combinatorially in terms of $b$.
It suffices to consider the chains of 
intertwiners that end at $E_b$ (i.e.,  
for $b=b^\circ$). 
Only non-invertible intertwiners may lead to reducibility. 
Moreover, $Rad=\{0\}$ implies that 
intertwiners of type $(T-t^{1/2})$ cannot appear 
when the chain goes from $E_b$ to $E_c$ avoiding 
{\em singular} intertwiners. 
Assuming that $\v'\subset \v$ is a proper DAHA
submodule, we proceed as follows. 
\smallskip

Given a chain of intertwiners,
let $E_b$ be the first Macdonald polynomial 
in $\v'$. Then 
the previous one is $E_a\not\in\v'\,$ and
the intertwiner between them can be only of type 
$(T+t^{-1/2})$.

We go from $E_b$ until the first {\em singular} 
intertwiner {\em if it exists} (not always).
Let $E_c$ be the last in this chain 
before this place for a reduced decomposition
extending that for $\pi_a$: $\la(\pi_a)\subset \la(\pi_c)$.
Then we apply the corresponding $T$ to $E_c$ 
and then go back following the sequence
$\la(\pi_c)\setminus \la(\pi_a)$ taken in the opposite
order, from the last root to the first.  

The resulting polynomial (the ``end" of this
chain) will belong to $\v'\cap \v(-a_\#)^\infty$.
The analysis shows that $Rad=\{0\}$ implies that it must be
proportional to $E_a$ with a nonzero coefficient of
proportionality. Therefore $E_a$ belongs to $\v'$,
a contradiction.
%
\vskip 0.2cm

This method, 
{\em reflection at the first singular place},
is actually of general nature and can 
be applied in various situations. For the first time in the
DAHA context it was used in \cite{CO} for $A_1$, where
it gives a complete decomposition of $\v$ including the
cases of roots of unity. Let us discuss some
combinatorial aspects of this method.
\smallskip
 
{\em Zigzags}.
In this proof, Lemma \ref{RADADE} is needed for managing
the situation when there are several intertwiners
of type $(T+t^{-1/2})$ in
$\la(\pi_c)\setminus \la(\pi_b)$.
For coinciding $t$ (as $k_{\sht}=\nu_{\lng}k_{\lng}$), 
this lemma is essentially on the combinatorics of the
sets $R_{ht}$ of roots $\al\in R_+$ with fixed
$ht\equal(\al,\rho^\vee)$. We introduce the  
{\em zigzags\,} alternating between 
$R_{ht}$ and $R_{ht+1}$ with the {\em links\,}
(as $t_{\lng}=t_{\sht}$)
corresponding to adding or subtracting
simple roots, quite a classical matter.
The claim is that {\em any maximal 
zigzag contains at least one endpoint
from $R_{ht}$.} 
 
This claim is actually
a combinatorial variant of the formula for
$\Pi_R$.
If the $k$\~parameters are arbitrary, then 
smaller subsets must be considered instead
of $R_{ht}$ and {\em links} become
somewhat more involved; however the  
claim about the endpoints holds. This approach
is expected to give a description of the cases 
when $\v/Rad$ is not irreducible. 

Generally, the technique
of intertwiners alone does not seem sufficient 
for decomposing $\v$ and becomes too combinatorial
even for managing the irreducibility of $\v$ subject
to $Rad=\{0\}$. We decided to omit the details
of the zigzag construction in this paper. 

\medskip
\subsection{Expected applications}\label{sec:EXPECTEDAP}
Decomposing the polynomial representation $\v$
is of key importance in the theory of DAHA and
is expected to have many applications. In the rational
case, the theory of $\v$ was started by Opdam; 
see \cite{DO,DJO}. From the viewpoint of
applications, the general $q,t$\~ case seems the
most fruitful, although the degenerate cases
have important applications too.

As a DAHA\~module, $\v$ is a 
{\em universal spherical representation}; 
this alone is sufficient to study it thoroughly.
Moreover, its identification
with the algebra of $X$\~polynomials makes $\v$
and all its constituents {\em commutative algebras}.
The quotients possess the unit, other constituents do not. 
There is also a natural projective action 
of $PSL_2(\Z)$  on {\em finite dimensional} 
constituents of $\v$ (maybe reducible) subject to 
certain technical restrictions, simple to control. 

The latter action exists because of the following
conceptual reason. The automorphism
$\tau_-$, one of two generators $\tau_\pm$
of the projective $PSL_2(\Z)$, is an outer
automorphism of DAHA formally corresponding 
to multiplication (conjugation)
by the $Y$\~Gaussian. It {\em always} acts in 
$\v$ and its constituents;
see Proposition \ref{TAUPSI}. However,
the generator $\tau_+$, represented
by the $X$\~Gaussian (an infinite Laurent series)
may act only in finite dimensional
constituents of $\v$ (the condition dim$<\infty$
is not always sufficient).
\smallskip

The expectations are that $\v$ and its quotients
serve quite a few examples of monoidal categories
(with tensoring), especially if the action 
of $PSL_2(\Z)$ is present there, and go well beyond
such examples. 

The celebrated Verlinde algebras are the key example.
They are interpreted in the DAHA theory as symmetric 
subalgebras of {\em perfect} quotients of $\v$  
as $q$ is a root of unity in the so-called
group case  $t=q$ (that is the simplest possible
setting in the DAHA theory apart from $t=1$).  

The multiplication in $\v$ leads to the {\em fusion} 
of integrable Kac\~ Moody modules. The DAHA\~action of 
$PSL_2(\Z)$ is nothing but the action of
 Verlinde operators $S,T$;
\ $T$ is multiplication be the $X$\~Gaussian, 
$S$ becomes the DAHA\~Fourier transform. 

Using the terminology from
\cite{C101}, symmetrizations
of {\em perfect representations\,} are natural 
generalizations of the  Verlinde algebras. However,  
the perfect representations
are well beyond the usual Verlinde algebras.
Let us comment on it.
\smallskip
 
First, one can take any 
$t$ instead of $t=q$ provided that
the corresponding perfect representation exists;
this means that the characters in the Verlinde theory 
will be replaced by the symmetric Macdonald polynomials
(treated as functions on at certain finite sets of points).

Second, the {\em non-symmetric} Macdonald polynomials
can be considered here instead of their symmetrization.
Furthermore, the {\em non-semisimple} Macdonald polynomials
can be taken, corresponding to {\em non-semisimple} 
counterparts of perfect representations.

Third, it is not necessary to assume that $q$ is a root of
unity (one of the main discoveries of the DAHA
approach); an important part of the theory of perfect
representations is for generic $q$ (then $k$
must be {\em singular}). They are finite dimensional
and has all key structures of Verlinde algebras.
\medskip

{\bf Possible relations.} The following is a sketch
of (some) known and expected applications 
of $\v$ and its constituents.

\smallskip
(a) Presumably all
{\em Verlinde-type algebras\,}, describing
``fusion" of integrable modules 
for the Kac-Moody algebras,
Virasoro algebras, $\w$\~algebras and similar objects,
are quotients or constituents of $\v$
(see confirmations in \cite{FHST},\cite{MT}). Moreover,
infinite dimensional constituents of $\v$
are expected to be connected with the theory  
at arbitrary Kac\~Moody central charge $c$. 
The most interesting case $|q|=1$ 
when $q$ is {\em not} a root of unity
(then $c\in \R$) may lead to the 
{\em Kac-Moody $L^2$\~theory}, generalizing 
the classical harmonic analysis on symmetric spaces.
Straight attempts to create such theory were unsuccessful.
\smallskip

(b) Similar expectations are for the tensor category
of {\em all} representations of Lusztig's quantum group
at roots of unity. The corresponding
Verlinde algebra describes the {\em reduced subcategory}
of this category; the equivalence with the definition
of Verlinde algebras in the Kac-Moody theory is due to 
Kazhdan, Lusztig \cite{KL2} and Finkelberg.
The first known examples of the {\em non-semisimple} Verlinde
algebras look very similar to what can be expected in the
so-called ``case of parallelogram". Note that the 
monoidal structure (fusion, tensoring) and the
action of $PSL_2(\Z)$ are generally difficult problems for 
the complete Lusztig's category.
\smallskip

(c) If an arbitrary {\em perfect representation}
is taken, then no categorical interpretation is
known and expected (there are no reasons
for integrality and
positivity of the structural constants). 
However, all such representations give
important examples of Fourier transform theories 
satisfying {\em all} standard  
classical properties. This line {\em directly} generalizes 
Fourier transforms associated
with the irreducible modules in the theory of
Weyl algebras (non-commutative tori) at roots of unity.
Among other applications, perfect representations are 
related to the generalized Macdonald eta-type identities, 
Gaussian sums \cite{C7}
and the so-called diagonal coinvariants \cite{Ha,Go,C29}.
\smallskip

(d) Semisimple {\em submodules\,} of $\v$ 
generalize the construction from \cite{FJMM}
of {\em ideals\,} in the ring of symmetric
polynomials of type $A_n$ linearly generated
by symmetric Macdonald polynomials.
These ideals are expected to be meaningful in the theory of
$\widehat{\mathfrak{gl}}_N$ and the corresponding $\w$-algebras,
presumably, via the duality from \cite{VV} and \cite{STU}.
They also give some kind of restriction maps
$\mathfrak{gl}_M\subset\mathfrak{gl}_N,$ although DAHA
generally do not have {\em straight} embeddings of this type
(unless the approach from \cite{BE} or similar methods 
involving completions of DAHA are used).
\smallskip

(e) The theory of $\v$ is connected
with the Plancherel formula on the {\em affine}
Hecke algebra due to Macdonald (the spherical
case), Matsumoto, Lusztig and many others,
describing the decomposition
of the regular representation of the affine
Hecke algebra. The regular representation of AHA
is interpreted as an induced module (depending 
on a generic weight) in the 
DAHA theory. Its spherical part 
is associated with the weight $-\rho_k$
and can be identified with $\v$. The relation to \cite{KL1,HO2} 
is via a new theory of Jantzen\~type filtration of DAHA
considered in the limit $q\to 0$; it may include
applications to square integrable and tempered
irreducible AHA\~representations.  
\smallskip

(f) In the rational case, {\em singular} $k$
(when the radical becomes nonzero) correspond to
singular multi-dimensional Bessel functions.
The most degenerate case is when
$\v_{rat}$ has a finite dimensional quotient. Such
quotients found important combinatorial applications,
for example, in \cite{Go}. They may be related to the
minimal conformal theories based on    
Virasoro-type algebras and their ``perturbations" (adding
additional parameters). The $q,t$\~case is connected with the
rational theory in many ways; for instance,
certain Verlinde algebras can be $q$\~deformed and then 
identified with their rational limits \cite{C7}. 
It may reflect various relations between
the Kac-Moody and Virasoro theories.
\smallskip


(g) The exponential map from \cite{C29}
can be used, in principle, to establish a
correspondence between the decomposition theory of
the polynomial representation for the rational
DAHA \cite{DJO,DO} and that in the $q,t$-case, 
although this approach (generally) requires analytic setting. 
Recent \cite{Et} is a step in this direction. 
Algebraically, the exponential map identifies
finite dimensional DAHA modules and their rational
degeneration. This map 
is connected with the localization functor,
the monodromy of a KZ\~type connection, from 
\cite{GGOR, VV1}. 
Using the exponential map and the KZ\~monodromy
for the polynomial and other infinite 
dimensional representations triggers interesting
problems and certainly must be studied thoroughly.
\smallskip

(h) The evaluation pairing of $\v$ plays the key 
role in the DAHA theory. For the Macdonald polynomials, 
it is given in terms of 
their values at certain ``shifted lattices" of points
(weights). These lattices of points are {\em exactly} those used 
in the theory of {\em interpolation polynomials\,} 
studied by Knop \~ Sahi, Okounkov \~ Olshantsky and others.
The definition requires the ``semisimplicity" of these sets
(directly related to the semisimplicity of $\v$). 
The consideration of the non-semisimple quotients of the
polynomial representation in this paper may lead to the 
theory of ``non-semisimple" interpolation polynomials. 
There are classical constructions that confirm that
such theory must exist.
\smallskip

(i) In the semisimple case, the polynomial representation
is completely governed by the intertwining operators;
the latter give a representation of the affine Weyl group
$\v$. The non-semisimple setting requires considering
intertwiners {\em together} with the $T$\~generators.
The corresponding formalism is connected with
the theory of virtual links in topology 
(when extended from $A$ to arbitrary root systems). 
We note that only the DAHA\~modules that do not involve 
non-invertible intertwiners (singular ones may appear) 
are {\em directly} connected with 
Kauffman's axioms. In a somewhat different (but related)
direction, the weight decompositions of semisimple DAHA modules 
may give something for the theory of ``categorization". 
\smallskip

(j) The {\em localization functor},
the monodromy map of a KZ\~type connection,
from \cite{GGOR} (the rational case) and
\cite{VV1} (the differential-trigonometric case) is
an important motivation for the direction of this paper.
For instance, the {\em affine exponents\,} are directly
connected with the non-semisimplicity of non-affine
Hecke algebras at roots of unity. Practically all problems 
discussed in the paper can be translated to the corresponding 
problems for Hecke algebras at roots of unity and/or for the
Schur algebras. Generally, the geometric theory
of Hecke algebras is better developed than that for DAHA.
However, in the DAHA theory there is a greater potential
of using relatively elementary (non-geometric tools)
like those developed in this paper.    
\smallskip

(k) Paper \cite{BE} on the parabolic induction
in the DAHA theory opens a systematic way for analyzing
induced representations (including $\v$). It does not
cover the $q,t$\~setting, but, presumably, can be 
generalized. The representations parabolically
induced from finite dimensional ones are in the focus of
\cite{BE}. Paper \cite{VV3} gives a classification
of finite dimensional modules in the spherical case for
generic $q$; it includes {\em perfect representations} 
and gives some non-semisimple ones. Before this paper,
a complete classification of finite dimensional modules was 
obtained only for $A$ in \cite{BEG} (all appeared perfect).
Understanding parallelism with the $p$\~adic theory is a 
natural challenge.


\medskip
{\bf Acknowledgements.}
The author is thankful to A.~Garsia, D.~Kazhdan,
E.~Opdam and N.~Wallach
for stimulating discussions.
Special thanks to P.~Etingof. The author
is grateful to RIMS, Kyoto University, for the invitations
in 2005, when the paper was started
and the author delivered a series of
lectures on the classification of semisimple DAHA modules
of type $A$, and in 2007, when this paper
was essentially completed. The author thanks M.~Kashiwara and
T.~Miwa for their hospitality and also thanks
IHES for the invitations in 2006-07. 

\medskip
\setcounter{equation}{0}
\section{Affine Weyl groups}
\setcounter{equation}{0}
Let $R=\{\al\}   \subset \R^n$ be a root system of type
$A,B,...,F,G$
with respect to a euclidean form $(z,z')$ on $\R^n
\ni z,z'$,
$W$ the {\dfont Weyl group} \index{Weyl group $W$}
generated by the reflections $s_\al$,
$R_{+}$ the set of positive  roots ($R_-=-R_+$)
corresponding to fixed simple 
roots $\al_1,...,\al_n,$
$\Ga$ the Dynkin diagram
with $\{\al_i, 1 \le i \le n\}$ as the vertices.

We will also use
the dual roots (coroots) and the dual root system:
$$R^\vee=\{\al^\vee =2\al/(\al,\al)\}.$$

The root lattice and the weight lattice are:
\begin{align}
& Q=\oplus^n_{i=1}\Z \al_i \subset P=\oplus^n_{i=1}\Z \om_i,
\notag
\end{align}
where $\{\om_i\}$ are fundamental weights:
$ (\om_i,\al_j^\vee)=\de_{ij}$ for the
simple coroots $\al_i^\vee.$

Replacing $\Z$ by $\Z_{\pm}=\{m\in\Z, \pm m\ge 0\}$ we obtain
$Q_\pm, P_\pm.$
Note that $Q\cap P_+\subset Q_+.$ Moreover, each $\om_j$ has all
nonzero coefficients (sometimes rational) when expressed
in terms of
$\{\al_i\}.$
Here and further see  \cite{Bo}.

The form will be normalized
by the condition  $(\al,\al)=2$ for the
{\em short} roots. This normalization coincides with that
from the tables in  \cite{Bo}  for the systems $A,C,D,E,G.$
Thus,

\centerline{
$\nu_\al\equal (\al,\al)/2$ can be either $1,$ or $\{1,2\},$ or
$\{1,3\}.$ }
\noindent
We will use the
notation $\nu_{\lng}$ for the long roots ($\nu_{\sht}=1$).

This normalization leads to the inclusions
$Q\subset Q^\vee,  P\subset P^\vee,$ where $P^\vee$ is
defined to be generated by
the fundamental coweights $\om_i^\vee.$

\smallskip
Let  $\vth\in R^\vee $ be the {\em maximal positive
coroot}. All simple coroots appear in its
decomposition in $R^\vee$. It also belongs to $R$, i.e. is a root,
because of the choice of normalization; so all simple roots
appear in its decomposition in $R$.

Also note that $2\ge (\vth,\al^\vee)\ge 0$ for $\al>0,$
$(\vth,\al^\vee)=2$ only for $\al=\vth,$ and
$s_{\vth}(\al)<0$ if $(\vth,\al)>0.$

See  \cite{Bo} to check that  $\vth$ considered
as a root
is maximal among all short positive roots of $R.$
It is also the least nonzero element in $Q_+^{+}=Q_+\cap
P_+=Q\cap P_+$
with respect to $Q_+.$

Setting
$\nu_i\ =\ \nu_{\al_i}, \
\nu_R\ = \{\nu_{\al}, \al\in R\},$ one has
\begin{align}\label{partialrho}
&\rho_\nu\equal (1/2)\sum_{\nu_{\al}=\nu} \al \ =
\ \sum_{\nu_i=\nu}  \om_i, \hbox{\ where\ } \al\in R_+,
\ \nu\in\nu_R.
\end{align}
Note that $(\rho_\nu,\al_i^\vee)=1$ as $\nu_i=\nu.$
We will call $\rho_\nu$ {\dfont partial} $\rho.$
\index{rho@rho $\rho$}
\index{raaaa@$\rho$}

\subsection{Affine roots}
The vectors $\ \tal=[\al,\nu_\al j] \in
\R^n\times \R \subset \R^{n+1}$
for $\al \in R, j \in \Z $ form the
\index{affine root system} {\dfont affine root system}
$\tR \supset R$ ($z\in \R^n$ are identified with $ [z,0]$).
We add $\al_0 \equal [-\vth,1]$ to the simple
 roots for the
{\dfont maximal short root}
\index{maximal short root $\vth$} $\vth$.
The corresponding set
$\tR_+$ of positive roots coincides with 
$R_+\cup \{[\al,\nu_\al j],\ \al\in R, \ j > 0\}$.

We will sometimes write $\tR=\widetilde{A}_n,
\widetilde{B}_n,\ldots, \widetilde{G}_2$ when dealing
with concrete affine root systems defined as above.

Any positive affine root $[\al,\nu_\al j]$
is a linear combinations
with non-negative integral coefficients of
$\{\al_i,\,0\le i\le n\}$.
Indeed, it is well known that 
$[\al^\vee,j]$ is such combination
in terms of  $\{\al_i^\vee,\, 1\le i\le n\}$
and $[-\vth,1]$ for the system of affine {\em coroots},
that is $\tR^\vee=\{ [\al^\vee,j],\, \al\in R,\, j\in\Z \}$.
Hence, $[-\al,\nu_\al j]=$  $\nu_\al[-\al^\vee,j]$ has
the required representation.

Note that the sum of the long roots is always long,
the sum of two short roots can be a long root only
if they are orthogonal to each other. This property
gives another justification of the claim 
that $\tR$ is a root system.

We complete the Dynkin diagram $\Ga$ of $R$
by $\al_0$ (by $-\vth$, to be more
exact); it is called {\dfont affine Dynkin diagram}
$\tGa$. One can obtain it from the
completed Dynkin diagram from \cite{Bo} for 
the {\em dual system}
$R^\vee$ by reversing all arrows. 
The number of laces between $\al_i$ and
$\al_j$ in $\tGa$ will be denoted by $m_{ij}.$

The set of
the indices of the images of $\al_0$ by all
the automorphisms of $\tGa$ will be denoted by $O$
($O=\{0\} \for E_8,F_4,G_2$). Let $O'=\{r\in O, r\neq 0\}$.
The elements $\om_r$ for $r\in O'$ are the so-called minuscule
weights: $(\om_r,\al^\vee)\le 1$ for
$\al \in R_+$.

Given $\tal=[\al,\nu_\al j]\in \tR,  \ b \in P$, let
\begin{align}
&s_{\tal}(\tz)\ =\  \tz-(z,\al^\vee)\tal,\
\ b'(\tz)\ =\ [z,\ze-(z,b)]
\label{ondon}
\end{align}
for $\tz=[z,\ze] \in \R^{n+1}$.

The \index{affine Weyl group $\widetilde{W}$}
{\dfont affine Weyl group} $\tW$ is generated by all $s_{\tal}$
(we write $\tW = \lan s_{\tal}, \tal\in \tR_+\ran)$. One can take
the simple reflections $s_i=s_{\al_i}\ (0 \le i \le n)$
as its
generators and introduce the corresponding notion of the
length. This group is
the semidirect product $W\lsmash Q'$ of
its subgroups $W=$ $\lan s_\al,
\al \in R_+\ran$ and $Q'=\{a', a\in Q\}$, where
\begin{align}
& \al'=\ s_{\al}s_{[\al,\,\nu_{\al}]}=\
s_{[-\al,\,\nu_\al]}s_{\al}\for
\al\in R.
\label{ondtwo}
\end{align}

\index{W@$\tW$,$\hW$:\  affine Weyl groups}
The
\index{extended Weyl group}
\index{affine Weyl group $\widehat{W}$}
{\dfont extended Weyl group} $ \hW$ generated by $W\and P'$
(instead of $Q'$) is isomorphic to $W\lsmash P'$:
\begin{align}
&(wb')([z,\ze])\ =\ [w(z),\ze-(z,b)] \for w\in W, b\in P.
\label{ondthr}
\end{align}
From now on,  $b$ and $b',$ $P$ and $P'$ will be identified.

Given $b\in P_+$, let $w^b_0$ be the longest element
in the subgroup $W_0^{b}\subset W$ of the elements
preserving $b$. This subgroup is generated by simple
reflections. We set
\begin{align}
&u_{b} = w_0w^b_0  \in  W,\ \pi_{b} =
b( u_{b})^{-1}
\ \in \ \hW, \  u_i= u_{\om_i},\pi_i=\pi_{\om_i},
\label{xwo}
\end{align}
where $w_0$ is the longest element in $W,$
$1\le i\le n.$

The elements $\pi_r\equal\pi_{\om_r}, r \in O'$ and
$\pi_0=\hbox{id}$ leave $\tGa$ invariant
and form a group denoted by $\Pi$,
 which is isomorphic to $P/Q$ by the natural
projection $\{\om_r \mapsto \pi_r\}$. As to $\{ u_r\}$,
they preserve the set $\{-\vth,\al_i, i>0\}$.
The relations $\pi_r(\al_0)= \al_r= ( u_r)^{-1}(-\vth)$
distinguish the
indices $r \in O'$. Moreover (see e.g., \cite{C2}):
\begin{align}
& \hW  = \Pi \lsmash \tW, \where
  \pi_rs_i\pi_r^{-1}  =  s_j \iif \pi_r(\al_i)=\al_j,\
 0\le j\le n.
\end{align}
\medskip
\subsection{The length 
\texorpdfstring{{\mathversion{bold}on $\hW$}}{}}
Setting
$\hw = \pi_r\tw \in \hW,\ \pi_r\in \Pi, \tw\in \tW,$
the length $l(\hw)$
is by definition the length of the reduced decomposition
$\tw= $ $s_{i_l}...s_{i_2} s_{i_1} $
in terms of the simple reflections
$s_i, 0\le i\le n.$ The number of  $s_{i}$
in this decomposition
such that $\nu_i=\nu$ is denoted by   $l_\nu(\hw).$

The \index{length $l(\hw)$} {\dfont length} can be
also defined as the
cardinality $|\la(\hw)|$
of the {\dfont $\la$\~set} of $\hw$\,:
\begin{align}\label{lasetdef}
&\la(\hw)\equal\tR_+\cap \hw^{-1}(\tR_-)=\{\tal\in \tR_+,\
\hw(\tal)\in \tR_-\},\
\hw\in \hW.
\end{align}
Respectively,
\begin{align}
&\la(\hw)=\cup_\nu\la_\nu(\hw),\
\la_\nu(\hw)\
\equal\ \{\tal\in \la(\hw),\nu({\tal})=\nu \}.
\label{xlambda}
\end{align}
Note that $\la(\hw)$ is closed with respect to positive linear
combinations. More exactly, if $\tal=u\tbe+v\tga\in \tR$ for
rational
$u,v>0$, then $\tal\in \la(\hw)$ if  $\tbe\in \la(\hw) \ni \tga$.
Vice versa, if $\la(\hw)\ni \tal=u\tbe+v\tga$ for
$\tbe,\tga\in \tR_+$ and rational $u,v>0$, then either $\tbe$ or
$\tga$
must belong to $\la(\hw)$. Also,
\begin{align}\label{laclosed}
&\tal=[\al,\nu_\al j]\in\la(\hw)\, \Leftrightarrow\,
\{\, [\al,\nu_\al i]\in \la(\hw)\notag\\
&\hbox{for\ all\ } 0\le i\le j\hbox{\ where\ }
i>0 \hbox{\ as\ } \al<0\,\}.
\end{align}
This property is obvious because $\hw([\al,\nu_\al i])=
\hw(\tal)+[0,\nu_\al (i-j)]$.
\smallskip

The coincidence with the previous definition
is based on the equivalence of the {\em length equality}

\begin{align}\label{ltutwa}
&(a)\ \ l_\nu(\hw\hu)=
l_\nu(\hw)+l_\nu(\hu)
\for \hw,\hu\in\hW
\end{align}
and the {\em cocycle relation}
\begin{align}
&  (b)\ \ \la_\nu(\hw\hu) = \la_\nu(\hu) \cup
\hu^{-1}(\la_\nu(\hw)),
\label{ltutw}
\end{align}
which, in its turn, is equivalent to
the {\em positivity condition}
\begin{align}\label{ltutwc}
& (c)\ \  \hu^{-1}(\la_\nu(\hw))
\subset \tR_+
\end{align}
and is also equivalent to the {\em embedding condition}
\begin{align}\label{ltutwd}
& (d)\ \  \la_\nu(\hu)\subset \la_\nu(\hw).
\end{align}

Formula (\ref{ltutw}) obviously includes 
the positivity condition (\ref{ltutwc}). It also
readily gives ($d$) and implies that
$$
\la_\nu(\hu) \cap \hu^{-1}(\la_\nu(\hw))\ =\
\hu^{-1}\bigr(\hu(\la_\nu(\hu)) \cap \la_\nu(\hw)\bigl)\
=\ \emptyset
$$ thanks to the general formula
$$
\la_\nu(\hw^{-1}) = -\hw(\la_\nu(\hw)).
$$
Thus it results in the equality $l_\nu(\hw\hu)=
l_\nu(\hw)+l_\nu(\hu)$ and we have the implications
$(a)\Leftarrow (b)\Rightarrow (c).$

The remaining implications $(a)\Rightarrow (b)\Leftarrow (c)$
are based on the
following simple general fact:
\begin{align}
&  \la_\nu(\hw\hu)\setminus \{ \la_\nu(\hw\hu)\cap
\la_\nu(\hu)\}
=\hu^{-1}(\la_\nu(\hw))\cap \tR_+ \hbox{ \ for \ any \ }
\hu,\hw.
\label{latutw}
\end{align}
For instance, the length equality ($a$) readily
implies ($b$) and ($d$) results in ($c$).
For the sake of completeness,  let us
deduce ($b$) from the positivity condition ($c$).
We follow \cite{C2}.

It suffices to check that  $\la_\nu(\hw\hu)\supset
\la_\nu(\hu).$
If there exists a (positive)  $\tal\in  \la_\nu(\hu)$
such that
$(\hw\hu)(\tal)\in \tR_+,$ then
$$
\hw(-\hu(\tal))\in \tR_-\ \Rightarrow\
-\hu(\tal)\in \la_\nu(\hw)\ \Rightarrow\
-\tal\in \hu^{-1}(\la_\nu(\hw)).
$$
We come to a contradiction with the positivity.
Hence $(a)\Leftrightarrow (b)\Leftrightarrow (c).$
\smallskip

Note that the embedding condition ($d$) readily gives
the following well-known fact. Let $R'$ be a root
subsystem $R'\subset\tR$ with the simple roots
that are simple in $\tR$ constituting
a connected subset of the affine Dynkin diagram
$\widetilde{\Gamma}$, $w_0'$ the greatest element
in the corresponding Weyl group $W'$. Then
\begin{align}\label{hwwoprime}
l(\hw)=l(w_0'\hw)+l(w_0') \Leftrightarrow
l(s_i\hw)=l(\hw)-1 \hbox{\ for\ all\ } \al_i\in R',
\end{align}
i.e., $\hw$ is divisible on the left by $w_0'$
in the sense of reduced decompositions
if and only if it is divisible by all $s_i$
for $\al_i\in R'$.
Indeed, it is equivalent to the ``inverse" statement
$$
\ \ \ \ \ \ \, l(\hw)=l(\hw w_0')+l(w_0')
\Leftrightarrow
l(\hw s_i)=l(\hw)-1 \hbox{\ for\ all\  }\al_i\in R',
$$
that directly follows from ($d$) and the inclusion
$\la(w_0')\subset\la(\hw)$. See also
Theorem \ref{INTRINLA} below.
\smallskip

Applying (\ref{ltutw}) to the reduced decomposition
$\hw=\pi_rs_{i_l}\cdots s_{i_2}s_{i_1},$
\begin{align}
\la(\hw) = \{\ &\tal^l=\tw^{-1}s_{i_l}(\al_{i_l}),\
\ldots,\ \tal^3=s_{i_1}s_{i_2}(\al_{i_3}),\notag\\
&\tal^2=s_{i_1}(\al_{i_2}),\ \tal^1=\al_{i_1}\  \}.
\label{tal}
\end{align}
It demonstrates directly that the cardinality
$l$ of the set $\la(\hw)$ equals $l(\hw).$
Cf. \cite{Hu},4.5.

This set can be introduced for non-reduced decompositions
as well. Let us denote it by $\tla(\hw)$ to differ from
$\la(\hw).$
It always contains $\la(\hw)$ and, moreover,
can be represented in the form
\begin{align}
&\tla(\hw)\ =\ \la(\hw)\, \cup\, \tla^+(\hw)\,
\cup\,-\tla^+(\hw),\label{tlaw}
\\ &\where \tla^+(\hw)\ =\
(\tR_+\cap\ \tla(\hw))\setminus \la(\hw).
\notag \end{align}
The coincidence  with $\la(\hw)$ is  for
reduced decompositions only.
\smallskip

Note that $\tla(\hw)$ depends on the choice of the
decomposition and it is actually a {\em sequence};
the roots in (\ref{tal}) are ordered naturally.
We will mainly
treat $\la(\hw)$ and $\tla(\hw)$ as {\em sequences\,}
in this paper,
for instance, when discussing the Bruhat ordering.
\smallskip

Let us consider the $\la$\~sets of the
reflections in $\tW$ and
check another standard property of the $\la$\~sets
(see, \cite{Bo} and 
\cite{Hu},4.6, Exchange Condition):
\begin{proposition}\label{HWSTAL}
For $\hw\in\hW$,
\begin{align}
&\la_\nu(\hw)\
=\ \{\tal>0, \ l_\nu( \hw s_{\tal}) \le l_\nu(\hw) \}.
\label{xlambda1}
\end{align}
\end{proposition}
{\em Proof.}
It suffices to consider  $\la(\hw).$ Also, we can the
inequality here
$l( \hw s_{\tal}) < l(\hw)$, since $l(s_{\tal})$
is odd (see below).
Thanks to (\ref{tal}), the set
$\{\tal>0\, \mid\, l( \hw s_{\tal}) < l(\hw) \}$ 
contains $\la(\hw).$
Obviously (\ref{xlambda1}) holds for $s_{\tal}=s_{\al_i}\ $
$(0\le i\le n).$

An arbitrary element $s_{\tal}$ for $\tal>0$
has a {\em reduced}
decomposition in the form:
$$
s_{\tal}=s_{i_1}s_{i_2}\cdots s_{i_p}s_m s_{i_p}\cdots
s_{i_2}s_{i_1},
\ i_{\bullet},m\ge 0.
$$
Indeed, $\al_i\in \la(s_{\tal})$ if and only if
$(\tal,\al_i)>0$ because
$\{\al_i\}$ is a minimal positive affine root. Given such $\al_i,$
the reflection $s_{\tal}$ is divisible by $s_i$ on the right
(i.e., has a reduced decomposition in the form $\cdots s_i$) and
$l(s_{\tbe})\le l(s_{\tal})$ for $\tbe\equal s_i(\tal).$
If $l(s_{\tbe})=l(s_{\tal}),$ then
$s_{\tbe}=s_is_{\tal}s_i=s_{\tbe}^{-1}$
is divisible by $s_i$ on the left or on the right,
which contradicts to
$(\tbe,\al_i)=-(\tal,\al_i)<0.$ Therefore
$l(s_{\tbe}) < l(s_{\tal})$
and we can proceed by induction.

We have also obtained that
the required decomposition can be started with an arbitrary
simple $\al_j\in \la(s_{\tal})$ taken as $\al_{i_1}$. Moreover,
the {\em sequence}
$\la(s_{\tal})$ can be expressed as follows
in terms of an (arbitrary) element $\tw\in\tW$ of 
minimal possible length such that
$\tw(\al)=\al_i$ for some $i\ge 0$:
\begin{align}\label{reflambda}
&\la(s_{\tal})=\{s_{\tal}(-\la(\tw))\}_{op}\,\cup\,\tal\,\cup\,
\la(\tw),
\end{align}
Given a sequence $\{\ \cdot\ \}$,
by $\{\ \cdot\ \}_{op}\ $, we mean the inversion of
its ordering.

If an arbitrary $\tw$ sending $\tal$ to a simple root is
taken here, then (\ref{reflambda}) holds for $\la(s_{\tal})$
considered as a {\em set} with possible cancelation of pairs
of opposite roots as in (\ref{tlaw}).
\smallskip

Now, we can prove (\ref{xlambda1}) by induction with respect
to the length $l(s_{\tal})$ of the roots $\tal$ such that
$l( \hw s_{\tal}) < l(\hw).$
As we already noticed, it holds for simple reflections.
Generally, 
let  $\al_i\in \la(s_{\tal}),\ \tbe\equal s_i(\tal)$,\ 
$s_{\tbe}=s_i s_{\tal} s_i,\ l(s_{\tbe})=l(s_{\tal})-2.$

If $\al_i\not\in \la(\hw),$ then
$$l(\hw s_{\tbe})=l((\hw s_i)(s_i s_{\tal} s_i))<
l(\hw)+1=l(\hw s_i),$$
and the induction statement gives that
$$\tbe\in \la(\hw s_i)\,\Rightarrow\, s_i(\tal)\in s_i(\la(\hw))
\,\Rightarrow\,\tal\in \la(\hw).
$$
Thus, we can assume that
$\al_i\in \la(\hw)$, equivalently, $\hw(\al_i)<0$, equivalently,
$l(\hw s_i)=l(\hw)-2.$ Since
$\al_i\in \la(s_{\tal})$ and $(\tal,\al_i)>0,$
we obtain:
$$ (\hw s_{\tal})(\al_i)=\hw(\al_i-2
\frac{(\tal,\al_i)}{(\tal,\tal)} \tal)<0
\hbox{\ provided\ } \hw(\tal)>0.
$$
If here $\hw(\tal)<0$, then $\tal\in \la(\hw)$,
which gives the desired.
Thus, it suffices to consider the case 
$(\hw s_{\tal})(\al_i)<0$, where 
$\al_i\in \la(\hw s_{\tal})$. 

In this case:
$$
l((\hw s_i)s_{\tbe})=l(\hw s_{\tal}s_i)\le
l(\hw s_{\tal})-1\le l(\hw)-2=l(\hw s_i).
$$
By induction, $\tbe\in \la(\hw s_i)$;
therefore $\tbe\in s_i(\la(\hw))$
and, finally, $\tal\in \la(\hw).$
\sq
\smallskip

The following is an immediate corollary of 
(\ref{reflambda}):
\begin{align} \label{talinla}
&\tal\in \la(\hw)\,\Leftrightarrow\, 
\la(s_{\tal}) =  \{\tbe,\,-s_{\tal}(\tbe)\,\mid\, 
s_{\tal}(\tbe)\in \tR_-\,,\,\tbe\in\la(\hw)\}.
\end{align}
\smallskip

In more detail, we have the following lemma.

\begin{lemma}\label{LEMREFLEC}
(i) Let $\tal=[\al,\nu_\al j]\in \tR_+$. 
Then $(\tal,\tbe)>0$
for each $\tbe\in \la(s_{\tal})$. Given 
$\tbe,\tga\in \la(s_{\tal})$, the sum 
$\tbe+\tga$ belongs to $\tR_+$ ($\Rightarrow$ 
$\tbe+\tga\in \la(s_{\tal})$)
if and only if $\nu_\al\ge \nu_\be=\nu_\ga$ and
$\tbe+\tga=[\al,\nu_\al j\,'\,]$.
In particular, $(\tbe,\tga)\ge 0$ unless the latter
condition holds.
 
(ii) Given a reduced decomposition of $s_{\tal}$,
its rank two Coxeter transformations of 
(consecutive) simple reflections 
are either $s_is_js_i\mapsto s_js_is_j$ with
the midpoints corresponding to $[\al,\nu_\al j\,'\,]
\in \la(s_{\tal})$ or $s_is_j\mapsto s_js_i$ 
otherwise. 

(iii) Let $\tbe$ be the first
root in  $\la(s_{\tal})$ (then it must be simple).
If $\tbe$
can be made a neighbor of $\tal$ in $\la(\tal)$ 
using the Coxeter transforms, then there exists
a reduced decomposition of $s_{\tal}$ such
that $\tbe$ is the first root in $\la(s_{\tal})$
and all roots $\tga$ between
$\tbe$ and $\tal$ are orthogonal to $\tbe$ unless they
satisfy 
$s_{\tbe}(\tga)=\tbe+\tga=[\al,\nu_\al j\,'\,],\,$
for some $j\,'<j$.
\end{lemma}
{\em Proof.}
The positivity  $(\tal,\tbe)>0$ is necessary 
(generally, not sufficient) for $\tbe $ to
belong to $\la(s_{\tal})$. If  $(\tbe,\tga)<0$,
then 
$s_{\tbe}(\tga)=\tga+\tbe\in \la(\hw)$ assuming that
$\nu_\be\ge \nu_\ga$. One has:
$$
s_{\tal}(\tbe+\tga)=
\tbe+\tga-\frac{(\tal,\tbe)}{\nu_\al}\tal-
\frac{(\tal,\tga)}{\nu_\al}\tal.
$$
However, the coefficient of $\tal$ must be $-(\tal,\tga)/\nu_\al$
unless the nonaffine components of $\tbe+\tga$ and 
$\tal$ coincide and $\nu_{\al}\ge\nu_\be=\nu_\ga$; 
we use that $s_{\tbe}(\tga)$ has the same length as $\tga$.

Claim (ii) follows from (i) because any sequence of roots
corresponding to a rank two Coxeter transformation contains
a pair of roots with the negative scalar product.
In (iii), we move $\tbe$
to the position next to $\tal$; all Coxeter transformations
we use must satisfy (ii). 
Then we can move $\tbe$ (back) to it is first position
(it is a simple root) using the Coxeter transforms in the
segment of $\la(s_{\tal})$ before $\tal$. 
\sq
\smallskip

The {\em sequence} $\la(\hw)=\{\tal^l,\ldots,\tal^1\}$, where
$l=l(\tw)$, and the element $\pi_r$
determine $\hw\in \hW$ uniquely:
\begin{align} \label{lambdainv}
&\al_{i_1}=\tal^1,\, \al_{i_2}=s^1(\tal^2), \ldots,
\al_{i_p}=s^1s^2\cdots s^{p-1}(\tal^p),\ldots\notag\\
&\al_{i_l}=s^1s^2\cdots s^{l-1}(\tal^l)\,,\where \notag\\
&s^p=s_{\tal^p},\ \hw=\pi_r s_{i_l}\cdots s_{i_1}=
\pi_rs^1\cdots s^l.
\end{align}
Notice the order of the reflections $s^p$ 
in the decomposition of $\hw$ is {\em inverse}.

A stronger fact holds.
{\em The $\la(\hw)$ considered as an unordered set
determines $\hw$ uniquely up to   
$\pi_r\in \Pi$ (on the left).} This statement
can be readily checked by induction with respect to $l=l(\hw)$.
Indeed, there exists at least one simple $\al_i\in \la(\hw)$
and {\em any} such $\al_i$ can be made $\al_{i_1}$; this means
that 
$l(\hw')=l-1$ for $\hw'=\hw s_i$.
Therefore the set $s_i(\la(\hw)\setminus \al_i)$ equals $\la(\hw')$
and determines $\hw'$ uniquely by the induction statement for $l-1$.

\medskip
\subsection{Reduction modulo 
\texorpdfstring{\mathversion{bold}$W$}{\em W}}
It generalizes the construction of the elements
$\pi_{b}$ for $b\in P_+.$ As a matter of fact,
the reduction modulo $W$ is a {\em formal} particular case
of a more general construction of the elements $\hw^\circ$ from 
Section \ref{sec:RightBruhat} below when $\tR^0=R$. This can be used
to obtain almost all claims of the following proposition, that is
from \cite{C4}.

\begin{proposition} \label{PIOM}
 Given $ b\in P$, there exists a unique decomposition
$b= \pi_b  u_b,$
$ u_b \in W$ satisfying one of the following equivalent conditions:

{(i) \  } $l(\pi_b)+l( u_b)\ =\ l(b)$ and
$l( u_b)$ is the greatest possible,

{(ii)\  }
$ \la(\pi_b)\cap R\ =\ \emptyset$.

The latter condition implies that
$l(\pi_b)+l(w)\ =\ l(\pi_b w)$
for any $w\in W.$ Besides, the relation $ u_b(b)
\equal b_-\in P_-=-P_+$
holds, which, in its turn,
determines $ u_b$ uniquely if one of the following equivalent
conditions is imposed:

{(iii) }
$l( u_b)$ is the smallest possible,

{(iv)\ }
if\, $\al\in \la( u_b)$  then $(\al,b)\neq 0$.
\end{proposition}
\qed

Note that the relation $l(\pi_b)+l(w)\ =\ l(\pi_b w)$
for any $w\in W$ is a special property of $\tR^0=R$
($\tR^0$ is from Section \ref{SUBSECTR0}).
Generally, $l(\hw^\circ\hu)\neq l(\hw^\circ)+l(\hu)$ for
$\hu\in \hW^0$, where
$\la(\hw^\circ)\cap \tR^0=\emptyset$.
\smallskip

Condition (ii) readily
gives a complete description of the set
$\pi_P=\{\pi_b, b\in P\}$, namely, 
only
$\,[\,\al<0,\,\nu_\al j>0\,]\,$
can appear in $\la(\pi_b)$ due to (\ref{laclosed}); see also
(\ref{lambpi}) below.
\smallskip

We note the following application of 
Theorem \ref{INTRINLA} below. A {\em sequence}
$$
\la\subset \tR_+[-]\equal\,\{[\al,\,\nu_\al j], \al\in R_-,j>0\}
$$
is in the form $\la(\pi_b)$  when and only when
the following three conditions hold:

(i) assuming that $\tal,\tbe,\tga=\tal+\tbe\in \tR_+[-]$,
if $\tal,\tbe\in \la$ then $\tga\in\la$ and $\tga$
appears between $\tal,\tbe$;
if $\tal\not\in\la$ then $\tbe$ belongs to $\la$
and appears in $\la$ before $\tga$;

(ii) if $\tal=[\al,\nu_\al j]\in \la$ then
$[\al,\nu_\al j\,'\,]\in \la$ as $j>j\,'>0$ and it appears in
$\la$ before $\tal$;

(iii) if $\tbe\in \la$ and $\tga=\tbe-[\al,\nu_\al j]\in 
\tR_+[-]$
for $\al\in R_+,\,j\ge 0$ then $\tga\in \la$ and it
appears before $\tbe$.

If $\la$ is treated as an unordered set, then
the existence of a representation
$\la=\la(\pi_b)$ for some $b$ is equivalent to 
(i+ii+iii) without the claims concerning
the ordering due to (\ref{lambpi})
\smallskip

Since $\pi_b = b u_b^{-1} = u_b^{-1} b_-, $
the set $\pi_P$ can be described explicitly in terms of $P_-$:
\begin{align}
&\pi_P = \{\, u^{-1} b_- \for b_-\in P_-\,,\ u\in W\notag\\
&\hbox{\ such\ that\ }
\al\in \la( u^{-1}) \Rightarrow (\al,b_-)\neq 0\, \}.
\label{lapiom}
\end{align}
Indeed, $\la(u^{-1}b_-)\,=\,
(-b_-)(\la(u^{-1}))\,\cup\, \la(b_-)$\,
and (\ref{lapiom}) is necessary and sufficient for
$\la(u^{-1}b_-)\cap R=\emptyset$. Note that using 
$\la(u_b^{-1})=-u_b(\la(u_b))$,
\begin{align}\label{lapiomu} 
&\al\in \la(u_b^{-1}) \Rightarrow\
(u_b^{-1}(\al),b)=
(\al,u_b(b))=(\al,b_-)\neq 0.
\end{align}
Actually, it suffices to check (\ref{lapiom}) for
simple roots $\al_i\in \la( u^{-1})$ only.
\smallskip

Using the longest element $w'_0$ in the centralizer
$W'_0$ of $b_-,$ the corresponding $u$ that can be taken
constitute the set
$$\{u\, \mid \, l(u^{-1}w'_0)\ =\ l(w'_0)+l(u^{-1})\}.$$
Their number is $|W|/|W'_0|.$
\smallskip

Note that Proposition \ref{PIOM},(ii) gives 
that any partial product
$\hu=s_{i_p}\cdots s_{i_1}$ for a reduced
decomposition $\pi_b=\pi_rs_{i_l}\cdots s_{i_1}$ 
belongs to $\pi_P$. Equivalently, if
$\hu$ satisfies $l(\pi_b)=l(\hu)+l(\pi_b\hu^{-1})$, 
then $\hu=\pi_c$ for $c=\hu\llb 0\rrb$ in the
notation from (\ref{afaction}).
\smallskip

For $\tal=[\al,\nu_\al j]\in \tR_+,$ one has:
\begin{align}
\la(b) = \{ \tal,\  &( b, \al^\vee )>j\ge 0 \iif \al\in R_+,
\label{xlambi}\\
&( b, \al^\vee )\ge j> 0 \iif \al\in R_-\},
\notag \\
\la(\pi_b) = \{ \tal,\ \al\in R_-,\
&( b_-, \al^\vee )>j> 0
\iif  u_b^{-1}(\al)\in R_+,
\label{lambpi} \\
&( b_-, \al^\vee )\ge j > 0 \iif
u_b^{-1}(\al)\in R_- \}, \notag \\
\la(\pi_b^{-1}) = \{ \tal\in \tR_+,\  -&(b,\al^\vee)>j\ge 0 \},
\label{lapimin}\\
\la(u_b)\ =\ \{ \al\in R_+,\ \ \ \,&(b,\al^\vee)> 0 \}.
\label{laumin}
\end{align}
For instance,
$l(b)=l(b_-)=-2(\rho^\vee,b_-)$ for
$2\rho^\vee=\sum_{\al>0}\al^\vee.$
\medskip
Note that nonempty
$\la(\pi_b^{-1})$ always contain nonaffine roots; thus
$\pi_b^{-1}$ can be represented in the form $\pi_c$ if and
only if $\pi_b=\pi_r$ for some $r\in O$: $\pi_P\cap\pi_P^{-1}=\Pi$.

Let us demonstrate how this formalism
works for calculating the set $\la(s_{\tbe})$ for
$\tbe=[-\be,l\nu_\be],$ where $\be\in R_+, l>0.$
This set can be determined directly from the 
definition of $\la(s_{\tbe})$, but it is useful to
establish a connection with $(-l\be)'=s_\be s_{\tbe}$.
One has:
\begin{align*}
&\la(s_\be)=s_\be(R_-)\cap R_+,\
s_{\tbe}(\la(s_\be))=\{[-\al,(\be,\al)l],\, \al\in \la(s_\be)\},\\
&\hbox{indeed,\ \ }
s_{\tbe}(\al)=\tal+(\al,\be^\vee)\tbe=s_\be(\al)+[0,(\al,\be)l]\\
&=-\al'+[0,(\al',\be)l] \for \al'=-s_\be(\al)\in -s_\be(\la(s_\be))=
\la(s_\be).
\end{align*}
Here $(\be,\al)$ must be greater than zero. Therefore
the above decomposition of $-l\be\in Q$ is {\em reduced} thanks to
the positivity condition ($c$) above, and we can apply
(\ref{xlambi}).

Note that the inequality $(\be,\al)>0$ is necessary for
$\al>0$ to belong to $\la(s_\be)$
but not sufficient; for $\be=\vth$ it is sufficient.
\smallskip

Using (\ref{xlambi}),
\begin{align*}
\la(-l\be)=\{\tal=[-\al,\nu_\al j]\, \mid\,
(\al^\vee,\be)l>j\ge 0 &\hbox{\ as\ }
\al\in R_-,\\
(\al^\vee,\be)l\ge j>0 &\hbox{\ as\ }
\al\in R_+\}.
\end{align*}
Finally,
\begin{align}\label{lamtbe}
\la(s_{\tbe})=&\{[-\al,\nu_\al j]\in \tR_+\ \mid\
(\al^\vee,\be)l>j\ge 0\}\ \cup\notag\\
&\{[-\al,\,(\be,\al)l]\ \mid\ \al>0<s_{\be}(\al),\, (\be,\al)>0 \}.
\end{align}

For example,
$\la(s_{\tbe})\cap R_+=\{\al'>0\, \mid\, (\al',\be)<0\}$; this set
is never empty unless $(\be,R_+)\ge 0$; the latter occurs, 
for instance, for $\be=\vth$.
\smallskip

Let us show explicitly how the transformation $-s_{\tbe}$ acts
in $\la(s_{\tbe})$:

(a) the subset
$\{[-\al,\,(\be,\al)l]\ \mid\ \al>0<s_{\be}(\al),\, (\be,\al)>0 \}$
maps exactly to
$\la(s_{\tbe})\cap R_+=\{[\al',0]\ \mid\
\al'>0,\,(\al',\be)<0\}$,

(b) the subset
$\{[-\al,\nu_\al j]\,\in\,\tR_+\, \mid\,
\ (\al^\vee,\be)l>j>0,\ \, \al>0<s_{\be}(\al)\}$ maps onto
$\ \ \{[-\al',\nu_{\al'} j]\in \tR_+\ \, \mid \ \,
((\al')^\vee,\be)l>j>0,\ \ \al'<0>s_{\be}(\al')\}$,

(c) the subset
$\{[-\al,\nu_\al j]\in \tR_+\, \mid\,
(\al^\vee,\be)l>j>0,\, \al>0>s_{\be}(\al)\}$ remains invariant
under the action of $-s_{\tbe}$.
\medskip

We will need later
the following \index{affine action $(wb)\llb z \rrb$}
{\dfont affine action} of $\hW$ on
$z \in \R^n$:
\begin{align}
& (wb)\llb z \rrb \ =\ w(b+z),\ w\in W, b\in P,\notag\\
& s_{\tal}\llb z\rrb\ =\ z - ((z,\al^\vee)+j)\al,
\ \tal=[\al,\nu_\al j]\in \tR.
\label{afaction}
\end{align}
For instance, $(b w)\llb 0\rrb=b$ for any $w\in W.$
The relation to the above action is given in terms of
the \index{affine pairing $([z,l], z'+d)$} {\dfont affine pairing}
$([z,l], z'+d)\equal (z,z')+l:$
\begin{align}
& (\hw([z,l]),\hw\llb z' \rrb+d) \ =\
([z,l], z'+d) \for \hw\in \hW,
\label{dform}
\end{align}
where we treat $(\,,\,+d)$ formally (one can add $d$
to $\R^{n+1}$ and extend $(\,,\,)$ correspondingly).

Introducing the
{\dfont basic affine Weyl chamber}
\begin{align}
&\CC^a\ =\ \bigcap_{i=0}^n \LL_{\al_i},\
\LL_{[\al,\nu_\al j]}=\{z\in \R^n,\
(z,\al)+j>0 \},
\notag \end{align}
we come to another interpretation of the $\la$\~sets:
\begin{align}
&\la_\nu(\hw)\ =\  \{\tal\in \tilde{R}_+, \,   \CC^a
\not\subset \hw\llb \LL_{\tal}\rrb, \, \nu_{\al}=\nu \}.
\label{lamaff}
\end{align}
Equivalently, taking a vector $\xi\in \CC^a$,
\begin{align}\label{lageomet}
\la(\hw)\,=\,\{\tal\in \tR\,\mid\,
(\tal^\vee,\xi+d)>0>(\tal^\vee,\xi'+d)\}
\end{align}
for $\xi'\in \hw^{-1}\llb \CC^a\rrb$.
Thus, we come to the following geometric description
of the $\la$\~sets (cf. Theorem \ref{INTRINLA},($a,b$) below).

\begin{proposition}\label{GEOMLA}
The $\la$\~sets are exactly those in the form
(\ref{lageomet}) for an arbitrary $\xi\in \CC^a$ and an 
arbitrary vector $\xi'$ inside
a certain (affine) Weyl chamber, i.e., provided
that $(\xi',\al^\vee)\not\in \nu_{\al}\Z$ for
all $\al\in R$. Given a generic segment in $\R^n$ 
from $\xi$ to $\xi'$, its consecutive
intersections with the affine root hyperplanes arrange such
set into a $\la$\~sequence.
\sq
\end{proposition}
\smallskip

Geometrically, $\Pi$ is the group of all elements of $\hW$
preserving $\CC^a$ with respect to the affine action.
Similarly, the elements $\pi_b^{-1}$ for $b\in P$ are exactly those
sending $\CC^a$ to
the basic nonaffine Weyl chamber
$\CC\equal\{z\in \R^n,$$ (z,\al_i)>0$ as $i>0\}.$
More generally,
given two finite sets of positive
affine roots $\{\tbe=[\be,\nu_\be i]\}$ and
$\{\tga=[\ga,\nu_\ga j]\},$
the closure of
the union of $\hw^{-1}\llb\CC^a\rrb$
over $\hw\in \hW$ such that
$\tbe\not\in \la(\hw)\ni \tga$ equals \\
\centerline{$
\{z\in \R^n,\  (z,\be)+i\ge 0,\ (z,\ga)+j\le 0
\ \hbox{\ for\ all\ } \tbe,\tga\}.$}
\smallskip

\rmk
Recall that for an arbitrary $z\in \R^n$, there
exists a unique $z_*=\tw \llb z \rrb\in \overline{\CC}^a$,
where $\tw\in \tW$, $\overline{\CC}^a$ is the closure
of $\CC^a$ in $\R^n$, and the centralizer of $z_*$ in
$\tW$ is generated by simple reflections.
For instance, $a_*=0$ for $a\in Q$, $b_*\in \{\om_r,r\in O\}$ for
$b\in P$.

The (classical) construction is as follows. We draw a line from
a small deformation $z^\ep$ of $z$ to an arbitrary
vector in $\CC^a$. It goes through
a chain of affine chambers 
and readily gives a reduced decomposition for
$\tw$ sending
$z^\ep$ to the required $z_*= \tw \llb z \rrb
\in \CC^a$.

If $\tu \llb z_* \rrb=z_*\in \overline{\CC}^a$, then  
$\tu \llb z_*^\ep\rrb$
remains in a small neighborhood of $z_*$ for a small
deformation $z_*^\ep\in \CC^a$ of $z_*$. Considering the line
from $z_*^\ep$ to $\tu \llb z_*^\ep\rrb$, we obtained
a decomposition of $\tu$, where all simple reflections
are in the hyperplanes sufficiently close to $z_*$. 
Thus all these hyperplanes are actually 
through $z_*$ and, therefore,
contain the {\em face} through $z_*$:
$
\{ z\,\mid\, (\tal^\vee,z_*+d)=0\Rightarrow
(\tal^\vee,z+d)=0 \for \tal\in \tR\}.
$\qed
\smallskip

\begin{lemma}\label{NOPIRB}
For an arbitrary $b\in P$, if $\pi_r\llb b\rrb=b$ for
$\pi_r\in \Pi$, then $r=0$, that is, $\pi_r=$id.
\end{lemma}
{\em Proof.} If $\pi_r\llb b\rrb=b$ then
$((-b)\pi_r b)\llb 0 \rrb =0$. The element $(-b)\pi_r b$
for $r\neq 0$
does not belong to $\tW$, which is normal in $\hW$.
However, 
$\hw\llb 0\rrb \Rightarrow \hw\in W$.

It is instructional to check this claim explicitly.
Introducing the map $i\mapsto \widehat{i}$
by $\al_{\,\widehat{i}\,}=\pi_r(\al_{i})$,
$\pi_r\llb b\rrb=b$ for $r\neq 0$ is equivalent to
\begin{align*}
&b=\sum_{j=1}^{n}n_j\om_j, \ n_{\,\widehat{i}\,}=n_i\in\Z
\for n_0\equal 1-(\vth,b), \ 0\le i\le n.
\end{align*}
Indeed, $\pi_r\llb b\rrb=b\, \Leftrightarrow \,
(\al_i^\vee,b+d)=(\al_{\,\widehat{i}\,}^\vee,b+d)$,
 that is
\begin{align*}
&(\al_{r'}^\vee,b)=1-(\vth,b)=(\al_r^\vee,b),
\where \widehat{r'}=0,\\
&(\al_j^\vee,b)=(\al_{\widehat{j}}^\vee,b) \for 0\neq j \neq r'.
\end{align*}
Setting $\vth=\sum_{j}\kappa_j\al_j^\vee$ and using that
$\kappa_r=1=\kappa_{r'\,}$,
\begin{align}\label{pirfixed}
b&=\sum_{j=1}^{n} n_j\om_j\ =\ \pi_r\llb b\rrb\in P
\Leftrightarrow\\
\{\ n_{\,\widehat{j}\,}&=n_j\in\Z \for 0\neq j\neq r'\,, \
n_{r}=n_{r'\,},\ n_r+\sum_{j=1}^n \kappa_j n_j=1\ \}.\notag
\end{align}

The sum $ n_r+\sum_j\kappa_j n_j$ is divisible by the order
of $\pi_r$ in the case of $A_n$, so there are no such $b$ for
this root system. Similarly,
they do not exist in the case of $B_{2l+1}\, (l>1)$
and for $\pi_r$ of order $4$ ($D_{2l+1},\, (l>2)$); indeed,
$\pi_r$ has no fixed nonaffine $\al_j^\vee$ in these cases
and the sum is divisible by $2$ .
In the remaining cases, $\pi_r\neq$id have fixed
nonaffine $\al_j^\vee$; the corresponding labels
$\kappa_j$ are $2$ for $B,C,D$ and  $3=|\Pi|$ for $E_6$,
therefore, the sum is
divisible by $2$ or $3$ and (\ref{pirfixed})
is impossible.
\sq

\smallskip
The element $b_{-}= u_b(b)$ is a unique element
from $P_{-}$ that belongs to the orbit $W(b)$.
Thus the equality   $c_-=b_- $ means that $b,c$
belong to the same orbit. We will also use
$b_{+} \equal w_0(b_-),$ a unique element in $W(b)\cap P_{+}.$
In terms of the elements $\pi_b,$
$$u_b\pi_b\ =\ b_-,\ \pi_b^{-1} u_b^{-1}\ =\ b^\varsigma_+,\ 
b^\varsigma=-w_0(b).$$

Note that $l(\pi_b w)=l(\pi_b)+l(w)$ for all $b\in P,\ w\in W.$
For instance,
\begin{align}
&l(b_- w)=l(b_-)+l(w),\ l(wb_+)=l(b_+)+l(w),
\label{lupiw}
\\
&l(u_b\pi_b w)=l(u_b)+l(\pi_b)+l(w) \for b\in P,\,
 w\in W.\notag
\end{align}
\smallskip

\subsection{Partial ordering on 
\texorpdfstring{\mathversion{bold}$P$}{\em P}}
It is necessary in the theory of nonsymmetric polynomials.
See \cite{O2,M4}. This ordering was also used in \cite{C2} in the
process of calculating the coefficients of $Y$\~operators.
The definition is as follows:
\begin{align}
&b \le c, c\ge b \for b, c\in P \iif c-b \in Q_+,
\label{order}
\\ &b \preceq c, c\succeq b \iif b_-< c_- \hbox{\ or\  }
\{b_-=c_- \hbox{\ and\ } b\le c\}.
\label{succ}
\end{align}
Recall that $b_-=c_- $ means that $b,c$
belong to the same $W$\~orbit.
We  write  $<,>,\prec, \succ$ respectively if $b \neq c$.

The following sets
\begin{align}
&\si(b)\equal \{c\in P, c\succeq b\},\
\si_*(b)\equal \{c\in P, c\succ b\}, \notag\\
&\si_-(b)\equal \si(b_-),\
\si_+(b)\equal \si_*(b_+)= \{c\in P, c_->b_-\}.
\label{cones}
\end{align}
are convex.
By {\dfont convex}, we mean that if
$ c, d= c+r\al\in \si$
for $\al\in R_+, r\in \Z_+$, then
\begin{align}
&\{c,\ c+\al,...,c+(r-1)\al,\ d\}\subset \si.
\label{convex}
\end{align}

The convexity of the intersections
$\si(b)\cap W(b), \si_*(b)\cap W(b)$
is by construction. For the sake of completeness,
let us check the convexity
of the sets $\si_{\pm}(b).$

Both sets are $W$\~invariant. Indeed,
$c_- > b_-$ if and only if  $b_+>w(c)>b_-$ for all
$w\in W.$
The set $\si_-(b)$ is the union of  $\si_+$ and the orbit $W(b).$
Here we use that $b_+\and b_-$ are the greatest and the least
elements of  $W(b)$ with respect to "$>$".
This is known (and can
be readily checked by the induction with respect to the
length - see e.g., \cite{C2}).

If the endpoints of (\ref{convex}) are between $b_+$ and $b_-$ ,
then
it is true for the orbits of all inner points even if  $w\in W$
changes the sign of $\al$ (and the order of the endpoints).
Also the elements from $\si(b)$
strictly between $c$ and $d$ (i.e.,
$c+q\al,\ $ $0<q<r$) belong to $\si_+(b).$ This gives the required.

The next propositions are essentially from \cite{C1}.
We will use the standard Bruhat ordering.
Given $\hw\in \hW$, the {\em standard Bruhat set}
$\b(\hw)$
is formed by $\hu$ obtained by striking out any number of
$\{s_{j}\}$ from a reduced decomposition of $\hw\in \hW$.
The set $\b(\hw)$ does not depend on the choice of the
reduced decompositions.

\begin{proposition }\label{BSTAL}
(i) Let $c=\hu\llb 0\rrb$,\, $b=\hw\llb 0\rrb$ and
$\hu\in \b(\hw)$. The latter means that
$\hu$ can be obtained by deleting
simple reflections from any reduced decomposition of $\hw$,
say, from the product of the reduced decompositions of $\pi_b$
and $w$ in the decomposition $\hw=\pi_b w$.
Then $c\succeq  b$ and  $b-c$ is a linear combination of the
non-affine components of the
corresponding roots from $\la(\hw^{-1})$; also,
$c=b$ if and only if
$s_{j}$ are deleted only from the reduced decomposition of $w$.

(ii) Letting  $b=s_i\llb c\rrb$
for $0\le i\le n,$ if the element $s_i\pi_c$ belongs
to $\{\pi_a,\,a\in P\}$ then it equals $\pi_b$. It happens
if and only if
$(\al_i,c+d)\neq 0.$ More precisely,
the following three conditions are equivalent:
\begin{align}
& \{c\succ b\}\Leftrightarrow  \{(\al_i^\vee,c+d)>0\}
\Leftrightarrow \{s_i\pi_c=\pi_b,\ l(\pi_b)=l(\pi_c)+1\}.
\label{aljb}
\end{align}
If the latter holds, then
$\la(\pi_b)=\pi_c^{-1}(\al_i)\cap\la(\pi_c)$, where
$$\pi_c^{-1}(\al_i)=u_c(\al_i)+[0,(c,\al_i)].$$
\qed
\end{proposition }
The following lemma from \cite{C4} completes (ii); it
describes the
case $(\al_i^\vee,c+d)=0.$

\begin{lemma}\label{VANISHL}
The condition $(\al_i,c+d)=0$ for $0\le i\le n$
equivalently, the condition $(\al_i,b+d)=0$, implies that
$u_c(\al_i)=\al_j$ as $i>0$ or
$u_c(-\vth)=\al_j$ as $i=0$ for a proper index $j>0.$
Given $c$ and $i$, the existence of such $\al_j$ and the equality 
$(\al_j,c_-)=0$ are equivalent to $(\al_i,c+d)=0$.
\end{lemma}
{\em Proof.} If $i>0$ then $\al\equal u_c(\al_i)>0$
and $(\al,c_-)=0.$
If $\al=\be+\ga$ for positive roots
$\be,\ga,$ then $(\be,c_-)=0=(\ga,c_-).$
Hence $\be'=u_c^{-1}(\be)> 0,$
$\ga'=u_c^{-1}(\ga)>0$
and therefore $\al_i=\be+\ga,$ which is impossible.
Thus $\al$ must be simple.

Similarly,  $(\vth,c)=1$ implies that
$\al\equal u_c(-\vth)>0, (\al,b_-)=-1$. Let
$\al=\be+\ga$, where
$\be>0<\ga.$ Since $\vth$ and $\al$ are short,
we can assume that at least one of them is short,
but it will follow automatically.
Then, transposing $\be\leftrightarrow \ga$
if necessary we obtain that
$(\be,c_-)=-1,\ (\ga,c_-)=0.$
Since $(\be,c_-)=\nu_\be(\be^\vee,c_-)$ we conclude that
$\nu_{\be}=1$ and $\be$ is short.
Hence $\be'=-u_c^{-1}(\be)>0$,
$ \ga'=-u_c^{-1}(\ga)<0$
and therefore $\vth=\be\,'\,+\ga\,'\,< \be'$. The latter results in
$\vth^\vee<(\be')^\vee,$ which is
impossible
since $\vth^\vee=\vth$ is the maximal positive root
in $R_+^\vee.$
\sq
\smallskip

Combining the lemma with (\ref{lambpi}), we come to the
following corollary:
\begin{corollary}\label{LAFORDOMINANT}
Let us take $c\in P$ such that $\la(\pi_c)$ contains 
$[-\al,\nu_\al j]$ for each $\al\in R_+$
where $j>0$ (then it holds for $j=1$).
It results in $-(c_-,\al_{i\,'\,})>0$ for all $n\ge i\,'\,>0$ and,
given $\hw\in \hW$,
the condition $l(\hw)+l(\pi_c)=l(\hw\pi_c)$
implies that $\hw\pi_c=\pi_b$ for some $b\in P$, 
having the same property as $\,c\,$ 
(due to $\la(\pi_c)\subset\la(\pi_b)$).
\end{corollary}
{\em Proof.}
Here the condition $-(c_-,\al_{i\,'\,})>0$ for all $i\,'\,>0$
(we can write $c_-\in P_{--}$) 
guarantees that $s_i\pi_c $ is in the form $\pi_b$
for any $i\ge 0$. 
This condition follows from the assumption
for $c$ (but is somewhat weaker). 
Provided $l(s_i\pi_c)=l(\pi_c)+1$,\ we obtain that 
$\la(s_i\pi_c)$=$\la(\pi_b)$ $\supset \la(\pi_c)$  
and can continue by induction with respect to $l(\pi_c)$. 

See also the next proposition,
that gives an explicit description
of the changes of the $\la$\~sets $\la(\pi_c)$
and $c_-$ 
upon the multiplication by simple reflections.
\sq
\smallskip

\begin{proposition}\label{BSTIL}
(i) Assuming (\ref{aljb}), let $i>0$. Then $b=s_i(c),\ $
$b_-=c_-$\, and
$u_b=u_c s_i.$ The set $\la(\pi_b)$ is obtained from
$\la(\pi_c)$ by adding $[\al,(c_-,\al)]$
for $\al=u_c(\al_i)$, i.e., by
replacing the inequality
$( c_-, \al^\vee )>j> 0$
in (\ref{lambpi}) with
$( c_-, \al^\vee )\ge j> 0.$
Here $\al\in R_-$ and $(c_-,\al^\vee)=(c,\al_i^\vee)>0.$

(ii) In the case $i=0,\ $ the following holds:
$b=\vth+s_{\vth}(c)$, the element $c$ is from $\si_+(b)$,
$b_-=c_- - u_c(\vth)\in P_-$ and $u_b=u_c s_{\vth}$.
For $\al= u_c(-\vth)=\al^\vee,$ the $\la$-inequality
$( c_-, \al)\ge j> 0$ is replaced
with $( c_-, \al)+2> j> 0$; respectively, the root
$[\al,(c_-,\al)+1]$ is added to $\la(\pi_c)$.
Here $\al=\al^\vee\in R_-$ and
$( c_-, \al)=-(c,\vth)\ge 0.$

(iii) For any $c\in P, r\in O',\ $ one has
$\pi_r\pi_c=\pi_b$ where $b=\pi_r\llb c\rrb.$
Respectively, $u_b=u_c u_r,\ b= \om_r+u_r^{-1}(c),\
b_-=c_- +u_c w_0(\om_{r}).$ In particular,
the latter weight always belongs to $P_-.$
\end{proposition}

\comment{
\smallskip
\subsection{Arrows in {\em P}}
We write $c\rightarrow b$ or $b\leftarrow c$ in the cases (i),(ii)
from the proposition above and
use the left-right arrow  $c\leftrightarrow b$ for (iii)
or when
$b$ and $c$ coincide.

By $c\rightarrow\!\rightarrow b,$  we mean that
$b$ can be obtained from $c$ by
a chain of (simple) right arrows. Respectively,
 $c\leftrightarrow\!\rightarrow b$ indicates that
$\leftrightarrow$
can be used in the chain.

Actually no greater than one
left-right arrow
is always sufficient and it can be placed right after $c.$
If such arrows are not involved , then this
ordering is obviously stronger than the standard Bruhat ordering
given by the procedure from Proposition \ref{BSTAL},(i),
that in its turn is stronger than "$\succ$".

If $l(\pi_b)=l(\pi_c)+l(\pi_b\pi_c^{-1})$ , then the reduced
decomposition of $\hw=\pi_b\pi_c^{-1}$ readily
produces a chain
of simple arrows from $c$ to $b.$
The number of right arrows is precisely
$l(\hw)$ (see(\ref{aljb})).
Only transforms of type (ii),(iii) will change the
$W$-orbits, adding negative short roots to the corresponding
$c_-$ for (ii) and
the weights in the form $w(\om_r)\ (w\in W)$ for (iii).

For instance,
let $c,b,b-c\in P_-.$ Then $\pi_c=c,\ \pi_b=b,\ $
$\hw=b-c$ is
of length $l(b)-l(c),$ and we  need to decompose $b-c.$
If $c=0,\ $ then $b-\om_r\in Q,$ the reduced decomposition of $b$
begins with $\pi_r\ (r\in O'),$ and
the first new $c_-$ is $w_0(\om_r).$
When there is no
$\pi_r$ and $b\in Q\cap P_-,$ the chain always
starts with  $-\vth.$

Let us now examine the arrows
$c\rightarrow b$ from the
viewpoint of the $W$-orbits.
\begin{proposition}\label{WORB}
(i) Given $c_-\in P_-$, any element in the form
$c_-+u(\om_r)\in P_- $ can be represented as $b_-$ for proper
$b\leftrightarrow c$. Respectively,
$c_-+u(\vth)\in P_-$
for $u\in W$ such that $u(\vth)\in R_-$
can be represented as $b_-$
for proper $b$ such that $W(c_-)\ni c\rightarrow b.$

(ii) In the case of $A,D,E,$
any element $b_-$ such that
$b_-\prec c_-$ (both belong to  $P_-$) can be obtained from
$c_-\in P_-$ using consecutive arrows
$W(c_-)\ni c\rightarrow b\in W(b_-).$ This cannot be true
for all root systems because only short roots
may be added to
$c_-$ using such a construction. \sq
\end{proposition}
}

\medskip
\section{Reduced decompositions}
We will discuss properties of the reduced decompositions
in connection with the corresponding $\la$\~sets and
$\la$\~sequences.
\smallskip

\subsection{Lambda--sequences}
Let us give an
intrinsic description of the {\em sequences\,}
$\la(\hw)$. See (\ref{lageomet}) for a
more geometric approach based on the
affine Weyl chambers.
\smallskip
\begin{maintheorem}\label{INTRINLA}
(i) A sequence $\la=\{\tal^l,\ldots \tal^1\}$ of pairwise distinct
roots from $\tR_+$ can be represented in the form
$\la(\hw)$ for a certain $\hw$
if and only if the following six conditions are satisfied:
\begin{align*}
&(a)\ \ \{\,\tal=\tal^q+\tal^r\not\in [0,\Z],\
\tal\in \tR_+\,\}\,\Rightarrow\, \tal\in \la\,;\\
&(b)\ \ \{\,\tal^p=\tal^q+\tal^r\,\}\, \Rightarrow\, 
\{\,p \hbox{\ is\ between\ }
q,r \hbox{\ in\ } \la\,\};\\
&(c)\ \ \{\,\tal^p=\tbe+\tga \for \tbe,\tga\in \tR_+\,\}\, 
\Rightarrow\,
\{\,\tbe\in \la \hbox{\ or\ } \tga\in \la\,\};\\
&(d)\ \ \{\tal^p=\tal^q+\tga \hbox{\ for\ }
\tR_+\ni\tga\not\in\la\}\,  \Rightarrow\, q<p\,;\\
&(e)\ \ \{\tal^p=[\al,\nu_\al j]\in \la,\,
[\al,\nu_\al j\,']\in \tR_+,\, 0\le j\,'\,<j\}\,
\Rightarrow\, [\al,\nu_\al j\,'\,]\in \la\,;\\
&(f)\ \ \{\tal^p=[\al,\nu_\al j],\ \tal^q=[\al,\nu_\al j\,'],
\ 0\le j\,'\,<j\}\,
\Rightarrow\, q<p\,.
\end{align*}

(ii) Assuming that ($a$),($c$),($e$) hold for a set $\la$, called
a {\dfont $\la$\~set},
there exists at least one its ordering, called a
{\dfont $\la$\~sequence},
satisfying conditions ($b$),($d$),($f$).
All such sequences are
in one-to-one correspondence with reduced decompositions of
$\tw\in \tW$ such that $\la(\tw)=\la$.

(iii) Conditions ($a,b,c,d$) imply conditions
\begin{align*}
&(a')\ \ \{\,\tal\,'\,=u\tal^q+v\tal^r\in \tR_+\,\}\, \Rightarrow\,
\tal\,'\,\in \la;\\
&(b')\ \ \{\,\tal^p=u\tal^q+v\tal^r\,\}\, \Rightarrow\, 
p \hbox{\ is\ between\ } q,r;\\
&(c')\ \ \{\,\tal^p=u\tbe+v\tga \for \tbe,\tga\in \tR_+\,\}\, 
\Rightarrow\,
\{\,\tbe\in \la \hbox{\ or\ } \tga\in \la\,\};\\
&(d')\ \ \{\,\tal^p=u\tal^q+v\tga \hbox{\ for\ }
\tR_+\ni\tga\not\in\la\,\}\,  \Rightarrow\, q<p,
\end{align*}
where $u,v$ are positive rational numbers. These conditions
coincide with ($a,b,c,d$) in the simply-laced case.
\end{maintheorem}
{\em Proof.}
The induction in $l$ will be used.
Claims (i,ii) are obvious as $l=1$.
We will establish claims (i,ii) for $l$
assuming that (i) holds for all $1\le l\,'\,<l$, 
that is the equivalence
($a,b,c,d,e,f$)$\,\Leftrightarrow\, \{\la=\la(\hw)\}$,
and that (ii) holds for such $l\,'\,$, 
that is the existence and description
of the $\la$\~sequences for any given $\la$\~set,

Given $\hw'$, the product $\hw's_i$ is reduced if and only if
$\al_i\not\in\la(\hw')=\{\tal^{l-1},\cdots,\tal^1\}$. Then
$\la(\hw's_i)=\{s_i(\la(\hw')),\al_i\}.$
Note that the last set is automatically
positive (belongs to $\tR_+$)
since $s_i(\tR_+)\cap\tR_-=-\al_i$. Indeed,
the decomposition of any $\tR_+\ni\tal\neq \al^i$
in terms of simple roots in $R_+$
and {\dfont imaginary roots\,}
$[0,\Z]$ contains either a simple root $\al_j\neq \al_i$ with
a positive coefficient or an imaginary root $[0,m]$ with $m>0$.

Let $\la(\hw)$ be the $\la$\~sequences for a reduced
decomposition $\hw=\pi_r s_{i_l}\cdots s_{i_1}$. Then
the set $\{s_{i_1}(\tal^l),\cdots, s_{i_1}(\tal^2)\}$
is $\la(\hw')$ for $\hw'=\hw s_{i_1}$
of length $l-1$, and we readily deduce ($a-f$) from the induction
assumption. The only if statement for $l$ is verified.
\smallskip

Let us check that $\la=\{\tal^p\}$ of length
$l$  satisfying ($a,b,c,d,e,f$) corresponds to a certain $\hw$.
First of all,
$\tal^1\in \la$ must be a simple root. Otherwise,
one can represent $\tal^p=\al_i+\tbe$ for a certain simple $\al_i$
where $\tbe$ is either from $\tR_+$ or is an
imaginary root $[0,\nu j]$ for some $j>0,\nu=\nu_{\tal^p}$.
However, such $\al_i$ must appear
in the sequence $\la$ before $\tal^1$, which is impossible.

If only ($a,c,e$) are imposed for $l'<l$,
then the same argument gives that $\la$ considered as
an unordered {\em set} contains at least
one simple root.
\smallskip

Using the notations $s^1=s_{\tal^1}$, $\tal'=s^1(\tal)$,
the roots $(\tal^p)'=s^1(\tal^p)$ are
all positive; see above.
Let us establish that
the sequence $\la\,'\,=s^1(\{\tal^l,\ldots, \tal^2\})$
satisfies ($a,b,c,d,e,f$); respectively, it satisfies 
($a,c,f$) if
$\la$ and $\la\,'\,$ are considered as unordered {\em sets\,}.

The claims ($a,b$) are obvious. Applying $s^1$ to any
$[\be,\nu_\be j]$
diminishes $j$ by a non-negative integer (by zero, if $s^1$
is nonaffine),
that gives ($e,f$).

As for ($c,d$), let $(\tal^p)'=\tbe\,'\,+\tga\,'\,$ 
for positive $\tbe\,'\,$ and $\tga\,'\,$; it suffices to 
consider only $\tga'=\tal^1$.
Then
$$
\tal^p=\tbe-\tal^1 \Rightarrow \tbe=\tal^p+\tal^1 \Rightarrow
\{\tbe=\tal^q, q<p\} \Rightarrow \tbe\,'\,\in \la\,'\,,
$$
where $\tbe=s^1(\tbe\,'\,)$. It gives ($c,d$) if $\tga\,'\,=\tal^1$.
Otherwise, both, $\tbe$ and $\tga$, are
positive, and ($c,d$) for $\la'$ follow from those for
$\la$.

Now we can use the induction assumption for $\la\,'\,$; it
gives that $\la\,'\,=\la(\hw')$ and we 
can go from $\hw'$ to $\hw=\hw's^1$ as above.
\smallskip

If only conditions ($a,c,e$) are imposed, then given {\em any}
ordering of the {\em set} $\la\,'\,$ of type ($b,d,f$),
the {\em sequence}
$\{s^1(\la\,'\,),\al^1\}$ satisfies ($b,d,f$) too. This 
is sufficient for establishing (ii).

The claims from (iii) are checked using the same induction
consideration.
\sq
\smallskip

Note that if $\tga=\sum_{m}\tal^{q_m}\in \tR_+$ where
$\tal^{q_m}\in \la(\hw)$, then, obviously, $\tga\in\la(\hw)$.
It formally results from Theorem
\ref{INTRINLA},(i) due to the known fact \cite{Bo}
that there exists a permutation of the indices $m$, 
$\tga=\tal^{q_1}+\tal^{q_2}+\ldots$, such that all
partial sums $\tal^{q_1}, \tal^{q_1}+\tal^{q_2}, \ldots$
are from $\la$.
\medskip

\subsection{Coxeter transformations}
We will prepare tools for studying transformations of
the reduced decompositions. The elementary ones are the
{\dfont Coxeter transformation} that are substitutions
$(\cdots s_is_js_i)\mapsto (\cdots s_js_is_j)$ in
reduced decompositions of $\hw\in \hW$ with the number
of factors $2,3,4,6$ as $\al_i$ and $\al_j$ are connected
by $m_{ij}=0,1,2,3$ laces in the affine Dynkin diagram.
They induce
{\em reversing the order} in the corresponding segments
(with $2,3,4,6$ roots)
of the sequence $\la(\hw)$. The corresponding roots form
a set identified with the set of positive roots
of type $A_1\times A_1$,$A_2$,$B_2$,$G_2$ respectively.
The action of Coxeter transformations in $\la$\~sets 
plays an important role in the paper.
\smallskip

\begin{proposition}\label{INTRINLAP}
(i) Given $\hw\in \hW$, the roots $\tal$ that
may appear in the beginning, the \underline{first roots}, of the
sequence
$\la(\hw)=\{\ldots,\tal\}$ for different reduced decompositions
of $\hw$ are exactly simple roots $\tal=\al_i\in\la(\hw)$.
The \underline{ last roots\,}, i.e.,
$\tal$ such that  $\la(\hw)=\{\tal,\ldots\}$ for a suitable
reduced decomposition
are exactly $\tal=-\hw^{-1}(\al_i)$ for
$\al_i\in\la(\hw^{-1})=-\hw(\la(\hw))$.

(ii) The \underline{last roots}
of $\la(\hw)$ are also exactly $\tal\in \la(\hw)$  satisfying
the following two conditions:

(a) $\tal \neq \tbe+\tga \hbox{\ for\ any\ two\ roots\ }
\tbe,\tga\in \la(\hw)$,

(b) $\tga \neq \tal+\tbe \hbox{\ for\ any\ } \tga\in \la(\hw),\
\tR_+\cup[0,\N]\ni \tbe\not\in \la(\hw)$.

(iii) Given a reduced decomposition of $\hw$,
if $\tga=u\tal+v\tbe \for \tal,\tbe,\tga\in \la(\hw),\ \Q\ni u,v>0$,
and $\tal,\tga$ are neighboring in
$\la(\hw)$, then proper Coxeter transformations in the
segment from $\tbe$ to $\tga$ ($\tga$ is excluded)
make $\tbe$ next to $\tga$, i.e., make the triple $\tbe,\tga,\tal$
a connected segment in the resulting $\la$\~sequence.
\end{proposition}
{\em Proof.} Claim (i) is simple and well known.
The demonstration
of (ii) is by induction; we assume that $\hw=\hw's_i$ as
$l(\hw)=l(\hw')+1$ and apply (ii) for $\hw'$. It is as 
follows.

Since $s_i(\tal\,'\,)>0$ for all 
$\tal\,'\,\in \tR_+\setminus \{\al_i\}$,
a root $\tal\in \la(\hw)$ satisfying ($a,b$) for $\hw$
equals
$s_i(\tal\,'\,)$ for $\tal\,'\,\in \la(\hw')$ satisfying ($a,b$) 
for $\hw'$ unless $\tal=\al_i$ ($\tbe\neq\al_i$ because
$\al_i\in \la(\hw$)). Thus $\tal\neq \al_i$ is the {\em last} 
for a certain reduced decomposition of $\hw'$ multiplied by
$s_i$ on the right.

We need only to check  that $\al_i$ can be made the {\em last} in
a sequence $\la(\hw)$ if
\begin{align}\label{tganeq}
\tga\neq \al_i+ \tbe \hbox{\ \,for\ any\ \,}
\tga\in \la(\hw)\not\ni \tbe\in \tR_+\cup [0,\N].
\end{align}

By induction, 
the last two roots of $\la(\hw)$ can be made
$\tde,\al_i,\ldots,$ 
($\la(\hw)=\{\tal,\al_i,\ldots\}$) for
a proper reduced decomposition of $\hw$. Here $\tde$ is the end
of $\la(\hw)$; it satisfies $s_i(\tde)>0$ since
$\tde=s_i(\tde\,'\,)$ for $\tde\,'\equal s_i(\tde)
\in \la(\hw')\subset \tR_+$.

Setting $\tde\,'\,=s_i(\tde)=\tde+m\al_i$, either
$m=0$ and we can simply transpose $\tde$ and $\al_i$ in $\la(\hw)$,
or $m>0$, or $m<0$. In the case $m>0$,
$\tde\,'\,\in \la(\hw)$ and it must be between $\tde$ and $\al_i$,
which is impossible (they are neighbors).
If $m<0$ and $\tde\,'\,\in \la(\hw)$,
then $\tde$ must be  between $\tde\,'\,$ and $\al_i$,
which is impossible too. Thus $m=-p<0$ and
$\tde\,'\,\not\in \la(\hw)$.

Finally, $\al_i+\tbe=\tde$ where the
root $\tbe\equal\tde\,'\,+(p-1)\al_i$ belongs to $\tR_+\cup [0,\N]$
(see \cite{Bo}).
If here $\tbe\in \la(\hw)$, then $\tde$
must be between $\tbe$ and $\al_i$, which is impossible because 
$\tde$ was the end of the sequence $\la(\hw)$. If 
$\tbe\not\in \la(\hw)$, then
it contradicts assumption (\ref{tganeq}); (ii) is checked.
\smallskip

Now let
$u\tal+v\tbe=\tga$ for the roots in $\la(\hw)$ as $u,v>0$.
We are going to make them consecutive
in a proper reduced decomposition. Note that $u,v>0$
is necessary and sufficient to make $\tga$ between $\tal$ and
$\tbe$.

One can suppose that $\tal=\al_i$
is the {\em first} in $\la(\hw)$.
Then $\tga$ is the {\em second} in this sequence
and the reduced decomposition reads as $\hw=\cdots s_j s_i$
for  $j$ such that $\tga=s_i(\al_j)=m\al_i+\al_j$ for $m\in \N$
and $\al_i,\al_j$ are connected in the affine Dynkin diagram.
Continuing, we can find a reduced decomposition of $\hw$
in the form $\hw=\hu \tv$, $\tv=\cdots s_i s_j s_i$,
$l(\hu \tv)=l(\hu)+l(\tv)$ where $\tv$ is the
longest possible product of $s_i$ and $s_j$, 
equivalently, $\al_i,\al_j\not\in \la(\tu)$.

Let as assume that $\tbe\in \tv^{-1}(\la(\tu))$.
Since $\tv(\tbe)$ is a (positive) linear combination of 
$\al_i,\al_j$, then either $\al_i$ or $\al_j$ must belong to
$\la(\tu)$, a contradiction. Hence $\tbe\in \la(\tv)$,
which proves (iii).
\sq
\smallskip

\comment{
(iv) Given a reduced decomposition of $\hw$,
if $\tga=\tal+\tbe$ for arbitrary
$\{\ldots\tbe,\ldots,\tga\,\ldots,\tal,\ldots\} \subset \la(\hw)$,
then proper Coxeter transformations inside the segment from
$\tal$ to $\tbe$ make these roots connected, i.e., result in
$\{\ldots\tbe,\tga\,\tal,\ldots\}$ provided
that these roots form the set of positive roots of
type $A_2$.

Generally, we do not suppose now that
$\tga$ and $\tal=\al_i$ are
neighbors in $\la(\hw)$. If there exists at least one simple
$\al_j\neq \al_i$ in $\la(\hw)$, then we can make $\al_j$ the
first in $\la(\hw)$ and, therefore, reduce the distance between
$\tal$ and $\tbe$ and proceed by induction unless
$\al_j=\tbe$. In the latter case, the roots $\tbe,\tga,\tal$
will be transformed to $\tbe,\tga,\tal$ by a chain of Coxeter
transformation. It is possible only if two of them
become neighbors in this process and we can use part (iii).

Similarly, we can assume that $\tbe$ a unique {\em end} 
of $\la(\hw)$ for all possible reduced decompositions of
$\hw$, i.e., is described by (i,ii). We will list all such
$\hw$ applying the following lemma, that will be also
used later in one of the main theorems.
\smallskip

\begin{lemma}\label{ONEENDS}
Given $\al_i$ for $0\le i\le$, there exists a unique element
$\tu^{i}\in \tW$ such that
\begin{align}\label{latuione}
\la(\tu^{i})=\{\tal\in \tR_+^0\,\mid\, \tal-\al_i\in
\widetilde{Q}_+= \oplus_{i=0}^n \Z_+\al_i\}
\end{align}
for all but $A_{n}$.
Explicitly, $\tu^{i}=u_i=w_0(w_0^{\om_i})$
\end{lemma}
\smallskip
}

\rmk
Claim (ii) is a demonstration that
Theorem \ref{INTRINLA} is sufficient for a {\em complete}
characterization of the {\em last roots\,} of the $\la$\~sequences.

Note that one can use here the plane geometric interpretation
of the reduced decompositions
from \cite{Ch5} for affine classical root systems
(see there the reference concerning the non-affine case
and the so-called reflection equation).
\sq
\smallskip

For $\la$ satisfying Theorem \ref{INTRINLA},
one can introduce {\dfont quasi-simple roots\,} $\tbe\in\la$ that
are not sums of the roots from $\la$
as in ($a$). Arbitrary $\tal\in\la$ are their
sums, but, generally, they are not linearly
independent vectors; for example, there are four quasi-simple
roots in $\la(w)$ for $w=(4231)\in \bS_4$ (the case of $A_3$).

It is possible to express
the set $\{\tbe\}$ of quasi-simple roots of $\tw$ in terms
of $\{\tbe\,'\,\}$ for $\tw'$ as
$\tw=\tw's_i$, $l(\tw)=l(\tw')+1$:
\begin{align}\label{tbesip}
&\{\tbe\}=\{\al_i\}\cup\{\,s_i(\tbe\,'\,)\, \mid\,
\tbe\,'\,+\al_i\not\in\la(\hw')\,\}.
\end{align}

\comment{
\begin{proposition}\label{QUASISIMPLE}
All roots in $\la(\hw)=\{\tal^l,\ldots,\tal^1\}$
are quasi-simple if and only if this sequence is non-decreasing,
i.e., either $(\tal^i,\tal^{i+1})=0$ or
$\tal^{i+1}-\tal^{i}$ is positive, a linear
combination of simple roots with non-negative coefficients
for all $l>i\ge 0$.
\end{proposition}
{\em Proof}. If the sequence $\la(\hw)$ is non-decreasing
then obviously all roots there are {\em quasi-simple}.
Let $\{\ldots,\tbe,\tal,\ldots\}$ be the first pair
in this sequence such that $\tbe-\tal$ is not positive.
It corresponds to
$\{\cdots,s_j,s_i,\cdots\}$ in the reduced decomposition
$\hw=\hu s_j s_i \tv$, where $\tal=\tv^{-1}(\al_i)$.\sq
}

\smallskip
The following variant of this definition (they coincide in
the simply-laced case) has applications to the
Bruhat ordering. We call
$\tbe\in\la$ a {\dfont pseudo-simple root} if $m\tbe$ is not a sum
of roots in $\la$ for any $m\in \N$.

\begin{proposition} \label{BRUHPSEUDO}
(i) Given a $\la$\~ sequence $\la=\la(\hw)$, the indices
$\{p\}$ of pseudo-simple
roots $\tal^p$ (see (\ref{lambdainv})) are exactly those
satisfying the condition $l(\hw')=l-1$ for $\hw'$ obtained
from $\hw$ by deleting $s_{i_p}$ in $\hw=\pi_r s_{i_l}\cdots
s_{i_1}$:
$$\hw'=\pi_r s_{i_l}\cdots s_{i_{p+1}}s_{i_{p-1}}\cdots s_{i_1}=
\hw s^p \for s^p=s_{\tal^p}$$.

(ii) The set $\{\tbe\}$ of pseudo-simple roots of
$\tw=\tw's_i$ such that $l(\tw)=l(\tw')+1$ is as follows.
For the set of pseudo-simple roots 
$\{\tbe\,'\,\}$ for $\tw'$,
$$
\{\tbe\}=\{\al_i\}\cup\{\,s_i(\tbe\,'\,)\, \mid\,
m\,'\,\tbe\,'\,+\al_i\not\in\la(\hw')\hbox{\ for \ any } 
m\,'\,\in \N\}.
$$

(iii)
Arbitrary roots from $\la(\hw)$ are linear combinations of
pseudo-simple roots with positive rational coefficients.

\end{proposition}
{\em Proof}. The condition $l(\hw')=l-1$ is equivalent to
the positivity  $s^p(\tal^q)\in \tR_+$ for $q=p+1,
\ldots,l$; see (\ref{tal}).
If $s^p(\tal^q)<0$ then $-\tbe=\tal^q-m\tal^p<0$
for $m=((\al^p)^\vee,\al^q)>0, \tbe>0$ and
$m\tal^p=\tbe+\tal^q$. However, $-\tbe$ cannot appear in
$\la$\~sets and has to 
coincide with $-\tal^r$ for some $r<p$;
see (\ref{tlaw}). Thus $m\tal^p=\tal^r+\tal^q$ and $\tal^p$ is
not pseudo-simple in $\la$.

Let us begin now with $\tal^p$ such that $m\tal^p=\tal^r+\tal^q$
for $q>p>r$. Considering a root subsystem of rank two containing
$\tal^p,\tal^r,\tal^q$ as positive roots,
the coefficient $m$ can be $1$ for $A_2$, $1,2$ for
$B_2$, $1,2,3$ for $G_2$. In either case, the root
$$s^p(\tal^q)=\tal^q-m\,'\,\tal^p=(m-m\,'\,)\tal^p-\tal^r \for
m\,'\,=((\al^p)^\vee,\al^q)
$$
is negative.

Claim (ii) is parallel to (\ref{tbesip}). Claim (iii)
follows from (ii). Indeed, $\al_i$ belongs to the set of
pseudo-simple roots in $\la(\hw)$. Also, given a pseudo-simple
root $\tbe\,'\,\in \la(\hw')$, the root $s_i(\tga\,'\,)$ will be 
pseudo-simple
in $\la(\hw)$ for the last $\tga\,'\,$ in the {\em sequence}
$\la(\hw')$ satisfying $\tga\,'\,\in u\tbe\,'\,+v\al_i$ as
$u,v\in \Q, u>0$.
\sq
\medskip

\subsection{Theorem about triples}
We are going to discuss a connection of the
non-{\em quasi-simple} roots and the Coxeter transformations;
the latter become Coxeter {\em permutations} in the context 
of $\la$\~sets.
The theorem below is expected to be 
an important tool for the classification
of semisimple representations of DAHA and similar questions.
It clearly demonstrates why dealing with the intertwining
operators for arbitrary root systems is significantly more 
difficult than in the $A_n$\~ theory (where much is known).
The classical theory of root systems \cite{Bo,Hu} is uniform
at level of the generators and relations. However, if the
``relations of relations" are considered the root systems
behave differently; the simplest are $A_n$ and the rank two
systems.
\smallskip

First of all,
positive roots $\al,\be$ from a rank two root system $R^2$
of type $A_2,B_2,G_2$ are simple if and only if
\begin{align}\label{ranktwosim}
&\tal+\tbe\in R^2 \hbox{\ and\ }
|\tal|\neq |\tbe| \for B_2, G_2,\where
|\tal|^2=2\nu_\al.
\end{align}
Note that  $\tal+\tbe$ is always a short root for
such $\tal,\tbe$.

Given a reduced decomposition of $\hw$,
the endpoints in $\la(\hw)$  
corresponding to (complete) Coxeter sub-products 
$(\cdots s_i s_j s_i)$ are such $\tal,\tbe$.
Vice versa, if $\tga=\tal+\tbe\in \tR_+$ and also
$|\tga|=|\tal|=|\tbe|$ as $|\tal|=|\tbe|$, 
then $\tal,\tbe$ come from (\ref{ranktwosim}) for
some $R^2$  unless $\tR$ is of type $\widetilde{G}_2$.


If $\tR$ is the affine system of type $\widetilde{G}_2$
and $|\tga|=|\tal|=|\tbe|$, then we need to assume additionally
that $\tal,\tbe,\tga$
are all short and there are no long roots among their
linear combinations, i.e., that they do not belong to 
any (finite) subsystem of $\tR$ of type $G_2$.
We call such triples $A_2$\~{\em pure-short}.
\smallskip

In the simply-laced case, we set (technically)
$lng=sht$, i.e., any conditions
that certain roots are long or short are disregarded as
$\tR=\widetilde{A},\widetilde{D},\widetilde{E}$ in the 
following theorem and below. Recall that we use the notation
$\tR=\widetilde{A}_n, \widetilde{B}_n,\ldots, \widetilde{G}_2$
for $\tR$ corresponding to $R=A_n,B_n,C_n$.

A connected part of the sequence 
$\la(\hw)$ isomorphic to the
sequence of all positive roots of type
$A_2, B_2$ or $G_2$ will be  called a 
{\dfont segment of rank two} in the next theorem.
Concerning {\em root subsystems} used in this theorem, 
they can be arbitrary; it suffices to suppose
that they are intersections of $\tR$ 
with $\Q$\~subspaces in $\Q[\tR]$.
\smallskip

\begin{maintheorem}\label{RANKTWO}
Given a reduced decomposition of $\hw\in \hW$,
let us assume that $\tal+\tbe=\tga$ for the roots
$\ldots,\tbe,\ldots,\tga,\ldots,\tal\ldots$
in $\la(\hw)$ ($\tal$ appears the first), where
only the following combinations of their lengths
are allowed:
\begin{align}
&(a)\, \hbox{lng}+\hbox{sht}=\hbox{sht} \hbox{ \ or\ }
\hbox{sht}+\hbox{lng}=\hbox{sht}
\ \,(\,\hbox{all\ systems\ }\tR\,),
\label{lngshtsht} \\
&(b)\, \hbox{lng}+\hbox{lng}=\hbox{lng}\ \,
 (\widetilde{B},\widetilde{F}_4) \hbox{\, \ or\ \,}
\hbox{sht}+\hbox{sht}=\hbox{sht} \ \,
(\widetilde{C},\widetilde{F}_4),
\label{shtshtsht} \\
&(c)\,\hbox{the\ roots\ \ } \tal,\tbe \hbox{\ \ are\ \ } 
A_2\hbox{\~pure-short\ \ when\ \ }
\tR= \widetilde{G}_2.
\label{puregtwo}
\end{align}

(i) Let $\,[\,\tbe,\tal\,]\,$ be a segment in $\la(\hw)$
from $\tal$ to $\tbe$. There exists $\tu\in \hW$ such that
$\tu(\tga)=\tga$ and the Coxeter transforms in 
$\,[\,\tu(\tbe),\tu(\tal)\,]\,$ 
can be used to make the triple 
$\tu(\tal),\tu(\tbe),\tu(\tga)=\tga$
part of a segment $L$ of rank $2$ in $\la(\hw)$.
Moreover, one can assume here that 
$\tu=s_{j_l}\cdots s_{j_1}$
and all consecutive products 
$\tu_m=s_{j_m}\cdots s_{j_1}$ ($\tu_{0}=$id)
leave $\tga$ invariant for $m=1,\ldots,l$ and 
$$
\tu_m(\tal)\in \la(\hw)\ni
\tu_m(\tbe),\ [\tu_m(\tal),\tu_m(\tbe)]\subset
[\tu_{m-1}(\tal),\tu_{m-1}(\tbe)], 
$$
where the Coxeter transforms can be used in
$[\tu_{m-1}(\tal),\tu_{m-1}(\tbe)]$ before 
finding the next $[\tu_{m}(\tal),\tu_{m}(\tbe)]$. 

One can take $\,\tu=$id \, in cases ($a$) for   
$\,\tR=\widetilde{A}_n,\,\widetilde{B}_n,\,
\widetilde{C}_n,\,\widetilde{F}_4,\,\widetilde{G}_2\, $ and 
also for the system $\widetilde{G}_2$ under ($c$).

(ii) For the triples of type ($b$),
such $L$ does not exists with $\tu=$id\, if (and only
if unless for $\widetilde{C}_{n\ge 4}$)
a root subsystem $R^3\subset \tR$ of type
$B_3$ or $C_3$ ($m=1,2$) can be found such that 
\begin{align}\label{rankthrees}
&\tbe=\ep_1+\ep_3,\, \tal=\ep_2-\ep_3,\
\ep_1-\ep_2\,\not\in\, [\,\tbe,\tal\,]\,
\not\ni\, m\ep_3, \hbox{\ or}\\
&\tbe=\ep_2+\ep_3,\, \tal=\ep_1-\ep_3,\
\ep_1-\ep_2\,\not\in\, [\,\tbe,\tal\,]\,
\not\ni\, m\ep_3\,;\label{rankthreess}
\end{align}
the notation is from \cite{Bo}, the positivity in
$R^3$ is induced from $\tR_+$. Equivalently,
the sequence $[\,\tbe,\tal\,]\cap R^3_+$ (with the
natural ordering) must be 
\begin{align}
&\{\,\ep_1+\ep_3,\, m\ep_1,\, \ep_2+\ep_3,\, \tga=\ep_1+\ep_2,\,
\ep_1-\ep_3,\,
m\ep_2,\, \ep_2-\ep_3\,\} \hbox{\ or\ }\label{rankthree}\\
&\{\,\ep_2+\ep_3,\, m\ep_2,\, \ep_1+\ep_3,\, \tga=\ep_1+\ep_2,\,
\ep_2-\ep_3,\,
m\ep_1,\, \ep_1-\ep_3\,\}\label{rankthreeo}
\end{align}
up to Coxeter transforms in $R^3$ and changing the order
of all roots in (\ref{rankthree}) to the opposite.
If $\tal$ is assumed simple in $R^3$, then
only (\ref{rankthrees}),(\ref{rankthree}) may occur.

If such $R^3$ exist, one
can still take $\tu=$id\, using Coxeter transforms
in the whole $\la(\hw)$ if the latter set
contains either $\ep_1-\ep_2$ or $m\ep_3$ for 
\underline{every} such $R^3$.

(iii) For the root system
$\tR$ of type $\widetilde{C}_{n\ge 4},
\widetilde{D}_{n\ge 4}\,$ or 
$\widetilde{E}_{6,7,8}$,
let as assume that $\tal$ is simple in $\tR_+$. Then
such $L$ does not exists with $\tu=$id\, if (and only
if unless for $\widetilde{C}$)
a root subsystem $R^4\subset \tR$ of type
$D_4$ can be found such that 
\begin{align}\label{rankthreesd}
&\tbe\,=\,\ep_1+\ep_3,\ \tga\,=\,\ep_1+\ep_2,\ 
\tal\,=\,\ep_2-\ep_3,\,
\\
&\{\,\ep_1-\ep_2,\,\ep_3-\ep_4,\,
\ep_3+\ep_4\,\}\cap [\,\tbe,\tal\,]\ =\ \emptyset\,;\notag
\end{align}
the notation is from \cite{Bo}, the positivity in
$R^4$ is induced from $\tR_+$. Equivalently,
the sequence $[\,\tbe,\tal\,]\cap R^4_+$  must be 
\begin{align}\label{rankthreed}
\{\,\tbe=\ep_1+\ep_3,\,&\ep_1-\ep_4,\, \ep_1+\ep_4,\, 
\ep_1-\ep_3, \tga=\ep_1+\ep_2,\\ 
&\ep_2+\ep_3,\, \ep_2+\ep_4,\, \ep_2-\ep_4,\, 
\tal=\ep_2-\ep_3\,\}\notag
\end{align}
up to Coxeter transforms in $R^4$. Equivalently,
$[\,\tbe,\tal\,]\cap R^4_+$ 
is the $\la$\~set of $s_{\th^4}$ in $R_+^4$ for
the maximal root $\th^4$.
Equivalently, $\tal=\al_2^4$, $\tga=\th^4$,
$\tbe=\th^4-\al_2^4$ and also $\al_2^4$ (a unique simple
root non-orthogonal to $\th^4$) is the only simple root
from $R^4$ in $[\,\tbe,\tal\,]$. 

In case of 
$\widetilde{C}_{n\ge 4}$ (when $\tGa$
contains a subdiagram of type  $D_4$),
either (\ref{rankthreesd})-(\ref{rankthreed})
must hold  or those from (ii) for $R^{3}$ 
if such $L$ does not exist with $\tu=$id.

If $\tal$ is not assumed simple in $\tR_+$, then
the condition is that $\tga=\th^4$ and 
$[\,\tbe,\tal\,]\cap R^4_+$ contains a unique simple
root from $R_+^4$. 
\end{maintheorem}
\medskip

\subsection{Discussion}
Before proving the theorem, let us discuss
some corollaries and general facts that we will need 
to clarify and verify claims (i,ii,iii).
\smallskip

We do not give the complete list of possibilities
for (iii) without the assumption that $\tal$ is simple
It is analogous to that in (ii).
Actually, it suffices to assume that it is simple in $R_+^4$
because one can always switch to considering 
$(\hw')^{-1}[\tbe,\tal]$
where $\hw'$ is the portion of the reduced decomposition of
$\hw$ {\em before} $\tal$. However, sometimes it is convenient
to avoid making $\tal$ simple; the next comment
can be readily extended to (iii).
\smallskip

\rmk
Note that (\ref{rankthree}) is the only case that may
occur if $\tal$ is assumed to be simple in $R_+^3$.
In this case, if $\ep_1-\ep_2$ or $m\ep_3$
appear in $\la(\hw)$, then it can happen only
{\em after} the last root in this sequence, $\ep_1+\ep_3$.
For instance, $m\ep_3$ must be after $m\ep_1$
since $\ep_3+(\ep_1-\ep_3)=\ep_1$ implies
that it is after $\ep_1+\ep_3$. 
As for (\ref{rankthreeo}), $m\ep_3$ (if present in $\la(\hw)$)
appears after the last root, $\ep_2+\ep_3$. Also, 
$\la(\hw)$ {\em must}
contain $\ep_1-\ep_2$ before $\ep_1-\ep_3$, since
$\ep_1-\ep_3=(\ep_1-\ep_2)+(\ep_2-\ep_3)$ and
$\ep_2-\ep_3$ is after $\ep_1-\ep_3$ in (\ref{rankthreeo}).

Similarly, if the inverse of (\ref{rankthree}),
\begin{align}\label{rankthreei}
&\{\,\ep_2-\ep_3,\, m\ep_2,\, \ep_1-\ep_3,\, \tga=\ep_1+\ep_2,\,
\ep_2+\ep_3,\,
m\ep_1,\, \ep_1+\ep_3\,\},
\end{align}
belongs to $\la(\hw)$, then the latter sequence {\em must}
contain $\ep_1-\ep_2$ and $m\ep_3$ (before $\ep_1+\ep_3$).
For instance, $(\ep_1-\ep_2)+(\ep_2+\ep_3)=\ep_1+\ep_3$
results in $\ep_1-\ep_2\in \la(\hw)$.
\sq
\medskip

{\bf Example of {\mathversion{bold}$E_6$}}.
The following example seems a good illustration
of (iii). Let $w=s_{\th}$ for the maximal root $\th=\om_2$ 
in $R$ of type $E_6$. Then $\al_2$ cannot be moved
in the triple 
$\{\be=\th-\al_2,\ga=\th,\al=\al_2\}$ from its first position 
and these roots  cannot be made neighboring in 
$[\be,\al]\in \la(s_{\th})$. Thus there must exist a 
$9$\~root set of type $D_4$ from (\ref{rankthreed}).
In notation from \cite{Bo}, the simple roots
from $R^4$ are as follows:
$\ep_1-\ep_2=\al_4,\, \ep_2-\ep_3=\al_2\,$ and 
$
\ep_3-\ep_4\,=\,\al_3+\al_4+\al_5,\ \,
\ep_3+\ep_4\,=\,\al_1+\al_3+\al_4+\al_5+\al_6.
$
The complete description is as follows
(we show the roots from $E_6$ under the
corresponding ones from the $9$\~set):
\vskip 0.1cm
\renewcommand{\baselinestretch}{0.1}. 
{\tiny
\begin{align*}
\ep_1&+\ep_3,&\ep_1&-\ep_4,&\ep_1&+\ep_4,&\ep_1&-\ep_3,& 
\ep_1&+\ep_2,&\ep_2&+\ep_3,&\ep_2&+\ep_4,&\ep_2&-\ep_4,&
\ep_2&-\ep_3:&\\ \\
12&321,& 11&211,& 01&210,& 12&221,& 12&321,& 
00&100,& 11&111,& 01&110,& 00&000.&\\
  &1   &   &1   &   &1   &   &1   &   &2   &   &1   &   
  &1   &   &1   &   &1   & 
\end{align*}\sq
}

\renewcommand{\baselinestretch}{1.}
\smallskip

{\bf Example of {\mathversion{bold}$F_4$}}.
Let $w=s_{\vth}$ for the {\em maximal short root}
 $\vth=\om_1$ in $R$ of type $F_4$. 
Then $\al_1$ cannot be moved from its first position
in the triple 
$\{\be=\vth-\al_1,\,\ga=\vth,\,\al=\al_1\}$. 
According to the theorem, there must exist a 
$7$\~root set of type $C_3$ from (\ref{rankthree}).
In notation from \cite{Bo}, it is as follows:
\renewcommand{\baselinestretch}{0.6}. 
{\small
\begin{align*}
\ep_1+&\ep_3,&2&\ep_1,&\ep_2+&\ep_3,&\ep_1+&\ep_2,&
\ep_1-&\ep_3,&2&\ep_2,&\ep_2-&\ep_3:&\\ \\ 
12&31,& 23&42,& 01&21,& 12&32,& 11&10,& 01&22, &00&01.&
\end{align*}
}
\sq
\renewcommand{\baselinestretch}{1.}\\
\medskip

Note, that  $\{\ep_i\}$ are the roots from
Theorem \ref{RANKTWO} in these two examples, 
i.e., they are the basic vectors used
in \cite{Bo} to describe $R^4$ and $R^3$ correspondingly
(not $\{\ep_i\}$ for $E_6$ and $F_4$).

One can readily find all {\em reflections} $s_{\tga}$
with the endpoints $\al=\al_i$,$\tbe=\tga-\al_i$ 
$\in \la(s_{\tga})$ that are {\em non-movable} inside 
$\la(s_{\tga})$ under the Coxeter transformations. The 
examples above can be generalized as follows.
\smallskip

\begin{proposition}\label{NONMOVREF}
(i) Let  us assume that the $\la$\~sequence
$\la(s_{\tga})$ for $\tga\in \tR_+$
has a unique beginning $\tga=\al_i$ for $0\le i\le n$.
Then $\tbe=-s_{\tga}(\al_i)$ is its unique end.
Provided that the lengths of $\tga$ and $\al_i$
coincide and the nonaffine components of these roots are not
proportional, $\tbe=\tga-\al_i$ and $\{\tbe,\tga,\tal\}$
form an $A_2$\~triple.

(ii) Let $\tga=[\ga,\nu_\ga j]$ for $\ga\in R_+,j\ge 0$.
Then the conditions from (i) are satisfied only for $i>0$
and if and only if

a) $(\ga,\al_i)>0$ for a unique $1\le i \le n$ ($\le 0$ otherwise),

b) moreover, $|\ga|=|\al_i|$ and  $\ga\neq \al_i$ for such $\al_i$.

\noindent
Here $j\ge 0$ can be arbitrary.

(iii) Let $\tga=[-\ga,\nu_\ga j]$ for $\ga\in R_+,j> 0$. 
Then (i) holds only for $i>0$ and if and only if
$\ga$ is either a maximal short root or
a maximal (long) root, $\vth'$ or $\th'$,
for a root subsystem $R'$ corresponding
to a connected subdiagram $\Ga'\in \Ga$ such that

a)  $\al_{0'}\not\in \Ga'$ for  
$\al_{0'}\in \Ga$ linked to $\al_0$ in the affine Dynkin diagram 
$\tGa$,

b)  there exists a unique $\Ga\ni \al_{m}\not\in \Ga'$ 
connected with $\Ga'$ by a link. 

\noindent
Here $i=m$ and $\vth'$ or $\th'$ are chosen to ensure that
$|\ga|=|\al_{m}|$ ($\Ga'$ must contain roots of length
$|\al_m|$); $j>0$ can be arbitrary.
\end{proposition}

{\em Proof.} Claim (ii) is straightforward. The uniqueness 
of $\al_i$ gives that $\al_0\not\in \la(s_{\tga})$, i.e.,
$i>0$. Indeed, 
$s_{\tga}(\al_0)=\al_0 +2\frac{(\ga,\vth)}{(\ga,\ga)}\ga$,
and the latter is a positive root since $(\ga,\vth)\ge 0$.
However, at least one simple root must sit in $\la(s_{\tga})$;
therefore, it can be only nonaffine.

Similarly, $(\ga,\al_{i'})\le 0$ for $i'=1,\ldots,n$ unless $i'=i$,
(when this inner product must be strictly positive). Otherwise, 
$\al_{i'}\in \la(s_{\tga})$, which contradicts the uniqueness.
The condition $|\ga|=|\al_i|$
here is necessary and sufficient to ensure that we really
deal with an $A_2$\~triple.
\smallskip

Assuming now that $\tga=[-\ga,\nu_\ga j]$ for $\ga>0,j>0$,
let us check that $i>0$ in (i). If $i=0$ then (i) holds if
and only if $\ga$ is short, $(\ga,\vth)=1$ and
$(\ga,\al_{i'})\ge 0$ for {\,\em all\,} $i'=1,\ldots, n$.
However, the latter implies that $\ga=\vth$, which
contradicts to  $(\ga,\vth)=1$. Thus $(\ga,\vth)=0$,
which holds if and only if the decomposition of
$\ga$ does not contain $\al_{0'}$ connected with $\al_0$
in $\tGa$ by a link.

Let $\Ga'\subset \Ga$ be the {\em support} of $\ga$ (the set of
all simple roots appearing in its decomposition). Then
any simple root $\al_{m}$ from (b) can be the {\em beginning} 
of the sequence $\la(s_{\tga})$ due to $(\ga,\al_{m})<0$, 
which implies that 
$$
s_{\tga}(\al_{m})=\al_{m}+
2\frac{(\al_{m},\ga)}{(\ga,\ga)}\tga\in \tR_-
\and \al_m\in \la(s_{\tga}). 
$$  
At least one such $\al_m$ exists. Thanks to the uniqueness,
one such $\al_m$ can exist. Moreover, 
$(\ga,\al_{i'})\ge 0$ must hold for all $\al_{i'}\in \Ga'$;
otherwise, there will be other simple roots in $\la(s_{\tga})$. 
These conditions imply
that either $\ga=\vth'$ as $\al_m$ is short or
 $\ga=\th'$ (the maximal positive root in $\Ga'$) as $\al_m$
is long, which results from $|\ga|=|\al_m|$.
Here we use that $\vth$ is a minimal positive root in $P_+$.
\sq
\smallskip

We can now describe all non-affine roots $\tga\in \tR_+$ such
that $\tga$ and the ends of $\la(s_{\tga})$ form a 
{\em non-movable} $A_2$\~triples; the examples of 
$E_6$ and $F_4$ considered above correspond to $\tga=\vth$.
Proposition \ref{NONMOVREF}, (iii) 
gives a complete description of such roots in the form
$\tga=[-\ga,\nu_\ga j]>0$ for $\ga>0$. Therefore we can 
restrict ourselves with $\tga=[\ga,\nu_\ga j]>0$. Moreover,
it suffices to assume that
$j=0$, since the answer is uniform with respect to
$j\ge 0$; see Proposition \ref{NONMOVREF}, (ii). We will consider
here the cases $E_6,F_4,B,C,D$; there are 
$7$ such $\ga>0$ for $E_7$ and $22$ for $E_8$. 
\vskip 0.1cm

The following are the lists of nonaffine roots $\ga$ such that
the ends $\be=\ga-\al_i, \al=\al_i$ of $\la(s_{\ga})$ are 
{\em non-movable} under
the Coxeter transforms within $\la(s_{\al})$ and  
$\{\be,\ga,\al\}$ form an $A_2$\~triple in the cases of 
$F_4,B,C$ (i.e., subject to  $|\al_i|=|\ga|$). The bar 
shows the place of the corresponding $\al_i$.
\medskip

{\bf The case of {\mathversion{bold}$E_6$}}.
The roots $\ga\in R_+$ with non-movable ends: 
\renewcommand{\baselinestretch}{0.6} 
{\small
\begin{align*}
01&\overline{2}10,& 1\overline{2}&210,& 01&2\overline{2}1,& 
12&\overline{3}21,& 12&321.         & \\
  &1              &              &1   &   &1              &  
  &1              &   &\overline{2} & \notag    
\end{align*}
}
\renewcommand{\baselinestretch}{1.}
The corresponding $\{\be=\ga-\al_i,\ga,\al_i\}$ automatically
form a triple.
\smallskip

{\bf The case of {\mathversion{bold}$F_4$}}.
The roots $\ga\in R_+$ with non-movable ends of $\la(s_\ga)$ and
subject to $|\al_i|=|\ga|$ are:
{\small
\begin{align*}
&01\overline{2}1,& &1\overline{2}20,& &12\overline{3}1,&
&123\overline{2},& &1\overline{3}42,& &\overline{2}342.&
\end{align*}
}
\vskip -0.3cm

{\bf The case of {\mathversion{bold}$B,C,D_n$}}.
Given $\al_i=\ep_i-\ep_{i+1}$, the corresponding
root $\ga$ equals  $\ep_{i-1}+\ep_{i}$ for
$i=2,\ldots,n-1$ (it is unique) provided
that $n\ge 3$ and $i<n-1\ge 3$ for $D_n$. 
The notation is from \cite{Bo}. 
\medskip

The {\em Coxeter sequences\,} of types $A_2,B_2,G_2$
in reduced decompositions of a given $\hw\in \hW$
can be naturally identified
with segments of $\la(\hw)$ isomorphic to the sequences of
positive roots of a rank two systems.
By a {\em Coxeter sequence} in a reduced
decompositions of $\hw$, we mean a
representation $\tw =\hw'\tv\hw''$ subject to
$l(\hw)=l(\hw')+l(\tv)+l(\hw'')$ such that
$\tv=(\cdots s_is_js_i)$
with $m_{ij}\ge 3$ factors. We will consider {\em all
possible reduced decompositions} of $\hw$ in the
next corollary, i.e., the Coxeter sequences are those
that can be made consecutive in at least one reduced
decomposition of $\hw$.
\smallskip

\begin{corollary}\label{RANKTWOCOR} 
(i) The triples satisfying one
of the conditions ($a$),($b$) or ($c$) and such that they can
be made consecutive in at least one reduced decomposition
of $\hw$ are in one-to-one correspondence
with  Coxeter sequences of type\\
$A_2$ in case ($a$) for the root systems $\widetilde{A},\widetilde{D},
\widetilde{E}$ or in cases ($b$),($c$),\\
$B_2$ in case ($a$) for the systems $\widetilde{B},\widetilde{C},
\widetilde{F}_4$, or\ \ $G_2$\, for ($a$) and $\widetilde{G}_2$.

(ii) Let $\tv\,(\tR^0_+)\subset \la(\hw)$
for $\tv\in \tW$ and 
the set of all positive roots
$\tR^0_+$ of a finite root subsystem $\tR^0$
generated over $\Z$ by a \underline{connected} subset of 
$\{\al_0,\cdots,\al_n\}$.
Then Coxeter transformations can be used to make
$\tv\,(\tR^0_+)$ a
(connected) segment in the corresponding $\la(\hw)$ 
if there are no triples $\tal,\tbe,\tga$
in $\la(\hw)\cap \tR^0_+$ that can be extended to
 a $7$\~set from (ii) or a $9$\~set from (iii)  
(upon the conjugation making $\tal$ simple).\\  
\sq
\end{corollary}
\medskip

\rmk
(i) The proof below is actually an algorithm of finding
such rank two segment(s) $L$. We will prove (i,ii,iii) by induction,
assuming that these claims hold for all
$$
\{\ldots,\tbe,\ldots,\tal,\ldots\}
\hbox{\ in\ any\ }\la(\hw)
\hbox{\ such\ that \ }
\ell[\tbe,\tal]=\ell'<\ell,
$$
where the {\em $\ell$\~length} $\ell[\tbe,\tal]$ equals the
number of roots in the segment $[\tbe,\tal]\subset \la(\hw)]$
including the endpoint. The algorithm below
diminishes $\ell$. The part of the theorem concerning
using $\tu\neq$id naturally emerges in this procedure too.
The last claim from (ii) about using the whole $\la(\hw)$
and its counterpart for (iii) (that was not formulated
explicitly) require a somewhat special consideration; 
it will be omitted.

(ii) Claim (i) of Corollary \ref{RANKTWOCOR} is a straight
application of (ii,iii) from the theorem. We will skip the
proof of claim (ii) of this corollary; it is not too 
important in this paper. Generally, there are quite a few 
situations when we can collect the roots from certain 
subsets of $\la(\hw)$ together
using Coxeter transforms. Claim (ii) gives an example.

Note that long roots that are sums of two short roots 
are excluded from the theorem; they exist for
$B,C,F$ if the short roots are orthogonal to each other.
Such pairs are not needed in Corollary \ref{RANKTWOCOR},(i) 
for the one-to-one correspondence 
with the Coxeter transformation of types $A_2,B_2,G_2$.

(iii) 
The theorem can be verified much simpler
for $\widetilde{A}_n$ and for the root systems 
$\widetilde{B}_n,\widetilde{C}_n$ where the triples are
{\em not} of type ($b$). Here either the planar 
interpretation can be used or the fact that the simple roots have
multiplicities no greater than $2$ in all positive roots.
We will discuss the latter ``numerical" approach in detail 
and give references concerning the planar interpretation.

Let us mention that the $7$\~sets from (ii) is somewhat simpler
to deal with than the $9$\~sets from (iii);
see Lemma \ref{LEMTHREE} and especially Lemma \ref{LEM3AB}.
Also, the case of $F_4$ is more involved than those of $B,C$
and requires a special consideration (some details will be
omitted).
\sq
\medskip

\section{Induction process}\label{sec:ProofTriples}
Let us first clarify property (\ref{rankthree})
for $\widetilde{B},
\widetilde{C},\widetilde{F}_4$ under ($b$).
\smallskip

\subsection{Subsystems 
\texorpdfstring{{\mathversion{bold}$B_3,C_3$}}{BC}}
We will describe all possible intersections of
$\la$\~sets $\la(\hw)$ with $R_+^3$ from (ii).
The notations from the previous section are used. 

\begin{lemma}\label{LEMTHREE}
Let the triple $\{\tal,\tbe,\tal+\tbe=\tga\}$ of
type ($b$) belong to the intersection
$\la^0$ of $\la(\hw)$ with  $R_+^3$, the 
set of the positive roots of type  $B_3$ or $C_3$ in $\tR$:
$$\{\, m\ep_i, \, 1\le i\le 3,\
\ep_i\pm\ep_j,\,  i<j\le 3,\}\,\subset\, \tR_+\hbox{\ as\ \ }
m=1,2.
$$
Assuming that $\tal$ is simple in $R_+^3$,
the ordered pair $\{\tbe,\tal\}$ coincides with
one of the following pairs:
\begin{align}
&(0):\{\ep_2-\ep_3,\ep_1-\ep_2\},\ \
(1):\{\ep_2+\ep_3,\ep_1-\ep_2\},\label{zeroone}\\
&(2):\{\ep_1-\ep_2,\ep_2-\ep_3\},\ \
(3):\{\ep_1+\ep_3, \ep_2-\ep_3\}.\label{twothree}
\end{align}
The second root (after $\ep_1-\ep_2$) in the
sequence $\la^0$ can be
$\ep_1-\ep_3$ and $m\ep_3$ in cases (0,1); the second root
(after $\ep_2-\ep_3$) can be
$\ep_1-\ep_3$ and $m\ep_2$ in cases (2,3).

Case (1a): (1) and the second root is $\ep_1-\ep_3$.
Then the third root is $m\ep_1$ or $\ep_2-\ep_3$ and
$$
\{\,\tbe=\ep_2+\ep_3,\,\ep_1+\ep_2,\, m\ep_3,\,\ep_1+\ep_3,\,m\ep_1,
\,\ep_1-\ep_3,\,\tal=\ep_1-\ep_2\,\}\,
\subset \,\la^0.
$$

Case (2a): (2) and the second root is $m\ep_2$. Then the
third root is $\ep_1-\ep_3$ or $\ep_2+\ep_3$ and
$$
\{\,\tbe=\ep_1-\ep_2,\,m\ep_1,\,\ep_1-\ep_3,\,\ep_1+\ep_2,\,
m\ep_2,\,\tal=\ep_2-\ep_3\,\}\,
\subset \,\la^0.
$$

Case (3a): (3) and $\ep_1-\ep_2\in \la^0$. Then
\begin{align*}
&\{\,\tbe=\ep_1+\ep_3,\, \ep_1+\ep_2,\, m\ep_1,\,
\ep_1-\ep_2,\, \ep_1-\ep_3,\, \tal=\ep_2-\ep_3\,\}\,
\subset\, \la^0.
\end{align*}
The second root in $\la^0$ can be $\ep_1-\ep_3$ or $m\ep_2$.

Case (3aa): (3a) and $m\ep_2\not\in \la^0$. Then $\la^0$
is precisely this set. These $6$ roots appear in this
very order in $\la^0$ modulo the Coxeter transforms in
this set.

Case (3b): (3) and $\ep_1-\ep_2\not\in \la^0$. Then
\begin{align*}
&\{\,\tbe=\ep_1+\ep_3,\, m\ep_1,\, \ep_2+\ep_3,\, \ep_1+\ep_2,\,
\ep_1-\ep_3,\,
m\ep_2,\, \tal=\ep_2-\ep_3\,\}\, \subset\, \la^0.
\end{align*}

Case (3bb): (3b) and
$m\ep_3\not\in \la(\hw)$. Then these $7$ roots
constitute $\la^0$ and they
appear in this very order in $\la^0$
modulo the Coxeter transforms in this set.
\end{lemma}
{\em Proof.}
Theorem \ref{INTRINLA} gives that
the intersection $\la^0$ of $\la(\hw)$ with the root system 
$\tR^0=R^3$ is a $\la$\~set
with respect to $R_+^3=R^3\cap\tR_+$. Using this fact,
all claims are straightforward.
For instance,
all possible orderings of $7$ roots from (3bb) are in one-to-one
correspondence
with reduced decomposition of $w=s_2 s_3 s_2 s_1 s_2 s_3 s_2$
in the Weyl group of type $BC_3$ in the notation from \cite{Bo}.
The Coxeter transformations in $\la^0$ are:
\begin{align*}
&[\ep_2+\ep_3,\, \ep_1+\ep_2,\,
\ep_1-\ep_3]\,\leftrightarrow\,
[\ep_1-\ep_3,\, \ep_1+\ep_2,\,
\ep_2+\ep_3],\\
&m\ep_1\,\leftrightarrow\,  \ep_2+\ep_3,\ \
\ep_1-\ep_3\,\leftrightarrow\, m\ep_2.
\end{align*}

Let us determine the set $\la^0$ and its ordering in cases
(3b,3bb) {\em directly} via Theorem \ref{INTRINLA};
it is an instructional exercise. The condition
$\ep_1-\ep_2\not\in\la(\hw)$ results in the following:
\begin{align*}
&(\ep_1-\ep_2)+(\ep_2+\ep_3)=\tbe\, \Rightarrow\, (\ep_2+\ep_3)
 \hbox{\ is \ before\ } \tbe,\\
&(\ep_1-\ep_2)+2\ep_2=\tga\, \Rightarrow\, 2\ep_2
\hbox{\ is\ before\ } \tga,\\
&(\ep_1-\ep_3)=(\ep_1-\ep_2)+\tal\,\Rightarrow\, (\ep_1-\ep_3)
\hbox{\  is\ after\ } \tal,\\
&2\ep_1=(\ep_1-\ep_2)+\tga\, \Rightarrow\, 2\ep_1
\hbox{\  is\ after\ } \tga.
\end{align*}

Using the condition $2\ep_3\not\in\la(\hw)$,
\begin{align*}
&(\ep_2+\ep_3)=\tal+2\ep_3\, \Rightarrow\, (\ep_2+\ep_3)
 \hbox{\  is\ after\ } \tal,\\
&(\ep_1-\ep_3)+2\ep_3=\tbe\, \Rightarrow\, (\ep_1-\ep_3)
\hbox{\  is\ before\ } \tbe.
\end{align*}
Then the following relations fix completely the order
of all $7$ roots between $\tal$ and $\tbe$
up to the Coxeter transformations in $\la^0$:
\begin{align*}
&2\ep_2=(\ep_2+\ep_3)+\tal\,\Rightarrow\, 2\ep_2
\hbox{\  is\ after\ } \tal, \hbox{\ since\ }
(\ep_2+\ep_3) \hbox{\ is \ after\ } \tal,\\
& 2\ep_1=(\ep_1-\ep_3)+\tbe\, \Rightarrow\, 2\ep_1
\hbox{\  is\ before\ } \tbe, \hbox{\ since\ }
(\ep_1-\ep_3) \hbox{\ is \ before\ } \tbe.
\end{align*}
\sq
\smallskip

Note that it is always possible to diminish
the distance between $\tal,\tbe$ {\em inside}
$R^3$ unless in case (3bb). One uses
the Coxeter transforms of type $B-C$ in cases
(1a) and (2a); the other cases are immediate.
It is of course a particular case of the Main Theorem.
The cases (0,1,2) can be described in complete
detail (similar to (3)) but we do not need it
in this paper. An analogous lemma exists for (iii);
it is useful when dealing practically with the 
{\em admissibility}.
\medskip

\subsection{Admissible triples}
We will begin with certain reductions and the consideration
of classical systems based on relatively straightforward
ways of diminishing $\ell$.
We will call $\{\tal+\tbe=\tga\in \la(\hw)\}$ satisfying
(ii) or (iii) when applicable an {\dfont admissible triple}.
This notion depends of course on choosing $\la(\hw)$, not
on the triple itself.
\smallskip

The following Lemma simplifies
dealing with the admissible triples of type (ii)
(its counterpart for (iii) will not be discussed).
 
\begin{lemma}\label{LEM3AB} In case (3a),
the triple $\{\ep_1-\ep_2,\ep_1-\ep_3, \ep_2-\ep_3\}$
is always admissible unless $\widetilde{R}$ is of
type $\widetilde{F}_4$.
\end{lemma}
{\em Proof.} Let us assume that
$$
\ep_1-\ep_2=\ep_1'+\ep_3',\,
\ep_1-\ep_3=\ep_1'+\ep_2',\,
\ep_2-\ep_3=\ep_2'-\ep_3'
$$
for $\ep_1',\ep_2',\ep_3'$ from $(R^3)'$
satisfying (3bb). The root system $\tR$ can be of
type $\widetilde{B}_n$ or $\widetilde{C}_n$.
Let us consider $\widetilde{B}_n$ for the sake of definiteness;
thus $m=1$. We will use that 
any {\em nonaffine} long root of type
$B$ can be {\em uniquely} represented as a sums of two 
{\em nonaffine} short roots. The last relation gives that
$$
\ep_2'=[0,x]-\ep_3,\,
\ep_3'=[0,x]-\ep_2\for x\in \N,\
\ep_1'=[0,-x]+\ep_1.
$$
The inequality $x>0$ is necessary  to make $\ep_2',\ep_3'$
positive (note that it can be insufficient depending on the
affine components of these roots).
Then $\ep_1'-\ep_2'=\ep_1+\ep_3+[0,-2x]$; it is
a positive root and does not belong to $\la(\hw)$ (due to (3bb)).
However, $\tbe=\ep_1+\ep_3$ belongs to $\la$ and 
therefore all {\em positive} roots
in the form $\tbe-[0,y]$ for $y\in \Z_+$ must also belong
to $\la(\hw)$, a contradiction.
\sq
\medskip

{\bf Minimality conditions.}
One can suppose that $\tal=\al_i$ is
the {\em first} root, and $\tbe$ is the {\em last} in $\la(\hw)$.
Note that if $\tga=\tal+\tbe$ is the second root
(or the last but one)
then Proposition \ref{INTRINLAP}, (iii) can be applied
and we may exclude such triples from the consideration.

Since $\tal$ is assumed simple in $R^3,R^4$, we do not
need the cases from (\ref{rankthreess}), 
(\ref{rankthreeo}) when $\tal=\ep_1-\ep_3$.
Also we do not need to 
consider possible inverting of all roots in 
(\ref{rankthree}). See 
Lemma \ref{LEMTHREE} (where $\tal$ was supposed
simple). 

If there exists at least one simple
$\al_j\neq \al_i$ in $\la(\hw)$, then we can make $\al_j$ the
first in $\la(\hw)$ and, therefore, reduce the distance between
$\tal$ and $\tbe$ and proceed by induction unless
$\al_j=\tbe$. In the latter case, the roots $\tbe,\tga,\tal$
will be transformed to $\tbe,\tga,\tal$ by a chain of Coxeter
transformation. Then two of these roots
become neighbors during this process somewhere; 
this case is governed by Proposition \ref{INTRINLAP}, (iii).
\smallskip

Generally, the theorem {\em for \underline{all} $\hw$ and all
reduced decompositions\,} is equivalent to the following
claim. {\em \underline{All} reduced decompositions
with the $\la$\~sets from $\tal$ to $\tbe$ containing
$\tga=\tal+\tbe$
such that this triple is \underline{admissible} (cannot be
included in the $7,9$-system from (ii,iii) if applicable)
have at least two \underline{first roots} or at least two
\underline{last roots}}.
\smallskip

\rem
The theorem guarantees that the cardinalities
of both sets, the {\em first roots\,} and {\em last roots\,},
will be greater than one, but considering the first roots here 
are technically sufficient (if known for all reduced 
decompositions) to proceed by induction.
Note that switching from $\hw$ to $\hw^{-1}$ and
reversing the order of the corresponding reflections and roots
gives that if the number of {\em first roots\,} is
(always) greater than one then the same holds for the number
of {\em last roots\,}. For instance, one can impose the
condition $|\tal|\ge|\tbe|$ when/if convenient.
\sq

Thus $\al_i$ can be assumed a {\em unique first root}
in the sequences $\la(\hw)$ for all reduced decompositions
of $\hw$.
This minimality constraint implies (but is not equivalent to)
the following conditions for $\tal=\al_i$:

(1) $\al_i$ belongs to any $\tde\in \la(\hw)$, i.e.,
$\tde-\al_i\in \widetilde{Q}_+\equal\oplus_{j=0}^n \Z_+\al_j $;

(2) for every $0\le j\neq i$, if $\tde-\al_j\in \tR_+$
then $\tde-\al_j\in \la(\hw)$.
\smallskip

Similarly, we can suppose that $\tbe$ is a {\em unique last root}
in $\la(\hw)$ for all possible reduced decompositions of
$\hw$; see Proposition \ref{INTRINLAP}, (ii). It is equivalent
to the conditions:
$$
\tde\neq \tbe+\tga\,'\,\hbox{\ as\ }
\tga\in \la(\hw)\not\ni \tga\,'\,,\
\tbe\neq \tde+\tga\,'\,\hbox{\ as\ }
\tga\in \la(\hw)\ni \tga\,'\,.
$$

When constructing the rank two segment
$L$ from the theorem by induction, we will also assume
that there are no {\em smaller \underline{admissible} triples\,}
$\tbe\,'\,+\al_i=\tga\,'\,$  in $\la(\hw)$ involving $\al_i$
different from $\tbe+\al_i=\tga$. Otherwise,
either we can move $\al_i$
from its first position, which is impossible (see above),
or $L$ is in the very beginning of $\la(\hw)$.
In the latter case, the Coxeter transformation in $L$
moves $\al_i$ from its first position.

Similarly, we will suppose that there are no smaller
admissible triples
$\tbe+\tal\,'\,=\tga\,'\,$  in $\la(\hw)$ involving $\tbe$
different from $\tbe+\al_i=\tga$.
\medskip

\subsection{Special cases}
We will begin with some special cases when the
theorem (and the corresponding algorithm for finding
$L$\~segments) are relatively simple.

{\bf The case of {\mathversion{bold}$\al_i=\al_0.$}} 
The induction argument in this case
requires only property (1) above.
Let $\al_i=\al_0=[-\vth,1]$; recall that $\vth$ is
the longest short root.
If $\tde=[\de,j\nu_\de]\in\la(\hw)$
as $\de>0$ then $\de=[\de,0]\in\la(\hw)$, but the latter
does not contain $\al_0$ in its decomposition. Hence,
$\tde\in \la(\hw)$ are always in the form
$\tde=[-\de,j]$ is for $\de>0$ when  $i=0$.

We set $\tbe=[-\be<0,j\nu_\be]$,
so $\tga=\tbe+\al_0=$ $[-\be-\vth,j\nu_\be+1]$.

In the simply-laced
case, $\be+\vth$ never belongs to $R_+$ due to the maximality
of $\vth$.

If there are two different root lengths, then
such $\tga$ can exist but must be long due to
the definition of $\vth$; therefore,
$\be$ must be short since \, lng $\pm$ lng = lng.
Therefore, it gives that $\tal+\tbe=\tga$ is in the form
\, sht+sht=lng;\,; however, this case was excluded from
the theorem (see ($b$)).
It concludes the (induction step in the) case $i=0$.

\smallskip
{\bf Systems {\mathversion{bold}$A,B,C$ under ($a$)}.}
Let $\tbe,\tal=\al_i,\tga=\al_i+\tbe$ and
the minimality conditions above are imposed.
Then a simple ``numerical" argument proves claims (i,ii)
under assumption ($a$) for the root systems $\widetilde{A},
\widetilde{B},\widetilde{C}$.
It is an instructional example; it also simplifies the
process of finding $L$ in these cases. The {\em general proof}
below will not use this approach (and will include these
special cases). The {\em admissibility conditions} from
(ii), (iii) are not applicable to these cases.

The reasoning below works in some other cases,
for instance, for ($c$). We will include the latter 
in the statement but omit the details.

\begin{lemma}\label{LEMABC}
Let us consider one of the following cases
from (\ref{lngshtsht}) and (\ref{puregtwo}):
\begin{align}\label{abcases}
&(a):\, \widetilde{A},\widetilde{B},\widetilde{C},
\widetilde{G}_2;  \ \,
(c):\, A_2-\hbox{pure-short\ roots\ for } \widetilde{G}_2.
\end{align}
Then (\ref{rankthree}) is not applicable and
one can always find the rank two segment
$L$ from the theorem with $\tu=$id.
\end{lemma}
{\em Proof.}
We may assume that $i>0$ and $|\al_i|\ge |\tbe|$. Thus
$\al_i$ is long  and $\tbe$
is short for $\widetilde{B},\widetilde{C},
\widetilde{F},\widetilde{G}$ in case ($a$).
We may also assume that
$\tga=\tbe+\al_i$ is a unique
admissible triple in $\la(\hw)$
involving $\al_i$ as the beginning and that
it is unique with $\tbe$ as the end.

The proof is mainly based on the minimality
assumption that all roots $\tde\in \la(\hw)$
must contain $\al_i$, i.e., the coefficients 
of the decomposition of $\tde-\al_i$ in terms of simple 
roots $\al_j\,(j\ge 0)$ are all non-negative.
In case ($b$), it is also necessary
to check that (\ref{rankthree}) always holds
under (\ref{abcases}); it is not difficult.
\smallskip

There are two possible subcases concerning $\tbe$.
\smallskip

{\em Subcase $\tbe=[\,-\be<0,\,j\nu_\be>0\,]$}. Then
$\tga=[-\ga,j\,\nu_{\ga}]=[\al_i-\be,j\,\nu_\be]$ for
$\ga=\be-\al_i\in R_+$. Indeed,
$\be$ and $\ga$ are of the same length because
either $\al_i$ is long or $\al_i,\be,\ga$ form a root system
of type $A_2$. The positivity of $\tga$ gives that
$\be-\al_i\in R_+$.

Using that $\tga=\tbe+\al_i$ is assumed to be a
unique such triple in $\la(\hw)$ 
involving $\al_i$,
we obtain that
$\tbe=[-\be,\nu_\be]$ for $\be>0$. I.e., one
can assume that $j=1$. The latter assumption and the
positivity of $\tbe$
implies that the multiplicity of $\al_0=[-\vth,1]$ in $\tbe$ 
is one for short $\be$ and two if it is long.

If $\be$ is short, then
$\vth-\be$ must contain $\al_i$ (in the decomposition
in terms of $\{\al_1,\ldots,\al_n\}$).
If $\be$ is long, then
$2\vth-\be$ must contain $\al_i$.
Combining it with
$\be-\al_i\in R_+$, we conclude that the multiplicity of
$\al_i$ in $\vth$ or $2\vth$ correspondingly must be at
least $2$. This is {\em impossible}. Thus
one can diminish the $\ell$\~length
and proceed by induction in the cases under
consideration. 
\smallskip

{\em Subcase $\tbe=[\,\be>0\,,\,j\nu_\be\ge 0\,]$.} 
Using the minimality, we can assume
that $j=0$ and $\tbe=\be$. Recall that
$\al_i$ is contained in $\be\in R_+$ and
$\tga=\ga=\al_i+\be$. In case ($a$), $\al_i$ is long, $\be$ is
short, and $\ga$ is short. Then the coefficient of $\al_i$
in $\ga$ can not be greater than $1$.  Therefore 
we can diminish the $\ell$\~length between $\tal$ and $\tbe$
and proceed by induction.

Note that formally this argument
can be used in the cases 
\begin{align*}
&\widetilde{B}(\hbox{sht}+\hbox{sht})
\hbox{\ or\ }
\widetilde{C} (\hbox{lng}+\hbox{lng}),
\end{align*}
that where excluded from the theorem since 
such combinations do not lead to $A_2$\~triples.
\sq
\medskip

\rmk
(i) In the cases $\widetilde{A},\widetilde{B},
\widetilde{C},\widetilde{D}$, one can
use the plane interpretation of
the reduced decompositions from \cite{Ch5}; it exists
for $G_2$ too. For instance, it is geometrically obvious 
that the roots 
$$
\tal=\ep_2-\ep_3,\ \tbe=\ep_1+\ep_3,\ \tga=\tal+\tbe=\ep_1+\ep_2,
$$
from (\ref{rankthree}) cannot be collected together.
Algebraically, the $7$\~set there corresponds to the 
reduced nonaffine decomposition
$w=s_2 s_3 s_2 s_1 s_2 s_3 s_2$, where
$s_1$ transposes $\ep_1$ and $\ep_2$, $s_2$ is the transposition
$\ep_2\leftrightarrow\ep_3$, and $s_3$
is the reflection $\ep_3\mapsto -\ep_3$. 
It is analogous for
the $9$\~root set of type $D_4$ from (\ref{rankthreed}).

(ii) Geometrically, one plots $3$ lines in the half-plane
$\{(x,y)\in \R^2,\, y>0\}$ subject to 
the reflection in the $x$\~axis. The intersection and
reflection points give a reduced decomposition
of $w$; the corresponding angles constitute $\la(w)$
(the reflection angle
must be multiplied by $2$ as $m=2$).
Here $\ep_i$ is  interpreted as
the (initial, before the reflection) angle of line $i$. 
They are counted with respect
to the intersection points with $x=b$ from the highest $y$ down
to the $x$\~axis for sufficiently large $b$ and then the
intersection and reflection points are examined as $a<x<b$,
provided that their ultimate transformation from $x=b$ to 
$x=a$ is $w$.

(iii) For $w$ from (i), the first two lines have ``almost"
coinciding nonzero angles $\ep_1,\ep_2$; $\ep_3$
is ``almost" zero and intersects the first two before and
after the reflection points. Generally, one takes $n$ lines
for $B_n,C_n$ and take $w$ for
$(n-1)$ ``almost" parallel non-horizontal lines and
one horizontal line. Then the first and the last angles (and
their sum) form a {\em non-admissible} triple;
the $\ell$\~distance between them cannot be diminished
by Coxeter transformations, which is obvious geometrically.
The theorem claims that this
example is, in a  sense, the only {\em obstacle} for collecting
triples in rank two segments. It is analogous for $D_n$.

(iv) The affine variant of this interpretation requires a portion
of $\R^2$ trapped between {\em two reflection lines\,}; it covers
the root system $C^\vee C$ including affine $B,C,D$.
The affine $A_n$\~system is described
by $n+1$ lines on a {\em cylinder} (with the periodic
$y$\~coordinate). This interpretation is helpful,
but the algebraic
approach is important even for the classical root systems;
the above simple proof for
$\widetilde{A},\widetilde{B},\widetilde{C}$ under ($a$)
is a clear demonstration.
Note that considering triples (there are many similar
situations) is the key for the technique of intertwiners.
\sq
\medskip

\subsection{Uniform construction}
Let  $\tbe,\tal=\al_i,\tga=\al_i+\tbe$ is an
{\em admissible triple}; one of the
conditions ($a$), ($b$), ($c$) must hold. We
impose the minimality conditions (1) and (2) above,  
assuming that $\al_i$
$\,(0\le i\le n)\,$ is a unique simple root in $\la(\hw)$ 

We will also assume later that neither $\tal$
nor $\tbe$
belong to {\em smaller admissible triples\,}.
The condition $|\al_i|\ge |\tbe|$ will be imposed
(unless stated otherwise); it can be
always provided using the inversion
$\hw\mapsto \hw^{-1}$ if necessary. It results in
$\, s_i(\tbe)=\tga,\ s_i(\tga)=\tbe$.
\smallskip

We will examine the \underline{second root}
in $\la(\hw)$ that
is $s_i(\al_j)=\al_j+\mu\al_i$
for $\hw=\cdots s_j s_i$\, as\, $0\le j\neq i$\, 
and $\mu\in \N$\ (if $\mu=0$ then $\al_j\in \la(\hw)$).

Note that $\mu>1$ only in the non-simply-laced case and
if and only if $\al_i$ is short and $\al_j$ is long.
Then $\tbe$ has to be short, due to $|\al_i|\ge |\tbe|$,
that implies in its turn that $\tga$ is short.

Let $\tal\,'\,\equal s_j(\al_i)=\al_i+\mu\,'\,\al_j$,
where $\mu\,'\,>1$ occurs only if $\al_j$ and $\al_i$
have different lengths and $|\al_i|>|\al_j|$.

Note that $\mu\mu\,'\,=m_{ij}=1,2,3$, where $m_{ij}$ is 
the number of
laces between $\al_i$ and $\al_j$ in the affine Dynkin
diagram, and at least one of $\mu,\mu\,'\,$ equals $1$,
namely, $\mu=1$ if $|\al_i|>|\al_j|$,
 $\mu\,'\,=1$ if $|\al_j|>|\al_i|$;  $\mu=1=\mu\,'\,$ if
$|\al_i|=|\al_j|.$

One has $\tbe\,'\,\equal s_j(\tbe)=\tbe+p\al_j$,
where $-\mu\,'\,\le p \le \mu\,'\,$.
If $|\tbe|=|\al_i|$ then $p=-\mu\,'\,,0,\mu\,'\,$; if
$|\tbe|<|\al_i|$ then $|\tbe|\le|\al_j|$ and $p=-1,0,1$.
\smallskip

If $p>0$, then $\tga\,'\,\equal s_j(\tga)=\tga+(\mu\,'\,+p)\al_j$.
If $|\tbe|=|\al_i|$, then  $|\tga|=|\al_i|$  and
$(\mu\,'\,+p)$ must be no greater than $1$ (as for $\al_i$),
which is impossible.
If these lengths are different, then $\tga$ has to be short
since \,sht+lng\, is\,  sht\, and $(\mu\,'\,+p)\ge 1$;
this is impossible too.
\smallskip

\begin{lemma}\label{LEMADMPRIME}
Let $m_{ij}=1$, which implies $\mu=1=\mu\,'\,$. Then
$\tal\,'\,=\al_i+\al_j\in \la(\hw)$ and
$\tbe\,'\,=s_j(\tbe)=\tbe+p\al_j$ for
$p\le 0$. In this case,
the triple $\{\tbe\,'\,,\tga\,'\,=\tal\,'\,+\tbe\,'\,,\tal\,'\,\}$
belongs to $\la(\hw)$. Moreover,
the triple $\{\tbe\,'\,,\tga\,'\,=\tal\,'\,+\tbe\,'\,,\tal\,'\,\}$
is admissible if\, $\{\tbe,\tga=\tal+\tbe,\tal\}$\, is
admissible. 
\end{lemma}
{\em Proof.}
Using $\mu=1=\mu\,'\,$, $\tal\,'\,=s_i(\al_j)\in \la(\hw)$.
Since $p\le 0$ (see above) and $\al_j\not\in\la(\hw)$,
$\tbe\,'\,=\tbe+p\al_j$
and $\tbe\,'\,$ has to be in $\la(\hw)$
(property (2) from the minimality conditions).
Therefore $\tga\,'\,$ belongs to this set too.

Let us suppose that the
triple $\tal\,'\,+\tbe\,'\,=\tga\,'\,$ can be extended to
the set of $7$ roots in the segment $[\tbe\,'\,,\tal\,'\,]$ of type
(\ref{rankthree}) as
$\tal\,'\,=\ep'_2-\ep'_3$, $\tbe\,'\,=\ep'_1+\ep'_3$, 
$\tga\,'\,=\ep'_1+\ep'_2$.

The set of roots $(R^3)'_+\equal\{m\ep'_i,\ep'_i\pm\ep'_j, i<j\}$
is positive in $\tR$ (by construction). Since $\tal\,'\,$ is simple
in $(R^3)'_+$,
the root $\al_j$ does not belong to $(R^3)'_+$; indeed,
$\tal=\tal\,'\,-\al_j>0$ in $\tR$ and therefore in $R^3$.
Therefore the image of  $(R^3)'_+$ under $s_j$ is positive
in $\tR$:
$$
R^3_+\equal\{m\ep_i,\ep_i\pm\ep_j, i<j\}
\subset \tR_+ \for \ep_k=s_j(\ep'_k), k=1,2,3.
$$
The roots 
$$
\tal\,=\,\ep_2-\ep_3,\  \tbe\,=\,\ep_1+\ep_3,\ \tga\,=\,\ep_1+\ep_2
$$
belong to $R^3_+$ by construction;
Lemma \ref{LEMTHREE} describes all possibilities.
Using Lemma \ref{LEM3AB} (and its counterpart for
$\widetilde{D}_n$ and $\widetilde{C}_{n\ge 4}$ with $R^4$
subsystems) we can conclude the proof for
the classical systems. The case of $\widetilde{F}_4$ 
requires a special consideration.
\smallskip
Without using this lemma we can proceed as follows.
\smallskip

We combine $(R^3)'$ and $R^3$ in a root
subsystem $\widehat{R}^4\subset \tR$ of type $B_4$ or $C_4$,
defined as an intersection
of $\tR$ with the $Q$\~span of $R^3$ and $\al_j$.
The latter root does not belong to $R^3$ or $(R^3)'$ and
is simple in $R^4$. The root $\al_j$, the set of roots
\begin{align}\label{r5prime}
&\{\tbe,\tbe\,'\,=\tbe+p\al_j,\tga,\tga\,'\,,\tal\,'\,=\al_i+\al_j,
\tal=\al_i\}\subset \la(\hw),
\end{align}
and also the images of the remaining $4$ roots
from the $7$\~set in $(R_+^3)'$ belong to $\widehat{R}_+^4$.
Generally, the latter system can be an affine extension
of $R^3$ with $\al_j$ being the affine simple root,
but this case can be excluded from the consideration.

Then we prove the lemma for $\widehat{R}^4$, i.e., deduce 
that $\{\tbe\,'\,,\tga\,'\,,\tal\,'\,\}$ is admissible from
the admissibility of $\{\tbe,\tga,\tal\}$.
Here we intersect $\la(\hw)$ with $R^4$; the
intersection is a $\la$\~set in the latter.
This method remains essentially the same in the 
simply-laced case (and for $\widetilde{C}_{n\ge 4}$), 
but the new subsystem will be
$\widehat{R}^5$ of type $D_5$ instead of $R^4$.

The logic of this approach becomes more transparent if
one first proves the Main Theorem for the nonaffine 
root systems of type $B_4$,$C_4$
and $D_5$. Then the initial
admissible triple $\{\tbe,\tga,\tal\}$ can be collected
in a rank two
$L$\~segment in $\widehat{R}^4$ (or $\widehat{R}^5$). 
However, $\al_j\not\in  \la(\hw)$ and
the order of roots in (\ref{r5prime}) must remain
the same if the order of $\{\tbe,\tga,\tal\}$ is
unchanged upon the Coxeter transforms. 
Therefore, in process of collecting 
$\{\tbe,\tga,\tal\}$ in a $L$\~segment of rank two
in $R^4$ (or $R^5$), 
we will automatically make 
$\{\tbe\,'\,,\tga\,'\,,\tal\,'\,\}$ collected together 
(in a $L$\~segment).
This gives the admissibility of 
$\{\tbe\,'\,,\tga\,'\,,\tal\,'\,\}$.
We will omit proving Main Theorem \ref{RANKTWO}
in the non-affine cases $B_4$,$C_4$,$D_5$,
and the details of the consideration of $\widetilde{F}_4$.
\sq
\comment{
\begin{align*}
&\tal\,'\,=s_{j}(\al_i)=
\al_i+\mu\,'\,\al_j, \ \tga\,'\,=s_j(\tga)=
s_j(\al_i+\tbe)=\tga+\mu\,'\,\al_j,\\
&\tga\,''\,=s_is_j(\tga)=\tga+(m_{ij}-1)\al_i+\mu\,'\,\al_j \for
m_{ij}=1,2,3.
\end{align*}}
\medskip

\subsection{\texorpdfstring{The case 
{\mathversion{bold}$m_{ij}>1$}}
{Unequal lengths}}
We will first consider the case when $s_j$ acts trivially
on $\tbe.$

{\bf Case {\mathversion{bold}$p=0$}}. Then $\tbe=\tbe\,'\,$
and $\al_i$ and $\tga$ must be of the same length 
for short $\al_j$ due to
$s_j(\tga-\al_i)=0$ and $(\tga,\al_j)=(\al_i,\al_j)$.
Recall that $s_i(\tde),s_j(\tde),s_is_j(\tde)$ are
positive for any
$\tde\in\la(\hw)$ after $s_j$, i.e., for all roots
excluding the first two roots, $\al_i$ and $\al_j+\mu\al_i$.
Note that the positivity of $s_i(\tde)$ and $s_j(\tde)$
readily follows from the fact that $\al_i\neq\tde\neq\al_j$.

If $m_{ij}=1$ then $\{\tbe,\tga\,'\,,\tal\,'\,\}$ is admissible 
due to Lemma \ref{LEMADMPRIME} and we can change
the position of $\tbe$
in the segment $[\tbe,\tal]\subset \la(\hw)$, which is
impossible due to the minimality conditions. 
Similarly, we can assume that 
$\al_i$ and $\tga$ are the same length for any $\al_j$;
otherwise the admissibility is automatic for 
unequal lengths (see case ($a$)). Thus,
it suffices to consider only the case $m_{ij}>1$, where
$\tR$ can be of types $\widetilde{B},\widetilde{C},
\widetilde{F}_4,\widetilde{G}_2$.
\smallskip

If $m_{ij}>1$ for these root systems, then the roots
$$
\al^1=\tbe=\ep_1-\ep_2,\
\al^2=\al_i=\ep_2-\ep_3,\ \al^3=\al_j=m\ep_3
$$
are simple roots in $\tR^0$, which is defined as
the intersection  of $\tR$ with the $\Z$\~span
of $\al_i,\al_j,\tga$ with the positivity with respect to
$\tR^0_+=\tR_+\cap\tR^0$. It is a {\em finite} root system
of rank $3$ unless in the $G_2$ case (we leave this
case to the reader). The set $\la(\hw)$ does not contain
$m\ep_3$ and therefore must contain $m\ep_2$. It is the
case (2a) of Lemma \ref{LEMTHREE}, the second root is
$m\ep_2$. One can inverse the order in the subset
$\{\ep_1-\ep_2, m\ep_1,\ep_1+\ep_2,m\ep_2\}$ in
$\la(\hw)$ using the induction statement in case ($a$)
of the theorem (here there is no problem with
(\ref{rankthree})).

After moving $m\ep_2$  from its first position inside the segment
$[\ep_1-\ep_2,m\ep_2]\in \la(\hw)$, the root
$\ep_1-\ep_3$ automatically becomes the second in $\la(\hw)$
(right after $\al_i=\ep_2-\ep_3$). Then part (iii) of
Proposition \ref{INTRINLAP}
guarantees that we can make the triple
$\{\ep_1-\ep_2,\ep_1-\ep_3,\ep_2-\ep_3\}$ a connected
sequence in $\la(\hw)$ for a proper reduced decomposition of
$\hw$. Then we can move $\al_i$ from its first position.
However the latter contradicts to the minimality conditions
and therefore concludes the consideration of the case $p=0$.
\smallskip

\rmk
The consideration of $p=0, m_{ij}>1$ for $\widetilde{B},
\widetilde{C},\widetilde{G}_2$
can be concluded following the special cases
$\widetilde{A},\widetilde{B},\widetilde{C}$ under
($a$) considered above.
The roots $\al^1=\tbe,\al^2=\al_i,\al^3=\al_j$
are {\em simple
roots\,} of $\tR^0$, 
the intersection  of $\tR$ with the $\Z$\~span
of $\al_i,\al_j,\tga$ with respect to
$\tR^0_+=\tR_+\cap\tR^0$. It is a {\em finite} root system
of rank $3$ unless in the $G_2$ case (we leave this
case to the reader). The intersection $\la(\hw)\cap \tR$
is a $\la$\~set in $\tR^0_+$.
For the sake of simplicity, let the triple and $\al_j$
be nonaffine. Then $\tga-2\al_i\ge 0$
and 
$s_is_j(\tga)$ 
$=\tga+(m_{ij}-1)\al_i+\mu\,'\,\al_j $ 
must contain $(m_{ij}+1)\al_i$, which is impossible.
Similar approach can be applied to $\widetilde{F}_4$.
\sq
\smallskip

{\em Thus, the case $p<0$ is sufficient to consider.}
\smallskip

To simplify the exposition we exclude the case 
$\{\,\mu\,'\,>1,\,p=-1\,\}$. It can be done since
$|\al_j| < |\al_i|$ in this case,, i.e., 
$\al_j$ is short. Therefore $|\tbe|=|\al_j|$ due to
$\tbe\,'\,=s_j(\tbe)=\tbe+p\al_j=\tbe+\al_j$. This gives
$|\tbe| < |\al_i|$. With this inequality, no admissibility
condition is needed and the proof of the theorem
is straightforward.

\comment{
Without using the inversion, one may argue
as follows (let us skip the case $G_2$).
First,
$$
\{\,\tbe,\, \tbe\,'\,=\tbe-\al_j, \,\tga,\, \tal\,'\,=\al_i+\al_j,\,
\al_i\,\}\,
\subset\,\la(\hw);
$$
here all roots but $\al_i$ are short.
Second, $(\tbe\,'\,,\al_i)=0$
and we can set
$$
\tbe\,'\,=\ep_1-\ep_2,\, \al_j=\ep_2-\ep_3,\, \al_i=2\ep_3
$$
for proper $\ep_1,\ep_2,\ep_3$, i.e., they generate a
root system $\tR^0\subset \tR$ of type $C_3$ and are
simple roots in $\tR^0_+=\tR^0\cap\tR_+$; see\cite{Bo}.
Third, the intersection of $\la(\hw)$ and $\tR^0_+$ is
a $\la$\~sequence and must contain (in the order of their
appearance):
\begin{align*}
\{\,&\tbe=\ep_1-\ep_3,\,\tga+\al_j=\ep_1+\ep_2,\,\tga+\tbe=
2\ep_1,\\
 &\tbe\,'\,=\ep_1-\ep_2,\,\tga=\ep_1+\ep_3,\, \tal\,'\,=
\ep_2+\ep_3,\,
\al_i=2\ep_3\,\}.
\end{align*}
Cf. Lemma \ref{TALTBE} below.
Finally, $\tga+\tbe-3\al_i>0$ in $\tR$ (see above)
and, therefore, in $\tR^0$,
due to the uniqueness of $\al_i$ in $\la(\hw)$, which is
impossible in $C_3$.
\sq
}
\smallskip

Thus, it suffices to assume that $|\tbe| = |\al_i|$
and $p=-\mu\,'\,$. We have:
\begin{align*}
&\tal\,'\,=s_{j}(\al_i)=
\al_i+\mu\,'\,\al_j,\ \, \tbe\,'\,=s_j(\tbe)=\tbe-\mu\,'\,\al_j,\ \,
s_j(\tga)=\tga
\end{align*}
in this case. Since $\tbe\,'\,$ is a positive root,
it must belong to $\la(\hw)$ (the minimality property 
(2) above).

The element $\tal\,'\,$ is a sum of positive roots
and therefore belongs to $\tR_+$.
However, generally, it can be apart from $\la(\hw)$.
The following lemma addresses this possibility.
\smallskip

\begin{lemma}\label{TALTBE}
We assume that  
$\al_i$,$\tga$, $\tbe$ are of equal lengths and
constitute an admissible triple in $\la(\hw)$. The
minimality conditions are imposed, in particular, $\al_i$
is a unique simple root in $\la(\hw)$; recall that 
$\al_i$ and
$s_i(\al_j)=\al_j+\mu\al_i$ are the first and
the second root in this set. The claim is that
$\tal\,'\,=s_j(\al_i)=\al_i+\mu\,'\,\al_j$ belongs to
$\la(\hw)$.
\end{lemma}
{\em Proof.} The inclusion $\tal\,'\,\in \la(\hw)$
can be wrong only when $\mu\neq \mu\,'\,$, i.e., as $m_{ij}>1$
(for $\widetilde{B},\widetilde{C},\widetilde{F}_4,
\widetilde{G}_2$). Also, either $\al_j$ must be long
and $\al_i,\tbe,\tga$ short, or it can be the other way round,
with short $\al_j$ and long $\al_i,\tbe,\tga$.
The latter variant is impossible for $\widetilde{G}_2$, and
both variants are completely analogous for
$\widetilde{B},\widetilde{C},\widetilde{F}_4$.
We will consider here only the first one, assuming that
$$
\mu\,'\,=1=-p,\ \mu>1,\ \tal\,'\,=\al_i+\al_j,\ \,
|\al_j|\,>\, |\al_i|\,=\,|\tbe|\,=\,|\tga|.
$$

Since $l(\hw)=l(\hw s_is_j)+l(s_js_i)$,
the following element must be a positive root in $\tR$:
$$
s_is_j(\tbe)=s_i(\tbe-\al_j)=
(\tbe+\al_i)-(\al_j+\mu\al_i)=\tbe\,'\,-(\mu-1)\al_i.
$$
Therefore $\tbe\,''\equal \tbe-(\al_i+\al_j)$ is a positive
root too; also $(\tbe\,'',\al_j)=0$.

Setting  $\, \al^1=\tbe\,'',\ \al^2=\al_i,\ \al^3=\al_j$,\,
let $\tR^0$ be a root system $\Z$\~generated by these
roots in $\tR$. Then they are simple roots in
$\tR^0_+=\tR^0\cap\tR_+$. The system $\tR^0$
is of type $C_3$ as $\mu=2$. In the case of $G_2$ ($\mu=3$),
one also has that $(\tbe\,'',\al_i)=0$. It implies that
$\tbe\,''$ can not be a (real) root. Therefore $\mu=3$ can
be excluded from further considerations; we will assume
that $\mu=2$.
\smallskip

Since $\la(\hw)\ni\tbe=\tbe\,''+(\al_i+\al_j)$, then
$\tbe\,''=\tbe\,'-\al_i\in \la(\hw)$ if 
$\al_i+\al_j\not\in \la(\hw)$ (i.e., if the claim
of the lemma does not hold). Then we proceed as follows.
\comment{
If $\tbe''=\tbe'-\al_i\in \la(\hw)$,
then $\tbe'=\tbe''+\al_i$ and $\al_i$ belongs to a
{\em smaller admissible triple},
different from the original one $\tga=\tbe+\al_i$.
Indeed, if $\tbe'=\tbe''+\al_i$ can be extended to
a $7$\~set in the form (\ref{rankthree}), the $m\ep_2$
there must coincide with $\al_j+\mu\al_j$, since there
exists a since such
set would belong to the same system $\tR^0$.
}
Setting
$$
\al^1=\tbe\,''=\be-\al_i-\al_j=\ep_1-\ep_2,\,
\al^2=\al_i=\ep_2-\ep_3,\,\al^3=\al_j=2\ep_3,
$$
the intersection $\la(\hw)\cap \tR^0_+$ is
a $\la$\~sequence in $\tR^0_+$, and it is not difficult
to check that it consists of the following roots in the
order of appearance:
\begin{align}\label{tbeep1ep3}
\{\,&\tbe=\ep_1+\ep_3,\,\tbe-\al_i-\al_j=\ep_1-\ep_2,\,
2\tbe-\al_j=2\ep_1,\, \tga=\ep_1+\ep_2,\\
 &\tbe-\al_j=\ep_1-\ep_3,\,\al_j+2\al_i=2\ep_2,\,
 \al_i=\ep_2-\ep_3\,\}.\notag
\end{align}
One can also use claim (ii) of the theorem
{\em for $\tR^0$ instead of $\tR$} here. It gives that 
$\al_i,\tbe$ can be transformed to a rank two segment $L^0$ 
of type $A_2$ inside 
$\la(\hw)\cap \tR^0_+$. Indeed,
$\tbe\,''$ is {\em simple} in $\tR^0_+$ and one can make
it the first root in this sequence instead of $\al_i$.
The set (\ref{tbeep1ep3})
is from case (2a) of Lemma \ref{LEMTHREE} and we can use
Lemma \ref{LEM3AB} to come
to a contradiction with the minimality conditions.

We obtained that $s_j(\al_i)=\al_i+\al_j\in \la(\hw)$
for long $\al_j$ 
and short $\al_i,\tbe,\tga$; 
the case of short $\al_j$ and long 
$\,\al_i,\tbe,\tga\,$ is analogous.
\sq
\medskip

\comment{
To make this final reduction clear, let us recalculated
the roots from $\la(\hw)$ considered above in terms of $\al^r$
in the notation from \cite{Bo} (the system $C_n$):
\begin{align*}
&\al_i=\al^2=\ep_2-\ep_3,
\al_j+2\al_i=2\al^2+\al^3=2\ep_2,
\tbe-\al_j=\al^1+\al^2=\ep_1-\ep_3,\\
&\tga=\al^1+2\al^2+\al^3=\ep_1+\ep_2,
\tbe=\al^1+\al^2+\al^3=\ep_1+\ep_3.
\end{align*}
Here the ordering is from $\la(\hw)$, although
$\al_j+2\al_i$ and $\tbe-\al_j$ can appear in the
opposite order since they are orthogonal to each other.
}

\subsection{Concluding the induction}
Thus, the case  $\mu=\mu\,'\,=1$ is sufficient for
obtaining a transformation to the
rank two $L$ from (i) or (ii,iii) containing
$\al_i,\tbe,\tga$. In this case,
the admissible triples constructed above 
for $\tga$ with the beginning at $\al_i+\al_j$ and 
end at $\tbe-\al_j$ lead to finding a simple
root in $\la(\hw)$ different from $\al_i$.
Let us demonstrate it.

We will omit the case of
$\widetilde{A}_n$ for the sake of uniformity; the theorem
has been already checked for this root system.
Then the
affine Dynkin diagram $\widetilde{\Gamma}$ is a tree.

Thus, $\mu=\mu\,'\,=1$ and
$\al_i+\al_j$ is the {\em second\,} root
after $\al_i$ in $\la(\hw)$. Using the admissible
triple constructed above, the induction hypothesis
makes it possible to move
$\al_i+\al_j$ from its second position in $\la(\hw)$
($\al_i$ remains
untouched). Therefore a simple 
$s_{k}\neq s_j$ with $\al_k$ of the same length as $\al_i$ 
must exist such that it can be made the 
{\em second} after $s_i$ upon a suitable transformation of 
the initial reduced
decomposition of $\hw$ (without touching $s_i$). 
Such $\al_{k}$
must be connected with $\al_i$ by a link in $\widetilde{\Ga}$
and therefore must be orthogonal to $\al_j\,$;
so $s_js_{k}=s_{k}s_j$.

We conclude that all {\em first roots\,} in $\la(\hw')$ 
for $\hw'=\hw s_i$ are pairwise orthogonal and there exists 
a reduced decomposition of $\hw$
in the form $\hw=\cdots s_{k} s_{j} s_i=\cdots s_{j} s_{k} s_i$.
Note that the corresponding $\la$\~set contains the following roots:
$$\tbe,\ldots ,\tbe-\al_j,\ldots,\tbe-\al_k,\ldots,\tga,\ldots,
\al_i+\al_{k},\al_i+\al_{j},\al_i;$$
here we can transpose $\al_i+\al_{k},\al_i+\al_{j}$
and the order of appearance of  $\tbe-\al_j,\tbe-\al_k$ can
be different.
\smallskip

Let us exclude for a while the cases $\widetilde{D},
\widetilde{E},\widetilde{B}$ (when $\widetilde{\Ga}$ 
is not a segment). Then
the number of {\em first roots\,} in $\la(\hw')$, is exactly $2$,
namely, they are $\al_j$ and $\al_k$.

The $\la$\~set of $\hw''=\hw s_i s_k$
contains an admissible triple that begins with $\al_j\,$;
explicitly, it is
$$\{s_k s_i(\tbe-\al_{j})=s_k(\tbe-\al_{j})=
\tbe-\al_j-\al_{k},\ldots, s_k s_i(\tga)=\tbe-\al_k,\ldots,
\al_j\}.
$$
Therefore it must have at least one {\em first
root} $\al_l\neq\al_j$. If $|\al_l|=|\al_i|$
then this root cannot be $\al_i$,
since otherwise $s_i s_k s_i=s_k s_i s_k$ and it would make $\al_k$
another {\em first root} of $\hw$. So it must be connected
by a link with $\al_k$ in $\widetilde{\Gamma}$. Indeed, otherwise
it would be the third {\em first root} in $\la(\hw')$, which
is impossible. This means that it cannot be connected
with $\al_j$ ($\widetilde{\Gamma}$ is a tree) and we have again 
the orthogonality of all {\em first roots\,} in $\la(\hw'')$.
The connection picture is 
$ \al_l \longleftrightarrow \al_k \longleftrightarrow \al_i
\longleftrightarrow \al_j$. 

\smallskip
We continue this process for $\hw'''=\hw s_i s_k s_l$;
the reflection $s_j$ can be made again 
the beginning of its reduced decomposition, respectively,
$\al_j$ can be assumed to be the first in 
$\la(\hw''')$. The latter sequence contains
the triple with $\al_j$ as its first endpoint; this triple
is the image of $\{\tbe-\al_j,\tga,\al_i+\al_j\}$
under $s_l s_k s_i$. Once again, another {\em first root}
$\al_m\neq\al_j$
must exist for $\hw'''$. Continuing to assume that
$|\al_l|=|\al_i|$, it cannot be $\al_i$ (otherwise we
could move $s_i$ through $s_l$ and use the Coxeter relation
for $s_i s_k s_i$) and it must be connected with $\al_l$ or
with $\al_k$. Hence, $\al_m$  cannot be connected with $\al_j$,
so it is orthogonal to $\al_j$. We can represent $\hw'''=
\cdots s_j s_m $ and then continue with $\hw''''=\hw'''s_m$.
\smallskip

Eventually, we come to an element $\hw s_i s_k s_l s_m\cdots$ of
length $3$ and with
the $\la$\~set that is a pure triple with $\al_j$ as the
beginning. However such a set does not have pairwise orthogonal
{\em first roots\,}. This contradiction proves
the existence of the rank two $L$ from (i,ii) of the theorem
(provided the admissibility). 
\smallskip

Recall, that it was proven 
for $\widetilde{\Ga}$ that is a segment and 
under the assumption
that we can always find the roots $\al_l,\al_m,\ldots$ of 
the same length as that for $\al_i$ (and in the
non-simply-laced case).

If a root in the latter sequence of {\em first roots}
has the length different from $|\al_i|$ 
(it can be only an endpoint of $\widehat{\Ga}$
unless in the $\widetilde{F}_4$\~case), then we use  
$k$ instead of  $j$ and add  
a simple root $\al_{l'}$ connected with
$\al_j$; then $\al_{m'}$ is connected with $\al_{l'}$
and so on. All their lengths will be coinciding
with $|\al_i|$ when moving in this direction 
(there is only one double link in $\widehat{\Ga}$).

A similar reasoning can be applied in the simply-laced case
when $\widetilde{\Ga}$ is not a segment (then one can, generally,
proceed using three directions). 

To complete the proof of the theorem we need to check that if
$\tal$ is before $\tbe$ in $\la(\hw)$ and this set
contains exactly the $7$ roots listed in (\ref{rankthree}) upon
the intersection with the $\Z$\~span of these seven roots, then

\noindent
1) all these $7$ roots belong to the segment  $[\tbe,\tal]$,

\noindent
2) their appearance  in $[\tbe,\tal]$ is as in (\ref{rankthree}).

The latter holds modulo the Coxeter transformations inside this set.
It readily results from Lemma \ref{LEMTHREE}, case (3).
The analogous claim for (\ref{rankthreed}) is equally simple. 
\sq
\smallskip

The following corollary, a variant of Lemma \ref{LEMREFLEC},
demonstrates how the theorem can be used in an important
particular case of reflections.

\begin{corollary}\label{CORTALPRIME}
Given $\tga=[\ga,\nu_\ga j]\in\tR_+$, let 
$\tal\in \la(s_{\tga})$ and
$\nu_\al=\nu_\ga$. Then $\tbe=-s_{\tga}(\tal)=\tga-\tal$
belongs to $\la(s_{\tga})$ 
due to (\ref{talinla}). 

(i) Given a reduced decomposition of $s_{\tga}$,
the triple  $\{\tbe,\tga,\tal\}$ for simple $\tal$
can be made consecutive in $\la(s_{\tga})$
if and only if $\tga\neq \th\,'\,$
in any root subsystem $R'\subset \tR$ of type
$B_3$,$C_3$ or $D_4$ containing this triple.
Here $\th\,'\,$ is the maximal positive root in 
$R'_+=R'\cap \tR_+$ unless for $C_3$; in the latter
case $\th\,'\,$ is
the maximal \underline{short} positive root.

(ii) Assuming that $\tal$ satisfies (i), let
$\tal\,'\,\in \la(s_{\tga})$ be a root between $\tal$ 
and $\tga$ such that $(\tal,\tal\,'\,)=0$
(see Lemma \ref{LEMREFLEC}). Then a pair
$\{\tal\,'\,,\tal\}$ can be made consecutive before
$\tga$ using rank two Coxeter transforms in 
$\la(s_{\tga})$.
These transforms are either in the form
$s_is_js_i\mapsto s_js_is_j$ with
the midpoints corresponding to 
$\tga\,'=[\ga,\nu_\ga j\,'\,]\in \la(s_{\tga})$ or 
in the form $s_is_j\mapsto s_js_i$ 
otherwise. 

(iii) Continuing, $\la(s_{\tga})$ can be transformed
to a sequence with $\tal$ as the beginning and such that the
only roots before $\tga$ non-orthogonal to $\tal$ in  $\la(s_{\tga})$
are $\tga\,'$ from (ii) for $j\,'\,<j$ and also 
$\tbe\,'\,=-s_{\tga\,'}(\tal)=\tga\,'-\tal\in \la(s_{\tga})$.
\end{corollary}  
{\em Proof.} 
For reflections $s_{\tga}$, the non-admissibility is
exactly as in (i). Indeed, we do not need to check
that $\tal$ is a unique simple root from
$R'_+=R'\cap \tR_+$ in $\la(s_{\tga})\cap R'$,
since all simple roots in $R'$
but $\tal$ are orthogonal to such $\tga$ .
For instance, the roots $\ep_1-\ep_2$ and $m\ep_3$ are
orthogonal to $\tga=\ep_1+\ep_2$ and therefore 
cannot belong to $\la(s_{\tga})$ in the case of
(\ref{rankthree}). 
Then claims (i,ii,iii) follow from the Main Theorem
and Lemma \ref{LEMREFLEC}). In (iii), we move $\tal$
to the position next to $\tga$; all roots it
``passes" have to be orthogonal to $\tal$ unless
they are in the form $\tga\,',\tbe\,'\,$. Then we
can move $\tal$ back to its first position.
\sq
\medskip

\section{Right Bruhat ordering}\label{sec:RightBruhat}
We will define the Bruhat ordering
on $\hW$ for the system $\tR$ relative to its  
root subsystem. The notations are from the previous
sections with a reservation about $\tR^0$.
\smallskip

\subsection{Basic properties}\label{SUBSECTR0}
By a {\dfont root subsystem} $\tR^0$, we will
mean the intersection of $\tR$ with a $\Z$\~lattice, namely,
$$
\tR^0=\tR\cap \tilde{\La}^0,\ \tilde{\La}^0\,\subset\,
\tilde{Q}\equal \sum_{\tal\in \tR}\Z\tal=[Q,\Z]\,;
$$
it is obviously a reduced root system in its own right.
{\em This meaning of $\tR^0$ will be fixed till the
end of the paper}.
According to \cite{Bo}, $\tR^0$ is a {\em subsystem} if and only
if it is {\em symmetric} ($-\tR^0=\tR^0$) and {\em closed},
i.e.,
$$
\tal+\tbe\in \tR\Rightarrow \tal+\tbe\in \tR^0 \for
\tal,\tbe\in \tR^0.
$$
For instance, all long roots always form a subsystem.
We refer to \cite{Bo} for the description of the
maximal (proper) subsystems in the nonaffine case.

A subsystem may be reducible, then it is a direct sum
of pairwise orthogonal
irreducible reduced root systems.
Almost all definitions and constructions for irreducible
systems can be extended to the reducible case;
we use them for reducible systems without comments.

Note that {\em not} all root systems that belong to $\tR$ are
root subsystems in the above sense. For instance,
in the notation from \cite{Bo},
the set $B_2^{\sht}$ of short roots in $B_2$, 
$$
B_2^{\sht}=\pm\{\al_1+\al_2,\ \al_2\}\ \subset\
\pm\{\al_1,\ \al_1+\al_2,\ \al_1+2\al_2,\ \al_2\},
$$
does not correspond to any $\La^0$. Indeed, 
otherwise it would contain $\al_1$, that is long.
\smallskip

The subset $B_2^{\lng}=\pm\{\al_1,\al_1+2\al_2\}$
of long roots in $B_2$ is
the intersection of $B_2$ with the lattice
$\La^0=$ $\{v\in Q, \, (v,\al_1)\in 4\Z\}$;
the normalization is as above: 
$(\al,\al)=2$ for {\em short} $\al$.

It is also a simplest example of a {\em root subsystem}
that is not an intersection of $B_2$ with any linear
$\Q$\~subspace in the $\Q$\~span of $Q$. Note that $\tR^0$
is an intersection of $\tR$
with a $\Q$\~linear subspace in $\Q^{n+1}$ if and only if
$\La\,'\,\equal\tilde{Q}/\tilde{\La}^0$ has no torsion.

The affine variant of this example is
the set $\tilde{B}_2^{\lng}$
of {\em affine} long roots in $[B_2,\Z]$; it is a root subsystem;
associated with $\tilde{\La}^0=$ $[\La^0,2\Z]$.
\smallskip

Thus, our {\em root subsystems\,}
are closed with respect to the integral
linear combinations whenever the results are roots, but
$\Q$\~linear combinations are not allowed here;
{\em this class is between
a wider class of all subsets that are root systems in
their own right and a narrower class of the intersections with
linear $\Q$\~subspaces}.
\medskip

{\bf Relations between {\mathversion{bold}$\tR$ and $\tR^0$}}. 
We are going to supply $\tR^0$ with the {\em induced}
systems of positive and simple roots and establish connections
of the $\la$\~sets for the corresponding
Weyl groups.

The subset of positive roots of $\tR^0$
is defined as $\tR_+^0=\tR_+\cap\tR^0$;
its minimal elements form simple roots $\{\al_j^0\}$
in $\tR_+^0$ and $\Z$\~generate the corresponding root lattice
$\tilde{Q}^0$. Note that the
intersection
$Q_+\cap \tilde{Q}^0$ is generally greater than $\tilde{Q}_+^0$
defined as the $\Z_+$\~span of the simple roots in $\tR_+^0$;
an example is $\tR^0=B_2^{\lng}$. The coincidence holds, for
instance, if every $\al_j^0$ contains a certain $\al_{j\,'\,}$
for $0\le j\,'\,\le n$  in its decomposition
such that distinct $j$ have distinct $j\,'\,$.
\smallskip

An important property of {\em root subsystems\,} is the compatibility
with taking the $\la$\~sets of the elements from
$\tW^0$ defined for  $\tR^0$ and those for
the main root system $\tR$.
It is based on the fact that intersections of
$\la$\~sets in $\tR$ with $\tR^0$ are $\la$\~sets with respect
to $\tR^0_+=\tR^0\cap \tR_+$ and the corresponding simple roots
in $\tR^0_+$. See Theorem \ref{INTRINLA}.

\begin{proposition} \label{BRUHATLA}
Let $\tR^0$ be a root subsystem of $\tR$ with a natural induced
subset $\tR_+^0=\tR_+\cap\tR^0$ of positive roots and the
corresponding set of simple roots in  $\tR_+^0$,
$\tW^0$ the Weyl group of $\tR^0$
generated by the reflections $s_{\tal}$ for $\tal\in \tR^0$ (the
simple reflection are sufficient), $\la^0 (\tu)$ the $\la$\~sets
defined within $\tR^0$ for $\tu\in \tW^0$.

Given $\hw\in \hW$, there exists a unique element $\hw|_0\in \tW^0$
such that $\la^0(\hw|_0)=\la(\hw)\cap \tR^0$. Here the
ordering of the roots is induced from that in $\la(\hw)$.
Explicitly, setting $g=|\la(\hw)\cap \tR^0|$,
\begin{align}\label{hw0s}
&\hw|_0=s^{p_1}\cdots s^{p_g}, \hbox{\ where\ }
\la(\hw)\cap \tR^0=\{\tal^{p_g},\ldots, \tal^{p_1}\},\, 
s^p=s_{\tal^p},\\
& \hw|_0\,=\,\cdots\,(s^{p_1}s^{p_2}s^{p_3}s^{p_2}s^{p_1})\,
(s^{p_1}s^{p_2}s^{p_1})\,(s^{p_1}) \hbox{\ is\ reduced\ 
in\ } \tW^0,\notag
\end{align}
where the elements in $(\,\cdot\,)$ are \underline{simple}
reflections in $\tW^0$. 
\end{proposition}
{\em Proof}. The conditions ($a,b,c,d$) from
Theorem \ref{INTRINLA} are compatible with the intersections
with $\tR^0$. Concerning ($c,d$), we use the following:
if the difference
of two roots from $\tR^0$ is a root in $\tR$, then it belongs to
$\tR^0$. As for (\ref{hw0s}), formula (\ref{lambdainv}) is applied.
\sq
\smallskip

Let $\tR^0\subset \tR$ be a {\em root subsystem}
of $\tR$. Given $\hw\in\hW$
and its reduced decomposition $\hw=\pi_rs_{i_l}\cdots s_{i_i}$,
the {\dfont right Bruhat set}
$\b^0(\hw)\subset \hW$ with respect to $\tR^0$ is formed by
the products $\hw'$ obtained from the decomposition of $\hw$
by striking out any  number of simple {\dfont right singular
reflections} $s_{i_p}$ satisfying by definition
$\tal^p\in \tR^0$ in the notation 
$\tal^p=s_{i_1}\cdots s_{i_{p-1}}(\al_p)$ from the previous
section (used in Proposition \ref{BRUHATLA}).

If the number of removed singular reflections is nonzero
then the notation will be $\b^0_o(\hw)$; thus
$\b^0(\hw)=\b^0_o(\hw)\cup \hw$.
\smallskip

\rmk
The {\em left Bruhat ordering}
with respect to $\tR^0$ is the $\hw\mapsto \hw^{-1}$ image
of the right Bruhat ordering. Explicitly, it is defined by
erasing (some of) $s_{i_p}$ satisfying the conditions
$\pi_rs_{i_l}\cdots s_{i_{p+1}}(\al_p)\in \tR^0$ for
reduced $\hw=\pi_rs_{i_l}\cdots s_{i_{1}}$.
Such simple reflections can be called {\em left singular}.

The left ordering is different from
the right Bruhat ordering unless $\tR^0$ is $\hW$\~invariant.
They obviously coincide for the usual Bruhat ordering on $\hW$,
when $\tR^0=\tR$. The example of the root subsystem
$\tilde{B}_2^{\lng}$ considered above satisfies the
$\hW$\~invariance condition and demonstrates that
the coincidence of the right and left Bruhat orderings
may occur for some non-trivial $\tR^0$. Obviously,
$\hW$\~invariant nonzero $\tR^0$ must be of maximal rank
due to the irreducibility of the action of $\hW$ in $\R^{n+1}$.
It is not difficult to describe all such cases, but we do
not need it here.

{\em Only the right Bruhat ordering will be used in this paper;
the word ``right" will be mainly omitted.}
\sq
\medskip

\subsection{Using Coxeter transforms}
The next theorem shows that the properties of the
right Bruhat ordering are similar to those of the usual one,
with a reservation that it is not invariant with respect
to the inversion $\hw\mapsto \hw^{-1}$.
\smallskip

Recall that the group generated
by all $s_{\tal}$ for $\tal\in \tR^0$ is denoted by $\tW^0$.
We put
\begin{align}\label{rightbruhat}
&\hw\,\ge_0\, \hw' \hbox{\ \ if\ \ } \hw'\in \b^0(\hw),
\hbox{\ \ respectively,\ \ }\\
&\hw\, >_0 \, \hw' \hbox{\ \ if\ \ } \hw'\in \b^0_o(\hw),
\hbox{\ i.e.,\ if\ } \hw'\neq \hw.\notag
\end{align}
\smallskip

\begin{theorem} \label{STRIKOUT}
Given a root subsystem $\tR^0\subset \tR$
and a reduced decomposition $\hw=\pi_rs_{i_l}\cdots s_{i_1}$
of $\hw\in \hW$,

(a) $\b^0(\hw)$ does not depend on the choice of the
reduced decomposition of $\hw$;

(b) $\b^0(\hw')\subset\b^0(\hw)$\  if\   $\hw'\in \b^0(\hw)$,\,
i.e., if the ordering $\,\ge_0\,$ is transitive;

(c) $\b^0(\hw')\subset\b^0(\hw s_{\tal^p})$, where $s_{i_p}$ is
the first (on the right) in the set of singular reflections deleted
when constructing $\hw'\in \b^0(\hw)$;

(d) $\cap_{\hw'}\b^0(\hw')=\{\hw^\circ\}$, where
$\hw^\circ$ is obtained from $\hw$ by crossing out all
singular $s_{i_p}$;

(e) $\hw^\circ$ is a unique element of minimal length in
the coset $\hw\tW^0$ for $\tW^0=$
$\lan s_{\tal}\,\mid\,\tal\in \tR^0\ran$;

(f) if\  $l(s_j \hw)<l(\hw)$\  for $0\le j\le n$
and also $\hw^{-1}(\al_j)\in \tR^0$, then
$s_j\hw\in \b^0_o(\hw)$.
\end{theorem}

{\em Proof.} Setting $s^p=s_{\tal^p},$ the
element $\hw'\in\b(\hw)$ obtained by striking out $g>0$
simple reflections $s_{i_p}$ for the indices
$p$ from the sequence
$\{p_g>,\ldots,>p_1\}$ equals $\hw'=\hw s^{p_g}\cdots s^{p_1}$.
This representation of $\hw'$ will be used constantly.
\smallskip

We will extend the definition of
$\b^0(\hw)$ to {\em possibly non-reduced},
decompositions of $\hw.$
The notation $\tilde{\b}^0(\hw)$ will be used; in this 
definition the {\em initial decomposition of $\hw$ must be
given}; this set may depend
on its choice. The decompositions,
of $\hw'$, possible non-reduced,
that are obtained from a given decomposition
of $\hw$ by striking out (some of) the corresponding 
singular simple reflections will be called {\em standard};
the elements from $\tilde{\b}^0(\hw)$ will be 
considered with the corresponding standard decompositions
(unless stated otherwise).  
We will generally distinguish the elements from
$\tilde{\b}^0(\hw)$ with coinciding
$\hw'$  if their {\em standard decompositions} are different.

If the decomposition of $\hw$ is fixed, then
given a standard decomposition of $\hw'\in \tilde{\b}^0_o$, its
$\tla$\~{\em sequence}, defined by formula (\ref{tal}), is
\begin{equation}\label{tlamset}
\tilde{\la}(\hw')= \{(s^{p_1}\cdots s^{p_g})(\tal^l),\ldots,
(s^{p_1}\cdots s^{p_h})(\tal^q),\ldots,\tal^1=\al_{i_1}\},
\end{equation}
where  $p_h$ is
the last $p$\~index such that $s_{i_p}$ is removed
before $s_{i_q}$. To be more exact, let $q$ run through
the set of all indices $\{l,\ldots,1\}$ apart from the
$p$\~sequence, then 
$p_h$ for $h=h(q)$ in (\ref{tlamset})
is the greatest index $p$ smaller than $q$.

We see that an arbitrary $\tbe\in\tla(\hw')$ is naturally
represented in the form
$\tbe=s^{p_1}\cdots s^{p_h}(\tal^q)$. Moreover,
$\tbe$ belongs to $\tR^0$ if and only if  $\tal^q\in \tR^0$;
$\tal^q$ is defined for the initial decomposition of $\hw$.
Indeed, the reflections $s^p$ and their products preserve $\tR^0$
for all $p\in \{p_g,\cdots,p_1\}$ by construction.


We see that, given an initial
decomposition of $\hw$, iterations (compositions) of
the deleting procedure in the
class of standard decompositions
do not give anything new; deleting singular 
simple reflections in
the standard decompositions (non-necessarily reduced)
of $\hw'\in \tilde{\b}^0(\hw)$ leads to $\hw''$ 
(with the resulting standard decompositions) from the 
same set $\tilde{\b}^0(\hw)$.

Such transitivity is based on the fact that 
singular simple reflections {\em remain singular} in the
corresponding $\hw'$ due to (\ref{tlamset}).
\smallskip

The following lemma extends ($a$) from
Theorem \ref{STRIKOUT} to the class of non-reduced
decompositions.

\begin{lemma}\label{LEMHCOX}
Homogeneous Coxeter transforms of
a given, possibly non-reduced, decomposition
of $\hw$ do not change the set $\tilde{\b}^0(\hw)$
considered simply as a set of elements in $\hW$ (i.e., without
the corresponding standard decompositions).
\end{lemma}
{\em Proof.}
From the view point of the
sequence $\tla(\hw)$ from (\ref{tal}), given a Coxeter
transform,
the corresponding consecutive roots $\tal^p$
constitute all positive roots
of a root system of rank two; this transform will
permute them changing their order in $\tla(\hw)$
to the inverse one.

The set $\tilde{\b}^0(\hw)$
obvious remains unchanged if no singular $\tal^p\, (\in \tR^0)\,$
are involved in this permutation. Otherwise,
the following configurations
of simple singular reflections may occur.
\smallskip

{\em First,} only {\em one} singular $\tal^p$ can be involved. We
examine removing the
corresponding simple singular reflection before and after the
Coxeter transformation of type $A_1\times A_1$, $A_2,B_2$, or $G_2$.

In the $A_2$\~case, the singular element can be
$\underline{s_i}s_{i+1}s_i$, $s_i\underline{s_{i+1}}s_i$,
or $s_is_{i+1}\underline{s_i}$; it is underlined.
The transformation $s_is_{i+1}s_i\mapsto s_{i+1}s_is_{i+1}$ in $\hw$
becomes respectively
$$
s_is_{i+1}\mapsto s_is_{i+1},\
s_{i+1}s_i\mapsto s_{i+1}s_i,\hbox{\ or\ } s_i^2\mapsto s_{i+1}^2.
$$
This is ``inside" $s_is_{i+1}s_i$ in
the decomposition of $\hw$.
Obviously the corresponding three $\hw'$ remain
unchanged. The other rank two cases are equally simple.

{\em Second,} at least {\em two} singular
simple reflections
can be involved in the transformation.
In the case of the Coxeter
relation of type $A_2$, {\em all} simple reflections have to be
singular since $\tR^0$ is stable under
addition and subtraction.
Therefore the right Bruhat sets for
$s_i s_{i+1}s_i$ and $s_{i+1}s_is_{i+1}$ are the usual Bruhat sets
(the whole $\bS_3$ generated by $s_i, s_{i+1}$ inside the 
decomposition of $\hw$) and obviously coincide.

Similarly, the reference to the standard Bruhat ordering
is sufficient in the remaining rank two cases, $A_1\times A_1$,
$B_2$ or $G_2$ provided that {\em all} involved
$\tal$ are from $\tR^0$.

{\em Third,} for $B_2,$ the configuration $\underline{s_i}s_{i+1}
\underline{s_i}s_{i+1}$ with exactly two singular reflections
may occur; cf. the example of $B_2^{\lng}$
considered above. It becomes
$s_{i+1}\underline{s_i}s_{i+1}
\underline{s_i}$ upon the transformation
$$
s_is_{i+1}s_is_{i+1}\mapsto
s_{i+1}s_is_{i+1}s_i.
$$
The elements $s_i$
($0,1$ or $2$ of them) are allowed to be removed from
the product $s_is_{i+1}s_is_{i+1}$ and its inverse.
The corresponding
right Bruhat sets (inside $\hw$) obviously coincide:
$$
\{s_is_{i+1}s_is_{i+1},\  s_is_{i+1}^2,\ s_{i+1}s_is_{i+1},\
\hbox{id}\}.
$$
Here $s_{i+1}$ may be underlined
instead of $\underline{s_i}$; it is completely analogous.

We can always treat the deleting procedure as right
multiplication by proper reflections; then the corresponding
elements $\hw'$ are
$$
\{\hw,\ \hw s_{\tal},\ \hw s_{\tbe},\ \hw s_{\tal}s_{\tbe}\},
$$
where $\tal,\tbe$ are orthogonal to each other.
The compatibility with the $B_2$\~transformation
simply means that this
$4$\~set remains unchanged if $s_{\tal}$ and $s_{\tbe}$
are transposed,
which is obvious because these reflection commute due to the
orthogonality of $\tal,\tbe$.
The same argument works for $G_2$ as exactly two
simple reflections are singular.
\smallskip

The last, the {\em fourth,} case is $G_2$ with three
singular reflections. They come from the $A_2$\~subsets of long or
short roots in the $G_2$\~system. It suffices  to check that
deleting $0,1$, $2$ or all three reflections
in the products $s_{\tal}s_{\tal+\tbe}s_{\tbe}$
and, respectively, in 
its inverse $s_{\tbe}s_{\tal+\tbe}s_{\tal}$
result in coinciding sets,
provided that $\tal+\tbe$ is a root.
It formally follows from the $A_2$\~consideration.

Directly, in terms of the simple reflections, we have either
$$
\underline{s_1}s_{2}\underline{s_1}s_{2}\underline{s_1}s_{2}
\ \ \mapsto\ \
s_2\underline{s_1}s_{2}\underline{s_1}s_{2}\underline{s_1}
$$
or the same products where $s_2$ is underlined
instead of $s_1$. Actually,
the latter corresponds to a subset of
all short roots in $G_2$ and $\Z$\~generates the whole system;
so it is not a {\em root subsystem} in our sense and this
case can be omitted.
The underlined reflections ($0,1,2$ or all of them)
are allowed to be removed from the products;
the right Bruhat sets coincide.

Similarly, the configuration for $G_2$ with two simple
singular reflections (the third case) is
$$
s_1s_{2}\underline{s_1}s_{2}s_1\underline{s_{2}} \hbox{\ \ or\ \ }
\underline{s_1}s_{2}s_1 \underline{s_{2}}s_1s_{2}.
$$
It leads to another, direct, justification of the coincidence 
for such configuration.
\sq
\smallskip

Lemma \ref{LEMHCOX}
obviously implies ($a$). Let us demonstrate that
the transitivity from ($b$) is a straightforward corollary
of this lemma too.

Generally, an arbitrary
reduced decompositions of $\hw'\in \b^0(\hw)$
can be obtained from a standard decomposition, maybe non-reduced,
of $\hw'$ using homogeneous Coxeter transformations,
removing the squares $s_i^2$, if any, then
applying homogeneous Coxeter transformations again, 
removing $s_i^2$ (if occur) and so on.

Respectively, the sequence
$\la(\hw')$ will be obtained from $\tla(\hw')$ by inverting
the order of roots in subsequences of type $A_1\times A_1$, $A_2$,
$B_2$ or $G_2$, respectively with $2,3,4,6$ neighboring roots,
and also by deleting neighboring pairs
$\{\tal,-\tal\}\subset\tla(\hw')$, if any.
Such deleting may diminish
$\tilde{\b^0}(\hw')$ if $s_i$ is singular;
the former procedure does not change this set as
has been already checked.

Thus $\b^0(\hw')\subset \tilde{\b}^0(\hw') \subset \b^0(\hw)$, 
which
gives ($b$). Claim ($c$) results from ($b$) by induction
with respect to the length $l=l(\hw)$. Indeed, 
if $s_{i_1}$ is removed when constructing $\hw'$, then
$\hw'\in \b^0(\pi_rs_{i_l}\cdots s_{i_2})$ because the
decomposition  $\pi_rs_{i_l}\cdots s_{i_2}$ remains reduced.
Otherwise, $\hw'\in \b^0(\pi_rs_{i_l}\cdots s_{i_2})s_{i_1}$.
In any case, we can proceed by induction.
\smallskip

Justification of ($d$) is similar to that of ($b$).
First, deleting {\em all} singular
simple reflections from a {\em standard} decomposition, possibly
non-reduced, of an arbitrary $\hw'\in \b^0(\hw)$
results in the same element
$\hw^\circ$; recall that the latter
is obtained from the initial reduced decomposition
of $\hw$ by deleting all singular simple reflections 
{\em at once}.
Second, the construction of $\hw^\circ$ is compatible
with the homogeneous Coxeter relations, which can be checked
following
Lemma \ref{LEMHCOX}. Third, this construction is compatible 
with removing the squares $s_i^2$ from the decompositions. 
It proves ($d$).
\smallskip

The characterization of $\hw^\circ$ from ($e$) is verified as follows.
The set $\la(\hw^\circ)$ contains no roots
from $\tR^0$,  because so does the set $\tla(\hw^\circ)$
constructed for a standard decomposition. Let us check
that this is a
defining property of $\hw^\circ$.

\begin{lemma}\label{UNIQWO}
There exists a unique element $\hw^*$ in the coset $\hw\tW^0$
such that $\la(\hw^*)\cap \tR^0=\emptyset$; it has minimal
possible length in this coset and $\hw^*=\hw^\circ$.
\end{lemma}
{\em Proof.}
If there are two such elements $\hw^\circ, \hw^*$, then
$\hw^*=\hw^\circ \hu$ for $\hu\in \tW^0$. One has
$\la(\hw^*)\subset$ $\hu^{-1}(\la(\hw^\circ))\cup \la(\hu)$ and
the former set is obtained from the latter by
removing the pairs $-\tal,\tal$. If $\hu\neq$id, then
$\la(\hu)$ contains at least one $\tal\in \tR^0$.
Such $\tal$ must disappear in $\la(\hw^*)$, which
is possible only if
$\hu^{-1}(\la(\hw^\circ))$ contains $-\tal$. However,
$\hu$ belongs to $\tW^0$ and
$\hu^{-1}(\la(\hw^\circ))\cap \tR^0=\emptyset$; we come to
a contradiction.

Elements  $\hw'$ of minimal length in $\tW^0$ cannot
have $\tal\in \tR^0$ in their $\la$\~sets $\la(\hw')$
because otherwise $l(\hw's_{\tal})<l(\hw')$ would hold.
Therefore, $\hw^\circ$ is of minimal length and unique
with this property.

Note that the existence of a unique element of minimal length
is essentially Corollary 3.4 from \cite{Dy},
well known for parabolic subgroups of Coxeter groups.
\sq
\smallskip

Let us check ($f$). Given a reduced decomposition
$\hw=\pi_rs_{i_l}\cdots s_{i_{1}}$,
$$
\al_j=\pi_rs_{i_l}\cdots s_{i_{p+1}}(\al_{i_p})\,=\,
\hw s_{i_1}\cdots s_{i_{p}}(\al_{i_p})
$$
for certain $p$. Then the element $s_j\hw$ is obtained from
$\hw$ by striking out $s_{i_p}$. Moreover,
$$
\tal^p=s_{i_1}\cdots s_{i_{p-1}}(\al_{i_p})=
s_{i_1}\cdots s_{i_{p-1}}s_{i_p}(-\al_{i_p})=-\hw^{-1}(\al_j)
\in \tR^0
$$
and $s_{i_p}$ has to be a singular simple reflection.
It completes the proof of the theorem.
\sq
\medskip

\subsection{Bruhat ordering \texorpdfstring{on 
{\mathversion{bold}$\tW^0$}}
{\em on tilde-W-0}}
The right Bruhat ordering can be
reduced to the standard Bruhat ordering in $\tW^0$,
that is established in the next proposition.
It readily gives another proof of the description
of the minimal elements  $\{ \hw^\circ\}$ from
Lemma \ref{UNIQWO}. This proposition gives a
simpler approach to the construction of $\hw^\circ$,
although the considerations based directly
on the Coxeter relations have their own advantages,
especially when the Bruhat construction
(deleting simple reflections) is extended to Hecke algebras.

For instance, it is important to determine
how the choice of a reduced decomposition of
$\hw\in \hW$ influences deleting (singular) $T_i$ in the
corresponding product $T_{\hw}$; here one can follow 
Lemma \ref{UNIQWO}. We note that {\em compatible
pairs of $R$\~matrices} to be discussed below are 
a certain formalization of this ``direct" approach.

Given $\hw=\pi_rs_{i_l}\cdots s_{i_{1}}$, let
$p\in \{p_g>p_{g-1}>,\ldots,p_j,\ldots,>p_{2}>p_1\}$
for the sequence of all singular
$1\le p_j\le l$ as $1\le j\le g$. We use the notation
\begin{align}\label{lambrelset}
&\tbe^j=\tal^{\,i_p}=\hw^{p}(\al_{i_p})
\hbox{\ for\ } p=p_j,\hbox{\ where\ } 
\hw^p\equal s_{i_1}\cdots s_{i_{p-1}},
\end{align}
and also will need $\hw_j$ obtained from $\hw^{p_j}$ by deleting
{\em all} singular simple reflections. We set
$\be_j=\hw_j(\al_{p_j})\in \tR^0,\, 1\le j\le g.$ Then
(cf. (\ref{hw0s}))
\begin{align} \label{lambdarelprodu}
&\hw_1=\hw^1,\ \hw_2=s_{i_1}\cdots s_{i_{p_1-1}}s_{i_{p_1+1}}\cdots
s_{i_{p_2-1}}=\hw^{p_2} s^{p_1},\notag\\
&\ldots,\,\hw_g=\hw^{p_g} s^{p_1}\cdots s^{p_{g-1}} 
\hbox{\ using\ } s^p=s_{\tal^p}\,;
\notag\\
&\be_{1}=\hw_1(\al_{p_1})=\tbe^1,\, \be_{2}=\hw_2(\al_{p_2}),
\ldots, \be_{g}=\hw_g(\al_{p_g}).
\end{align}

\begin{proposition}\label{BRUHTW0}
(i) The minimal element $\hw^\circ$ from
($d,e$) of Theorem \ref{STRIKOUT} equals
$\hw\,\hw|_0^{-1}$ for the unique element
$\hw|_0\in \tW^0$ such that $\la^0(\hw|_0)=\la(\hw)\cap \tR^0$
constructed in (\ref{hw0s}).

(ii) Explicitly, $\la^0(\hw|_0)=\{\tbe^j, j=1,\ldots, g\}$,
\begin{align} \label{lambdainvrel}
&\hw|_0=s_{\be_g}\cdots s_{\be_1}= s_{\tbe^1}\cdots s_{\tbe^g},
\end{align}
where all roots $\be_j$ are simple in $\tR^0_+=\tR^0\cap\tR_+$,
and the first product is a reduced decomposition of $\hw|_0$
in terms of the simple reflections of $\tW^0$ defined with respect
to $\tR^0_+.$

(iii) The right Bruhat sets can be expressed in terms of the
standard Bruhat sets $\b(\tu;\,\tR^0_+)\subset \tW^0$
for $\tu\in\tW^0$ (defined
for the root system $\tR^0_+$ and the corresponding set of
simple roots there):
$$\b^0(\hw)\ =\ \hw^\circ\, \b\,(\hw|_0\,;\,\tR^0_+)
\hbox{\ \ for\ arbitrary\ } \hw\in \hW.
$$
\end{proposition}
{\em Proof}.
The formula for $\hw^\circ$ in terms
of $\hw|_0$ readily follows from Proposition \ref{BRUHATLA}.

Given a reduced decomposition
of $\hw$, the element $\hw|_0$ is naturally
represented as the product $s_{\tbe^1}\cdots s_{\tbe^g}$
for $\{\tbe_j\}=\tR^0\cap \la(\hw)$ and coincides with
$s_{\be_g}\cdots s_{\be_1}$ in (\ref{lambdainvrel}) by construction. 

Recall that these two products
correspond to moving simple singular reflections in the
decomposition of $\hw$ to the right beginning
with the {\em last} singular reflection, $s_{p_g}$, and 
beginning with the {\em first} 
singular reflection, that is $s_{p_1}$,
respectively. Notice that $s_{\be_j}$ appear in this
product in the order inverse to the order of $s_{\tbe^j}$.

Using Theorem \ref{INTRINLA}, we conclude that
$\{\tbe_j\}$ is a $\la$\~set in $\tR^0_+$ and, therefore,
it is exactly $\la^0(\hw|_0)$; hence,
$\hw|_0=s_{\be_g}\cdots s_{\be_1}$ is a reduced decomposition
in terms of simple reflections. It completes (ii).

Any increasing sequence of {\em simple} reflection in the
decomposition of $\hw|_0$ is associated with a decreasing
sequence of positive roots 
$\{\tal^{p_j}\}\subset\tR^0\cap \la(\hw)$
(the order becomes inverse) and, therefore, corresponds to
an increasing sequence of singular simple reflections
in the initial decomposition of $\hw$.
Deleting the former sequences of singular reflections
corresponds to deleting the latter sequences, and the other way
round. It gives (iii).
\sq
\medskip

\section{Double Hecke algebras}
\setcounter{equation}{0}
By  $m,$ we denote the least natural number
such that  $(P,P)=(1/m)\Z.$  Thus
$m=2 \for D_{2k},\ m=1 \for B_{2k} \and C_{k},$
otherwise $m=|\Pi|$.

The double affine Hecke algebra depends
on the parameters
$q, t_\nu,\, \nu\in \{\nu_\al\}.$ It will be defined
over the ring
$$
\Q_{q,t}\equal\Q[q^{\pm 1/m},t^{\pm 1/2}]
$$
formed by
polynomials in terms of $q^{\pm 1/m}$ and
$\{t_\nu^{\pm 1/2} \}.$
We set
\begin{align}
&   t_{\tal} = t_{\al}=t_{\nu_\al},\ t_i = t_{\al_i},\
q_{\tal}=q^{\nu_\al},\ q_i=q^{\nu_{\al_i}},\notag\\
&\where \tal=[\al,\nu_\al j] \in \tR,\ 0\le i\le n.
\label{taljx}
\end{align}

It will be convenient to use the parameters
$\{k_\nu\}$ together with  $\{t_\nu \},$ setting
$$
t_\al=t_\nu=q_\al^{k_\nu} \for \nu=\nu_\al, \and
\rho_k=(1/2)\sum_{\al>0} k_\al \al.
$$

For instance, by $q^{(\rho_k,\al)}$, we mean
$\prod_{\nu\in\nu_R}t_\nu^{((\rho_\nu)^\vee,\al)}$;
here $\al\in R$, $(\rho_\nu)^\vee=\rho_\nu/\nu$, and 
this product contains
only {\em integral} powers of $t_{\sht}$ and $t_{\lng}.$
Note that $(\rho_k,\al_i^\vee)=k_i=k_{\al_i}=
((\rho_k)^\vee,\al_i)$ for $i>0$; 
$(\rho_k)^\vee\equal \sum k_{\nu}(\rho_\nu)^\vee$. Also,
$(\rho_k,b_+)=-(\rho_k,b_-)$ for $b_+=w_0(b_-)$ (see above).

For pairwise commutative $X_1,\ldots,X_n,$
\begin{align}
& X_{\tb}\ =\ \prod_{i=1}^nX_i^{l_i} q^{ j}
\iif \tb=[b,j],\ \hw(X_{\tb})\ =\ X_{\hw(\tb)}.
\label{Xdex}
\\
&\hbox{where\ } b=\sum_{i=1}^n l_i \om_i\in P,\ j \in
\frac{1}{ m}\Z,\ \hw\in \hW.
\notag \end{align}
For instance, $X_0\equal X_{\al_0}=qX_\vth^{-1}$.

Later $Y_{\tb}=Y_b q^{-j}$ will be needed. Note the
negative sign of $j$. For instance,
$Y_0\equal Y_{\al_0}=q^{-1}Y_\vth^{-1}$.

We set $(\tilde{b},\tilde{c})=(b,c)$ ignoring the affine extensions.

\smallskip
\subsection{Main Definition}
We will use that $\pi_r^{-1}$ is $\pi_{r^*}$ and
$u_r^{-1}$ is $u_{r^*}$
for $r^*\in O\ ,$  $u_r=\pi_r^{-1}\om_r.$
The reflection $^*$ is
induced by an involution of the nonaffine Dynkin diagram
$\Gamma.$

\begin{definition}
The  double  affine Hecke algebra $\HH\ $
is generated over $ \Q_{ q,t}$ by
the elements $\{ T_i,\ 0\le i\le n\}$,
pairwise commutative $\{X_b, \ b\in P\}$ satisfying
(\ref{Xdex}),
and the group $\Pi,$ where the following relations are imposed:

(o)\ \  $ (T_i-t_i^{1/2})(T_i+t_i^{-1/2})\ =\
0,\ 0\ \le\ i\ \le\ n$;

(i)\ \ \ $ T_iT_jT_i...\ =\ T_jT_iT_j...,\ m_{ij}$
factors on each side;

(ii)\ \   $ \pi_rT_i\pi_r^{-1}\ =\ T_j \iif
\pi_r(\al_i)=\al_j$;

(iii)\  $T_iX_b T_i\ =\ X_b X_{\al_i}^{-1} \iif
(b,\al^\vee_i)=1,\
0 \le i\le  n$;

(iv)\ $T_iX_b\ =\ X_b T_i$ if $(b,\al^\vee_i)=0
\for 0 \le i\le  n$;

(v)\ \ $\pi_rX_b \pi_r^{-1}\ =\ X_{\pi_r(b)}\ =\
X_{ u^{-1}_r(b)}
 q^{(\om_{r^*},b)},\  r\in O'$.
\label{double}
\end{definition}
\qed

Given $\tw \in \tW, r\in O,\ $ the product
\begin{align}
&T_{\pi_r\tw}\equal \pi_r\prod_{k=1}^l T_{i_k},\where
\tw=\prod_{k=1}^l s_{i_k},
l=l(\tw),
\label{Twx}
\end{align}
does not depend on the choice of the reduced decomposition
(because $\{T\}$ satisfy the same ``braid'' relations
as $\{s\}$ do).
Moreover,
\begin{align}
&T_{\hv}T_{\hw}\ =\ T_{\hv\hw}\  \hbox{ whenever}\
 l(\hv\hw)=l(\hv)+l(\hw) \for
\hv,\hw \in \hW. \label{TTx}
\end{align}
In particular, we arrive at the pairwise
commutative elements
\begin{align}
& Y_{b}\ =\  \prod_{i=1}^nY_i^{l_i} \iif
b=\sum_{i=1}^n l_i\om_i\in P,\where
Y_i\equal T_{\om_i},
\label{Ybx}
\end{align}
satisfying the relations
\begin{align}\label{TYTL}
&T^{-1}_iY_b T^{-1}_i\ =\ Y_b Y_{\al_i}^{-1} \iif
(b,\al^\vee_i)=1,
\notag\\
& T_iY_b\ =\ Y_b T_i \iif (b,\al^\vee_i)=0,
 \ 1 \le i\le  n.
\end{align}
\smallskip
The origin of this construction is due to Bernstein and
Zelevinsky. The definition (\ref{Ybx}) and
relations (\ref{TYTL}) are due to
Lusztig (see, e.g.,\cite{L});
in the DAHA context, they are from \cite{C15}.

\subsection{Automorphisms}
The following maps can be uniquely extended to
automorphisms of
$\HH\ $(see \cite{C15},\cite{C4},\cite{C12}):
\begin{align}
\vep:\ &X_i \mapsto Y_i,\   Y_i \mapsto X_i,\
 T_i \mapsto T_i^{-1}\,(i\ge 1),\,
t_\nu \mapsto t_\nu^{-1},\,
 q\mapsto  q^{-1},\label{veppol}\\
\tau_+:\  &X_b \mapsto X_b, \ Y_r \mapsto
X_rY_r q^{-\frac{(\om_r,\om_r)}{2}},\
\pi_r \mapsto q^{-(\om_r,\om_r)}X_r\pi_r,
\notag\\
\tau_+:\ &Y_\vth \mapsto q^{-1}\,X_\vth T_0^{-1}
T_{s_\vth},\, T_0\mapsto q^{-1}\,X_\vth T_0^{-1},
\and
\label{taux}\\
\tau_-\ &\equal  \vep\tau_+\vep,\ \
\si\equal \tau_+\tau_-^{-1}\tau_+\ =\
\tau_-^{-1}\tau_+\tau_-^{-1}= \vep\si^{-1}\vep,
\label{tauminax}
\end{align}
where $r\in O'.$
They fix
$T_i\,(i\ge 1).$ The $\tau_{\pm},\si$
fix $\ t_\nu,\ q$
and fractional powers of $t_\nu,\ q.$

Note that $\tau_-$ acts trivially
on $\{T_0,\pi_r,Y_b\},$ $\si$ sends $X_b$ to $Y_b^{-1}.$

In the definition of $\tau_\pm$ and $\si,$
we need to add $q^{\pm 1/(2m)}$ to
$\Q_{q,t}.$
These automorphisms actually act in the central extension
of the {\em elliptic
braid group} defined by the relations of $\HH\, $
where the quadratic relation is dropped and
fractional powers of $q$ are treated as central
elements.

The elements $\tau_\pm$ generate the projective $PSL(2,\Z),$
which is isomorphic to the braid group $B_3$ due to Steinberg.
Adding $\vep,$
we obtain a projective action of $PGL(2,\Z).$
\medskip

The following {\em anti-involutions\,} are of key importance
for the theory of the polynomial representation:
\begin{align}
&X^\star=X^{-1},Y^\star=Y^{-1}, T_i^\star=T_i^{-1},\
 q^\star=q^{-1},t_\nu^*=t_\nu^{-1},\label{star}\\
&\phi\equal \vep\,\star=\star\,\vep:\
X_b\mapsto Y_b^{-1},\
T_i\mapsto T_i\ (1\le i\le n),\label{starphi}
\end{align}
where the latter preserves $q,t_\nu$ and their fractional
powers.

We will also use
$\eta\equal \vep\si.$ It conjugates
$q,t$ and is uniquely defined
from the relations
\begin{align}
&\eta:\ T_i\mapsto T_i^{-1},\ X_b\mapsto X_b^{-1},\
\pi_r\mapsto \pi_r,\notag\\
&\where 0\le i\le n,\ b\in P,\ r\in O'.
\label{etatxpi}
\end{align}

\smallskip
\subsection{Intertwining operators}
The {\em $X$-intertwiners\,} (see \cite{C1})
are introduced as follows:
\begin{align}
&\Phi_i\ =\
T_i + (t_i^{1/2}-t_i^{-1/2})
(X_{\al_i}-1)^{-1} \for 0\le i\le n,\notag\\
&\Phi_i^\diamond\ =\ \Phi_i (X_{\al_i}-1)\ =\
T_i (X_{\al_i}-1)+ t_i^{1/2}-t_i^{-1/2},
\notag\\
& G_i=\Phi_i (\phi_i)^{-1},\
\phi_i= t_i^{1/2} +
(t_i^{1/2} -t_i^{-1/2})(X_{\al_i}-1)^{-1}.
\label{Phix}
\end{align}

They belong to $\HH\ $ extended by
the rational functions in terms of $\{X\}$. The $G$ are called
the {\dfont normalized intertwiners\,}. The elements
$G_i, 0\le i\le n,\ r\in O',$
satisfy the same relations
as $\{s_i,\pi_r\}$ do, so the map
\begin{align}
\hw\mapsto G_{\hw}\ =\ \pi_r G_{i_l}\cdots G_{i_1},
\where \hw=\pi_r s_{i_l}\cdots s_{i_1}\in \hW,
\label{Phiprodx}
\end{align}
is a  well defined homomorphism  from $\hW.$

The intertwining property is
$$
G_{\hw} X_b G_{\hw}^{-1}=X_{\hw(b)}
\where X_{[b,j]}= X_b q^{j}.
$$
We will refer to $G_{\hw}$ as intertwiners too;
$\Phi_i,\Phi_i^\diamond, G_i,\pi_r$ will be called
{\em simple intertwiners.}

As to $\Phi_i$ and $\Phi_i^\diamond$, they
satisfy the  homogeneous Coxeter relations
and those with $\Pi_r.$ So we may set
$\Phi_{\hw} =$ $\pi_r \Phi_{i_l}\cdots \Phi_{i_1}$ for the reduced
decompositions; similarly,
$$\Phi_{\hw}^\diamond =\pi_r \Phi_{i_l}^\diamond\cdots
\Phi_{i_1}^\diamond=
\Phi_{\hw}(X_{\tal^l}-1)\cdots(X_{\tal^1}-1)
$$
in the notation from (\ref{tal}):
$\tal^p=s_{i_1}\cdots s_{i_{p-1}}(\al_{i_p})$.
They intertwine $X$ in the same way as $G$ do.

The $\Phi_w$ for $w\in W$ are well known
in the theory of affine Hecke algebras.
The affine intertwiners are
the key tool in the theory of semisimple
and spherical representations of DAHA, including
applications to the Macdonald polynomials and
the Harish-Chandra\~ Opdam spherical transform.
\smallskip

\subsection{Intertwiners and 
\texorpdfstring{{\mathversion{bold}$T_i$}}{T}}
Concerning the relation of $\Phi$ to $\{T_i\}$, the
following holds:
\begin{align}\label{tiphitj}
&T_i\Phi_{\hw}=\Phi_{\hw}T_j  \hbox{ \ if\ }
s_i\hw=\hw s_j,\ l(\hw s_j)=l(\hw)+1=l(s_i\hw)
\end{align}
for $i,j\ge 0$.
These conditions are equivalent to
$s_j(\la(\hw))=\la(\hw)$ and
$\hw(\al_j)=\al_i$.

An explicit description of such $\hw$ and the justification
are straightforward.
Indeed, there must exist a reduced decomposition of
$\hw$ that begins with $\hu$ such that $\hu s_j=s_{j\,'\,}\hu$
is an elementary Coxeter transformation
($j\,'\,=j$ unless it is of type $A_2$). One has:
\begin{align}\label{tphit}
T_{j\,'\,}\Phi_{\hu}&=(\Phi_{j\,'\,}-
(t_{j\,'\,}^{1/2}-t_{j\,'\,}^{-1/2})/
(X_{\al_{j\,'\,}}-1))\Phi_{\hu}\\
&=\Phi_{\hu}\Phi_j-\Phi_{\hu}
(t_{j}^{1/2}-t_j^{-1/2})/(X_{\al_{j}}-1))=
\Phi_{\hu}T_j.\notag
\end{align}
Then $s_i\hw'=\hw's_{j'}$ for $\hw=\hw'\hu$ and we can continue by
induction.

Let us generalize (\ref{tiphitj})
using the right Bruhat ordering from
Section \ref{sec:RightBruhat}
associated with an arbitrary {\em root subsystem} $\tR^0$ of $\tR$.

\begin{theorem}\label{PHIBRUHAT}
Given a reduced decomposition
$\hw =\pi_r s_{i_l}\cdots s_{i_1}$, we define
$\widetilde{\Phi}_{\hw}$ by replacing $\Phi_{i_p}$ with
$T_{i_p}$ when $p$ is {\dfont singular}, i.e.,
$\tal^p=s_{i_1}\cdots s_{i_{p-1}}(\al_{i_p})\in$
$\tR^0$; the $\widetilde{\Phi}_{\hw}^\diamond$  and
$\widetilde{G}_{\hw}$ are defined correspondingly.

(a) The element $\widetilde{\Phi}_{\hw}$ remains unchanged under
homogeneous Coxeter transformations of the reduced decomposition
of $\hw$ if the $A_2$\~transforms with singular middle $i_p$ and
non-singular edges are avoided,
that are as follows
\begin{align}\label{atwosing}
&s_{i\pm 1}s_{i}s_{i\pm 1}\mapsto s_{i}s_{i\pm 1}s_{i},
\hbox{\ \, for\ \, }  i=i_p,\ \, i_{p+1}=i\pm 1=i_{p-1},\\
&\hbox{where\ \, } \tal^p\in \tR^0 \hbox{\ \, and\ \,}
\tal^{p-1}\not\in \tR^0\not\ni \tal^{p+1}. \notag
\end{align}

(b) Similarly, $\widetilde{\Phi}_{\hw}^\diamond$ are invariant under
homogeneous Coxeter transformations apart from (\ref{atwosing}).
This restriction is not needed for
$\widetilde{G}_{\hw}$, i.e.,\,
the latter does not depend on the choice of the
reduced decomposition for arbitrary homogeneous Coxeter
transformations if the normalized intertwiners are used.

(c) If $\widetilde{\Phi}^{2}$ is obtained from
$\widetilde{\Phi}^1=\widetilde{\Phi}_{\hw}$
when transformations (\ref{atwosing}) are
allowed, then the difference $\widetilde{\Phi}^{1}\,-\,$
$\widetilde{\Phi}^{2}$
is a linear combination of the terms in the form
$\widetilde{\Phi}_{\hw'}P(X)$ for rational functions $P(X)$ from
\begin{align}\label{pnonsing}
&
\Z[q_\nu,\, (t_\nu^{1/2}-t_\nu^{-1/2}),\,
X_{\be},\, (1-X_{\tal})^{-1}],\
\be\in R, \tal\in \tR\setminus \tR^0,
\end{align}
where $\widetilde{\Phi}_{\hw'}$ are
defined for some reduced decompositions of
$\hw'\in\b_o^0(\hw)$ from Theorem \ref{STRIKOUT}.

(d) The linear space generated by $\widetilde{\Phi}_{\hw'}P(X)$
for $\hw'\in \b^0(\hw)$ and
$$
P(X)\in \Q_{q,t}[X_{b},\, (1-X_{\tal})^{-1}\,\mid\,
b\in B,\, \tal\not\in \tR^0]
$$
is a left module over
$\Q_{q,t}[X_{b},\, (1-X_{\tal})^{-1}\,\mid\,
b\in B,\, \tal\not\in \hw^{-1}(\tR^0)]$.
\end{theorem}
{\it Proof}.
All Coxeter rank 2 transformations but the exceptional one from
(\ref{atwosing}) are checked similar
to (\ref{tphit}); one follows the proof of Theorem \ref{STRIKOUT}.

The cases of $\widetilde{\Phi}^\diamond$ and $\widetilde{G}$
are completely analogous.
The unitary property $G_i^2=1$ guarantees
that the transformations (\ref{atwosing}) do not change
$\widetilde{G}_{\hw}$. We arrive at ($a,b$).
\smallskip

Let the transformation be $s_{i\pm 1}\underline{s_{i}}s_{i\pm 1}
\mapsto
s_{i}\underline{s_{i\pm 1}}s_{i}$, where we underline the singular
reflection and $i=i_p,\ i_{p-1}=i\pm 1= i_{p+1}$
inside a decomposition of $\hw$. Then, setting
$\al=\al_i,\, \be=\al_{i\pm 1},\, t=t_i,$
\begin{align}\label{phitphi}
&\widetilde{\Phi}^1={\Phi}_{i\pm 1}T_i{\Phi}_{i\pm 1},\,
\widetilde{\Phi}^2={\Phi}_{i}T_{i\pm 1}{\Phi}_{i},\notag\\
&\widetilde{\Phi}^1-\widetilde{\Phi}^2=
\frac{(t^{1/2}-t^{-1/2})^3\, (X_{\be}^{-1}-X_{\al}^{-1}) }
{(2-X_{\al}-X_{\al}^{-1})(2-X_{\be}-X_{\be}^{-1})}.
\end{align}
The denominators of (\ref{phitphi})
become products of $(X_{\tal}-1)$ for $\tal\not\in \tR^0$
when moved to the right all the way. We move them
through non-singular $\Phi_j$ using the intertwining
property and through singular $T_j$ using directly 
the relations (\ref{tixi}).
In either case, the denominators are
conjugated by the corresponding $s_j$.

The product $\widetilde{\Phi}^1$
where $\Phi_{i\pm 1}T_i{\Phi}_{i\pm 1}$ is omitted
becomes 
$\widetilde{\Phi}_{\hw'}$ defined for the
{\em standard}, possibly non-reduced, decomposition of
$\hw'$ obtained from $\hw$ by deleting 
$s_{i\pm 1} s_i s_{i\pm 1}$.
If this decomposition is non-reduced, then
one needs either the quadratic relations for $T_j^2$ for
singular $s_j$ or the formula for $\Phi_j^2$, {\em an
$X$\~function},
for non-singular $s_j$. Using these relations may be preceded by
homogeneous Coxeter transformations. It results in additional
terms of the same type (\ref{pnonsing}). Eventually, we will come
to $\Phi_{\hw'}$ for reduced decompositions of $\hw'\in \b^0(\hw)$,
as required in ($c$).

Another transformation that may occur here is
``deleting" singular $T_j$ due to (\ref{tixi}) 
in process of moving $X$\~functions to the right,
maybe followed by further reductions described above.
The $X$\~function here will be either 
$\Phi_j^2$ or those due to (\ref{tixi});
more generally, one can take here
arbitrary rational functions of $X$ with the
denominators that become from (\ref{pnonsing}) when 
moved to the right. This proves ($c$) and ($d$).
\sq
\smallskip

Note that when $(X_{\al_i}-X_{\al_{i\pm 1\,})})$ is moved to 
the right it becomes divisible by $(X_{\tal^p}-1)$ for $i=i_p$ 
in the absence of $\{T_j\}$ in the corresponding
portion of $\widetilde{\Phi}_{\hw}$.
Otherwise,
we will need to move this difference through $T_j$ for singular
indices $j=i_q\, (q<p)$, which may destroy the divisibility.
More generally, we obtain the following.

\begin{corollary}\label{SINGCOM}
Given a reduced decomposition $\hw=\pi_r s_l\cdots s_1$ and
a singular $s_{i_p}$,
if $s^p\, s^q=s^q\, s^p$ for all singular $i_q$ with $q<p$,
where the notation $s^q=s_{\tal^q}$ from
Theorem \ref{STRIKOUT} is used, then
$\widetilde{\Phi}^1-\widetilde{\Phi}^2$ is divisible by
$(X_{\tal^p}-1)$ on the right.
\sq
\end{corollary}
\smallskip

\subsection{Compatible {\em r}--matrices}
Generalizing the construction of $\widetilde{\Phi}$,
we come to the notion of the
compatibility of {\em $\r$\~matrices\,}. It extends the key
property of the (quantum) $\r$\~matrix
$\{\r_{\tal},\,\tal\in \tR_+\}$
defined in \cite{Ch5}, which is as follows. The products
$\r_{\tw}=\r_{\tal^l}\cdots\r_{\tal^1}$ must 
depend only on $\tw$
and must not depend on the
{\em reduced} decompositions $\tw=s_{i_l}\cdots s_{i_1}$ for
any $\tw\in \tW$. Here $\r_{\tal}$ are elements
in some algebra $\mathfrak{R}$ with a unit.
Equivalently, the {\em cocycle relation} must be satisfied:
\begin{align}\label{cocyc}
&\r_{\tw} \r_{\tu}={}^{\tu^{-1}}(\r_{\tw})\r_{\tu}, \hbox{\ where\ }
\tu^{-1}(\ldots\r_{\tal}\ldots)=\ldots \r_{\tu^{-1}(\tal)}\ldots,
\end{align}
provided that $l(\tu\tw)=l(\tu)+l(\tw)$.
\smallskip

If \, $\{\,{\r}_{\tw},\, \tw\in \tW\,\}$ is extended
to $\hW\ni \hw=\pi_r\tw$ by the formula
$\r_{\hw}=\p_r\r_{\tw}$ for a homomorphism
$\Pi\ni \pi_r\mapsto \p_r\in \mathfrak{R}$
and (\ref{cocyc}) holds,
then $\r_{\hw}$ is called a $\hW$\~{\em extension}
(or a $\hW^\flat$\~extension for $\Pi^\flat$ instead of
$\Pi$).
It simply means the commutativity relations:
$\p_r{\r}_{\tal}={\r}_{\tal}\p_r$.
\smallskip

Given a {\em root subsystem} $\tR^0$,
an $\r$\~matrix $\{\r^0_{\tal}\}$  defined for $\tal\in \tR_+^0$ 
is called {\dfont completely compatible} with $\r$ if 
the following {\em compatibility property} holds.
The elements
$\widetilde{\r}_{\tw}$ obtained from  $\r_{\tw}$
by replacing $\r_{\tal}\mapsto \r^0_{\tal}$
whenever $\tal\in \tR^0_+$
have to remain unchanged when the
(homogeneous) Coxeter transformations
are applied to reduced decompositions of $\tw$, i.e.,
they must depend only on $\tw$. 
Note that the complete compatibility
does not result in (\ref{cocyc}) for
$\widetilde{\r}_{\tw}$ unless $\tR^0$ is $\tW$\~invariant.

Respectively, if
$\p_r\widetilde{\r}_{\tal}=\widetilde{\r}_{\tal}\p_r$,
then \, $\{\,\widetilde{\r}_{\hw},\, \hw\in \hW\,\}$
is called a {\em $\hW$\~extension}
of\, $\{\,\widetilde{\r}_{\tw},\, \tw\in \tW\,\}$.

In the context of $\r$\~matrices,
the Coxeter transformations are inversions 
($w_0$\~permutations) of the consecutive  elements 
$\r_{\tal^p}$ for the segments $\{\tal^p\}\subset \la(\hw)$ 
that can be identified 
with consecutive positive roots of rank $2$,
i.e., with those of types $A_1\times A_1$, $A_2,\, B_2,$
or $G_2$.

An $\r$\~matrix $\r^0$ is
called {\dfont partially compatible} with $\r$
if the transformations from
(\ref{atwosing}) are excluded, that are transformations
of type $A_2$ such that the second (middle)
root belongs to $\tR^0$ but the first and the third
(the edges) do not.
\smallskip

{\bf Rank two relations.}
First of all,
all rank two relations from \cite{Ch5}
for $\r$ and 
$\widetilde{\r}$ (individually) must
be satisfied. Then one needs to check the following
rank two {\em compatibility relations\,}:
$$
(A_1 \times A_1):\ \widetilde{\r}_{\alpha }
\r_{\beta}= \r_{\beta}\widetilde{\r}_{\alpha },
$$
\begin{align}\label{specialatwo}
(A_2,\ \hbox{general}):\
&\widetilde{\r}_{\alpha } \r_{\alpha + \beta }
\r_{ \beta } =
\r_{ \beta }\r_{\alpha + \beta }\widetilde{\r}_{\alpha },\notag\\
(A_2,\ \hbox{special}):\
&\r_{\alpha } \widetilde{\r}_{\alpha + \beta } \r_{ \beta } =
\r_{ \beta }\widetilde{\r}_{\alpha + \beta }\r_{\alpha },
\end{align}
where the {\em special relation} is needed to ensure the
{\em complete compatibility} of $\r$ and $\widetilde{\r}$
(otherwise, the compatibility will be only {\em partial}), and
\begin{align*}
(B_2,\ 1\sim):\ &\widetilde{\r}_{\alpha } \r_{\alpha + \beta }
\r_{\alpha + 2\beta }\r_{ \beta } =  \r_{ \beta }
\r_{\alpha + 2\beta }\r_{\alpha + \beta }\widetilde{\r}_{\alpha }\\
&\r_{\alpha } \widetilde{\r}_{\alpha + \beta }
\r_{\alpha + 2\beta }\r_{ \beta } =  \r_{ \beta }
\r_{\alpha + 2\beta }\widetilde{\r}_{\alpha + \beta }\r_{\alpha }\\
&\r_{\alpha } \r_{\alpha + \beta }
\widetilde{\r}_{\alpha + 2\beta }\r_{ \beta } =  \r_{ \beta }
\widetilde{\r}_{\alpha + 2\beta }\r_{\alpha + \beta }\r_{\alpha }\\
&\r_{\alpha } \r_{\alpha + \beta }
\r_{\alpha + 2\beta }\widetilde{\r}_{ \beta } =  \widetilde{\r}_{
\beta }
\r_{\alpha + 2\beta }\r_{\alpha + \beta }\r_{\alpha },
\end{align*}
\begin{align*}
(B_2,\ 2\sim):\ &\widetilde{\r}_{\alpha } \r_{\alpha + \beta }
\widetilde{\r}_{\alpha + 2\beta }\r_{ \beta } =  \r_{ \beta }
\widetilde{\r}_{\alpha + 2\beta }\r_{\alpha + \beta }
\widetilde{\r}_{\alpha }\\
&\r_{\alpha } \widetilde{\r}_{\alpha + \beta }
\r_{\alpha + 2\beta }\widetilde{\r}_{ \beta } =  \widetilde{\r}_{
\beta }
\r_{\alpha + 2\beta }\widetilde{\r}_{\alpha + \beta }\r_{\alpha },
\end{align*}
and also the corresponding compatibility relations for $G_2$,
obtained from the ones for $\r$,
\begin{align*}
&\r_{\alpha } \r_{\alpha + \beta }
\r_{2\alpha + 3\beta }\r_{\alpha + 2\beta }
\r_{\alpha + 3\beta }\r_{ \beta } =
\r_{ \beta } \r_{\alpha + 3\beta }
\r_{\alpha + 2\beta }\r_{2\alpha + 3\beta }
\r_{\alpha + \beta }\r_{\alpha },
\end{align*}
by replacing $\r\mapsto \widetilde{\r}$ for: \ 
1) {\em one} of the roots ($1\sim$) from the set
$$
\al, \ \alpha + \beta,\ 2\alpha + 3\beta,
\ \alpha + 2\beta,\ \alpha + 3\beta,\ \beta,
$$
2) those in the pairs $\{\alpha, \alpha + 2\beta\}$
$\{2\alpha + 3\beta,\beta\}$
of orthogonal roots ($2\sim$) and\ 3) for those in the
$A_2$\~subsystems ($3\sim$):
$$
\{\al,2\alpha + 3\beta,\alpha + 3\beta \} \hbox{\ \ or\ \ }
\{\alpha + \beta,\alpha + 2\beta,\beta \}.
$$
Here the roots involved constitute all positive roots in the
corresponding rank two subsystem taken clockwise and then
counterclockwise. See the proof of Lemma \ref{LEMHCOX}.
The roots $\al,\be$ generating a rank two system
are {\em affine} from $\tR$.
\smallskip

We mention that not all compatibility relations from these
lists can appear for a given pair $\tR^0\subset R$
due to our definition of a {\em root subsystem}.
The relations we give are ``generic"; not all of them
must be checked for concrete $\tR^0\subset R$.

For instance, the
subsets of {\em all} short roots do not form
a {\em root subsystem} for $B_2$ or $G_2$, as well
as the pairs of orthogonal roots for $G_2$;
such roots $\Z$\~generate the whole rank two system.
The corresponding compatibility $\r$\~relations have to be 
omitted if we stick to the definition of the {\em root subsystem}
$\tR^0$ from this paper.
\medskip

\subsection{Examples}
Let us assume that 
there exists a homomorphism 
$\tW\to \mathfrak{R}^*$ sending  
$s_{\tal}\mapsto \s_{\tal}$  and
$\s_{\tal}(\r_{\tbe})\s_{\tal}=\r_{s_{\tal}(\tbe)}$.
The simplest abstract example of {\em complete} (respectively,
{\em partial}) compatibility
is $\{\r^0_{\tal}=\s_{\tal},\, \tal\in \tR^0\}$ for {\em unitary}
(respectively, {\em arbitrary}) $\r$. 
By {\em unitary}, we mean that $\{\r_{\tal}\}$ are
defined for all $\tal\in\tR$ and satisfy the conditions
$\r_{\tal}\r_{-\tal}=1$.
\smallskip

The standard Hecke-type examples of $\r$\~matrices
and the compatibility are as follows.

Let us assume that $T_i$ and $\tW$ act in a vector space $V$ and
$T_{\tal}\equal\tw T_i \tw^{-1}$ does not depend on the choice
of $\tw$ such that $\tw(\al_i)=\tal$ for an arbitrary
$\tal\in \tR_+$. Then $\{\r_{\tal}=s_{\tal}T_{\tal}\}$
is a {\em constant} $\r$\~matrix.

An example of such $V$ is the polynomial
representation (see below)
restricted to the affine Hecke algebra generated by
$\{T_i,\, i\ge 0\}$ with
the usual action of $\tW$ on Laurent polynomials.
In the abstract theory of quantum
Knizhnik\~Zamolodchikov equations, the action of $\tw$ through
the projection $\tW\to W$ is used, i.e., with
$T_{[\al,\nu_{\al}j]}$ depending only on $\al\in R$ 
(see \cite{Ch5}). 

Another (well-known) example is for the root system of
type $A_n$ only. We set $T_{\al}=T_{ij}$ for $\al=\ep_i-\ep_j$
and the so-called Baxter\~Jimbo matrices $T_{ij}$
in $V_N\otimes V_N$ in the $i,j$\~components of the
tensor power $V_N^{\otimes (n+1)}$ of the
fundamental $N$\~dimensional representation $V_N$ of $GL_N$.
\smallskip

{\bf Using intertwiners.}
Let $\La_{a}\, (a\in Q)\,$ be multiplicative {\em scalar}
variables: they satisfy all the formulas
for $X_a, X_{\tal}$ and are supposed to commute with
$\{T_{\tal}\}$ defined via $T_{\tal}\equal\tw T_i \tw^{-1}$ 
as above.
Here $Q$ is the root lattice, but it can be replaced
by $P$ (or $B$) if the extended case is considered.

The main properties of the intertwiners $\Phi$
and $G$ can be reformulated (generalized) as follows. 
The collections
$\r$ and $\g$ given by
\begin{align}\label{rgmatrix}
&\r_{\tal}\equal s_{\tal}T_{\tal}+
(t_{\al}^{1/2}-t_{\al}^{-1/2})/(\La_{\tal}-1),\
\tal\in \tR_+, \\
&\g_{\tal}\equal \r_{\tal}
(t_{\al}^{1/2}+(t_{\al}^{1/2}-t_{\al}^{-1/2})/(\La_{\tal}-1))^{-1},
\notag
\end{align}
are {\em non-constant} End$(V)$\~valued $\r$\~matrices.
Their extensions to $\hW$ are given by $\pi_r\mapsto 1$.

Note that one can use the formulas
$T_{\tal}=\tw T_i \tw^{-1}$ to define $T_{\tal}$ for all $\tal$,
non-necessarily positive. Then
$\g$ is unitary, $\g_{\tal}\g_{-\tal}=1$, $\r$ is not.
\smallskip

The (algebraically non-trivial) fact that $\r,\g$ are
$\r$\~matrices can be deduced from the consideration
of the intertwiners $\Phi,\,G$
in the representation of $\HH^\flat$ induced from a {\em generic}
character
$\La:\,X_a\mapsto \La_a$ of $\Q_{q,t}[X]$;
they are applied to the cyclic vector. The space
of this induced representation
is naturally isomorphic to the group
algebra of the affine Hecke algebra generated by $\{T_i,\, i\ge 0\}$.
Note that these $\r$\~matrices are functions with the values
in any $\HH^\flat$\~modules $V$ with an action of $\tW$
such that the elements $T_{\tal}=\tw T_i \tw^{-1}$ are well defined
(do not depend on the choice of $\tw$ such that 
$\tw(\al_i)=\tal$).

Theorem \ref{PHIBRUHAT} readily gives that the collection
$\{\g^0_{\tal}\equal s_{\tal}T_{\tal}=\r^0_{\tal}\}$ is
an $\r$\~matrix that is 
{\em completely compatible} with $\g$ and also
{\em partially compatible} with $\r$ from
(\ref{rgmatrix}) for any root subsystem
$\tR^0$ provided that
$$
\La_{\tal}=1 \Longrightarrow \tal\in \tR^0.
$$
We must also impose
the inequality $\La_{\tal}\neq t_{\al}^{\pm 1}$ for $\g$ whenever
$\tal\not\in\tR^0$.

We note that the case of $\g$ here is related to
Proposition 8.11 from \cite{L} (in the affine case).
The considerations from Section 6 from this paper are similar 
to our ``abstract" approach to the relative Bruhat ordering and 
{\em completely} compatible $\r$\~matrices. 
\smallskip

\rmk
(i) There is a variant of this construction when $\g$ is always used
instead of $\r$ as $\La_{\tal}\neq 1,t_{\al}^{\pm 1}$,
i.e., we make $\widetilde{\r}$ as unitary as possible.
It gives an example of {\em partially compatible triple} 
$\{\g^0,\r,\g\}$.

(ii) The {\em trivial} $\{\g^0_{\tal}=1,\ \tal\in\tR^0\}$ is
{\em completely compatible} with $\g$; the corresponding
$\widetilde{\g}_{\hw}$ equals $\g_{\hw^\circ}$ for the minimal
element $\hw^\circ$ from Theorem \ref{STRIKOUT},($d$).
Thus the {\em complete compatibility}, equivalently,
the fact that $\widetilde{\g}_{\hw}$ does not depend
on the choice of the reduced decomposition of $\hw$, is
not obvious even in this (simplest) case.

(iii)
The $\La$\~factors that appear when the Coxeter transformations
of type (\ref{atwosing}) are used in the products
for $\widetilde{\r}$ are
scalar and easy to control; so the difference between
$\r$ and $\g$ is simply a matter of
normalization. It becomes much more involved
when the intertwiners $\widetilde{\Phi}_{\hw}$ are used
to define $\widetilde{\r}$. 
Replacing $\Phi$ by the normalized intertwiners $G$ is
possible only in DAHA\~modules where non-invertible
intertwiners do not appear.
\sq
\smallskip

{\bf A connection with virtual links}.
If $\tR^0$ is not involved 
and one can uses $\g^0=\{s_{\tal}T_{\tal}\}$
instead of $\g$ for any roots, then the
rank two relations for  $\{\g^0,\g\}$ 
are essentially the axioms
of the {\em virtual links\,}  generalized to
arbitrary reduced affine root systems. 
Here $\g^0$ serves usual links, $\g$ stays for virtual links. 
In the DAHA theory, the combinations like
$\g_{12}^0\g_{13}^0\g_{23}$ (with exactly two $\g^0$)
cannot appear. These ones also do not satisfy the
triangle (Reidemeister III) property in the theory
of virtual links.

Recall that $\g$ is
unitary and the transformations
(\ref{atwosing}) are allowed. In the context
of the $\r$\~matrices, (\ref{atwosing}) is
the {\em special relation} of type $A_2$
from (\ref{specialatwo}) with
$\g^0$ is in the middle and $\g$ at the edges.
If $\r$ is used instead of $\g$, then
(\ref{specialatwo}) does not hold, so
the general DAHA\~theory, generally,
requires more sophisticated approach.

In a more abstract way, the algebra formally
generated by DAHA and the intertwining operators
(defined for $X$ {\em or} for $Y$, normalized or
not) seems the key object. It is challenging to 
understand its ``topological meaning".
We also note that the subcategory 
of  semisimple DAHA\~modules
may have something to do with the invariants
of virtual links and the ``categorification". 
\medskip

\section{Polynomial representation}
\setcounter{equation}{0}
From now on we will switch from $\HH\ $ to an
{\dfont intermediate subalgebra} $\HH^\flat\subset\HH\ $
with $P$ replaced by a lattice $B$ between $Q$ and $P$
(see \cite{C12}). Respectively, $\Pi$ is
changed to the preimage $\Pi^\flat$ of $B/Q$ in $\Pi.$ Generally,
there can be two different
lattices $B_x$ and $B_y$ for $X$ and $Y.$
We consider only $B_x=B=B_y$ in the paper;
respectively, $a,b\in B$ in $X_a,Y_b.$

We also set $\hW^{\flat}=B\cdot W\subset \hW,\ $
and replace $m$ by the least $\tilde{m}\in \N$ such that
$\tilde{m}(B,B)\subset \Z$ in the definition of the
$\Q_{q,t}.$

The relations (iii,iv) from Definition \ref{double} (cf.
 \cite{L}) read as:
\begin{align}
&T_iX_b -X_{s_i(b)}T_i\ =\
(t_i^{1/2}-t_i^{-1/2})\frac{s_i(X_b)-X_b}{
X_{\al_i}-1},\ 0 \le i\le  n.
\label{tixi}
\end{align}
Replacing $X_b$ by $Y_b^{-1}$
we obtain the dual "$T-Y$" relations.
\smallskip

Note that $\HH^\flat$ and the polynomial representations
(and their rational and trigonometric
degenerations)  
are actually defined over $\Z$ extended by the parameters of DAHA.
However the field $\Q_{q,t}$ will be sufficient in this paper.
\smallskip

\subsection{Definitions}
The Demazure-Lusztig operators
are as follows:
\begin{align}
&T_i\  = \  t_i ^{1/2} s_i\ +\
(t_i^{1/2}-t_i^{-1/2})(X_{\al_i}-1)^{-1}(s_i-1),
\ 0\le i\le n;
\label{Demazx}
\end{align}
they obviously preserve $\Q[q,t^{\pm 1/2}][X_b]$.
We note that only the formula for $T_0$ involves $q$:
\begin{align}
&T_0\  = \ t_0^{1/2}s_0\ +\ (t_0^{1/2}-t_0^{-1/2})
(X_0 -1)^{-1}(s_0-1),\hbox{\ where\ }\notag\\
&X_0=qX_\vth^{-1},\
s_0(X_b)\ =\ X_bX_{\vth}^{-(b,\vth)}
 q^{(b,\vth)},\
\al_0=[-\vth,1].
\end{align}

The map sending $ T_j$ to the corresponding operator from
(\ref{Demazx}), $X_b$ to $X_b$
(see (\ref{Xdex})) and 
$\pi_r\mapsto \pi_r$ induces a
$ \Q_{ q,t}$\~linear
homomorphism from $\HH^\flat\, $ to the algebra of
linear endomorphisms
of $\Q_{ q,t}[X]$.
This $\HH^\flat\,$-module is faithful
and remains faithful when   $q,t$ take
any nonzero complex values assuming that
$q$ is not a root of unity.
It will be called the
{\dfont polynomial representation};
the notation is
$$
\v\equal \Q_{q,t}[X_b]\ =\ \Q_{q,t}[X_b, b\in B].
$$

The images of the $Y_b$ are called the
{\dfont difference Dunkl operators\,}. To be more exact,
they must be called difference-trigonometric Dunkl operators,
because there are also
difference-rational Dunkl operators.

The polynomial representation
is the $\HH^\flat\,$-module induced from the one dimensional
representation $T_i\mapsto t_i^{1/2},\,$ $Y_i\mapsto t_i^{1/2}$
of the affine Hecke subalgebra $\h_Y^\flat=\lan T_i,Y_b\ran.$
Here the following PBW-type theorem is used.
For arbitrary nonzero $q,t,$ any element $H \in \HH^\flat\, $
has a unique decomposition in the form
\begin{align}
&H =\sum_{w\in W }\,  g_{w}\, f_w\, T_w,\
g_{w} \in \Q_{q,t}[X_b],\ f_{w} \in \Q_{q,t}[Y_b].
\label{hatdecx}
\end{align}
\smallskip

{\bf Invariant form.} The definition is in terms of
the {\dfont truncated theta function}
\begin{align}
&\mu\ = \mu^{(k)}=\prod_{\al \in R_+}
\prod_{i=0}^\infty \frac{(1-X_\al q_\al^{i})
(1-X_\al^{-1}q_\al^{i+1})
}{
(1-X_\al t_\al q_\al^{i})
(1-X_\al^{-1}t_\al^{}q_\al^{i+1})}.\
\label{mu}
\end{align}
It is considered as a Laurent series with the
coefficients in  $\Q[t_\nu][[q_\nu]]$ for
$\nu\in \nu_R.$

The constant term of a Laurent series $f(X)$ will be denoted
by $\langle  f \rangle.$ Then
\begin{align}
&\langle\mu\rangle\ =\ \prod_{\al \in R_+}
\prod_{i=1}^{\infty} \frac{ (1- q^{(\rho_k,\al)+i\nu_\al})^2
}{
(1-t_\al q^{(\rho_k,\al)+i\nu_\al})
(1-t_\al^{-1}q^{(\rho_k,\al)+i\nu_\al})
}.
\label{consterm}
\end{align}
Recall that $q^{(z,\al)}=q_\al^{(z,\al^\vee)},
\ t_\al=q_\al^{k_\al},$ so we can set here
$q^{(\rho_k,\al)+i\nu_\al}=q_\al^{(\rho_k,\al^\vee)+i}$.
This formula is equivalent to the Macdonald constant term conjecture 
from \cite{M1}, proved in complete generality in \cite{C2}.

Let $\mu_\circ\equal \mu/\langle \mu \rangle$.
The coefficients of the Laurent series $\mu_\circ$ are
from the field of rationals
$\Q(q,t)\equal\Q(q_\nu,t_\nu), \where \nu\in \nu_R.$
Note that  $\mu_\circ^*\ =\ \mu_\circ$ for the involution
$$
 X_b^*\ =\  X_{-b},\ t^*\ =\ t^{-1},\ q^*\ =\ q^{-1}.
$$
This involution is the restriction
of the anti-involution $\,\star\,$
from (\ref{star}) to $X$-polynomials (and Laurent series).
These two properties of $\mu_\circ$
can be directly seen from the difference relations
for $\mu.$
\smallskip

\begin{proposition}
Setting,
\begin{align}
&\langle f,g\rangle_\circ\ \equal\ \langle
\mu_\circ f\ {g}^*\rangle\ =\
\langle g,f\rangle_\circ^* \for
f,g \in \Q(q,t)[X],
\label{innerpro}
\end{align}
the polynomial representation is $\star$-unitary:
\begin{align}
&\langle H(f),g\rangle_\circ\  = \
\langle f, H^\star(g)\rangle_\circ \for H\in \HH,
\ f\in \Q_{q,t}[X].
\label{staru}
\end{align}
\end{proposition}
\smallskip

\subsection{Macdonald polynomials}
There are two equivalent definitions of the
{\dfont nonsymmetric Macdonald polynomials\,}, denoted by
$E_b(X) = E_b^{(k)}$ for $b\in B$. They belong to
$\Q(q,t)[X_a,a\in B]$ and,
using the pairing $\lan\ ,\ \ran$,
can be introduced by means of
the conditions
\begin{align}
&E_b-X_b\ \in\ \oplus_{c\succ b}\Q(q,t) X_c,\
\langle E_b, X_{c}\rangle_\circ = 0 \for B \ni c\succ b.
\label{macd}
\end{align}
They are well defined
because the  pairing
is nondegenerate (for generic $q,t$) and form a
basis in $\Q(q,t)[X_b]$.

This definition is due to Macdonald (for
$k_{\sht}=k_{\lng}\in \Z_+ $), who extended
Opdam's nonsymmetric polynomials introduced
in the differential case in \cite{O2}
(Opdam mentions Heckman's unpublished lectures in \cite{O2}).
The general (reduced) case was considered in \cite{C4}.

Another approach is based on the $Y$\~operators.
We continue using the same notation $X,Y,T$
for these operators acting in the polynomial
representation.

\begin{proposition}
The polynomials $\{E_b, b\in B\}$
are unique (up to proportionality) eigenfunctions of
the operators $\{L_f\equal f(Y_1,\ldots, Y_n),
f\in \Q[X]\}$
acting in $\Q_{q,t}[X]:$
\begin{align}
&L_{f}(E_b)\ =\ f(q^{-b_\#})E_b\, \hbox{\ for\ }\,
b_\#\equal b- u_b^{-1}(\rho_k),
\label{Yone} \\
& X_a(q^{b})\ =\
q^{(a,b)}\ ,\hbox{where\ } a,b\in B,\
 u_b=\pi_b^{-1}b, \label{xaonb}
\end{align}
$u_b$ is from Proposition \ref{PIOM},\
$b_\#=\pi_b(\!(-\rho_k)\!)$.
\label{YONE}
\end{proposition}
\qed

The coefficients of the Macdonald
polynomials are rational functions in terms of $q_\nu,t_\nu$
(here either approach can be used).

\medskip
{\bf Symmetric polynomials}.
Following Proposition \ref{YONE}, the
{\dfont symmetric Macdonald polynomials\,} $\ P_b=P_b^{(k)}\ $
can be introduced as
eigenfunctions of the $W$\~invariant
operators $L_f=$ $f(Y_1,\ldots,Y_n)$
defined for symmetric, i.e.,  $W$\~invariant,
polynomials $f$ as follows:
\begin{align}
L_{f}&(P_{b_-})=f(q^{-b_-+\rho_k}) P_{b_-},\ \
b_-\in B_-,\notag\\
&P_{b_-}=\sum_{b\in W(b_-)}X_b \mod
\oplus_{c_-\succ b_-}\Q(q,t) X_c,\
\label{Lf}
\end{align}
Here it suffices to take the
{\dfont monomial symmetric functions\,}, namely,
$f_b=\sum_{c\in W(b)}X_c$ for $b\in B_-.$
For minuscule $b=-\om_r$ and $b=-\vth,$
the restrictions of the operators
$L_{f_b}$ to symmetric polynomials
are the
{\dfont Macdonald operators\,} from \cite{M2,M3}.

The $P$\~polynomials
are pairwise orthogonal with respect to
$\langle\ ,\ \rangle_\circ$
as well as $\{E\}.$
Since they are $W$\~invariant,  $\mu$ can be replaced
by the symmetric measure\~function due to Macdonald:
\begin{align}
&\de\ = \de^{(k)}=\prod_{\al \in R_+}
\prod_{i=0}^\infty \frac{(1-X_\al q_\al^{i})
(1-X_\al^{-1}q_\al^{i})
}{
(1-X_\al t_\al q_\al^{i})
(1-X_\al^{-1}t_\al^{}q_\al^{i})}.\
\label{delp}
\end{align}
The corresponding pairing remains $\ast$-hermitian because
the function
$\de_\circ$ is $\ast$\~invariant.

These polynomials were introduced in \cite{M2,M3}.
They were used for the first
time in Kadell's unpublished work (classical root systems).
In the one-dimensional case, they are due to Rogers.

The connection between $E$ and $P$ is as follows
\begin{align}
&P_{b_-}\ =\ \P_{b_+} E_{b_+}, \ b_-\in B_-,\
b_+=w_0(b_-),
\notag\\
&\P_{b_+}\equal\sum_{c\in W(b_+)}
\prod_\nu t_\nu^{l_{\nu}(w_c)/2} T_{w_c}, \where
\label{symmetr}
\end{align}
$w_c\in W$ is the element of the least length such that
$c=w_c(b_+)$;
the notation $\P$ will be used if the summation is over all $w$. 
See \cite{O2,M4,C4}. 

\smallskip

\subsection{Using intertwiners}
The $Y$\~intertwiners serve as creation operators
in the theory of nonsymmetric Macdonald polynomials.
They are defined as follows:
\begin{align}\label{tauintery}
\Psi_i=\tau_+(T_i)+
\frac{t_i^{1/2}-t_i^{-1/2}}{Y_{\al_i}^{-1}-1},\ i\ge 0,\ \
P_r=\tau_+(\pi_r),\ r\in O',\\
\Psi_{\hw}=P_r\Psi_{i_l}\ldots\Psi_{i_1}\hbox {\ for\
reduced\ decompositions\ }
\hw=\pi_rs_{i_l}\ldots s_{i_1}.\notag
\end{align}
Recall that
$\tau_+(T_0)=X_0^{-1}T_0^{-1}=q^{-1}\,X_\vth T_0^{-1},$
$\tau_+(\pi_r)=q^{-(\om_r,\om_r)/2}X_r\pi_r.$
The elements $\Psi_{\hw}$
are the images $\tau_-\si(\Phi_{\hw})$
of the $X$\~intertwiners $\Phi_{\hw}.$

Indeed, thanks to the relations $\si(X_b)=Y_b^{-1},$
$\tau_-(Y_b)=Y_b:$
\begin{align*}
&\tau_-\si(\Phi_i X_b)=\tau_-\si(X_{c} \Phi_i)\
\Rightarrow\
\tau_-\si(\Phi_i)(Y_b^{-1})=Y_c^{-1}\tau_-\si(\Phi_i)
\end{align*}
for $c=s_i(b).$
However,
$$\tau_-\si(\Phi_i)=\tau_+(T_i)+
(t_i^{1/2}-t_i^{-1/2})(Y_{\al_i}^{-1}-1)^{-1}
$$
because $\tau_-\si=\tau_+\tau_-^{-1}$ and $\tau_-$
preserves $T_i.$ Similarly,
$\tau_-\si(\pi_r)$ $=\tau_+(\pi_r).$
We obtain that
$$
\Psi_{\hw}Y_b=Y_{\hw(b)}\Psi_{\hw}, \hbox{\ where\ }
Y_{[b,j]}\equal Y_b q^{-j}.
$$
The same property holds for
the normalized intertwiners
$$F_{\hw}\equal\tau_-\si(G_{\hw}),
$$ which are {\em unitary},
i.e., induce a homomorphism from $\hW$ to a proper localization
of $\HH^\flat.$

The automorphism $\tau_-$ acts in $\v$ and commutes with
the $Y$-operators.
The following proposition describes its action on the
$\Psi$-intertwiners.

\begin{proposition}\label{TAUPSI}
(i) For generic $q,t$ (or for arbitrary
$q,t$ provided that the polynomial
$E_b$ for $b\in B$ is well defined),
\begin{align}\label{taumineb}
\tau_-(E_b)=q^{-\frac{(b_-\,,\,b_-)}{2}+(b_-\,,\,\rho_k)}\,E_b
\for P_-\ni b_-\in W(b).
\end{align}

(ii) For any $q,t$ and $Y_0=q^{-1}Y_{\vth}^{-1},$
\begin{align}\label{tauminpsi}
\tau_-(T_i)&=T_i\, (i>0),\ \ \, \tau_-(\tau_+(T_0))=
\tau_+(T_0)^{-1}Y_0,\\
\tau_-(\Psi_i)&=\Psi_i (i>0),\ \
\tau_-(\Psi_0)=\Psi_0 Y_0=Y_0^{-1} \Psi_0,\notag \\
\tau_-(\tau_+(\pi_r))&=q^{(\om_r,\om_r)/2}Y_r\tau_+(\pi_r)
=q^{-(\om_r,\om_r)/2}\tau_+(\pi_r)Y_{r^*}^{-1}.\notag
\end{align}
\end{proposition}
{\em Proof.} These claims can be checked using
directly the definitions.
One can also identify $\tau_-$ in $\v$  with
the operator of multiplication by
$\widetilde{\tau_-}=\ga_y(1)^{-1}\,\ga_y$ for the
{\em $Y$-Gaussian} defined as follows:
\begin{align}\label{gausery}
&\ga_y\equal \sum_{b\in B} q^{(b,b)/2}Y_b,\
\ga_y(1)=\sum_{b\in B} q^{(b,b)/2}q^{(b,\rho_k)}.
\end{align}
Here we assume that $0<q<1$ and use that
$\v$ is a union of finite dimensional spaces
preserved by the $Y$-operators.

The irreducibility of $\v$ for generic $q,t$ is sufficient
to conclude that $\widetilde{\tau_-}=\tau_-\,$ because
the conjugation by $\widetilde{\tau_-}$ induces
$\tau_-$ in $\HH^\flat.$
This implies their coincidence
for generic $q$ (apart from roots of unity) and $t_\nu$,
which results in (\ref{tauminpsi}) for {\em arbitrary} $q,t.$
\sq

Setting
\begin{align}
\Psi_i^b=\Psi_i(q^{b_{\#}})=&
\tau_+(T_i) + (t_i^{1/2}-t_i^{-1/2})
(X_{\al_i}(q^{b_{\#}})-1)^{-1},
\label{Phijb}
\\
F_i^b=(\Psi_i\psi_i^{-1})(q^{b_{\#}})=&
\frac{ \tau_+(T_i) + (t_i^{1/2}-t_i^{-1/2})(X_{\al_i}
(q^{b_{\#}})-1)^{-1}
}{
t_i^{1/2} + (t_i^{1/2} -t_i^{-1/2})(X_{\al_i}
(q^{b_{\#}})-1)^{-1} },
\label{Gjb}
\end{align}
we come to the Main Theorem 5.1 from \cite{C1}.

\begin{theorem}
Given $c\in B,\ 0\le i\le n$ such that
$(\al_i, c+d)> 0,$
\begin{align}
&E_{b}q^{-(b,b)/2} \ =\ t_i^{1/2} \Psi_i^c(E_c)
 q^{-(c,c)/2}
\for b= s_i\llb c\rrb.
\label{Phieb}
\end{align}
If  $(\al_i, c+d)=0,$ then
\begin{align}
&\tau_+(T_i) (E_c) \ =\ t_i^{1/2} E_c, \ 0\le i\le n,
\label{Tjeco}
\end{align}
which  results in the relations $s_i(E_c)=E_c$ as
$i>0$.  For $b=\pi_r\llb c\rrb,$ where
the indices $\,r\,$ are from $O',$
\begin{align}
&q^{-(b,b)/2+(c,c)/2}E_b\ =\ \tau_+(\pi_r)(E_c)\ =\
X_{\om_r}q^{-(\om_r,\om_r)/2}\pi_r(E_c).
\label{pireb}
\end{align}
Also $\tau_+(\pi_r)(E_c)\neq E_c$ for $\pi_r\neq$id,
since $\pi_r\llb c\rrb\neq c$ for any $c\in B$
due to Lemma \ref{NOPIRB}.
\sq
\label{PHIEB}
\end{theorem}
\smallskip

Using the standard Bruhat ordering and the set $\b(\hw)$ from
Proposition \ref{BSTAL}, the theorem results in the following.
\begin{corollary} Given a reduced decomposition
$\pi_b=\pi_r s_{i_l}\cdots s_{i_1}$ and $\hw'\in \b(\pi_b)$,
we define $T_{\hw'}$ for the corresponding, possibly non-reduced,
decomposition of $\hw'$. Then
\begin{align}
&\tau_+(T_{\hw'})(1) \in \oplus_{\pi_{b'}\in \b(\pi_b)}
\Q_{q,t}\,E_{b'}
\subset \oplus_{a\succeq b} \Q_{q,t}\,X_{a}.
\label{Tjecoy}
\end{align}
\end{corollary}
{\em Proof.}
We use (\ref{tauintery}), that is 
$$
\Psi_i=\tau_+(T_i)+
\frac{t_i^{1/2}-t_i^{-1/2}}{Y_{\al_i}^{-1}-1},\
P_r=\tau_+(\pi_r),
$$
and that $\{E_c\}$ are eigenfunctions of the $Y$\~operators.
\sq
\smallskip

We can now renormalize the $E$\~polynomials as follows:
\begin{align}
&\hE_{b}\ \equal\ \tau_+(G_r G_{i_l}^{c_l}\ldots
G_{i_1}^{c_1})(1),
\where
\label{ehatb}
 \\
&c_1=0, c_2=s_{i_1}\llb c_1\rrb,\ldots,
c_l=s_{i_l}\llb c_{l-1}\rrb\,
\for \pi_b= \pi_r s_{i_l}\ldots s_{i_1}.
\notag \end{align}
These polynomials
do not depend on the particular choice of
the decomposition of $\pi_b$ (not necessarily
reduced) and are proportional to $E_b$ for all
$b\in B$:
\begin{align}
E_{b}q^{-(b,b)/2} \ =&\ \prod_{1\le p\le l}
\bigl(t_{i_p}^{1/2}\phi_{i_p}(q^{c_p})\bigr)\
\hE_b\notag\\
=&\ \prod_{[\al,j]\in \la\,'\,(\pi_b)}
\Bigl(
\frac{
1- q_\al^{j}t_\al X_\al(q^{\rho_k})
 }{
1- q_\al^{j} X_\al(q^{\rho_k})
}
\Bigr)\ \hE_b.
\label{ebebhat}
\end{align}
The polynomials $\hE_b$ are
directly connected with the spherical polynomials.
\smallskip

\subsection{Spherical polynomials}
The following renormalization
of the $E$-polynomials is of major importance in
the Fourier analysis (see \cite{C4}):
\begin{align}
\e_b\ \equal&\ E_b(X)(E_b(q^{-\rho_k}))^{-1},\where
\notag\\
E_{b}(q^{-\rho_k}) \ =&\ q^{(\rho_k,b_-)}
\prod_{[\al,j]\in \la\,'\,(\pi_b)}
\Bigl(
\frac{
1- q_\al^{j}t_\al X_\al(q^{\rho_k})
 }{
1- q_\al^{j}X_\al(q^{\rho_k})
}
\Bigr),
\label{ebebs}\\
\la\,'\,(\pi_b)\ =&\
\{[\al,j]\ |\  [-\al,\nu_\al j]\in \la(\pi_b)\}.
\label{jbseto}
\end{align}
We call them
{\dfont spherical polynomials\,}.
Explicitly
(see (\ref{lambpi})),
\begin{align}
\la\,'\,(\pi_b)\ =& \{[\al,j]\, \mid\,\al\in R_+,
\label{jbset}
\\
&-( b_-, \al^\vee )>j> 0
\iif  u_b^{-1}(\al)\in R_-,\notag\\
&-( b_-, \al^\vee )\ge j > 0 \iif
u_b^{-1}(\al)\in R_+ \}.
\notag \end{align}
Formula (\ref{ebebs}) is the Macdonald
{\dfont evaluation conjecture}
in the nonsymmetric variant from \cite{C4}.
See \cite{C3} for the symmetric evaluation conjecture.

Note that one has to consider only long $\al$ (resp., short)
if $k_{\sht}=0$ (resp.,
$k_{\lng}=0$) in the $\la\,'\,$-set.
All formulas involving $\la$ or $\la\,'\,$ have to be modified
correspondingly in such case.

We have the following {\dfont duality formula}
for $b,c\in B\, :$
\begin{align}
&\e_b(q^{c_{\#}})\ =\ \e_c(q^{b_{\#}}),\
b_\# = b- u_b^{-1}(\rho_k),
\label{ebdual}
\end{align}
that is the main justification of the definition of $\e_b$.

Combining (\ref{ebebhat}) and (\ref{ebebs}),
we conclude that
\begin{align}
&\e_{b}\ =\ q^{(-\rho_k+b_-,b_-)/2} \hE_b.
\label{ebebeq}
\end{align}
See also \cite{C1}.
\smallskip

The proof of the duality formula
(\ref{ebdual}) is based on
the anti-involution $\phi$ from
(\ref{starphi}):
$$
\phi:
X_b\mapsto Y_b^{-1}\mapsto X_b,\
T_i\mapsto T_i\ (1\le i\le n),\ q\mapsto q,t\mapsto t.
$$
Following \cite{C3,C4} (see also (\ref{valofh}) below),
we define the {\dfont evaluation pairing}
for $f,g\in \v$,

\begin{align}
 &\{f,g\}\ =\ \{L_{\imath(f)}(g(X))\}\ =\
 \{L_{\imath(f)}(g(X))\}(q^{-\rho_k}),
\label{symfg}
\\
&\imath(X_b)\ =\ X_{-b}\ =\ X_b^{-1},\
\imath(z)\ =\  z \for
z\in \Q_{q,t}\, ,
\notag \end{align}
where $L_f$ is from Proposition \ref{YONE}.
This pairing is symmetric and induces $\phi$ in $\HH^\flat$.
\smallskip

As an application of (\ref{ebebeq}),
we obtain that, given $b\in B,$
the spherical polynomial
$\e_b$
is well defined for $q,t\in \C^*$ if
\begin{align}
&\prod_{[\al,j]\in \la\,'\,(\pi_b)}
\bigl(
1- q_\al^{j}t_\al X_\al(q^{\rho_k})\bigr)\ \neq\ 0.
\label{ebebnze}
\end{align}

Similarly (see \cite{C1}, Corollary 5.3),
the polynomial $E_b$ exists if
$$
\prod_{[\al,j]\in \la\,'\,(\pi_b)}
\bigl(
1- q_\al^{j}X_\al(q^{\rho_k})\bigr)\neq 0.
$$
If $b\in B_-$ and the latter inequality holds for
$b_+=w_0(b)\in B_+,$
then the symmetric polynomials
$P_b$ is well defined.
\medskip

{\bf Proof of the evaluation formula.}
Another application of the intertwiners is a different approach
to the evaluation formula (\ref{ebebs}); it is especially
important when $q,t$ are arbitrary (non-generic)
and $\v$ is not supposed semisimple.
Note that the other two known proofs of
the evaluation formula are based respectively
on the duality and the technique of shift operators (see
\cite{C3,C4}). The following proposition readily gives
(\ref{ebebs}).

\begin{proposition}\label{COREBEC}
(i) Let us assume that $E_c'$ is a $Y$\~eigenvector
satisfying (\ref{Yone}) for $c_\#$.
We will introduce $E_b'$ using either
formula (\ref{Phieb}) for $b=s_i(\!(c)\!)$, where  $i\ge 0$
and $l(\pi_b)=1+l(\pi_c)$,
or using (\ref{pireb}) for
$b=\pi_r(\!(c)\!).$
Then $E_b'$ satisfies (\ref{Yone}) for $b_\#$ and
\begin{align}
&q^{-(\rho_k,b_-)}E_{b}'(q^{-\rho_k})  =
q^{-(\rho_k,c_-)}E_{c}'(q^{-\rho_k})
\frac{
1- q_\al^{(\tal^\vee\,,\,c_- +d)}\,t_\al X_\al(q^{\rho_k})
 }{
1- q_\al^{(\tal^\vee\,,\,c_- +d)}\, X_\al(q^{\rho_k})
},
\label{ebebc}\\
&\tal=u_c(\al_i),\ \ \,
\al=u_c(\al_i) \hbox{\ for\ } i>0,
\ \ \, \al=u_c(-\vth) \hbox{\ for\ } i=0,\notag\\
& q^{-(\rho_k,b_-)}E_{b}'(q^{-\rho_k})  =
q^{-(\rho_k,c_-)}E_{c}'(q^{-\rho_k}) \for b=\pi_r(\!(c)\!).
\label{ebebpi}
\end{align}

(ii) More generally, let $\widetilde{E}=\Psi_i(E)$
for $E\in\v$, $\,i\ge 0$, assuming that $(Y_{\al_i}-1)^{-1}(E)$
is well defined.
Then
\begin{align}\label{eprimerho}
&\widetilde{E}(q^{-\rho_k}) =
((\frac{t_i^{1/2}Y_{\al_i}^{-1}-t_i^{-1/2}}
{Y_{\al_i}^{-1}-1})(E))(q^{-\rho_k}) \for i>0,\\
&\widetilde{E}(q^{-\rho_k}) =
((\frac{t_0^{1/2}Y_{\al_0}^{-1}-t_0^{-1/2}}
{Y_{\al_0}^{-1}-1}Y_\vth^{-1})(E))(q^{-\rho_k}) \for i=0\notag.
\end{align}
If $\widetilde{E}=P_r(E)\ $, where $P_r=\tau_+(\pi_r)$, then
$\ \widetilde{E}(q^{-\rho_k}) = (Y_r^{-1}(E))(q^{-\rho_k})$.
\end{proposition}

{\em Proof.} We will move $\Psi_i^c$ in
$\{1,E_b'\}=E_{b}'(q^{-\rho_k})$  to the first component
using (\ref{etatxpi}) and then back. In more detail,
we need the following
relations:
\begin{align}
&\phi(\tau_+(T_i))=(\star\cdot\tau_+\cdot\eta)(T_i)=
\tau_+(T_i),\ \,\phi(P_r)=X_r u_r^{-1}.\label{phitaut}
\end{align}
Recall that $\tau_+(T_i)=T_i$ for $i>0$ and
$\tau_+(T_0)=X_0^{-1}T_0^{-1}$ for $X_0=qX_{\vth}^{-1}$.

The $X$\~rational functions $\phi(\Psi_i^c)(1)$ and
$\phi(P_r)(1)$ can be readily calculated. Then we move them
back to the second component, replacing $X$ by $Y^{-1}$ and
applying the result to $E_c'$, which is a $Y$\~eigenvector.
\smallskip

Claim (ii) is a natural generalization of (i). 
First, using
(\ref{phitaut}),
\begin{align}
&\phi(\Psi_i)=\phi(\tau_+(T_i)+
\frac{t_i^{1/2}-t_i^{-1/2}}{Y_{\al_i}^{-1}-1})=
\tau_+(T_i)+\frac{t_i^{1/2}-t_i^{-1/2}}{X_{\al_i}-1}.
\label{phipsi}
\end{align}
Second, the images of the $\Psi$\~intertwiners
in the polynomial representation are
\begin{align}\label{phiPsi}
&\phi(\Psi_i)(f)=\frac{t_i^{1/2}X_{\al_i}-t^{-1/2}}
{X_{\al_i}-1}s_i^x(f),\ \phi(P_r)(f)=\pi_r^x(f),\where \\
&s_i^x(f)=s_i(f)\hbox{\ for\ } i>0,\,
s_0^x(f)=X_\vth s_\vth(f),\ \pi_r^x(f)=X_r u_r^{-1}(f).\notag
\end{align}
Let us give an explicit demonstration for $i=0$:
\begin{align*}
&\phi(\Psi_0)=X_0^{-1}(\frac{t_0^{1/2}X_0-t_0^{-1/2}}{X_0-1}s_0)\\
-&X_0^{-1}\frac{(t_0^{1/2}-t_0^{-1/2})X_0}{X_0-1}
+\frac{t_0^{1/2}-t_0^{-1/2}}{X_0-1}\\
=&\frac{t_0^{1/2}X_0-t_0^{-1/2}}{X_0-1}X_0^{-1}s_0=
\frac{t_0^{1/2}X_0-t_0^{-1/2}}{X_0-1}s_0X_0.
\end{align*}

Third, let $\widetilde{E}=\Psi_i(E)$. We set
$\de_{i0}=1,\,0$ respectively for $i=0$ or $i>0$. Then
the evaluation $\widetilde{E}(q^{-\rho_k})$ equals
\begin{align*}
&\{(\frac{t_i^{1/2}X_{\al_i}-t_i^{-1/2}}
{X_{\al_i}-1}s_i^x)(1),E\}=
\{(\frac{t_i^{1/2}X_{\al_i}-t_i^{-1/2}}
{X_{\al_i}-1}X_\vth^{\de_{i0}})(1),E\}\\
=&\{1,(\frac{t_i^{1/2}Y_{\al_i}^{-1}-t_i^{-1/2}}
{Y_{\al_i}^{-1}-1}Y_\vth^{-\de_{i0}})(E)\} =
(\frac{t_i^{1/2}Y_{\al_i}^{-1}-t_i^{-1/2}}
{Y_{\al_i}^{-1}-1}Y_\vth^{-\de_{i0}})(E)(q^{-\rho_k}).
\end{align*}
The evaluation of $\widetilde{E}=P_r(E)$ is calculated similarly.
\sq
\smallskip

A variant of the same construction gives that
\begin{align}\label{twerho}
&(T_w(E))(q^{-\rho_k})\ =\ \prod_\nu t_\nu^{l_{\nu}(w)/2}\,
E(q^{-\rho_k})
 \for E\in \v,\ w\in W,
\end{align}
where $l(w)$ is the length of {\em nonaffine} $w$.

This formula can be readily used to deduce the formula for
$P_{b}(q^{-\rho_k})$
for the symmetric Macdonald polynomials $\,P_b\, (b\in B_-)\,$
from that for the
nonsymmetric polynomials. See (\ref{symmetr}).
\medskip

\comment{
%
%


The paper is mainly devoted to the irreducibility
of the polynomial representation of the Double
affine Hecke algebra, DAHA, for arbitrary 
reduced root systems and generic ``central charge" $q$. 
Part I is devoted to combinatorial and general 
algebraic aspects of the theory.               

The technique of intertwiners from \cite{C1} in a
non-semisimple variant is the main tool.
It is important for the
decomposition of the polynomial representation
in terms of irreducible DAHA modules and for its
weight decomposition. We focus on the principal
aspects of the technique of intertwiners and discuss
only basic (and instructional) applications.
Generally, it is more efficient to combine 
the intertwiners with other approaches, to be 
considered in further papers.
\smallskip

There are several methods that can be used now in the
study of the polynomial representation of DAHA and
its degenerations. Certainly the localization functor
(the KZ\~monodromy) from \cite{GGOR,VV1}
must be mentioned, as well as   
the geometric methods of \cite{VV3} 
and the parabolic induction from recent \cite{BE}.
The technique of intertwiners provides constructive
(relatively elementary)  
tools for managing the irreducibility of the
polynomial representation and its constituents 
for any $q,t$ based on combinatorics of affine Weyl groups.
    
\smallskip
An important general objective of this technique
is finding a counterpart of the classical theory of 
{\em highest vectors} for DAHA (and AHA),
complementary to the geometric method of \cite{KL1}.
It involves difficult combinatorial problems 
and is known only for type $A$ and in some cases of small ranks.
However, the geometric DAHA methods are far from simple too
and explicit theory of DAHA modules is needed in quite a few
applications (see Section \ref{sec:EXPECTEDAP}).     
\medskip

{\bf Main constructions.}
The {\em polynomial representation}, denoted by $\v$ 
in the paper, is well known to be irreducible and semisimple
for generic values of the DAHA-parameters $q$ and $t=q^k.$
It becomes reducible either when $q$ is a root
of unity or for generic $q$ and special $t.$

We perform a complete analysis of the
irreducibility and semisimplicity of $\v$ 
for generic $q.$ Another application
is a construction of the canonical semisimple
submodule in $\v$ generalizing that of type $A$
from \cite{FJMM} (the symmetric variant) 
and \cite{Ka} (the non-symmetric case). 
 
\smallskip

We begin with a description of all {\em singular} 
$t=q^k$ making (by definition) 
the radical of the evaluation pairing 
nonzero. The answer is given in terms of the principal 
values (at $t^{-\rho}$) of the nonsymmetric Macdonald
polynomials \cite{C3,C4,C100} 
for the weights sufficiently large to ensure the
existence of these polynomials. 
The same answer can be obtain using a generalization 
of the method from \cite{O1,O3} in the 
rational case based on the {\em shift-operator}. 
The latter is used to calculate the principal value  of
the $t$\~discriminant $D_Y$ applied to the 
$t$\~discriminant $D_X$ in $\v$, 
where $X,Y$ are the generators of DAHA.   
The approach via the Macdonald polynomials
has no rational counterpart.

The {\em evaluation pairing} is defined as follows: 
$$
\{E,F\}\,=\, E(Y^{-1})(F(X))(t^{-\rho}),\ E,F \in \v\,;
$$
for instance, $D_Y(D_X)(t^{-\rho})=\{D_{X^{-1}},D_X\}$.
In the simply-laced case,
the radical $Rad$ of this pairing is zero if 
and only $\v$ is irreducible. 
This equivalence becomes more subtle in the 
non-simply-laced case, as well as the formula for 
$D_Y(D_X)(t^{-\rho})\,$.
The $q,t$\~theory provides the best (and direct)
method for calculating this formula, including 
managing the rational case through the limiting 
procedure from DAHA to its rational degeneration.
\smallskip

The main objects of this paper are
the {\em chains of the intertwiners\,} and  
{\em non-semisimple Macdonald polynomials\,} 
(defined via such chains).
We define a system of subspaces in $\v$
with the Macdonald polynomials as top elements
in a punctured neighborhood of a singular $t$\,
and then extend this construction to singular $t$.
\vfill
\eject
}

\section{Shift-operator} \label{SHIFTOP}
\setcounter{equation}{0}
In this section, we take $B=P$.
By $[H]_\dagger$, we mean the restriction of the
image of the element $H$ acting in the polynomial
representation $\v$ to the subspace $\v^W$ of
$W$\~invariant Laurent polynomials. The image
of $[H]_\dagger$ is in $\v^W$ if $H$ is
$W$\~{\em invariant}, i.e., belongs to
$$\HH^W\equal\{A\in \HH\, \mid\, T_iA=AT_i,\for i>0\}.
$$

We fix a subset $v\subset \nu_R$ and introduce the
{\dfont shift operator}:
\begin{align}
&\s_v^t \ =\
(\x_v^t) ^{-1}\y_v^t,\
\S_v^{ q,t}= [\s_v^t]_\dagger=
\ (\x_v^t )^{-1}[\y_v^t]_{\dagger},
\label{shiftoper}\\
&\x_v^t \equal \prod_{\al\in R_+}^{\nu_\al\in v}
((t_\al X_{\al})^{1/2}-
(t_\al X_{\al})^{-1/2}),\notag\\
&\y_v^t  = \prod_{\al\in R_+}^{\nu_\al\in v}
(t_\al Y_{\al}^{-1})^{1/2}-
(t_\al Y_{\al}^{-1})^{-1/2}).\notag
\end{align}
Note that the
elements $\x_{v}^t,\ \y_{v}^t $
belong respectively to $\Z_t[X_b],\
\Z_t [Y_b]$ for $b\in P.$

\subsection{Shift--formula}
The following proposition is from \cite{C2},
\cite{C12}. The $W$\~invariant operators
$L_f^{q,t}$ for $W$\~invariant $f$  and the
{\em symmetric} Macdonald polynomials $P_b^{q,t}(X)$
defined in (\ref{Lf}) are used.

\begin{proposition}\label{SHIFTOPER}
The operators $\s_{v}^t$ and
its restriction $\S_v^{ q,t}$
to $ \Q_{ q,t}[X]^W$ are $W$-invariant and
preserve the latter space. Provided that $t_{\nu}=1$
whenever
$\nu\not\in v,$
\begin{align}\label{shiGL}
& \S_v^{ q,t}L_f^{ q,t}\ =\ L_f^{ q, tq_v} \S_v^{ q,t}
\for f\in \C[X]^W,
\end{align}
where $ tq_v=\{t_\nu q^\nu , t_{\nu\,'\,}\}$
for $\nu\in v\not\ni \nu\,'\,,$
\begin{align}\label{SPg}
&\S_v^{ q,t} (P_{b}^{ q,t})= g_v^{ q,t}(b)
P_{b+\rho_v}^{ q, tq_v}, \for\\
&g_v^{ q,t}(b)\equal
\prod_{\al\in R_+,\nu_\al\in v}
(X_\al(q^{(\rho_k-b)/2}) -
t_\al X_\al(q^{(b-\rho_k)/2})), \notag\\
\end{align}
where $\rho_v=\sum_{\nu\in v}\rho_\nu,
\ P_c=0 \for c\not\in P_-.$
\end{proposition}
\smallskip

We will also need the element
$$ \overline{\y}_v^t \equal \prod_{\al\in R_+,\nu_a\in v}
((t_\al Y_{\al})^{1/2}-
(t_\al Y_{\al})^{-1/2});$$
it is $\y$ with $Y_\al^{\pm1/2}$ replaced by
$Y_\al^{\mp1/2}$ (or $\x$ with $Y$ instead of $X$). 

The product $\overline{\y}_v^t \y_v^t$ is obviously a
$W$-invariant polynomial. So we can apply the formulas
in (\ref{Lf}):
\begin{align}\label{bargg}
&(\overline{\y}_v^t {\y}_v^t) (P_b^{q,t}(X))=
\overline{g}_v^{q,t}(b)g_v^{q,t}(b) \for\\
&\overline{g}_v^{ q,t}(b)\equal
\prod_{\al\in R_+,\nu_\al\in v}
(t_\al^{-1}X_\al(q^{(b-\rho_k)/2}) -
 X_\al(q^{(\rho_k-b)/2})) \notag\\
=&\prod_{\al\in R_+,\nu_\al\in v}
 t_\al^{-1}q^{(\al,(b-\rho_k))/2}
\prod_{\al\in R_+,\nu_\al\in v}
(1- t_\al q^{(\al,\rho_k-b)})\notag\\
=&q^{(\rho_v,(b-\rho_k))}
\prod_{\nu\in v}t_\nu^{-\hbox{\tiny\kapp}_\nu}
\prod_{\al\in R_+,\nu_\al\in v}
(1- t_\al q^{(\al,\rho_k-b)}),\notag
\end{align}
where $\kapp_\nu$ is the number of
positive roots with $\nu_\al=\nu.$

Note the formula
\begin{align*}
&\overline{g}_v^{q,t}(b)g_v^{q,t}(b)=
\prod_{\al\in R,\nu_\al\in v}
(t_\al^{1/2}X_\al(q^{(b-\rho_k)/2}) -
t_\al^{-1/2} X_\al(q^{(\rho_k-b)/2})),
\end{align*}
where the product is over {\em all} roots.

We are going to employ
the duality and the
Macdonald (symmetric) evaluation conjecture,
the formula for the ``principle value" of the
Macdonald polynomials proved in \cite{C12} (in the $A$-case,
both were verified by Koornwinder).
They read as:
\begin{align}\label{PbPc}
&P_b(q^{c-\rho_k})P_c(q^{-\rho_k})\ =\
P_c(q^{b-\rho_k})P_b(q^{-\rho_k}), \ b,c\in P_-,\\
P_b(q^{-\rho_k})& = P_b(q^{+\rho_k}) =
X_b(q^{\rho_k})\prod_{\al\in R_+}
\prod_{ j=1}^{-(\al^\vee,b)}
\Bigl(\frac{1- q_\al^{j-1}t_\al X_\al(q^{\rho_k})}
{1- q_\al^{j-1}X_\al(q^{\rho_k})}\Bigr).\notag
\end{align}

As a matter of fact, the second formula follows from the first
via the Pieri rules. See \cite{C12}. The derivation
of the evaluation formula from the Pieri rules in the
case of $A_n$ is due to Koornwinder.
\smallskip

{\bf Norm-formulas via the evaluation ones}.
It is important to note
that the duality also results in the norm-formulas,
including the celebrated Macdonald constant term conjecture.

The simplest known deduction of the norm-formulas
is based on the intertwining
operators acting on the {\em spherical polynomials\,}
that are the {\em nonsymmetric} Macdonald
polynomials taking the value $1$ at $q^{-\rho_k}$.
The intertwiners preserve the latter normalization up
to certain simple factors. They act changing the norms 
of $E$\~polynomials, but this is simple to control. 
Thus, using intertwiners we naturally
come to the norm-formulas in the {\em spherical
normalization}, i.e., for spherical polynomials.

Now one can use the evaluation
formulas (principle value formulas) to get the
norm-formulas in the standard monic normalization of
the leading terms of the $E$\~polynomials.
The last step is reproving the {\em symmetric} 
Macdonald norm conjecture
(stated in \cite{M1} and proven in \cite{C2});
we apply the symmetrization and 
use the formulas for  $\{P_{b_-}\}$
in terms of $\{E_b\}$.

This procedure gives a transparent deduction
of the norm formulas from the evaluation ones. 
However it does
not clarify why the norms appear the products 
{\em very much
similar} to those in the definition of the Macdonald measure.

The ``conceptual" reason for this coincidence 
(which is a result of a straight calculation in
the above deduction) is as follows. 

The above calculation in terms of the intertwining operators
is actually {\em equivalent} to the calculation of 
the {\em Fourier transform},
acting from $\v$ to the {\em $\HH^\flat$\~representation in 
delta-functions}. It
sends the spherical polynomials exactly to the corresponding
{\em delta-functions\,} defined with respect to the discretization
of the Macdonald measure (with the natural normalization).
It was established in \cite{C1} (see also \cite{C4,C101}).
Therefore the norms of the
spherical polynomials practically coincide with the
{\em values\,} of the Macdonald measure at the corresponding 
weights. This gives a complete clarification of
the constant term conjecture and concludes (the corresponding
part of) \cite{M1}.

\medskip
\subsection{Formula for 
\texorpdfstring{{\mathversion{bold}$\overline{\y}(\x)$}}{Y(X)}}
Later on, the shift operator will be used only
when the set $v$ is the whole $\nu_R.$
The suffix $v$ will be dropped in the formulas,
the indices $\nu$ are arbitrary from
$\nu_R.$
For instance,
$tq_v=tq=\{t_\nu q_\nu,\, \nu\in\nu_R\}.$
We will also use $\rho_\nu^\vee=(1/2)\sum_{\nu_\al=\nu}\al^\vee:$
\begin{align}\label{trhovee}
&t^{\rho^\vee}=\prod_{\nu\in \nu_R}t_\nu ^{\rho_\nu^\vee},\
t^{(\al,\rho^\vee)}=\prod_{\nu\in \nu_R}t_\nu ^{(\al,\rho_\nu^\vee)}
=q^{(\al,\rho_k)},
\\
&t^{\kapp}=\prod_{\nu\in \nu_R}t_\nu ^{\hbox{\tiny\kapp}_\nu} \for
\kapp_\nu=|\{\al\in R_+,\nu_\al=\nu\}|.\notag
\end{align}

\begin{maintheorem} \label{YOFXEV}
Let $v=\nu_R$ and  $b=-\rho.$ Using the notation
$q^{(\rho_k-b)}=(tq)^{\rho^\vee}=
\prod_{\nu}(t_\nu q_\nu)^{\rho_\nu^\vee},$ and
$$
(tq)^{(\,\cdot\,,\rho^\vee)}=
\prod_\nu (t_\nu q_\nu)^{(\,\cdot\,,\rho_\nu^\vee)},
\hbox{\ e.g.,\ } (tq)^{(\rho,\rho^\vee)}=
\prod_\nu (t_\nu q_\nu)^{(\rho,\rho_\nu^\vee)}=
q^{(\rho,\rho+\rho_k)},
$$
we come to the relation
\begin{align}\label{yofx}
&\overline{\y}^{t}(\x^{t})=C_{q,t}P_{b}^{q,t} \for
C_{q,t}=\overline{g}^{q,t}(b)\\
=&\prod_{\al\in R_+}
(t_{\al}^{-1}(tq)^{-(\al,\rho^\vee)/2}-
(tq)^{(\al,\rho^\vee)/2})\notag\\
=&t^{-\kapp}(tq)^{(\rho,\rho^\vee)}
\prod_{\al\in R_+}
(1-t_\al(tq)^{(\al,\rho^\vee)})\notag.
\end{align}
Applying the evaluation formula for $P^{q,t}_{-\rho}(X),$
\begin{align}\label{yofxeval}
&\overline{\y}^{t}(\x^{t})(q^{-\rho_k}) =
C_{q,t} P^{q,t}_{-\rho}(q^{-\rho_k})
=t^{-\kapp}(t^2 q)^{-(\rho,\,\rho^\vee)}\Pi_{\tR},\\
&\Pi_{\tR}\equal\prod_{\al\in R_+}
\Bigl( (1- t_\al(tq)^{(\al,\,\rho^\vee)})
\prod_{ j=1}^{(\al^\vee,\,\rho)}
\frac{
(1- q_\al^{j-1}t_\al t^{(\al,\rho^\vee)})}
{(1- q_\al^{j-1}t^{(\al,\rho^\vee)})}\Bigr)\notag\\
=
&\prod_{\al\in R_+}
\Bigl( (1- q_\al^{k_\al+(\al^\vee,\,\rho+\rho_k)})
\prod_{ j=1}^{(\al^\vee,\,\rho)}
\frac{
(1- q_\al^{j-1+k_\al+(\al^\vee,\,\rho_k)})}
{(1- q_\al^{j-1+(\al^\vee,\,\rho_k)})}\Bigr),\notag\\
&\hbox{where\ \ }
t^{-\kapp}(t^2 q)^{-(\rho,\,\rho^\vee)}\ =\ 
q^{-\sum_{\nu}\,\nu\, k_\nu\, \hbox{\tiny{\kapp}}_\nu}\,
q^{-(\rho,\,\rho+2\rho_k)}.\notag
\end{align}
\end{maintheorem}
{\em Proof.} The choice $b=-\rho$  makes
$\S^{q,t}(P_b^{q,t})$ a constant, which equals to
$g^{q,t}(b).$ Therefore
$\y^{t}(P_b^{q,t})=g^{q,t}(b)\x^{t}$ and
$$
\overline{\y}^{t}\y^{t}(P_b^{q,t})=
g^{q,t}(b)\overline{\y}^{t}(\x^{t}).$$
Thus $\overline{\y}^{t}(\x^{t})=(g^{t})^{-1}
\overline{\y}^{t}\y^{t}(P_b^{q,t})=$
\begin{align*}
=&\ \ \overline{g}^{q,t}(b)g^{q,t}(b)(g^{q,t}(b))^{-1}P_b^{q,t}=
\overline{g}^{q,t}(b)P_b^{q,t}\\
=&\prod_{\al\in R_+}
(t_\al^{-1}X_\al(q^{(b-\rho_k)/2}) -
 X_\al(q^{(\rho_k-b)/2}))P_b^{q,t}\\
=&\prod_{\al\in R_+}
( t_\al^{-1}q^{(\al,b-\rho_k)/2} -
 q^{(\al,\rho_k-b)/2} ) P_b^{q,t}\\
=&\prod_{\al\in R_+}
 t_\al^{-1}q^{(\al,(b-\rho_k))/2}
\prod_{\al\in R_+}
(1- t_\al q^{(\al,\rho_k-b)}) P_b^{q,t}\\
=&\prod_{\nu\in \nu_R}
t_\nu^{-\hbox{\tiny\kapp}_\nu}
(t_\nu q_\nu)^{-(\rho,\rho_\nu^\vee)}
\prod_{\al\in R_+}
\Bigl( (1-t_\al \prod_{\nu\in \nu_R}
(t_\nu q_\nu)^{(\al,\rho_\nu^\vee)}\Bigr)P_b^{q,t}.
\end{align*}
The ``principle value"
of $P_b$ as $b=-\rho$ equals
\begin{align*}
&P_{-\rho}(q^{-\rho_k}) =
q^{-(\rho,\rho_k)}\prod_{\al\in R_+}
\prod_{ j=1}^{(\al^\vee,\rho)}
\Bigl(\frac{1- q_\al^{j-1}t_\al q_\al^{(\al^\vee,\rho_k})}
{1- q_\al^{j-1} q_\al^{(\al^\vee,\rho_k})}\Bigr)\\
=\ &(\prod_{\nu\in \nu_R}
t_\nu^{-(\rho,\rho_\nu^\vee)})
\prod_{\al\in R_+}
\prod_{ j=1}^{(\al^\vee,\rho)}
\frac{1- q_\al^{j-1}t_\al
\prod_{\nu\in \nu_R} t_\nu^{(\al,\rho_\nu^\vee)}}
{1- q_\al^{j-1}
\prod_{\nu\in \nu_R} t_\nu^{(\al,\rho_\nu^\vee)}}\\
=\ &t^{-(\rho,\rho^\vee)}
\prod_{\al\in R_+}
\prod_{ j=1}^{(\al^\vee,\rho)}
\frac{1- q_\al^{j-1}t_\al
t^{(\al,\rho^\vee)}}
{1- q_\al^{j-1}
t^{(\al,\rho^\vee)}}.
\end{align*}

Combining this formula with the previous one,
we obtain:
\begin{align*}
&C_{q,t}\cdot P_{-\rho}(q^{-\rho_k})=
\prod_{\nu\in \nu_R}
t_\nu^{-\hbox{\tiny\kapp}_\nu}
(t_\nu q_\nu)^{-(\rho,\rho_\nu^\vee)}
\prod_{\al\in R_+}
(1-t_\al\prod_{\nu\in \nu_R}
(t_\nu q_\nu)^{(\al,\rho_\nu^\vee)})\\
\cdot&
\prod_{\nu\in \nu_R}
t_\nu^{-(\rho,\rho_\nu^\vee)}
\prod_{\al\in R_+}
\prod_{ j=1}^{(\al^\vee,\rho)}
\Bigl(\frac
{1- q_\al^{j-1}t_\al
\prod_{\nu\in \nu_R} t_\nu^{(\al,\rho_\nu^\vee)}}
{1- q_\al^{j-1}
\prod_{\nu\in \nu_R} t_\nu^{(\al,\rho_\nu^\vee)}}\Bigr)\\
=&t^{-\kapp}
(t^2 q)^{-(\rho,\rho^\vee)}
\prod_{\al\in R_+}
(1- t_\al(tq)^{(\al,\rho^\vee)})
\prod_{\al\in R_+}\prod_{ j=1}^{(\al^\vee,\rho)}
\Bigl(\frac{
(1- q_\al^{j-1}t_\al t^{(\al,\rho^\vee)})}
{(1- q_\al^{j-1}t^{(\al,\rho^\vee)})}\Bigr).
\end{align*}
\sq
\smallskip

Generalizing the principle value
for the $X$\~polynomials, that is at $q^{-\rho_k}$,
we set
\begin{align}\label{valofh}
\{H\}\equal
H(1)(q^{-\rho_k}) \hbox{\ \ for \ an\ operator\ } H
\hbox{\ acting\ in \ \ }\v.
\end{align}

Applying formulas (3.5)-(3.7) from Key Lemma 3.3 from
\cite{C3}, we come to
the following connection of the values of the symmetric
Macdonald polynomials at $q^{-\rho_k}$ as $k\mapsto k+1.$
Here the evaluation formula is used too.

\begin{corollary} \label{YXPB}
Let $P=P_b^{q,t}$ for $b\in P_-$  be the Macdonald
polynomial for $t$ and $P_{b+\rho}^{q,tq}$ the one for $tq$;
it is zero if $b+\rho\not\in P_-.$
Then
\begin{align}\label{kconfin}
&\{\overline{\y}^{t}\x^t\}P_{b+\rho}^{q,tq}(q^{-\rho_{k+1}})=
\overline{g}^{q,t}(b)P_b^{q,t}(q^{-\rho_k})\\
=&\prod_{\al\in R_+}
(t_\al^{-1}X_\al(q^{(b-\rho_k)/2}) -
 X_\al(q^{(\rho_k-b)/2}))P_{b}^{q,t}(q^{-\rho_k}), \notag
\end{align}
where
$\{\overline{\y}^{t}\x^t\}=
\overline{\y}^{t}(\x^t)(q^{-\rho_k}).$ Equivalently,
\begin{align}\label{kconfinnx}
&P_{-\rho}^{q,t}(q^{-\rho_k})P_{b+\rho}^{q,tq}(q^{-\rho_{k+1}})\\
&=\
\prod_{\al\in R_+}\frac{
t_\al^{-1}X_\al(q^{(-\rho_k+b)/2}) -
 X_\al(q^{(\rho_k-b)/2})}
{
t_\al^{-1}X_\al(q^{(-\rho_k-\rho)/2}) -
 X_\al(q^{(\rho_k+\rho)/2})}\,
P_b^{q,t}(q^{-\rho_k})\notag\\
&=\ q^{(b+\rho\,,\,\rho)}
\prod_{\al\in R_+}
\frac{1- q_\al^{k_\al+(\al^\vee,\,\rho_k-b)}}
{1- q_\al^{k_\al+(\al^\vee,\,\rho_k+\rho)}}\,P_b^{q,t}(q^{-\rho_k}).
\notag
\end{align}
 \sq
\end{corollary}
\smallskip

\medskip
\subsection{Rational limit}
Setting $k_{\sht}=k=\nu_{\lng} k_{\lng},\ 
t_{\sht}=t=t_{\lng}$, the {\em Coxeter exponents\,}
are defined from the formula

\begin{align}\label{tdegr}
\prod_{\al\in R_+} \frac{1- t^{1+(\al,\,\rho^\vee)}}
{1-t^{(\al,\,\rho^\vee)}}=
\prod_{i=1}^n \frac{1-t^{\,m_i+1}}
{1-t},\ \ 2\rho^\vee=\sum_{\al\in R_+}\,\al^\vee.
\end{align}

In the simply-laced case, 
the products $\Pi_{\tR}$ 
and  
$C_{q,t} P^{q,t}_{-\rho}(q^{-\rho_k})$
can be expressed in terms of the
Coxeter exponents, Namely, $\Pi_{\tR}=$
\begin{align}
\prod_{\al\in R_+}
(1- &t(tq)^{(\al,\rho)})
\prod_{ j=1}^{(\al,\rho)}
\Bigl(\frac{
(1- q^{j-1} t^{1+(\al,\rho)})}
{(1- q^{j-1}t^{(\al,\rho)})}\Bigr)=
\prod_{i=1}^n \frac{\prod_{j=0}^{m_i} (1-q^{j} t^{\,m_i+1})}
{1-t},\notag\\
&C_{q,t} P^{q,t}_{-\rho}(q^{-\rho_k})
=t^{- | R_+ |}
(t^2 q)^{-(\rho,\rho)} \prod_{i=1}^n
\frac{\prod_{j=0}^{m_i} (1-q^{j} t^{\,m_i+1})}
{1-t}. \label{yofxdeg}
\end{align}
\smallskip

\rmk
(i) Here, technically, the  product $\prod_{\al\in R_+}
(1- t(tq)^{(\al,\rho)})$ ensures
the cancelation of the binomials in the denominator; it is not
needed  in (\ref{tdegr}). The process
of such cancelation is not quite immediate even in the
simply-laced case. It requires a combinatorial reformulation
of (\ref{tdegr}) in terms of the sequence of the numbers
of positive roots $\al$ of given  $(\al,\rho^\vee),$
the {\em height\,}; cf. Lemma \ref{COMBR01}.

(ii) In our approach, the
right-hand side automatically contains
$|R_+|$ extra factors in the numerator, which corresponds to
the classical formula $m_1+\ldots m_n=|R_+|$. 
The latter requires (simple) calculating the leading terms 
in (\ref{tdegr}). Another famous formula
$(m_i+1)\cdots(m_n+1)=|W|$, which is immediate from 
the Poincar\'e polynomial interpretation of (\ref{tdegr}),
becomes the formula for the product of {\em rational exponents}, 
a counterpart of (\ref{yofxdeg}) for the rational DAHA.

(iii) It is not difficult to calculate $\Pi_{\tR}$ in the
non-simply-laced case under $k_{\sht}=k=\nu_{\lng} k_{\lng}$.
The structure of the formula is practically the same as in
(\ref{yofxdeg}), but the indices $j$ 
satisfy some non-trivial
divisibility conditions and can become greater
than $m_i$. We note that there is no significant
simplification of $\Pi_{\tR}$ under the
substitution $k_{\sht}=k_{\lng}$ (without $\nu$).
\sq
\medskip

{\bf Poincar\'e polynomial}.
Let us recall the generalization of (\ref{tdegr})
to the case of two different $k$\~parameters from \cite{Ma0}
and its relation to the Poincar\'e polynomial of $W$.
We use the notion of the partial length $l_\nu(w)$
from (\ref{xlambda}), where $\nu\in \nu_R$. One has:
\begin{align}\label{tdegra}
&\Pi_R\equal
\sum_{w\in W}\prod_\nu\,t_{\nu\,}^{\,l_\nu(w)}\, =\,
\prod_{\al\in R_+} \frac{1- q_{\al}^{k_\al+(\al^\vee,\,\rho_k)}}
{1-q_{\al}^{(\al^\vee,\,\rho_k)}}\,.
\end{align}
See \cite{C2} for the proof of this formula via the 
$r$\~matrices and its applications to the Macdonald 
norm conjecture.
\smallskip

Let us list the counterparts of (\ref{tdegr}) in the
non-simply-laced cases. Actually, they can be readily
obtained from the formulas for the affine
exponents considered below by picking the binomials
without $q^j$ for $j>0$, i.e., the binomials given
entirely in terms 
of $t_{\sht}=q^{k_{\sht}}$ and 
$t_{\lng}=q_{\lng}^{k_{\lng}}$. In the $B-C$ cases, 

\smallskip
\begin{align}\label{tdegrb}
&B_n:\ \ \Pi_R\ =\ 
\prod_{m=0}^{n-1}\, \Bigl(\frac{1-t_{\lng}^{m+1}} {1-t_{\lng}}\Bigr)\,
\Bigl(\frac{1-t_{\lng}^{2m}t_{\sht}^2} 
{1-t_{\lng}^{m}\,t_{\sht\,}}\Bigr),\\ 
&C_n:\ \ \Pi_R\ =\ 
\prod_{m=0}^{n-1}\, \Bigl(\frac{1-t_{\sht}^{m+1}} {1-t_{\sht}}\Bigr)\,
\Bigl(\frac{1-t_{\sht}^{2j}t_{\lng}^2} 
{1-t_{\sht}^{m}\,t_{\lng\,}}\Bigr).\ \label{tdegrc}
\end{align}
These formulas are $t_{\lng}\leftrightarrow t_{\sht}$\~dual
to each other,
which readily follows from the interpretation of $\Pi_R$
in terms of the Poincar\`e polynomial. The same symmetry
holds in the formulas below. In the case of $G_2$,
\begin{align}\label{tdegrg}
&G_2:\ \ \Pi_R\ =\ 
\frac{(1-t_{\lng}^{2})\, ({1-t_{\sht}^{2}})\,
({1-t_{\lng}^{3}t_{\sht}^{3}})} 
{(1-t_{\lng}^{\ })\, ({1-t_{\sht}^{\ }})\,
({1-t_{\lng}^{\ }t_{\sht}^{\ }})}\,. 
\end{align}

In the case of $F_4$,  
\begin{align}\label{tdegrf}
F_4:&\ \ \Pi_R\ =\ \frac
{(1-t_{\lng}^{2})\, ({1-t_{\lng}^{3}})\,
(1-t_{\sht}^{2})\, ({1-t_{\sht}^{3}})}
{(1-t_{\lng}^{\ })\, ({1-t_{\lng}^{\ }})\,
(1-t_{\sht}^{\ })\, ({1-t_{\sht}^{\ }})}\,\\
&\times\ \frac
{({1-t_{\lng}^{4}t_{\sht}^{2}}) ({1-t_{\lng}^{2}t_{\sht}^{4}}) 
({1-t_{\lng}^{4}t_{\sht}^{4}}) ({1-t_{\lng}^{6}t_{\sht}^{6}})} 
{({1-t_{\lng}^{2}t_{\sht}^{\ }}) ({1-t_{\lng}^{\ }t_{\sht}^{2}}) 
({1-t_{\lng}^{\ }t_{\sht}^{\ }}) ({1-t_{\lng}^{3}t_{\sht}^{3}})}\,. 
\notag
\end{align}
\medskip 

{\bf Rational limit}.
In the rational limit $q=e^h, h\to 0,$
(\ref{yofxeval}) results in the formula from
Theorem 9.8 \cite{O3} (the simply-laced case) and
from Theorem 4.11 \cite{DJO}, where the root systems
$B_n,F_4$ and $I_2(2m)$ where considered.
In the rational setting, the limiting formulas for
$B$ and $C$ are obtained form each other by the
transposition $k_{\lng}\leftrightarrow k_{\sht}$;
$G_2=I_2(6)$. 

For the sake of concreteness,
let us describe here the limiting procedure in the simply-laced
case only. Let $Y_b=e^{-\sqrt{h}\,\widetilde{y}_b\,},$
$X_b=e^{\sqrt{h}x_b}.$
Then the $h$\~leading term of
$\widetilde{y}_b$ becomes the {\dfont Dunkl operator}
\cite{Du}:
$$
y_b=\partial_b+\sum_{\al\in R_+}\frac{k(b,\al)}{x_\al}(1-s_\al)
\for \partial_b(x_c)=(b,c),
$$
acting in the polynomial representation $\C[x_b].$ The
limit of the shift operator $\s$ is
$$\prod_{\al\in R_+}x_\al^{-1}\prod_{\al\in R_+}y_\al.$$

The $h$\~leading
term of the expression
$\overline{\y}(\x)$ is
$$
(-h)^{|R_+|}(\prod_{\al\in R_+}y_\al)((\prod_{\al\in R_+}x_\al)),
$$
where we apply the $y$\~operator to the $x$\~polynomial.
It is a {\em constant}, so the ``principle value", which is
taken at zero in the rational theory, simply coincides
with it.

The leading term of the
product in the right-hand side of (\ref{yofxdeg})
can be readily calculated; it is 
$$
(-h)^{| R_+ |}|W|
\prod_{i=1}^n \prod_{j=1}^{m_i}(j+ (m_i+1) k).
$$
We use here that
the product
$\prod_{i=1}^{n}(1-t^{m_i+1})/(1-t)$ tends to $|W|.$

We arrive at the following formula due to Opdam.
\begin{corollary}\label{XYOPDAM}
For $y_b=
\partial_b+\sum_{\al\in R_+}\frac{k(b,\al)}{x_\al}(1-s_\al),$
\begin{align}\label{yofxrat}
&(\prod_{\al\in R_+}y_\al)((\prod_{\al\in R_+}x_\al))=
|W|\prod_{i=1}^n \prod_{j=1}^{m_i}(j+ (m_i+1) k),
\end{align}
\sq
\end{corollary}
\smallskip

\rmk
Formula (\ref{yofxrat}) was proved in \cite{O1} in the
crystallographic case and then extended in \cite{O3,DJO}
to groups generated by complex reflections; its relation
to the reducibility of the rational polynomial representation
was established in the latter paper. See also paper \cite{DO}.
We obtain it as a rational limit of the $q,t$\~formula;
the rational limit is an entirely algebraic procedure. 

Our approach via the 
{\em affine exponents\,} results in a direct 
link to the defining product formula for the 
{\em Coxeter exponents\,} in (\ref{tdegr}); see also 
(\ref{tdegrb},\ref{tdegrc},\ref{tdegrg},\ref{tdegrf}).
This approach establishes a relation
with the $p$\~adic theory (Macdonald-Matsumoto) corresponding
to the limit $q\to 0$. Recall that $\Pi_{\tR}$ becomes
the Poincar\'e polynomial $\Pi_R$ under the latter limit;
$\Pi_R$ appears virtually in all constructions of the 
$p$\~adic theory of spherical functions. 

Opdam used  the differential-{\em trigonometric} theory in 
\cite{O1}; the rational theory alone appeared insufficient
for (\ref{yofxrat}). Its direct justification is known only
in the $A$\~case (Dunkl, Hanlon).

Another approach (originated in \cite{O3}) is based on 
paper \cite{GU} devoted to the semisimplicity of the 
(non-affine) Hecke algebras. It works for arbitrary finite 
Coxeter groups, that is beyond the frameworks of our method. 

Reversing, in 
a sense, \cite{O3} (or using \cite{GGOR}), we can now apply
the affine exponents to {\em reprove} the formula 
from \cite{GU} in the crystallographic case.  
Applications of the affine exponents 
are expected for the semisimplicity of the 
Ariki\~Koike\~Cherednik algebras and for similar algebras.
It is related to the $q\leftrightarrow q'$\~duality of
$\Pi_{\tR}$, which will not be discussed here.
\sq

\smallskip
\subsection{Non-simply-laced case}
Recall that we use the normalization $\nu_{\sht}=1$
for $\nu_{\al}=(\al,\al)/2$ and $q_\al=q^{\nu_\al}$.
We will begin with the formula in the case of
the root system $B_n;$ $\ep_m$ are from the $B$\~table of
\cite{Bo}, $(\ep_l,\ep_m)=2\de_{lm}$. Then
\begin{align*}
R_+=&\{\hbox{sht\ \,}:\ \ep_m= \al_m+\ldots+\al_n, \
m=1,\ldots,n\} \and \\
&\{\hbox{lng}_-:\ \ep_l-\ep_m=\al_{l}+\ldots \al_{m-1}\},
\hbox{ \ as\ }
 n\ge m>l\ge 1,\\
&\{\hbox{lng}_+:\ \ep_l+\ep_m=
\al_{l}+\ldots \al_{m-1}+2(\al_m+\ldots+\al_{n})\}\,;\\
(\ep_m^\vee, \rho_k)&\ = \ k_{\sht}+2(n-m)k_{\lng},\
\rho_1=\rho_{\sht}=\sum_{m=1}^n\ep_m/2\,,
\\
((\ep_l-\ep_m)^\vee,\rho_k)&\ =\ (m-l)k_{\lng},\ \ \rho_2=
\rho_{\lng}=\sum_{m=1}^n (n-m)\ep_m,
\\
((\ep_l+\ep_m)^\vee, \rho_k)&\ =\  k_{\sht}+
(2n-m-l)k_{\lng},\ \rho_k=k_{\sht}\rho_1+k_{\lng}\rho_2.
\end{align*}
Let us separate the terms with $j=1$ constituting the part
of (\ref{yofxeval}) without $q$. Recall that
the factor $t^{-\kapp}q^{-(\rho,\,\rho+2\rho_k)}$ is
disregarded in the definition of $\Pi_{\tR}$:
\begin{align}\label{yofxevalb}
\Pi_{\tR}\equal&\prod_{\al\in R_+}
\Bigl( (1- q_\al^{k_\al+(\al^\vee,\,\rho+\rho_k)})
\prod_{ j=1}^{(\al^\vee,\,\rho)}
\frac{
(1- q_\al^{j-1+k_\al+(\al^\vee,\,\rho_k)})}
{(1- q_\al^{j-1+(\al^\vee,\,\rho_k)})}\Bigr)\notag\\
=\ &\frac{Q'_1}{Q'_2}\
\prod_{m=0}^{n-1}\ (\frac{1-q^{2(m+1)k_{\lng}}}{1-q^{2k_{\lng}}})\
(\frac{1-q^{2(2mk_{\lng}+k_{\sht})} }
{1-q^{2mk_{\lng}+k_{\sht}}}).
\end{align}
The numerator of the product part ($t$\~pure part) 
of the latter formula 
is obviously divisible by the denominator. It is much less obvious
for $Q'_1/Q'_2$ (this follows of course from the construction 
in terms of the shift operator). Let $\widetilde{Q}_{1}$ and
$\widetilde{Q}_{2}$ be the
{\em reduced numerator and denominator} upon reducing
the {\em coinciding} factors but without the divisions, 
$\widetilde{Q}'_{1,2}$ the sub-products where all
$t$\~pure binomials are removed. Then

\begin{align}\label{yofbfinal}
\widetilde{Q}'_1\, &=\,(1-q^{2k_{\sht}+1})
\prod_{m=2}^n\prod_{j=1}^{m-1} (1-q^{2mk_{\lng}+2j})\\
&\times\, \prod_{m\in M}
\prod_{j=1,2,\ldots}^{i_m+\sigma_mj<2+2m}\,(1-q^{2mk_{\lng}+
2k_{\sht}+i_m+\sigma_mj}),\notag\\
M\ =\ &\{1,2,\ldots,2[n/2],2[n/2]+2,2[n/2]+4,\ldots,2n-2\},\notag
\end{align}
the $j$\~step $\sigma_m$ is $1$ for $m=2,4,...,2[(n+1)/2]+2$,
and $2$ otherwise; the shift $i_m$ is $-1$ for odd $m$
and zero otherwise.

Let us give the list of exponents in the binomials
$(1-q^{\{\cdot\}})$ that appear in $\widetilde{Q}'_1$ for $n=3$:

\begin{align}\label{bthreenum}
&
(1 + 2 k_{\sht} );\ (2 + 4 k_{\lng} );\
(2 + 6 k_{\lng} ),\ (4 + 6 k_{\lng} );\\
&
(1 + 2 k_{\lng} + 2 k_{\sht} ),\
(3 + 2 k_{\lng} + 2 k_{\sht} );\notag\\
&(1 + 4 k_{\lng} + 2 k_{\sht} ),\
(2 + 4 k_{\lng} + 2 k_{\sht} ),\notag\\
&(3 + 4 k_{\lng} + 2 k_{\sht} ),\
  (4 + 4 k_{\lng} + 2 k_{\sht} ),
\  (5 + 4 k_{\lng} + 2 k_{\sht} );\notag\\
&(2 + 8 k_{\lng} + 2 k_{\sht} ),\
(4 + 8 k_{\lng} + 2 k_{\sht} ),\notag\\
&(6 + 8 k_{\lng} + 2 k_{\sht} ),\
(8 + 8 k_{\lng} + 2 k_{\sht} ).\notag
\end{align}

The formula for the {\em reduced} denominator
$\widetilde{Q}'_2$ reads as follows:
\begin{align}\label{yofbdenom}
&\widetilde{Q}'_2\ =\
\prod_{m=1}^{n-1}\prod_{j=1}^{2m} (1-q^{k_{\sht}+2mk_{\lng}+j}).
\end{align}
The number of terms here is $n(n-1)$.

In the case of $n=3$, the  $q$\~exponents of $\widetilde{Q}'_2$ are

\begin{align}\label{bthreeden}
&(1 + 2k_{\lng} +k_{\sht} )\, ,
(2 + 2k_{\lng} +k_{\sht} );\ \ (1 + 4k_{\lng} +k_{\sht} ),\\
&(2 + 4k_{\lng} +k_{\sht} )\, ,\,
(3 + 4k_{\lng} +k_{\sht} )\, ,\,
(4 + 4k_{\lng} +k_{\sht} ).\notag
\end{align}
The number of factors in $\widetilde{Q}'_2$ is
always smaller that in $\widetilde{Q}'_1$
by the number of positive roots $|R_+|$
($9$ in the case of $B_3$).
\medskip

\section{Affine exponents}
Continuing the previous section,
we come to a general notion of the {\dfont affine exponents\,}
that are the exponents of the terms $(1-q^{\{\cdot\}})$
in the {\em reduced denominator}  $\widetilde{Q}_1$ and
{\em reduced numerator}  $\widetilde{Q}_2$ 
of $\Pi_{\tR}$ from (\ref{yofxevalb}).
Affine exponents can
be {\em positive} (from the numerator) or 
{\em negative} (from the denominator); 
(\ref{bthreenum}) and (\ref{bthreeden})
give such exponents without the pure $k$\~terms
for $\widetilde{B}_3$. The {\em reduction}
is simple removing
coinciding terms in the numerator and denominator without
performing any actual divisions. 

In contrast to the {\em classical exponents\,}\,
(Coxeter exponents),
the reduced denominators $\widetilde{Q}_2$ can be nontrivial
if $R$ is not simply-laced; in the simply-laced case,
only the factors $(1-t)$  appear in $\widetilde{Q}_2$.
\smallskip 

Always, as in  (\ref{tdegra}),
the numerator $\widetilde{Q}_1$ is divisible by the
denominator $\widetilde{Q}_2$,
which is not quite simple to justify without the
approach based on the shift operator (in the
non-simply-laced case). As in the classical case,
the total divisibility guaranties that
the coefficients of the $q,t$\~expansion  
of $\Pi_{\tR}$
in terms of the {\em $q,t$\~powers} 
$q_\nu^a t_{\sht}^b t_{\lng}^c$ for $a,b,c\in \Z_+$
will be from $\Z$.  
Indeed, the final formula will be a product of terms
in the form $(1-A^m)/(1-A)$, where $A$ is a 
$q,t$\~power, and also terms $(1-A)$ coming from the rational
exponents. The latter terms, generally,
result in negative coefficients;
the coefficients would be all from $\Z_+$ 
if the rational exponents were disregarded. 
The coefficients
of the classical product $\Pi_R$ are from $\Z_+$.
\medskip

\subsection{Rational exponents}
In this paper, {\dfont rational exponents}\,
play an important role. Their list is
obtained from the list of affine exponents 
{\em in the numerator}\, $\widetilde{Q}_1$ by 
removing the binomials that are divisible by some 
binomials in the denominator  $\widetilde{Q}_2$.
Only such exponents can lead to the {\em $k$\~roots\,} 
of the product (\ref{yofxevalb})
provided that $q$ is not a root of unity and imposing\ 
$
q^a\prod_\nu t_\nu^{b_\nu}=1\ \Rightarrow\ 
a+\sum_\nu\nu k_\nu b_\nu=0.
$\ 
More formally,
rational exponents are the affine {\em positive} (from the
numerator) exponents, such that the corresponding
binomials are not divisible by any binomials in the
denominator. There is one-to-one correspondence
between the non-rational exponents and the 
{\em negative} exponents (from the denominator), 
which makes this definition meaningful and also ensures that
{\em the number of rational exponents is always $|R_+|$.}

\medskip
\rmk
(i) Concerning ``integrality" of the
$q$\~expansion of $\Pi_{\tR}$,
we would like to discuss
the formula for the graded multiplicities of the adjoint
representation $\mathfrak{g}$ of the simple Lie 
algebra $\mathfrak{g}$ associated with $R$ in
the exterior algebra $\Lambda\mathfrak{g}.$

It was conjectured by Joseph in the simply-laced case that
\begin{align}\label{Joseph}
&\sum_{m\ge 0}[\,\Lambda^m\mathfrak{g}\,:\,\mathfrak{g}\,]\,q^m=
(1+q^{-1})\prod_{i=1}^{n-1}
(1+q^{2m_i+1})\sum_{i=1}^n q^{2m_i},
\end{align}
a variant/generalization of a series of formulas on the
structure of $\Lambda\mathfrak{g}$ due 
to Kostant and others. Bazlov extended it to 
the non-simply-laced case (which is
nontrivial) and proved it in \cite{Baz} by using the
$Y$\~operators. 

(ii) He calculates the ``non-affine" part
of the expansion of $\mu$ and interprets a proper sum 
of the corresponding coefficients {\em for a special value
of $t$} as (\ref{Joseph}).  See also \cite{Ion},
where the method is somewhat different and additional
information can be found. 

We do not expect \cite{Baz},\cite{Ion} to be directly
connected with our affine (or rational) exponents. 
Recall that the main part of $\Pi_{\tR}$ comes from the
principle value $P_{-\rho}(q^{-\rho_k})$, the
{\em $q,t$\~dimension of \,$V_\rho$}\,;\, the latter
does not appear in these two papers. However, the 
coefficients of the polynomial
$P_{-\rho}$ are connected with the coefficients calculated
in \cite{Baz},\cite{Ion}.
It may be an indication that these papers can be extended
and some variants of (\ref{Joseph}) may be used to interpret
$P_{-\rho}(q^{-\rho_k})$ and $\Pi_{\tR}$. 
\smallskip

(iii) Some of the {\em rational exponents}\ 
appear in \cite{Baz} and \cite{Ion}. It is not surprising 
because the $k$\~zeros of rational exponents
are the {\em singular $k$\~parameters}\,, 
the values of $k$ when $\v$ has a nontrivial radical
$Rad$ of the evaluation pairing $\{\, ,\,\}$.
To be exact, such zeros give the cases
when the $t$\~discriminant $\x^t$,
the $t$\~anti-invariant for the non-affine Hecke subalgebra
of $\HH^\flat$ from \ref{shiftoper}, belongs to $Rad$ (which
makes it nonzero).
If $k_\nu$ are replaced by the sets $k_\nu+\Z_+$, then 
the zeros of affine exponents describe {\em all} singular $k$. 

This interpretation of the rational exponents
is in the case when all multiplicative
relations among $q,t_{\nu}$ come from the 
$\Z$-relations among $1,\nu k_\nu$;
generally, the {\em affine exponents}\ must be used,
say, $t$ can be $\ze q^{s/r}$ for  
$\ze=\sqrt[r]{1}\neq 1$.

Therefore the rational or affine exponents are
inevitable  in formulas that somehow require the 
non-degeneracy of the Macdonald polynomials and related 
structures.
\sq

\medskip
In the simply-laced case, as we know,
the {\em rational exponents}\ 
are all positive {\em affine exponents}\ 
(from the numerator)
with nonzero integer components $j$. They are 
$\{j+k(m_i+1),\,1\le j\le m_i\}$.
For $\widetilde{B}_n, \widetilde{C}_n,
\widetilde{F}_4, \widetilde{G}_2$,
the lists are given in the following theorem.

Due to their definition, the 
product of the {\em affine rational exponents}\ 
(as they are, without $1-q^{\{\cdot\}}$)
is proportional to the limit of (\ref{yofxevalb})
as $q\to 1$. So this product describes
$(\prod_{\al\in R_+}y_\al)((\prod_{\al\in R_+}x_\al))$,
which generalizes Corollary \ref{XYOPDAM} to arbitrary,
possibly non-simply-laced, root systems.

The products of affine rational exponents 
in the non-simply-laced cases 
are proportional to those from \cite{DJO},Theorem 4.11. 
Recall that we process the product (\ref{yofxevalb})
combinatorially, by canceling coinciding factors as one  
does with the classical $\Pi_R$. The ``non-divisible" 
binomials (from the numerator) 
lead to \, {\em affine rational exponents}\,.  
\smallskip

\comment{
(i)
Upon identifying proportional 
rational exponents, 
they can be re-introduced entirely in terms of the rational 
DAHA theory. Namely, we can define them as binomial 
factors in Opdam's 
formula for 
$(\prod_{\al\in R_+}y_\al)((\prod_{\al\in R_+}x_\al))$;
see (\ref{XYOPDAM}).

The {\em multiple singular $k$ in the
rational theory}, the $k$\~zeros of these factors,
generally do not appear in the $q,t$\~theory.
Due to $1-q^{\{\cdot\}}$, proportional rational exponents come
from {\em distinct} factors in $\Pi_{\tR}$ in our approach
(unless for $D_{even}$\,!). We note that the nontrivial
multiplicities are important in the theory of Bernstein-Sato
polynomials; see \cite{O1} for the connection. The interpretation
of the rational exponents as the reducibility condition for
the {\em rational} polynomial representation \cite{DJO} 
does not involve their multiplicities. The
non-trivial multiplicities
create technical difficulties for Opdam's approach.
}

\rmk
(i) The multiplicities of the rational exponents in the
limit are the only lost information if
the rational exponents are defined entirely
within the rational theory, i.e., without the above
approach via the rational affine exponents in the $q,t$\~theory.
The nontrivial multiplicities
create technical difficulties for Opdam's approach. 
Due to the presence of $1-q^{\{\cdot\}}$ in the affine formulas,
proportional rational exponents come
from {\em distinct} factors in $\Pi_{\tR}$ in our approach
(unless for $D_{even}$\,!). 
The cases of nontrivial multiplicities in the 
limit are as follows. 

For $\widetilde{B}_n, n\ge 4$,
such pairs are $\{m k_{\lng}+j\,,\, lmk_{\lng}+lj\}$
as $2\le lm\le n, j/m\not\in \Z$. For $\widetilde{C}_n$,
$k_{\lng}$ must be replaced here by $k_{\sht}$. 
Also, the pairs of the rational exponents
$\{1 + 2k_{\lng} + 2k_{\sht}\, ,\, 3 + 6k_{\lng} + 6k_{\sht}\},
\{3 + 2k_{\lng} + 2k_{\sht}\, ,\, 9 + 6k_{\lng} + 6k_{\sht}\}$
for $\widetilde{F}_4$ from
(\ref{redfdegrees}) cannot be distinguished under the
rational limit. 

These are
the only cases when the multiplicities occur in the rational
limit for generic $k_{\sht},k_{\lng}$ (in the non-simply-laced case).
If $k_{\sht}=\nu_{\lng}k_{\lng}$ (for instance, for
$A,D,E$), then there are many multiple exponents in the rational
limit.  

We stick to the ``combinatorial" (via $q,t$) definition of the
affine rational exponents in this paper. 
\smallskip

(ii) 
Technically, the  reduction of coinciding factors is 
somewhat simpler among the affine rational exponents,  
but the case of ``divisible" (non-rational) exponents is not very
different; see below. Note that the duality between 
the rational exponents of $\widetilde{B}$
and $\widetilde{C}$ from the theorem does not hold for the
affine exponents (its certain affine counterpart exists
but is of more sophisticated nature). 

However, Theorem \ref{OTHERAFFINE} below shows
that the strict duality holds for {\em another} 
affine extension of $R$, namely, for 
$R\times \Z=\{[\al,j]\}$. 

Under the rational limit, the rational
$B\leftrightarrow C$\~ duality simply reflects
the fact that the rational Dunkl operators are 
$k_{\sht}\leftrightarrow k_{\lng}$\~coinciding for 
$B$ and $C$. Indeed, in the rational case, 
the construction of these operators does not 
depend on the choice of the affine extension of $R$.
\smallskip

(iii) 
We define in this paper the affine roots in the form
$[\al,\nu_\al j]$ (with the $\nu_\al$\~factors) 
and introduce DAHA using the {\em coinciding} affine systems 
$\{\tR,\tR\}$ for both, $X$ and $Y$. This 
makes $\HH\,$ invariant under the action of the 
automorphisms $\tau_{\pm},\si$ and simplifies
other considerations, like the theory of the
evaluation pairing $\{\cdot,\cdot\}.$
Almost all results of this paper can be transferred to
the case of the affine root system 
$\widehat{R}\equal R\times \Z$ and, moreover, the
for combination $\tR,\widehat{R}$ taken for $X$ and for $Y$ 
in the definition of $\HH\, $ and related structures. Note 
that DAHA and the $q,t$\~shift-operator
from \cite{C2}) were defined for the combination
$\{\widehat{R},\tR\}$ of affine extensions. Such choice 
is exactly needed for the affine duality discussed below.

There is also flexibility with choosing the lattices
$Q\subset B\subset P$ for $X$ and the lattice between
$Q^\vee$ and $P^\vee$ for the generators $Y$. However, the
choice of the lattices does not affect the definition 
of the shift-operator.
\sq
\smallskip

\begin{theorem}\label{BPRIMDEG}
(i) The number of rational exponents is always $|R_+|$
and they are all simple (unless for $D_{even}$).
In the $\widetilde{B}_n$\~case, the list of the corresponding
powers in the binomials $(1-q^{\,\{\, \cdot \ \}})$ is
\begin{align}\label{redbdegrees}
&\{2k_{\sht}+1\},\ \{2mk_{\lng}+2j\,,\ 2\le m\le n,\, 0<j<m\},\\
&\{2k_{\sht}+2mk_{\lng}+2j+1,\ \ 1\le m<n,\ \, 0\le j\le m\}.
\notag
\end{align}
Up to proportionality,
the rational exponents
in the $\widetilde{C}_n$\~case are obtained from (\ref{redbdegrees})
by the transposition $k_{\lng}\leftrightarrow k_{\sht}$. Explicitly:
\begin{align}\label{redbdegreesc}
&\{2k_{\lng}+1\},\ \{mk_{\sht}+j\,,\ 2\le m\le n,\, 0<j<m\},\\
&\{2k_{\lng}+2mk_{\sht}+2j+1,\, 1\le m<n,\, 0\le j\le m\}.\notag
\end{align}

(ii) In the case of $\widetilde{G}_2$, the list is:
\begin{align}\label{redgdegrees}
\{
&(1 + 2 k_{\lng}),\ (1 + 2 k_{\sht}),
\ (1 + 3 k_{\lng} + 3 k_{\sht}),
\\
&(2 + 3 k_{\lng} + 3 k_{\sht}),\ (4 + 3 k_{\lng} + 3 k_{\sht}),\
(5 + 3 k_{\lng} + 3 k_{\sht})\}.\notag
\end{align}

(iii) In the case of $\widetilde{F}_4$, the list is:
\begin{align}\label{redfdegrees}
\{
&(2 + 4k_{\lng}),\ (2 + 6k_{\lng}),\ (4 + 6k_{\lng}),\\
&(1 + 2k_{\sht}),\ (1 + 3k_{\sht}),\ (2 + 3k_{\sht}),\notag\\
&(1 + 2k_{\lng} + 2k_{\sht}), (3 + 2k_{\lng} + 2k_{\sht}),
\notag\\
&(2 + 8k_{\lng} + 4k_{\sht}), (6 + 8k_{\lng} + 4k_{\sht}),
(10 + 8k_{\lng} + 4k_{\sht}),
\notag\\
&(1 + 2k_{\lng} + 4k_{\sht}), (3 + 2k_{\lng} + 4k_{\sht}),
\ (5 + 2k_{\lng} + 4k_{\sht}),
\notag\\
&(1 + 4k_{\lng} + 4k_{\sht}), (3 + 4k_{\lng} + 4k_{\sht}),
\notag\\
&\ \ \ \ \ \ \ \ \ \ \ \ \ \ \ \ \ \ \ \ \ \ \  
 (5 + 4k_{\lng} + 4k_{\sht}),\ (7 + 4k_{\lng} + 4k_{\sht}),
\notag\\
&(1 + 6k_{\lng} + 6k_{\sht}), (3 + 6k_{\lng} + 6k_{\sht}),
\ (5 + 6k_{\lng} + 6k_{\sht}),
\notag\\
&(7 + 6k_{\lng} + 6k_{\sht}), (9 + 6k_{\lng} + 6k_{\sht}),
(11 + 6k_{\lng} + 6k_{\sht})\}.\notag
\end{align}
\end{theorem}
\sq
\medskip

\subsection{The case of
\texorpdfstring{{\mathversion{bold}$\widetilde{C}_n$}}
{C}}
Now the inner product is $(\ep_l,\ep_m)=\de_{lm}$ and
\begin{align*}
R_+\ =\ \{\ 2\ep_m, \
m=1&,\ldots,n,\and \ep_l\pm\ep_m,\
 n\ge m>l\ge 1\ \},\\
((2\ep_m)^\vee, \rho_k)\ =& \ k_{\lng}+(n-m)k_{\sht},\
\rho_2=\rho_{\lng}=\sum_{m=1}^n\ep_m,
\\
((\ep_l-\ep_m)^\vee,\rho_k)\ =&\ (m-l)k_{\sht},\ \ \rho_1=
\rho_{\sht}=\sum_{m=1}^n (n-m)\ep_m,
\\
((\ep_l+\ep_m)^\vee, \rho_k)\ =&\  2k_{\lng}+
(2n-m-l)k_{\sht},\ \rho_k=k_{\sht}\rho_1+k_{\lng}\rho_2.
\end{align*}

The complete formula for the product $\Pi_{\widetilde{C}_n}$
reads as:
\begin{align}\label{yofxevalc}
&\prod_{\al\in R_+}
\Bigl( (1- q_\al^{k_\al+(\al^\vee,\,\rho+\rho_k)})
\prod_{ j=1}^{(\al^\vee,\,\rho)}
\frac{
(1- q_\al^{j-1+k_\al+(\al^\vee,\,\rho_k)})}
{(1- q_\al^{j-1+(\al^\vee,\,\rho_k)})}\Bigr)\\
=
&\frac{\prod_{m=2}^{n}\,\prod_{j=0}^{m-1}\ (1-q^{mk_{\sht}+j})}
{(1-q^{k_{\sht}})^{n-1}\,
\prod_{m=0}^{n-1}\prod_{0\le j\le m+1}^
{\hbox{\tiny even\ }j\hbox{\tiny \ for\ even\ }m}\
(1-q^{2k_{\lng}+mk_{\sht}+j})}\notag\\
\times&
\prod_{m=[\frac{n+1}{2}]}^{n-1}\,\prod_{1\le j\le 2m+1}^
{j\hbox{\tiny \,mod\,}2=1}
\ (1-q^{2k_{\lng}+2mk_{\sht}+j})\notag\\
\times&
\prod_{m=0}^{n-1}\,\prod_{0\le j\le 2m+2}^{j\hbox{\tiny\,mod\,}2=0}
\ (1-q^{4k_{\lng}+2mk_{\sht}+j}).\notag
\end{align}

We did not separate the terms that come from
the Poincar\'e polynomial $\Pi_R$ and the terms
that contain $q$; see (\ref{tdegrc}).
Generally, this product looks
somewhat simpler than the one for $\widetilde{B}_n$.
We do not have any statements about their connection.

Let us give the formula for the number $N_-(C_n)$
of {\em negative affine exponents} (i.e., in the denominator):
$$
N_-(C_n)=n-2+[\frac{n+1}{2}][\frac{n+3}{2}]/2+[\frac{n+2}{2}]^2.
$$
Respectively, $N_+(C_n)=|R_+|+ N_-(C_n)=n^2+N_-(C_n)$. 
These formulas for $B_n$ are very different:
$N_-(B_n)=n^2+n+1$; see (\ref{yofbdenom}).
\medskip

\subsection{The case of
\texorpdfstring{{\mathversion{bold}$\widetilde{G}_2$}}
{G}}
We will give the complete formula for the affine
exponents in the $\widetilde{G}_2$\~case.
The following is the list of the $k$\~heights
$(\al^\vee,\rho_k)$ in the notation of \cite{Bo}:
\begin{align}\label{rootsg2}
&\hbox{short\ } \ \al_1,\ \ \ \ \ \ \ \
(\al^\vee,\rho_k)= k_{\sht},\notag\\
&\hbox{short\ }\ \al_1+\ \al_2,\
(\al^\vee,\rho_k)= 3k_{\lng}+ k_{\sht},\notag\\
&\hbox{short\ } 2\al_1+\ \al_2, \
(\al^\vee,\rho_k)= 3k_{\lng}+ 2k_{\sht},\notag\\
&\hbox{long\ \ } \ \al_2,\ \ \ \ \ \ \ \ \
(\al^\vee,\rho_k)= k_{\lng},\notag\\
&\hbox{long\ \ } 3\al_1+\ \al_2 ,\
(\al^\vee,\rho_k)= k_{\lng}+ k_{\sht},\notag\\
&\hbox{long\ \ } 3\al_1+2\al_2,\
(\al^\vee,\rho_k)= 2k_{\lng}+ k_{\sht},\notag\\
&\rho_k=(2\al_1+\al_2)k_{\sht}+(3\al_1+2\al_2)k_{\lng}.
\end{align}
Recall that the normalization is $(\al_1,\al_1)=2$.
\smallskip

The $12$ {\em affine exponents}\ in the numerator are:
\begin{align}\label{g2affinen}
&(6 k_{\lng}),\ (3 + 6 k_{\lng} ),\ (2 k_{\sht} ),
\ (1 + 2 k_{\sht} ),\\
&  (1 + 3 k_{\lng}+ 3 k_{\sht} ),\ (2 + 3 k_{\lng}+ 3 k_{\sht} ),\
(4 + 3 k_{\lng}+ 3 k_{\sht} ),\notag\\
&  (5 + 3 k_{\lng}+ 3 k_{\sht} ),\ (9 k_{\lng}+ 3 k_{\sht} ),\
(3 + 9 k_{\lng}+ 3 k_{\sht} ),\notag\\
&  (6 + 9 k_{\lng}+ 3 k_{\sht} ),\ (9 + 9 k_{\lng}+ 3 k_{\sht} ).
\notag
\end{align}

The $6$ {\em affine exponents} \ in the denominator are:
\begin{align}\label{g2affined}
&\{(3 k_{\lng}),\ (k_{\sht} ),\ (3 k_{\lng}+ k_{\sht} ),\
(1 + 3 k_{\lng}+ k_{\sht} ),\\
&  (2 + 3 k_{\lng}+ k_{\sht} ),\ (3 + 3 k_{\lng}+ k_{\sht} ) \}.
\notag
\end{align}
\medskip

\subsection{The case of
\texorpdfstring{{\mathversion{bold}$\widetilde{F}_4$}}
{F}}
In the notation of \cite{Bo}, the roots are
\begin{align}\label{rootsf4}
\hbox{lng:\ } &1000,0100,1100,0120,1120,1220,\notag\\
              &0122,1122,1222,1242,1342,2342,\notag\\
\hbox{sht:\ } &0010,0001,0110,0011,0111,1110,\notag\\
              &1111,0121,1121,1221,1231,1232;\notag\\
(\,abcd^\vee\,,\,\rho_k\,&)\ =\ 
l_{\nu,\lng}(a+b)k_{\lng}+l_{\nu,\sht}^{-1}(c+d)k_{\sht},\\
\nu=\nu_{abcd},\  l&_{\lng,\sht}=2=l_{\sht,\lng},\  
l_{\lng,\lng}=1=l_{\sht,\sht}.\notag
\end{align}

Let us give the list of the {\em corresponding}
$(\al^\vee,\rho_k):$
{\small
\begin{align*}
&\hbox{\underline{lng}:\ \, }
k_{\lng},\, k_{\lng},\, 2 k_{\lng},\, k_{\lng} + k_{\sht},\, 
2 k_{\lng} + k_{\sht},\, 3 k_{\lng} + k_{\sht},\, 
k_{\lng} + 2 k_{\sht},\\ 
& 2 k_{\lng} + 2 k_{\sht},\, 3 k_{\lng} + 2 k_{\sht},\, 
3 k_{\lng} + 3 k_{\sht},\, 
4 k_{\lng} + 3 k_{\sht},\, 5 k_{\lng} + 3 k_{\sht};\\
&\hbox{\underline{sht}:\ \, }
k_{\sht},\, k_{\sht},\, 2 k_{\lng} + k_{\sht},\, 2 k_{\sht},\, 
2 k_{\lng} + 2 k_{\sht},\, 4 k_{\lng} + k_{\sht},\, 
4 k_{\lng} + 2 k_{\sht},\\
& 2 k_{\lng} + 3 k_{\sht},\, 4 k_{\lng} + 3 k_{\sht},\, 
6 k_{\lng} + 3 k_{\sht},\, 
6 k_{\lng} + 4 k_{\sht},\, 6 k_{\lng} + 5 k_{\sht}.
\end{align*}}
\smallskip

The $47$ {\em affine exponents} \ in the numerator:
\begin{align}\label{f4affinen}
&(4k_{\lng} ),\  (2 + 4k_{\lng} ),\  (6k_{\lng} ),\  (2 + 6k_{\lng} ),
 (4 + 6k_{\lng} ),\\
&(2k_{\sht} ),\  (1 + 2k_{\sht} ),\  (3k_{\sht} ),\ (1 + 3k_{\sht} ),\
 (2 + 3k_{\sht} ),
\notag\\
&(1 + 2k_{\lng} + 4k_{\sht} ),\
(3 + 2k_{\lng} + 4k_{\sht} ),\ (5 + 2k_{\lng} + 4k_{\sht} ),
\notag\\
&(8k_{\lng} + 2k_{\sht} ),\ (2 +8k_{\lng} + 2k_{\sht} ),\ldots,\
(8 +8k_{\lng} + 2k_{\sht} ),
\notag\\
&(4k_{\lng} + 4k_{\sht} ),\ 
(1 + 4k_{\lng} + 4k_{\sht} ),\ldots,\ (7 + 4k_{\lng} + 4k_{\sht} ),
\notag\\
&(8k_{\lng} + 4k_{\sht} ),\
(2 + 8k_{\lng} + 4k_{\sht} ),\ldots,(10 + 8k_{\lng} + 4k_{\sht} ),
\notag\\
&(1 + 6k_{\lng} + 6k_{\sht} ), 
(3 + 6k_{\lng} + 6k_{\sht} ),\ldots, (11 + 6k_{\lng} + 6k_{\sht} ),
\notag\\
&(12k_{\lng}+6k_{\sht}),\,
(2+12k_{\lng}+6k_{\sht}),\ldots,\,
(16+12k_{\lng}+6k_{\sht}).\notag
\end{align}

Note that the integer step here is $1$ 
for $2k_{\sht}$,$3k_{\sht}$ and $(4k_{\lng}+4k_{\sht})$ or $2$
otherwise.
\smallskip

The $23$ {\em affine exponents} \ in the denominator:
\begin{align}\label{f4affined}
&(2k_{\lng} ),(2k_{\lng}),(k_{\sht}),(k_{\sht}),(2k_{\lng} +k_{\sht} ),
(1 +2k_{\lng} +k_{\sht} ),\\
&(2+2k_{\lng} +k_{\sht} ),\  
(2k_{\lng} + 2k_{\sht} ),\ (2 + 2k_{\lng} + 2k_{\sht} ),
\notag\\
&(4k_{\lng} + k_{\sht} ), (1 + 4k_{\lng} + k_{\sht} ),\ldots,
(4 + 4k_{\lng} + k_{\sht} ),\notag\\
&(6k_{\lng} + 3k_{\sht} ),\ (1 + 6k_{\lng} + 3k_{\sht} ),\ \ldots,\
(8 + 6k_{\lng} + 3k_{\sht} ).\notag
\end{align}
\medskip

\subsection{Affine duality} 
Without going into detail
let us formulate the theorem about the variant of
$\Pi_{\tR}$ for {\em another choice of the affine extension}
of $R$; let 
$\widehat{R}\equal \{[\al\in R,j\in\Z]\}$. For such
affine system, the factors $\nu_\al$
do not appear in the integer components of the affine roots
and, respectively, we will not need $q_\nu$ in the formula
for $\Pi_{\widehat{R}}$.

There is no change in the simply-laced case, so only
the non-simply-laced systems are sufficient to 
consider. 
We come to the following variant of the previous construction.

The corresponding DAHA and the shift-operator
will be now for the pair of affine root system
$\{\widehat{R},\tR\}$ used respectively for $X,Y$ as 
in \cite{C2}. The connection with the radical of $\v$ 
is very much the same as it is for $\{\tR,\tR\}$, which is
the only case considered in this paper (Theorem \ref{RADZERO}).
The duality below is a variant of the Langlands duality
and is directly related to the Fourier 
transform for DAHA of type $\{\widehat{R},\tR\}$.

For this combination of affine extensions, the composition
$\y\circ \x$ is invariant with respect to the 
DAHA\~Fourier transform (treated as an abstract anti-involution);
this leads to the required duality. 
\smallskip 

The rational DAHA do not depend on the particular
choice of the pair of affine root systems and 
{\em rational exponents\,} remain the same for any 
such choices {\em up to proportionality}. The corresponding
binomials in $\Pi_{\widehat{R}}$ are  
$(1-q^{\{\,rat\ exp\, \}})$; upon the substitutions 
$t_\al=q^{k_\al}$,
they can be different from those in $\Pi_{\tR}$,
where we set $t_\al=q^{\nu_\al k_\al}$. Note that
the binomials in terms of $t_{\lng},t_{\sht}$ only
(without nontrivial powers of $q$) remain the same as 
for $\Pi_{\tR}$\,; the Poincar\'e polynomial is unchanged.

\smallskip

\begin{theorem} \label{OTHERAFFINE}
(i) In the cases 
$\widehat{B}_n,\widehat{C}_n,\widehat{F}_4,
\widehat{G}_2$, the product
\begin{align}\label{yofxevalnew}
&\Pi_{\widehat{R}}\equal\prod_{\al\in R_+}
\Bigl( (1- q\,^{k_\al+(\al^\vee,\,\rho+\rho_k)})
\prod_{ j=1}^{(\al^\vee,\,\rho)}
\frac{
(1- q\,^{j-1+k_\al+(\al^\vee,\,\rho_k)})}
{(1- q\,^{j-1+(\al^\vee,\,\rho_k)})}\Bigr)
\end{align}
is a regular function with the corresponding 
rational exponents that are \underline{proportional}\, 
to those from Theorem \ref{BPRIMDEG} for $\Pi_{\tR}$. The  
\underline{affine duality}\, holds\,:
$$\Pi_{\widehat{R^{}}}\,(\,q,\,k_{\lng},\,k_{\sht}\,)\ =\ 
\Pi_{\widehat{R^\vee}}\,(\,q,\,k_{\sht},\,k_{\lng}\,).$$

(ii) The $\Pi_{\widehat{R}\,}$\~product for $\widehat{B}_n $ 
equals 
\begin{align}\label{newbcaffine}
&\frac{\prod_{m=2}^n \prod_{j=0}^{m-1}(1-q^{mk_{\lng}+j})\,
\prod_{m=0}^{n-1}\prod_{j=0}^{2m+1}(1-q^{2k_{\sht}+2mk_{\lng}+j})}
{ (1-q^{k_{\lng}})^{n-1}\, 
\prod_{m=0}^{n-1}\prod_{j=0}^{m}(1-q^{k_{\sht}+mk_{\lng}+j})}
\end{align}
with $\frac{n^2+3n-2}{2}$ \underline{negative} affine exponents.
The proportional pairs of rational exponents here are
$\{mk_{\lng}+j\,,\,lmk_{\lng}+lj\}$ as $j/m\not\in \Z$ and
$lm<n,m>1$.

(iii)
The affine exponents for $\widehat{G}_2$ are
\begin{align}\label{newg2affine}
\hbox{positive:\ }&
(2 k_{\lng}),\, (2 k_{\sht}),\, (1 + 2 k_{\lng}),\, (1 + 2 k_{\sht}),\\
&\{\,(j + 3 k_{\lng} + 3 k_{\sht}),\ j=0,1,2,3,4,5\},\notag\\
\hbox{negative:\ }&
(k_{\lng}),\ (k_{\sht}),\ (k_{\lng}+k_{\sht}),\ 
(1 + k_{\lng}+k_{\sht}),\notag
\end{align}
with $\ 10/4\ $ positive/negative exponents. All rational
exponents are pairwise not proportional.

(iv)
The affine exponents for  $\widehat{F}_4$
are $(lk_{\lng}+sk_{\sht}+j)$, where $j=0,1,\ldots,l+s-1$
and
\begin{align}\label{newf4affine}
\hbox{positive\,:}&\ 
[ls]\in\{ [20],[02],[30],[03],[24],[42],[44],[66]\},\\
\hbox{negative\,:}&\ 
[ls]\in\{ [10],[01],[10],[01],[12],[21],[11],[33]\},\notag
\end{align}
with $\ 42/18\ $ positive/negative exponents. The proportional
pairs of rational exponents 
are $\{[44]+2,[66]+3\}$ and $\{[44]+6,[66]+9\}.$
\sq
\end{theorem}
\medskip

\section{The chain of intertwiners}
\setcounter{equation}{0}
We are going to develop the technique of intertwiners
aiming at decomposing the polynomial representation $\v$
when the action of the $Y$\~operators is non-semisimple.
We will start with some
basic properties of the generalized eigenvectors.

Recall that
$c_\#= c- u_c^{-1}(\rho_k)=\pi_c(\!(-\rho_k)\!)$
for $c\in B,$ where
the affine action $(\!(\cdot)\!)$ from (\ref{afaction})
is used. For instance, $-0_\#=\rho_k$.
See (\ref{Yone}) and (\ref{macd}).

\smallskip
\subsection{Generalized eigenvectors}
This section is for arbitrary nonzero $q,t$ including the case
when $q$ are roots of unity. Almost all statements will require the
constraint $t_\nu\neq 1$.
For instance, we will
constantly use that
$$
q_\al^{(\tal^\vee,c_-+d)-(\al^\vee, \rho_k)}\,\neq\, t_\al^{\pm 1}
\hbox {\ if\ }  q_\al^{(\tal^\vee,c_-+d)-(\al^\vee, \rho_k)}=1.
$$
Sometimes $t_\nu\neq \pm 1$ will be needed; later  
it  will be imposed permanently.
There will be also quite a few claims that hold for
generic $q$ only; this condition will be stated explicitly.

The  space of the {\dfont generalized eigenvectors\,}
$\v^\infty(\xi)\subset\v$ corresponding to a given weight
$\xi=-c_\#$ is as follows:
\begin{align}\label{vyyi}
&\v^s(\xi)\equal
\{v\in \v\ |\ (Y_a-q^{(a,\xi)})^s(v)=0\},\\
&\v(\xi)\ =\ \v^1(\xi),\ \, \v^\infty(\xi)=
\cup_{s>0}\, \v^s(\xi),\notag
\end{align}
for all $a\in B.$
\smallskip

{\em The symbols $q^\xi$ and the weights $\xi$ are always 
identified if the corresponding characters 
$a\mapsto q^{(a,\xi)}$ coincide for $a\in B.$}
\smallskip

The spaces $\v(-c_\#)^{\infty}$ are finite dimensional for
generic $q$ and infinite dimensional when $q$ is a root of unity.
The dimension of the space $\v(-c_\#)^{\infty}$ equals the
number of $b$ satisfying (\ref{vyyi}).
Given an {\em arbitrary one-parametric deformation}
${q}^\natural,\{{t}^\natural_\nu\}$ of $q,t$  which makes
the polynomial representation semisimple, the space
$\v(-c_\#)^{\infty}$ is the limit $\{{q}^\natural\to q,\,
{t}^\natural\to t\}$
of the linear space that is the direct sum
$$
\oplus\, \Q_{{q}^\natural,{t}^\natural}{E}^\natural_b
\hbox{\ such\ that\ }
q^{-c_\#}=q^{-b_\#}
$$
for the Macdonald polynomials ${E}^\natural_b$ defined for
${q}^\natural,{t}^\natural.$

The limit is defined in a sense of vector
bundles over curves, to be exact, over a small one-dimensional
disc with the center at $q,t$.
We will call this limit {\em flat limit}.
The corresponding space is the linear span of the limits
of all linear combinations of ${E}^\natural_b$ divided by proper
powers of the deformation parameter. Obviously,
$\v=\oplus_c \v(-c_\#)^{\infty}$ for pairwise different
$q^{-c_\#}$ (i.e., if the characters $\,q^{-(a,\,c_\#)}\,$ are
different for $a\in B$).

\smallskip

The space $\v(-c_\#)$ always contains at least one Macdonald
polynomial $E_{c}.$ Indeed, one can take $c=c^\circ$
assuming that
\begin{align}\label{primarydef}
&q^{-c_\#}=q^{-c^\circ_\#} \hbox{\ and\ }
q^{-b_\#}\neq q^{-c^\circ_\#}\hbox{\ for\ all\ } b
\hbox{\ such\ that\ }
B\ni b\succ c^\circ.
\end{align}
Under this assumption, there is a
unique $Y$\~eigenvector
of weight $-c_\#$ in the space
$\oplus_{b\succeq c^\circ}\Q_{q,t} X_b;$ it is proportional to
$E_{c^\circ}.$ We will call such $c^\circ$
{\dfont primary} elements; they exist for any $c.$

We are going to apply the simple intertwiners
$P_r,\Psi_i$ from (\ref{tauintery})
to the spaces $\v(-c_\#)$
following (\ref{Phieb}).
The $P_r$ are always well defined and invertible.
The intertwiner $\Psi_i$
is not well defined in $\v(-c_\#)$ if and only if
$\Psi_i^c$ from (\ref{Phijb}), an element in
the non-affine Hecke algebra $\H,$
is infinity. The latter occurs exactly when

\begin{align}\label{interinfty}
&q_\al^{(\tal^\vee,c_- +d)-(\al^\vee, \rho_k)}=1,
\end{align}
where we set
$$
\tal=u_c(\al_i) \and \al=u_c(\al_i) \for i>0,\
\al=u_c(-\vth) \for i=0. 
$$

Sometimes it will be convenient
to renormalize the $\Psi$\~intertwiners as follows:
\begin{align}\label{Psiprime}
&\Psi_i^\diamond=\Psi_i\cdot (Y_{\al_i}^{-1}-1),\ P_r^\diamond=P_r,\
0\le i\le n,\, r\in O.
\end{align}
Given $\hw\in \hW,$ the corresponding $\Psi^\diamond$\~intertwiner
$\Psi_{\hw}^\diamond$ maps $\v^s(-c_\#)$ to $\v^s(-b_\#)$
for $b=\hw(\!(c)\!).$ In contrast to $\Psi_{\hw},$ the
intertwiners $\Psi_{\hw}^\diamond$ are
always well defined, although they can be non-invertible
and identically zero in some $\v(-c_\#).$
\smallskip

\subsection{The 
\texorpdfstring{{\mathversion{bold}$\widetilde{E}$}--polynomials}
{E-tilde-polynomials}}
\label{sec:tildeE}
Let $\widetilde{V}_{0}=\Q_{q,t}1\subset \v.$
Given a reduced decomposition
$\pi_c=\pi_r s_{i_l}\cdots s_{i_1}$ for $c\in B,$
we define $\widetilde{V}_c$
by induction as follows.


(a) Let
$\widetilde{V}_b=\Psi_i (\widetilde{V}_c)=
\Psi_i^\diamond(\widetilde{V}_c)$
for $b=s_i(\!(c)\!)$
unless $\Psi_i^c\in \H$ is infinity.

(b) If $\Psi_i^c\in \H$ is infinity,
i.e., (\ref{interinfty}) holds, then
$\widetilde{V}_b\equal
\widetilde{V}_c+ \tau_+(T_i)(\widetilde{V}_c).$

(c) For $\pi_r\in \Pi^\flat=\{\pi_r,\, \om_r\in B\}$,\ \,let
$\widetilde{V}_b\equal P_r(\widetilde{V}_c)$,\ where
$b=\pi_r(\!(c)\!)$.
\smallskip

Note that $(\al_i,c+d)\neq 0$ if and only if
$s_i\pi_c$ is represented in the form $\pi_b.$
This always holds for ($b$). Indeed, using
$\tal=u_c(\al_i),$ one gets $(\al_i,c+d)=0\, \Rightarrow$
$(\tal,c_- +d)=0$; the latter is impossible due to
(\ref{interinfty}).

Given a reduced decomposition,
$\pi_c=\pi_rs_{i_l}\cdots s_{i_1},$ we set
\begin{align}\label{hatefinal}
&\widetilde{\Psi}_{i_p}\ =\ \Psi_{i_p}\ \, \hbox{\ if\ }\, \
q_\al^{(\tal^\vee,b_-+d)-(\al^\vee, \rho_k)}\neq 1\\
&\hbox{for\ \,} b = s_{i_{p-1}}\cdots s_{i_1}(\!(0)\!),\ \,
\tal=u_b(\al_{i_p}),\notag\\
&\widetilde{\Psi}_{i_p}=\tau_+(T_{i_p}) \hbox{\ if\ }
q_\al^{(\tal^\vee,b_-+d)-(\al^\vee, \rho_k)}=1.
\label{hatefinall}
\end{align}
Equivalently, the elements
$s_{i_p}$ from (\ref{hatefinall}) are {\dfont singular},
that is $\tal^p\in \tR^0$; see the definition of $\tR^0$
in (\ref{troeq}) below.
The intertwiner
$\widetilde{\Psi}_{i_p}$ will be
called singular too.

Given $c\in B$, we define the {\dfont non-semisimple polynomial}
$\widetilde{E}_c\in \widetilde{V}_c$ to be proportional to
$$
P_r\widetilde{\Psi}_{i_l} \cdots
\widetilde{\Psi}_{i_p}\cdots \widetilde{\Psi}_{i_1}(1).
$$
\smallskip

Note that always $\widetilde{\Psi}_{i_p}(\widetilde{E}_c)=0$
if $(\al_{i_p}^\vee,c+d)=0,$ i.e., when
$s_i\pi_c\not\in \pi_B.$ It may happen only for
{\em non-singular} $s_{i_p}$, i.e., when (\ref{hatefinal})
holds.
In contrast to the semisimple case,
the specializations $\Psi_i^c$  
are sufficient only as dim$\widetilde{V}_c=1.$
Generally, the whole $\Psi_i$ in terms of $Y$
from (\ref{tauintery}) have to be involved. 

The definition of $\widetilde{E}_c$ can be naturally
extended to {\em possibly non-reduced} decompositions of $\pi_c$.
In this case, the notation will be $\widetilde{E}^\dag_c$;
we call these polynomials {\dfont non-semisimple non-reduced}.

We will mainly need  $\widetilde{E}$\~polynomials and
$\widetilde{V}$\~spaces for reduced decompositions.
Note that these polynomials are defined so far
only up to proportionality and {\em depend on the
choice of the decomposition} of $\pi_c.$
The $\widetilde{V}$\~spaces {\em depend on the particular
choice of the reduced decomposition} too but their dependence
can be controlled, as we will see later.

Given a decomposition of an {\em arbitrary} $\hw\in \hW$,
one can introduce $\widetilde{E}_{c}$ and $\widetilde{E}^\dag_{c}$
\, for\, $c=\hw(\!(0)\!)$\, in the same manner;
this definition coincides with the previous one
if $\hw=\pi_c$ is taken.


Generalizing, this construction can be originated
at a given Macdonald
polynomial $E_c$ instead of $E_0=1.$ For
instance, one may begin with {\em primary} $c=c^\circ$
from (\ref{primarydef}), when
$q^{-a_\#}\neq q^{-c^\circ_\#}$ for all  $a\succ c^\circ.$
Given $c$, we define $\widetilde{E}_b$ for a decomposition of $\hw$
such that $\hw(\!(c)\!)=b$. The same notation $\widetilde{E}_b$
will be used we will explicitly mention the {\dfont chain origin}
if necessary. The reduced decompositions of $\hw=\pi_b\pi_c^{-1}$
will be mainly needed, i.e., those satisfying
$l(\hw)=l(\pi_b)-l(\pi_c)$;
otherwise $\dag$ will be added.

We call the resulting
sequence $\{\widetilde{E}_b,\ldots,\widetilde{E}_c=E_c\}$
a {\dfont chain originated at} $E_c$.
The same terminology will be
used for the chains of $\widetilde{V}$\~spaces
$\{\widetilde{V}_b,\ldots,
\widetilde{V}_c=\Q_{q,t}\widetilde{E}_c\}$.
\smallskip

Concerning the relation to the spaces of generalized vectors,
it is obvious that $\widetilde{V}_c\subset \v(-c_\#)^\infty.$
The main advantage of $\widetilde{V}_c$ versus  $\v(-c_\#)^\infty$ is
that the former are defined ``locally", following a decomposition
of $\pi_c$. The spaces of generalized eigenvectors
are defined ``globally". Also, generally,
$\widetilde{V}_c$ are smaller than $ \v(-c_\#)^\infty.$
For instance, the spaces
$ \v(-c_\#)^\infty$ are infinite dimensional
when $q$ is a root of unity; the spaces $\widetilde{V}_c$
are always finite dimensional. Here a natural challenge
is to make the spaces  $\{V_c\}$  
generating the whole $\v$ (and containing $E_c$ in the
semisimple case) as small as possible.
\smallskip

\subsection{The spaces 
\texorpdfstring{{\mathversion{bold}$\v_c$}}{Vc}}
A demerit
of  $\widetilde{V}_c$ is that they are not
limits of any natural spaces in
semisimple deformations of $\v$ and may depend
on the reduced decompositions. We will introduce a
somewhat greater system of
finite dimensional spaces $\v_c$ that are such
limits and depend only on the corresponding $c$. 
Their pullbacks $\v_c^\natural$ will be  defined in terms of
the right Bruhat ordering from Section \ref{sec:RightBruhat}
and Theorem \ref{STRIKOUT}.
The sequence of necessary definitions is as follows.

{\em First,} let us introduce a {\em root subsystem}
$\tR^0\subset \tR$
using the relation (\ref{interinfty}) for $c=0$:
\begin{align}\label{troeq}
&\tR^0\equal\{\tal=[\al,\nu_\al j]\in \tR\, \mid \,
q^{\nu_\al j-(\al, \rho_k)}=q_\al^{j-(\al^\vee, \rho_k)}=1\},\\
&\hbox{explicitly,\ \ }
q_\al^j\,\prod_{\nu\in\nu_R}\,t_\nu^{-(\al^\vee,\rho_\nu)}
=1 \for \rho_\nu=\frac{1}{2}
\sum_{\al>0,\nu_\al=\nu}\al\notag.
\end{align}
Recall that the latter product contains
only integral powers of $t_{\sht}$ and $t_{\lng}$ since
$(\rho_k,\al_i^\vee)=k_i=k_{\al_i}$ for $i>0$.

Obviously $\tR^0$ is a {\em root subsystem} of $\tR$, it
is closed with respect to the addition and subtraction 
(if the result is
a root). Note that, generally, it is not true that
$\tR^0$ is an intersection of the $\Q$\~span
$\Q\lan\tR^0\ran$ and $\tR$. Indeed, let $0<q<1$ and
\begin{align}\label{tzeq}
&t_\nu=\ze_\nu\, q_\al^{k_{\al}},\ k_\nu\in (1/N_{\nu})\Z,\
\ze_\nu^{N_\nu}=1 \for N_\nu\in \N.
\end{align}
We may assume that $\ze_\nu$ is a primitive $N_\nu'$-th
root of unity for $N_\nu'\,\mid\, N_\nu$.
Then (\ref{troeq}) is equivalent to
\begin{align}\label{tzeqeq}
&j-(\al^\vee, \rho_k)=0 \and
N_\nu'\, \mid\, (\al^\vee,\rho) \for \nu=\nu_\al.
\end{align}

The divisibility conditions $N_\nu'\, \mid\, (\al^\vee,\rho)$
may not be compatible with taking fractional linear combinations
of $\tal=[\al,\nu_\al j]\in \tR^0$.
Concerning {\em root subsystems},
see Section \ref{sec:RightBruhat}.
In the absence of the roots of unity $\ze_\nu$, when $N_\nu'=1$,
the first equation in (\ref{tzeqeq}) has a solution
if and only if $(\al^\vee,\rho_k)$ is in $\Z_+$ 
(in $\in \N$ as $\al<0$); generally, this
involves divisibility conditions too and may be
incompatible with fractional linear combinations.
\smallskip

{\em Second}, given $c\in B,$ we take $\pi_b\in \b^0(\pi_c),$ i.e.,
$\pi_b$ is obtained from a reduced decomposition
$\pi_c=\pi_rs_{i_l}\cdots s_{i_1}$ by striking out some
of the simple reflections $s_{i_p}$ such that
$\tal^p=s_{i_1}\cdots s_{i_{p-1}}(\al_p)$ belong to $\tR^0$.
Then we pick certain {\em reduced} decompositions of the elements
$\pi_b$  and construct the corresponding
polynomials $\widetilde{E}_b$. Recall that
$\widetilde{E}_{b}=\widetilde{\Psi}_{\pi_b}(1)$.

The {\dfont little generalized
eigenspace} $\v_c$ can be defined now as
the linear span of $\widetilde{E}_b$ introduced above
for the elements $b\in B$.
\smallskip

{\em Third,} $\v_c$ in the semisimple case 
is constructed as follows.
For an arbitrary one-parametric deformation
${q}^\natural,\{{t}^\natural_\nu\}$ of $q,t$ that makes
the polynomial representation semisimple, the space
${\v}^\natural_c$ is defined as the direct sum
$$
\oplus_{\pi_b\in\b^0(\pi_c)}\,\Q_{{q}^\natural,\,{t}^\natural}\,
{E}^\natural_b
$$
for the Macdonald polynomials ${E}^\natural_b$ defined for
the (generic) parameters ${q}^\natural,\,{t}^\natural.$

Due to the definition of the $E$\~polynomials and
thanks to Proposition \ref{BSTAL},
$$
{\v}^\natural_c\,\subset\,
(\oplus_{a\succ c}\Q_{{q}^\natural,\,{t}^\natural}X_a)
\,\oplus\,\Q_{{q}^\natural,\,{t}^\natural}X_c,
$$
and the projection onto the leading monomial $X_c$
is nonzero here.

Theorem \ref{STRIKOUT} guarantees that ${\v}^\natural_c$
does not depend on the choice of the reduced decomposition
of $\pi_c$.

\begin{maintheorem} \label{VINVHAT}
(i) The {\underline flat limit}
$\{{q}^\natural\to q,\, {t}^\natural\to t\}$
of the space
${\v}^\natural_c$ does not depend on the choice of the
one-parametric deformation and coincides with $\v_c$;
this limit is defined in a sense of vector bundles
over a disc at $q,t.$ In particular, $\v_c$ does not depend
on the choice of the reduced decompositions of the elements
$\pi_b$ needed in the definition of $\widetilde{E}_b$.

(ii) Given a reduced decomposition of
$\pi_c$, the polynomials $\widetilde{E}_c$ belong to
$\oplus_{a\succeq c}\Q_{q,t} X_a$
and have a nonzero leading component; from now on they
will be normalized by the relation
\begin{align}\label{ecxfilt}
&\widetilde{E}_c-X_c\in
\oplus_{a\succ c}\Q_{q,t} X_a.
\end{align}
The polynomials
$\{\widetilde{E}_b\, (\pi_b\in \b^0(\pi_c))\}$
are linearly independent and form a basis of $\v_c$ for any
choices of the reduced decompositions of $\pi_b$.

(iii) The space $\v_b$ defined for reduced $\pi_b=s_i\pi_c$, i.e.,
when $l(\pi_b)=l(\pi_c)+1$,
satisfies the property
$\v_c+ \tau_+(T_i)(\v_c)\subset \v_b$\, if
(\ref{interinfty}) holds (case ($b$) above);
$\widetilde{V}_b=$ $\widetilde{V}_c+ \tau_+(T_i)(\widetilde{V}_c)$
for such $s_i$.
In particular, the space $ \widetilde{V}_c$ belongs to
$\v_c$ for any $c\in B$; also dim$\,\v_b>$dim$\,\v_c$ and
dim$\,\widetilde{V}_b>$dim$\,\widetilde{V}_c$.

(iv) The dimension of $\v_b$ for a reduced $\pi_b=s_i\pi_c$
decreases if and only if\,
$\Psi_i^c=\tau_+(T_i)-t_i^{1/2}$,\, i.e., when
\begin{align}\label{hatefinalminus}
&q_\al^{(\tal^\vee,c_-+d)-(\al^\vee, \rho_k)+k_\al}= 1,
\end{align}
and also there exists $\pi_a\in \b^0_o(\pi_c)$
such that $s_i\pi_a\not\in \pi_B$ for
a certain $a\in B$, equivalently, $(\al_i,a+d)=0$.

(v) For invertible $\Psi_i^c$,
dim$\,\v_b$=dim$\,\v_c$ and
dim$\,\widetilde{V}_b$=dim$\,\widetilde{V}_c$.
When $t_i\neq -1$ and 
$\Psi_i^c=\tau_+(T_i)+t_i^{-1/2}$, i.e. when 
\begin{align}\label{hatefinalplus}
&q_\al^{(\tal^\vee,c_-+d)-(\al^\vee, \rho_k)-k_\al}= 1,
\end{align}
dim$\,\v_b>$dim$\,\v_c$  if and only if
there exists $\hw'\in \b^0_o(\pi_c)$
such that $\hw'=s_i\pi_a, l(\hw')=l(\pi_a)+1$ for
a certain $a\in B$ and $\hw'\not\in \pi_B$.
\end{maintheorem}

{\em Proof.} Let $\pi_c=\pi_rs_{i_l}\cdots s_{i_1}$ be a reduced
decomposition. For any $q',\,t'$,
one can introduce the polynomials
$$\widetilde{E}^{q',\,t'}_c=
P_r\widetilde{\Psi}_{i_l} \cdots
\widetilde{\Psi}_{i_p}\cdots \widetilde{\Psi}_{i_1}(1),
$$
where $\widetilde{\Psi}_{i_p}$ is either $\Psi_{i_p}$ or
$\tau_+(T_{i_p})$ according to (\ref{hatefinal}), namely,
$$
\hbox{where\, either\ \,}
\tal^p=s_{i_1}\cdots s_{i_{p-1}}(\al_{i_p})\not\in \tR^0
\hbox{\ \, or\ \,} \tal^p\in \tR^0.
$$

Using the definition of $\Psi_i$,
$$
\Psi_{i_p}-(t_{i_p}^{1/2}-t_{i_p}^{-1/2})
(Y_{\al_{i_p}}^{-1}-1)=\tau_+(T_{i_p}).
$$
These formulas express $\tau_+(T_{i_p})$
in terms of the $Y$\~intertwiners $\Psi$; they hold for
$\widetilde{E}^{q^\natural,\,t^\natural}_c$ or for any
$\widetilde{E}^{q',\,t'}_b$, where
$\pi_b\in \b^0(\pi_c)$. Since the Macdonald polynomials
are $Y$\~eigenvectors, we readily obtain that the polynomials
$\widetilde{E}^{q^\natural,\,t^\natural}_b$ belong
to ${\v}^\natural_c$ for such $b.$
The limits of these polynomials are exactly
$\widetilde{E}_b$ by construction.

We see that
\begin{align}\label{ecfilt}
&\widetilde{E}_c \in \oplus_{a\succeq c}\Q_{q,t} X_a,
\end{align}
and that the projection of $\widetilde{E}_c$ onto
$\Q_{q,t} X_c$ is nonzero. The next lemma makes
this fact explicit and also shows that
the latter claim  formally follows from (\ref{ecfilt}).

\begin{lemma}\label{PSIHATE}
For $b,c,\hw$ satisfying $\pi_b=\hw\pi_c$ subject to
$l(\pi_b)=l(\hw)+l(\pi_c)$, i.e., if
$\pi_b=\hw\pi_c$ is reduced,
the monomial $X_b$ appears in
$\widetilde{\Psi}_{\hw}(\widetilde{E}_c)$
with a nonzero coefficient.
\end{lemma}

{\em Proof.}
It suffices to consider
$\hw=s_i.$ Let $(\al_i^\vee,c+d)>0.$
It is straight to check (see below) that, given $i>0$, the
non-affine simple intertwiners
$\Psi_i$ transfers
the space from (\ref{ecfilt}) to that for
$s_i(\!(c)\!)=s_i(c)$, because so does $T_i.$
However this argument is not applicable
to $\Psi_0$; this is the main part of the claim.
\smallskip

{\em First,} we assume that the intertwiner
$\Psi_i$ is not infinity, i.e., $\Psi_i$ is non-singular. 
There are two 
subcases, $i>0$ and $i=0$.
Let us begin with $i>0.$
Then $\widetilde{\Psi}_i(X_c)$ satisfies
(\ref{macdgenn}) for $b$ and
$$
\Psi_i(\widetilde{E}_c)-
t_i^{-1/2}X_b\in \oplus_{a\succ b}\Q_{q,t} X_a.
$$

When $i=0,$ we
follow the end of the proof of Theorem 5.1 from \cite{C1}
(see also Theorem \ref{PHIEB} above)
and come to the same relation with the coefficient
$t_i^{-1/2}q^{(c,c)/2-(b,b)/2}$ instead of $t_i^{-1/2}.$
Here it is sufficient to know that $\Psi_i(\widetilde{E}_c)$ is
a linear combination of $X_b$ modulo lower terms; then we
only need to calculate the coefficient of $X_b$ in
$(\tau_+(T_0))(X_c).$
\smallskip

{\em Second,} let $\Psi_i$ be infinite (singular).
The consideration
is analogous. For $i>0,$ the element $T_{i}$ maps
$$
\oplus_{a\succeq c}\Q_{q,t} X_a\mapsto
\oplus_{a\succeq b}\Q_{q,t} X_a.
$$
For $i=0,$ we
use again that the unwanted terms $X_a$ will not appear when
$\widetilde{\Psi}_i$ is applied; then the coefficient
of $X_c$ in $\tau_+(T_0)(X_b)$ is calculated.\\
\sq

The consideration of the leading terms readily
gives that the polynomials
$\widetilde{E}^{q^\natural,t^\natural}_b$ for
$\pi_b\in \b^0(\pi_c)$ form a basis of ${\v}^\natural_c$
and
$\{\widetilde{E}_b\}$ is a basis
of the limit of  ${\v}^\natural_c$. Therefore the latter
does not depend on the choice of the one-parametric
deformation and \ $\lim
 {\v}^\natural_c=$ $\v_c$.

As a by-product,
we obtain that the latter space does not
depend on the choice of the reduced decomposition of $\pi_c$.
As we will see, the interpretation of $\v_c$ as a limit is not
necessary here; this fact follows directly
from Theorem \ref{PHIBRUHAT}.

Thus claims (i) and (ii) are checked; (iii) results from
Theorem \ref{STRIKOUT},($e$). Claims (iv,v) are straightforward;
we use that (\ref{hatefinalplus}) and (\ref{hatefinalminus})
are invariant under the action of $\tW^0$ and crossing out
singular reflections in $\pi_c$.

For instance, $\v_b$ can become smaller than $\v_c$ for  
reduced $\pi_b=s_i\pi_c$ 
if and only if there exists $\pi_a\in \b^0_o(\pi_c)$
such that $s_i\pi_a\not\in \pi_B$ for
a certain $a\in B$, equivalently,  $(\al_i,a+d)=0$, equivalently,
$\pm\pi_a^{-1}(\al_i)$ is a simple non-affine root. The latter
means that (\ref{hatefinalminus}) holds.
\sq
\smallskip

\subsection{Further properties}
Let us begin with some applications to the
Macdonald polynomials.
Given $c\in B$ and a reduced decomposition
$\pi_c=\pi_rs_{i_l}\cdots s_{i_1}$,
we have the following embeddings:
\begin{align}\label{vspacesall}
&\widetilde{V}_c\subset\v_c\subset \v(-c_\#)^\infty.
\end{align}
The dimension of $\v(-c_\#)^\infty$ equals one if and only if
the coset $\pi_c\hW^\flat[-\rho_k]$ does not contain the elements
in the form $\pi_b$ for $b\neq c$,
where
$$
\hW^\flat[\xi]\equal
\{\,\hw\in \hW^\flat\, \mid\, q^{\hw\llb \xi\rrb}=q^{\xi}\,\}.
$$
In this case, the eigenvalue $q^{-c_\#}$ is
$Y$\~simple and the Macdonald polynomial $E_c$ exists.
All such $c$ are {\em primary}, but, generally, the set
of primary $c=c^\circ$ is broader; see (\ref{primarydef}).

Note that if  $\pi_c\tW^0$ does not contain the elements
$\pi_b$ for $b\neq c$ for $\tW^0=\lan s_{\tal}\,\mid\,
\tal\in \tR^0\ran$ (see Proposition \ref{BRUHATLA}),
then dim$\v_c=1$. Since $b\succ c$ for such $b$, primary $c=c^\circ$
automatically satisfy this condition. Also
$\la(\pi_c)\cap \tR^0=\emptyset \Rightarrow$ dim$\v_c=1$,
i.e., the dimension is one if there are no singular $s_{i_p}$
in reduced decompositions of $\pi_c$. 

Generally,
dim$\,\v_c=1$ {\em if and only if} there are no elements in
the form of $\pi_b$ in the set $\b^0_o(\pi_c)$.
Any $\v_c$ contains at least one Macdonald polynomial;
indeed, $E_{c^\circ}$ for
{\em any} primary $c^\circ$ (there can be several for
a given element $c$) can be taken. 

The equality dim$\v_c=1$ guarantees
that $E_c$ exists and does not depend on the presentation
of $\v$ as a limit of the generic
semisimple polynomial representation ${\v}^\natural.$
The equality dim$\widetilde{V}_c=1$ gives that the Macdonald
polynomial $E_c$ exists, i.e., is a $Y$\~eigenfunction,
have the required structure of its monomials and a
nonzero leading term, however it may depend
on the choice of the limiting procedure.

If the limit of
$\widetilde{E}^{q^\natural,\,t^\natural}_c$ exists only for
a {\em certain choice} of the deformation parameter then
this limit is an  eigenvector from the $Y$\~eigenspace
$\v_c\cap \v(-c_\#)$. As a matter of fact, the whole space
$\v_c$ with the filtration corresponding to $>_0$ is
such a limit if we switch here to 
more general understanding of the limiting
procedure (involving induced $\HH^\flat$\~modules).

Thus:
\begin{align}\label{dimoneimplications}
&\hbox{dim}\widetilde{V}_c=1\,\Leftarrow\,\hbox{dim}\v_c=1
\,\Leftarrow\,
c=c^o\,\Leftarrow\, \hbox{dim}\v(-c_\#)^\infty=1,\\
&\hbox{dim}\v_c=1\ \Leftrightarrow\
\hbox{\ existence\ and\ \underline{total}
\ uniquiness\ of\ } E_c.\notag
\end{align}
\smallskip

Let us check that dim$\v(-c_\#)^\infty=1$ for sufficiently
big $c\in B$ if $q$ is not a root of unity and, therefore,
$\hW^\flat[-\rho_k]$ is finite.
Recall that
\begin{align*}
&\pi_c=cu_c^{-1}=u_c^{-1}c_- \for
c_-=u_c(c)\in B_-\,,\, u\in W,\\
&c_\#= c- u_c^{-1}(\rho_k)=\pi_c(\!(-\rho_k)\!),\
-0_\#=\rho_k,
\end{align*}
where
the affine action $(\!(\cdot)\!)$ from (\ref{afaction})
is used.

The following condition is sufficient for dim$\v(-c_\#)^\infty=1$:
\begin{align}\label{cbhw}
&(c_- +[\hw]_b,\al_i)<0 \hbox{\ for\ all\ } i>0 \hbox{\ and \ all\ }
\hw\in \hW^\flat[\rho_k],
\end{align}
where we use the decomposition $\hw\equal[\hw]_b[\hw]_u$ with
$[\hw]_b\in B,\, [\hw]_u\in W.$
This condition means that $c_-+[\hw]_b$ is {\em anti-dominant}
for all $\hw\in\hW^\flat[\rho_k]$, which always holds
for {\em sufficiently large} $b$\, provided
that $\hW^\flat[\rho_k]$ is finite, equivalently, $q$ is not
a root of unity. Let us check it.

First, $[\hw]_u\neq$ id for
$\hW^\flat[\rho_k]\ni\hw\neq$ id; otherwise $[\hw]_b+\rho_k=\rho_k$
and $[\hw]_b=0$. Second,
$$
\pi_c\hw=(u_c^{-1}[\hw]_u)([\hw]_u^{-1}(c_- +[\hw]_b)),\
\pi_{c}\hw\in \pi_B \Leftrightarrow [\hw]_u^{-1}(c_- +[\hw]_b)
\in B_-.
$$
The latter and relations (\ref{cbhw}) imply that $[\hw]_u=$id.

Replacing here the
{\em complete stabilizer} $\hW^\flat[-\rho_k]$ by its subgroup
$\tW^0$, we obtain the following {\em criterion} for
dim$\v_c=1$ provided that $q$ is not a root of unity:

\begin{align}\label{vcdimo}
&(c_- +[\tw]_b,\al_i)<0 \hbox{\ \,for\ all\ \,} i>0 
\hbox {\,\ and\ all\ \,}
\tw\in \tW^0,
\end{align}
where $\tw=[\tw]_b[\tw]_u$ is defined as above.

\medskip
\section{The structure of 
\texorpdfstring{{\mathversion{normal}$\v_c$}}{Vc}}
The following proposition
provides exact tools for
calculating the spaces $\v(-c_\#)^\infty$
and $\v_c$. 
\smallskip

\subsection{Multiplication in 
\texorpdfstring{{\mathversion{bold}$\pi_B$}}{pi-B}}
We continue using the decomposition
$\hw=[\hw]_b[\hw]_u$ with $[\hw]_u\in W$ and
$[\hw]_b\in B$ and other notations from the previous 
section.
\smallskip

\begin{proposition}\label{HWALBE}
(i) The element $\pi_c\hw$ for $c\in B,\hw\in \hW^\flat$
can be represented in the form $\pi_b$ if and only if
the following three conditions hold for \underline{every}
$\al\in R_+$,
\begin{align*}
(a)\ \ \ \, \ \hw(\al)&\ \not\in\ -[R_+,Z_+]\,,\ \, equivalently,\\
&\la(u)\cap \la(\hw)=\emptyset, \hbox{\ where\ } u=[\hw]_u,\\
(b)\hbox{\ if\ \ } \hw(\al)&\ =\ \,
[-\be,j_\circ\nu_\be] \for \be>0,j_\circ> 0,\\
&\hbox{then\ }[-\be,j_\circ\nu_\be]\not\in\la(\pi_c),\\
(c)\hbox{\ if\ \ } \hw(\al)&\ =\ \,
[\be,-j_\circ\nu_\be] \for \be>0,j_\circ> 0,\\
&\hbox{then\ } [-\be,j_\circ \nu_\be]\in\la(\pi_c).\\
\end{align*}
\noindent
Here $\nu_{\be}=\nu_{\al}$ in ($b,c$). Imposing
($a$), the roots $\al\in R_+$ satisfying ($b$) constitute
all $\la(u)\ni\al\not\in\la(\hw)$, ($c$) describes all
$\la(\hw)\ni \al\not\in \la(u)$.

(ii) Let $q$ be not a root of unity. Then the elements
$\,\hw\in\hW^\flat[-\rho_k]$ have pairwise distinct
\,$W$\~projections\, $u=[\hw]_u$; also,
the elements $[z,j]$ from the
$\Z$\~span $\tilde{Q}^0$ of $\tR^0$ have pairwise distinct
$z$\~components. For instance,
$u=$id$\, \, \Leftrightarrow\, \hw=$id and,
given $\hw\in\hW^\flat[-\rho_k]$, the component
$j_\circ\nu_\be$ in ($b$) or ($c$) can be uniquely 
determined from the relations
$$
(b'):\, [\be+\al,-j_\circ\nu_\be]\in \tilde{Q}^0,\ \
(c'):\, [\be-\al,-j_\circ\nu_\be]\in \tilde{Q}^0.
$$
Condition ($a$) holds for any $\hw\in\hW^\flat[-\rho_k]$ if
\begin{align}\label{knegative}
&(\Z_+k_{\sht}+Z_+k_{\lng})\cap\Z_+=\{0\}.
\end{align}

(iii) If $q$ is not a root of unity and ($a$) holds for
$\hw=bu\in\hW^\flat[-\rho_k]$,
then $\pi_c\hw$
is \underline{not} in the form $\pi_b$ if there exists
\underline{at least one}
$\al\in R_+$ such that $\be=-u(\al)>0$ and also
\begin{align}\label{suffbig}
&[-\be,j\nu_\be]\in \la(\pi_c) \for j>0
\hbox{\ satisfying\ }
[\be+\al,-j\nu_\be]\in \tilde{Q}^0.
\end{align}

Let $\hW^\flat[-\rho_k]_a$ be the set of elements
$\hw$ from $\hW^\flat[-\rho_k]$ under condition ($a$) and
$[\,\hW^\flat[-\rho_k]_a\,]_u$ the set of their $W$\~projections.
If $\la(\pi_c)$
contains at least one $[-\be,j\nu_\be]$ for every
$\be=-u(\al)$ satisfying (\ref{suffbig}),
where $u\in [\,\hW^\flat[-\rho_k]_a\,]_u$, then
dim$\,\v(-c_\#)^\infty=1$.

(iv) Given $\hw=bu\in \hW^\flat[-\rho_k]_a$,
let us assume that
(iii) does not hold, i.e.,
no $\al$ exist satisfying (\ref{suffbig})
for $\hw$ and $u=[\hw]_u$.
Then $\pi_c\hw$ is \underline{not}
from $\pi_B$ \underline{if and only} if
there exists \underline{at least one} $\al\in R_+$
such that $\be=u(\al)>0$ and
\begin{align}\label{suffsmall}
&[-\be,j\nu_\be]\not\in \la(\pi_c) \for j>0
\hbox{\ satisfying\ } [\be-\al,-j\nu_\be]\in \tilde{Q}^0.
\end{align}

When  $\pi_c=$id, dim$\,\v(-0_\#)^\infty=1$ if and only if
$\al$ satisfying (\ref{suffsmall})
exist for the $W$\~projection $u=[\hw]_u$ of every
element $\hw\in \hW^\flat[-\rho_k]_a$.
\end{proposition}

{\it Proof}. If $\hw(\al)=[-\be,-j\nu_b]$ for $\be>0,j\ge 0$,
then $\al\in\la(\hw)$ but
$-\al\not\in\hw^{-1}(\la(\pi_c)$ because
$-\al= \hw^{-1}([\be,j_\circ\nu_\be]).$ Thus $\al$ will
appear in $\la(\pi_c\hw)$ and the latter set cannot be in
the form $\la(\pi_b)$ for any $b$.
Recall that $\la(\pi_c\hw)$
is obtained from $\hw^{-1}(\la(\pi_c)\cup \la(\pi_c)$
by removing all pairs $\{\tal,-\tal\}$.

Condition ($b$) describes $\al\in R_+$
that may appear in $\la(\pi_c\hw)$ because of $\hw^{-1}(\pi_c)$;
here $\hw^{-1}([-\be,j_\circ\nu_\be])=\al$
and $[-\be,j_\circ\nu_\be]$ must not be from $\la(\pi_c)$.

Condition ($c$) gives $\al$ from $\la(\hw)$; here
$\hw^{-1}([-\be,j_\circ\nu_\be])=-\al$ and
$\al$ will {\em not} appear in $\la(\pi_c\hw)$ only if
$[-\be,j_\circ\nu_\be]\in \la(\pi_b)$.

Note that if $j>j_\circ$ under ($c$), then
$\hw^{-1}([-\be,j\nu_\be])$ becomes a positive root with the
negative non-affine component $\al$; such roots can appear
in $\pi_b$.

The equivalence from ($a$) and the
interpretation of ($b,c$) under\, $\la(\hw)\cap\la(u)=\emptyset$\,
follow directly from\ 
$\hw(\al)=bu(\al)=[u(\al),-(b,u(\al))]$.

As for (ii), given
$u\in [\hW^\flat[-\rho_k]]_u= W\cap(\hW^\flat[-\rho_k]\,B)$,
there exists a unique
$\hw\in \hW^\flat[-\rho_k]$ such that $\hw=bu$ for a certain
$b\in B$, since $q$ is not a root of unity. It is analogous
for $\tilde{Q}^0$.

Let us assume that $(\Z_+k_{\sht}+Z_+k_{\lng})\cap\Z_+=\{0\}$.
For instance, the conditions $\Re(k_\nu)<0$ for all $\nu$
are sufficient to impose. 
If ($a$) does not hold for such $k$,  i.e.,
$\hw(\al)\not= [-\be,-j\nu_\be]$ for
$\{\al>0,\be>0,j\ge 0\}$,
then
$$(\al,-\rho_k)=([-\be,-j\nu_\be],-\rho_k+d)\,=\,
(\be,\rho_k)-j\nu_\be
$$
and $(\al+\be,\rho_k)=j \nu_\be\ge 0$, which is impossible;
see (\ref{troeq}).

Claims (iii) and (iv) correspond respectively to
cases ($b$) and ($c$) for $\hw \in \hW^\flat[-\rho_k]$.
\sq
\smallskip

%

Claim (iii) from the proposition gives that
the eigenspace $\v(-c_\#)^\infty$ is one-dimensional
for ``sufficiently big" $c$, (iv) gives that dim$\,\v_c=1$
for ``sufficiently small" $c$ (see below).
We will mainly need (iii). 
\smallskip

\rmk
We note that (i) can be used in the theory of
Schubert manifolds of the affine Grassmanian
defined for the maximal parahoric subgroup in
the corresponding loop group. \sq
\smallskip

Let us switch to the spaces $\v_c$,
which requires changing $\hW^\flat[-\rho_k]$ to $\tW^0$.
This reduction reduces the number of
possibilities for $\tu$ in applications. Recall that
$\tW^0$ is defined entirely in terms of $\tR^0$ from
(\ref{troeq}) and is, generally,
simpler to control then the whole centralizer
$\hW^\flat[-\rho_k]$ of $\rho_k$ in $\hW^\flat$. Also, 
the Bruhat ordering  $\pi_c\tu\,<_0\, \pi_c$ is simpler
for such $\tu$; see
Theorem \ref{STRIKOUT} and Proposition \ref{BRUHTW0},(iii).

\begin{corollary}\label{NONDECR}
We employ the proposition as \,$q$\, is not a root of unity
switching to the set $\tW^0_a=\tW^0\cap\hW^\flat[-\rho_k]_a$.
If for \underline{every} $u=[\tu]_u$ for $\tu\in\tW^0_a]_u$
such that $\,\pi_c\tu <_0 \pi_c\, $, the set
$\la(\pi_c)$ contains at least one $[-\be,j\nu_\be]$ with
$\be=-u(\al)>0$ satisfying (\ref{suffbig}),
then dim$\,\v_c=1$.

(i) Let us assume that $\la(\pi_c)$  contains the set
$\tR^1_+[-]=\tR^1\cap \tR_+[-]$, where
\begin{align}\label{critrootsp}
&\tR_+[-]\equal\{\tal=[-\al,j\nu_\al],\,\al>0,j>0\} \hbox{\ and\ }\\
&\tR^1 \equal\{\,\tal=[-\al,j\nu_\al]\in \tR\, \mid \,
q_\al^{k_\al+ j+(\al^\vee,\rho_k))}=1\, \}.\label{critrootone}
\end{align}
Let $j=0$ in (\ref{critrootone}) occur for simple $\al$ only,
for instance, this always holds as $t_\nu$ are not roots of unity.
This condition implies that $[-\be,j\nu_\be]$ from ($b$) belongs to
$\tR^1_+[-]$ if and only if $\al=\al_i$ for some $i>0$.
Then dim$\,\v_c=1$ and, moreover, $\la(\pi_c)$ contains $\tR^0$.

(ii) Continuing to impose that $j=0$ in (\ref{critrootone})
holds for simple $\al$ only, let $S$ be a subset of
$\{1,2,\ldots,n\}$. The root lattice $Q$ for $R$ will be
used. We assume now that
\begin{align}\label{dimvonep}
& q^{2\rho_k^u}\,\not\in\, q^{Q}, \hbox{\ where\ }
2\rho_k^u\,\equal\,\sum_{\al\in\la(u)}\,k_{\al}\al,
\hbox{\ for \ all\ elements } \\
& u\in W_1^S\equal\{w\in W \mid
\la(u)\cap \{\al_1,\ldots,\al_n\}=\al_i,\ i\in S\}.\notag
\end{align}
Then dim$\,\v_{c'}=1$ and, moreover, $\la(\pi_{c'})$
contains $\tR^0$
for any $\pi_{c'}$ such that
$\la(\pi_{c'})$ contains all roots
$[-\be,j\nu_\be]\in\tR^1_+[-]$ unless $\be$ is 
in the form $-u(\al_i)$\ for $i\in S$ and for\ 
$u\in W_1^i\subset W_1^S.$
\end{corollary}

{\em Proof.}
The first claim is from
Proposition \ref{HWALBE}, part (iii), where 
$\tu\in\tW^0$ is taken.
The statement that $\la(\pi_c)$ contains $\tR^0$
formally follows from the relation dim$\,\v_{b}=1$ if
the latter holds for {\em all}
$b$ such that $\la(\pi_c)\subset\la(\pi_b)$. Indeed,
if $\tR^0_+\not\subset \la(\pi_c)$ then there exists
$b$ such that $l(\pi_b)=l(\pi_b\pi_c^{-1})+l(\pi_c)$
and $\la(\pi_b)\setminus \la(\pi_c)$ contains {\em singular}
$\tal$. Therefore, dim$\,\v_{b}>1$ for such $b$ with a singular
{\em last root} in $\la(\pi_b)$.

Let us check (i). The second condition concerning $j=0$
readily gives that
$\tu(\al)=[-\be,j\nu_\be]\in\tR^1_+[-]$ under ($b$) (from
the proposition) holds
for $\al\in R_+$ and that $\tu\in \hW[\-\rho_k]^\flat$ if and
only if $\al=-u^{-1}(\be)$ is a simple root.

Let $\tu\in \tW^0_a$. Then there exists at
least one simple $\al_i\in \la(u)$. Recall that ($a$) from
Proposition \ref{HWALBE} is imposed (see the definition
of $\tW^0_a$).
Following (\ref{suffbig}) from (iii)
(corresponding to case ($b$) from
this proposition), we take $-\be=u(\al_i)$ for $u=[\tu]_u$
and extend it to $[-\be,j\nu_\be]\in\tR^1_+[-]$.
Since $\la(\pi_c)$ contains {\em all} elements from
$\tR^1_+[-]$ by assumption, (i) is verified.

The demonstration of (ii) is a variant of the same argument.
If $\la(u)$ contains at least two simple roots, then we can
proceed as in (i); recall that
$W_1^S$ is the set of $w\in W$ such that
$\la(u)$ contains exactly one simple root $\al_i$ for
certain $i\in S$.
The elements $\pi_{c'}\tu$ for $\tu$
with $u\not\in W_1^S$ cannot be represented in the form
$\pi_b$ by assumption. Concerning $u\in W_1^S$, the
existence of its extension $\tu=bu\in \hW^\flat[-\rho_k]$ 
is {\em equivalent}
to the ``$q$\~integrality" of $2\rho_k^u=\rho_k-u(\rho_k)$,
namely, to the relation
\begin{align}\label{utotu}
&q^{\rho_k-u(\rho_k)}\,=\,q^{2\rho_k^u}\,\in\, q^{Q},
\end{align}
which is not allowed in (ii) for $i\in S$.
\sq
\medskip

\subsection{The semisimple submodule}
Let us assume that $q$ and $\{t_\nu\}$
are not roots of unity and, also, 
$q^a\prod_\nu t_\nu^{b_\nu}=1$ implies
$a+\sum_\nu\nu k_\nu b_\nu=0.$ 

 Then
(\ref{critrootone}), (\ref{dimvonep}) read respectively as
\begin{align}\label{critrootoneu}
&\tR^1=\{\,\tal=[-\al,j\nu_\al]\in \tR\,
\mid \, k_\al+ j+(\al^\vee,\rho_k)=0\, \},\\
&2\rho_k^u\ =\ \sum_{\al\in\la(u)}\,k_{\al}\al\,\not\in\, Q
\hbox{\ for\ all \ } u\in W_1.
\label{tworhou}
\end{align}
Here $2\rho_k^u$  belongs to $Q$ if and only if
$\tu=bu\in \hW^\flat[-\rho_k]$ exists; $\tu$ is a unique
pullback of such $u$ (if it exists).

We assume that $q,t$ satisfy this condition in the
next theorem to make its statement more transparent;
generally, the ``$q$\~integrality" is sufficient to use
instead of the condition for $-\be$ below, namely,
$$q^{-\be+\nu_\be\Z_+}\cap q^{\tR_+^1[-]+\tR_+^0[-]}
=\emptyset.$$
\smallskip

\begin{theorem} \label{FEIGIN}
Continuing part (i) of Corollary \ref{NONDECR},
we impose the inclusions  $\tR_+^1\subset\la(\pi_c)$,
which result in
dim$\, \v_c=1$ and $\tR_+^0\subset\la(\pi_c)$. Let
$\v_{ss}$ be a linear space with a basis
$E_c=\widetilde{E}_c$. We also
impose (\ref{dimvonep}) from (ii) for $u\in W_1^i$
such that $-\be=u(\al_i)<0$ is \underline{not} a sum
of the \underline{non-affine projections}
of the roots from $\tR_+^1[-]$ and $\tR_+^0[-]$. Respectively,
$S$ will be the set of $i$ when such $-\be$, let us call them
\underline{indecomposable}, exist for at least one $u\in W_1^i$.

Then the space $\v_{ss\,}$ is a $Y$\~semisimple
$\HH^\flat$\~submodule of $\v$, which is irreducible if and
only if $\tR_+^{-1}\subset\la(\pi_c)$ for all $c$ corresponding
to $E_c\in \v_{ss\,}$, where
\begin{align}\label{critrootoneuminus}
&\tR^{-1}=\{\,\tal=[-\al,j\nu_\al]\in \tR\,
\mid \, -k_\al+ j+(\al^\vee,\rho_k)=0\, \}.
\end{align}
\end{theorem}

{\em Proof.} We have already checked the semisimplicity
of $\v_{ss\,}$. The only possibility to get an element apart
from $\v_{ss\,}$ when applying the generators of $\HH^\flat$
to $E_c\in \v_{ss\,}$ is when $X_{\al_j}(q^{c_\#})=t_j$ for some
$j\ge 0$. Indeed, it can happen only when
$s_j\pi_c$ is {\em not} a reduced decomposition
(then $s_j\pi_c$ can be represented in the form $\pi_{c'}$) with
$s_j$ corresponding to a non-invertible intertwiner $\Psi_j^c$.
The latter cannot be $\tau_+(T_j)-t_j^{1/2}$ because in this
case the condition $\tR^1_+[-]\subset \la(\pi_c)$ remains
unchanged for $\pi_{c'}=s_j\pi_c$ (if the latter belongs to $\pi_B$).
Thus it must be $\tau_+(T_j)+t_j^{-1/2}$, equivalently,
$X_{\al_j}(q^{c_\#})=t_j$.

Since dim\,$\v_c=1$, we can assume that $s_j$ is
the last simple reflection in a reduced decomposition of
$\pi_c$ and the {\em last root} $[-\be,j\nu_\be]$ in
$\la(\pi_c)$ belongs to $\tR_+^1$. If dim\,$\v_{c'}=1$ for
(non-reduced) $\pi_{c'}=s_j\pi_c$ above, then
$(\tau_+(T_j)+t_j^{-1/2})(E_c)=0$ and
$\tau_+(T_j)(E_{c})\in \v_{ss\,}$, i.e.,
we stay within $\v_{ss\,}.$

Generally, $\v_{c'}$ can be of dimension
greater than $1$ because one root from $\tR_+[-]$, namely
$[-\be,j\nu_\be]$, may be missing in $\la(\pi_{c'})$.
The component $-\be$ of this root is {\em indecomposable}
(see above). Indeed, if the root
$[-\be,\nu_\be j]\in \tR_+[-]$ is the {\em last root}
in $\la(\pi_{c})$, then it is {\em not} a sum of two roots
in this set, which contains $\tR_+^1[-]$ by construction
and $\tR_+^0$ due to Corollary \ref{NONDECR},(i).

Now let us follow part (ii) of this corollary. We use that
dim$\, V_{c'}$=1 if for {\em any} $u\in [\tW^0]_u$
we can find {\em at least one} $\al\in R_+$ satisfying
(\ref{suffbig}). The only problem with finding such
$\al$ may occur when
$\la(u)$ contains precisely one {\em simple} root $\al_i$.
In this case, dim$=1$ if $\la(\pi_{c'})$ contains $\tu(\al_i)$
assuming that the extension of \,$u$ to $\tu=ub\in\tW^0$
exists. However only one root, $[-\be,j\nu_\be]$,
from $\tR^1_+[-]$ can be absent in $\la(\pi_{c'}).$
Therefore, dim$\, V_{c'}=1$ unless $[-\be,j\nu_\be]=\tu(\al_i)$. We
do not allow the latter simply imposing the condition that such $u$
do not have extensions in $\tW^0$.

Due to the semisimplicity of $\v_{ss\,}$, its irreducibility
is {\em equivalent} to the absence of non-invertible intertwiners
in {\em reduced} decompositions of $\pi_b$ after a certain
$\pi_c$ corresponding to $E_c\in \v_{ss\,}$.
\sq
\smallskip

The condition $\tR_+^{1}[-]\subset \la(\pi_c)\Rightarrow
\tR_+^{-1}[-]\subset \la(\pi_c)$, which ensures the irreducibility
of $\v_{ss\,},$ holds almost always,
for instance, in the simply-laced case.
See Lemma \ref{ZFROMM} below.

When $k_{\sht}=k=\nu_{\lng}k_{\lng}$, 
only the roots
with the smallest
possible $(\be,\rho^\vee)$ satisfying $k((\be,\rho^\vee)+1)\in \N$
can be the {\em last roots} in $\la(\pi_c)$ provided that
this set contains
$\tR^1_+[-]$. The constraint ($a$) from Proposition
\ref{HWALBE} simply means that
$k<0$ in this case (only
rational $k<0$ are sufficient to consider);
otherwise, $\tR^0=\emptyset$ and the polynomial representation
is semisimple irreducible. Note that
$$\tR^1_+ \equal \tR^1\cap\tR_+=
\tR^1_+[-]\,\cup\, \{\al_i,i>0\} \hbox{\ since\ } k<0.$$
\smallskip

We mention that there are other ways to justify that $\v_{ss\,}$
is a $\HH^\flat$\~submodule (especially, for equal $k$);
however, always a reduced decomposition of the
elements $\pi_c$ can be found for $E_c$ appearing in $\v_{ss\,}$
such that dim$\,\v_{c'}$=1 occurs for $\pi_{c'}$
``before" $\pi_c$ (in this decomposition).
\medskip

{\bf The case of {\mathversion{bold}$A_n$}}. 
Let connect our $\v_{ss\,}$ in the case of $A_n$ 
with the construction from \cite{Ka},  
which extends that from \cite{FJMM}. In the latter paper, a
symmetric variant of $\v_{ss\,}$ was defined (for $GL_{n+1}$).
Namely, for an arbitrary negative rational
$k$ (with the denominator no greater than $n+1$), a
set of weights was found such that
the corresponding symmetric Macdonald polynomials exist
and linearly generate the space closed with respect to
the multiplication. The authors note that the symmetric
Macdonald polynomials actually exist for a bigger set of
weights. As a matter of fact, this remark is closely
connected with our approach; we prove that $\v_{ss\,}$ is a
DAHA\~submodule using that $\v_{c'}$ become one-dimensional
{\em before} $c'$ reach the set $\{c\}$ 
corresponding to $\{E_c\}$ linearly generated  $\v_{ss\,}$.

Paper \cite{Ka} contains a statement
equivalent to our one in the $A_n$\~case. We note that the
technique of intertwiners and the analysis of $\v_{ss\,}$
is significantly  simpler for $A_{n}$ than for other
root systems. In \cite{Ka}, $\v_{ss\,}$  was a part of the
conjectural decomposition of $\v$ in terms of irreducible
modules (the Kasatani
conjecture) justified in \cite{En} for $r\neq 2$
using the localization functor for the degenerate DAHA of
type $A$.

Let $t=q^k, k=-s/r$ and $(r,s)=1, n+1\ge r>1$ provided that
$q^at^b=1$ for $a,b\in \Q$ implies 
$a+kb=0.$ 

Using the notation from \cite{Bo}, $\al_i=\ep_i-\ep_{i+1},$
$$
2\rho=n\ep_1+(n-2)\ep_2+\ldots+(n-2(j-1))\ep_j+\ldots+
(-n)\ep_{n+1}.
$$
Given $1\le i\le n$, let us determine all permutations
$u\in W_1^i$ that can be lifted to elements $ub$ from
$\hw[-\rho_k]$. Such permutations can be described as
partitions of the segment $[1,\ldots,n]$ in terms
of $2m>2$ consecutive (connected) {\em segments\,}
$L_p=[p',p'']$ with $|L_p|=p''-p'+1$ elements:
\begin{align} \label{lsegments}
&L_1,L_2,L_3,\ldots,i,\ldots,L_{2m-2}, L_{2m-1},L_{2m}\ni n,\\
&\hbox{where\ }|L_p|>0, |L_p|=0\hbox{\ mod\ } r
\hbox{\ for\ } p=2,3,\ldots,2m-1,
\notag\\
&\hbox{\ and\ }
|L_1|+|L_3|+|L_5|+\ldots+|L_{2m-1}|=i,\ L_p=[p',p'']
\notag
\end{align}
Note that only $L_1,L_{2m}$ are allowed to be empty.
The corresponding permutation $u$ is
$$
\{L_1,L_3,L_5,\ldots,L_{2m-3},L_{2m-1},L_2,L_4,\ldots,
L_{2m-2},L_{2m}\}.
$$

This $u$ sends $\al_i=\ep_i-\ep_{i+1}$ to
$-\be=\ep_x-\ep_y$ for $x=(2m-1)'', y=2'$ and
$(\be,\rho)=((n-2(y-1))-(n-2(x-1)))/2=x-y$; here $x>y$.
The relation
$(\be+\al_i,\rho_k)\in \Z$ (see ($b,b'$)
from Proposition \ref{HWALBE}), 
which guaranties that
$-\be$ is a projection of the element from $\tR^1_+[-]$,
gives that the cardinalities of the
segments between $L_2$ and $L_{2m-1}$ are divisible by $r$
and  $x-y+1$ is divisible by $r$.

Since $m>1$, $x-y+1$ must be at least $2r$ that means that
$-\be$ is always decomposable in the sense of
Theorem \ref{FEIGIN}. Thus $\v_{ss\,}$ is always an $\HH\,$\~
submodule in the $A$\~case when $r>1$.
\smallskip

{\em A counterexample for $A_1$}.
The simplest example when $\v_{ss\,}$ is {\em not} a
$\HH\,$\~submodule
of $\v$ is for $\tR=\widetilde{A}_1$ and
$k=-m\in -\N$; the relation dim$\,\v_{c'}>1$ 
does not hold in this case.
Here $u=s_1$ can be lifted to
$\tu=[-\al_1,m]\in \tW^0\simeq S_2$
and $[-\be,\nu_\be j]=\tu(\al_1)=[-\al_1,2m]$. The latter
root is the only element of $\tR^1_+[-]$. The condition
$k+ j+(\be,\rho_k)=0$
from the definition of $\tR^1_+[-]$ in (\ref{critrootoneu})
reads as: $-m+2m-m(\al_i,\rho)=0$.

The space
$\v_{c}=\widetilde{V}_c$ becomes one-dimensional
for $\pi_c=(s_1s_0)^m$; here
$\la(\pi_c)=\{[-\al_1,2m],\ldots,[-\al_1,1]\}$ and the
last intertwiner in the chain is $\tau_+(T_1)-t^{1/2}$. 
The spaces $\v_{c'}$
are two-dimensional for previous 
$\pi_{c'}=\cdots s_0s_1s_0$, i.e., when the number $l$ of simple
reflection in this product satisfies $2m-1\ge l\ge m$;
here $\la(\pi_{c'})=\{[-\al_1,l],\ldots,[-\al_1,1]\}$.

The polynomial representation is irreducible but not
semisimple for such $k$. Note that there are no {\em
positive} affine roots $\tal$ in this case such that
$-k+(\tal,-\rho_k+d)=0$; the roots $-\al_1$ and $[\al_1,-2m]$
satisfying this condition are negative. This means that
$\tR^{-1}_+[-]$ is empty and 
$\tau_+(T_i)+t^{-1/2}$ will not appear in the chains of
intertwiners in $\v$. 
\sq
\medskip

\subsection{Eigenspaces for
\texorpdfstring{{\mathversion{bold}$\pi_c\tw$}}
{(pi-c)\t{w}}}
The following is another application of Proposition 
\ref{HWALBE}, which is important for {\em exact} calculating 
dim$\, \v_c$ and dim$\, \v[-c_\#]^\infty$.
For instance, it can be applied to reinstate Theorem
\ref{FEIGIN}.
 
We continue using the notation
$\tw=[\tw]_b [\tw]_u$. Let
$\la^1_R(\tw)$ be
the subset of {\em simple non-affine} roots in $\la(\tw)$,
that is $\la(\tw)\cap R^1$, where we set
 $R^1\equal\{\al_i,1\le i\le n\}$. 
When $q$ is not a root of unity,
\begin{align}\label{laonetw}
&\la^1_R(\tw)=\{\al_i,i>0\,\mid\, [\tw]_u(\al_i)>0, 
[\tw]_u(\al_i)\neq \al_{i\,'\,},\, i\,'\,>0\}.
\end{align}

The set  $\tR^{\pm 1}$ can be either
as in Theorem \ref{FEIGIN} under
$$
q^a\prod_\nu t_\nu^{b_\nu}=1 \Rightarrow
a+\sum_\nu\nu k_\nu b_\nu=0,
$$ 
or as in Proposition \ref{HWALBE} without this assumption.

We constantly use that $\tW^0$ and
$\hW^\flat[-\rho_k]$ preserve $\tR^{\pm 1}$.
For instance, $\tw(R^1\setminus \la_R^1(\tw))\subset \tR^1$
and $-\tw(\la_R^1(\tw))\subset \tR^{-1}$ in the Main Theorem,(i) 
below.

\begin{maintheorem}\label{HWALSIMPLE}
(i) Given $\tw\in \tW^0$, we impose condition ($a$)
from (i) in Proposition \ref{HWALBE}
for every $\al\in R_+$; it is \underline{necessary} to ensure
$\pi_c\tw\not\in\pi_B$ (recall that (a) results from 
(\ref{knegative})).
Setting $\la_R^1=\la_R^1(\tw)$,
\begin{align}\label{pictwsimple}
&\pi_c\tw\in \pi_B\,\Leftrightarrow\, 
\{-\tw(\la_R^1)\subset \la(\pi_c),\,
\tw(R^1\setminus \la_R^1)\cap\la(\pi_c)=\emptyset\}. 
\end{align}
Let us take $\tw$ such that $\pi_c\tw\in \b^0_o(\pi_c)$,
i.e., $\pi_c\tw$ must be  
obtained from a reduced decomposition
of $\pi_c$ by deleting some singular reflections. 
Then 
\begin{align}\label{pictwsimapv}
&\v_c=\sum_b \Q_{q,t}\,\widetilde{E}_b \hbox{\ for\ }
\pi_b=\pi_c\tw \in \b^0_o(\pi_c).
\end{align}

(ii) Let us consider a reflection $\tw=s_{\tal}$ in (i) for
$\tal=[-\al,\nu_\al j]\in \tR_+^0$, where $\al>0,j>0$.
Then ($a$) holds, $[\tw]_u=s_{\al}$ and  
$$\la_R^1(\tw)\ =\ \{\al_{i\,'\,}\,\mid\, (\al_{i\,'\,},\al)<0,\, 
{i\,'\,}>0\}.
$$ 
Let $\pi_c\tw$ be the result of deleting \underline{one} singular
simple reflection $s_p$ in a reduced decomposition of $\pi_c$.
Then $\pi_c\tw=\pi_b$, equivalently 
$\widetilde{E}_b\in \v_c$, if and only if
\begin{align}\label{talpicb}
s_{\tal}(\al_{i\,'\,})=\al_{i\,'\,}+
2\frac{(\al,\al_{i\,'\,})}{(\al,\al)}\tal\, \not\in\,
\la(\pi_c) \hbox{\ \,for\ all\  } (\al_{i\,'\,},\al)>0;
\end{align}
here the condition  
$\{-\tw(\la_R^1)\subset \la(\pi_c)\}$ holds
automatically. If 
$s_{\tal}(\al_{i\,'\,})\in\la(\pi_c)$ for  $i\,'\,>0$ satisfying
(\ref{talpicb}), then $\pi_c s_{\tal}\not\in \pi_B$ 
and this root can appear in the sequence $\la(\pi_c)$
only {after} (the root corresponding to) $s_p$;
For instance, (\ref{talpicb}) holds if $s_p$
is the last in the decomposition of $\hw$. 

(iii) If $(\al_{i\,'\,},\al)<0$ for $i\,'\,>0$, then 
$-s_{\tal}(\al_{i\,'\,})\in \tR^{-1}_+[-]$ appears in
$\la(\pi_c)$ before $s_p$. 
For arbitrary $\tal$ (not associated with $s_p$), the conditions 
$\la(\pi_c)\ni$ $-s_{\tal}(\al_{i\,'\,})\in \tR_+^{-1}[-]$ 
and (\ref{talpicb})
are necessary and sufficient to ensure that 
$\pi_c s_{\tal}\in \pi_B$.  More generally,
provided ($a$) for $\hw$ as in (i), 
\begin{align}\label{pictwsimapinfty}
&\v(-c_\#)^\infty=
\sum_b \Q_{q,t}\,\widetilde{E}_b \hbox{\ for\ }
\pi_b=\pi_c\hw,\, \hw\in \hW^\flat[-\rho_k],
\hbox{\ where}\\
\pi_c\hw\in &\pi_B\,\Leftrightarrow\,
\{-\hw(\la_R^1(\hw))\subset \la(\pi_c),\,
\hw(R^1 \setminus \la_R^1(\hw))\cap\la(\pi_c)=
\emptyset\}.\notag 
\end{align}
If $\hw$ does not satisfy ($a$) from Proposition \ref{HWALBE}
or one of the latter conditions in $\{\ ,\ \}$ does not hold 
for \underline{every} $\hw\in \hW^\flat[-\rho_k]$ then
dim$\,\v(-c_\#)^\infty=1$; vice versa, the condition
dim$=1$ implies that all three must hold.
\end{maintheorem}
{\em Proof}.  It is a version of Proposition \ref{HWALBE}
where we use that $\hv\not\in\pi_B$
implies that $\la(\hv)$ contains at least one
{\em simple} non-affine root. i.e., a root in the form
$\al_{i\,'\,}$ for $i\,'\,>0$. We use (\ref{talinla}) to check
that $\{-\tw(\la_R^1)\subset \la(\pi_c)\}$ holds
automatically in this case.

The justification of claims about
the position of $s_{\tal}(\al_{i\,'\,})$ from (ii) and
$-s_{\tal}(\al_{i\,'\,})$ from (iii) respectively
after and before $s_p$ is as follows.
Let $s_{\tal}(\al_{i\,'\,})$ 
$=\al_j+2\frac{(\al,\al_{i\,'\,})}{(\al,\al)}\tal\ $
assuming that $i\,'\,>0,\,(\al,\al_{i\,'\,})>0$. Then 
$s_{\tal}(\al_{i\,'\,})$ is a sum
of $\al_{i\,'\,}$ and $\tal\in\la(\pi_c)$ 
with positive coefficients. Since $\al_{i\,'\,}\not\in\la(\pi_c)$,
$\tal$ must appear before  $s_{\tal}(\al_{i\,'\,})$ for any
reduced decomposition of $\pi_c$. See Theorem \ref{INTRINLA}.
The claim from (iii) is analogous. \\
\sq
\medskip 

{\bf Spaces {\mathversion{bold}$\widetilde{V}_c$}.}
Given an arbitrary, {\em non-necessarily reduced},
decomposition
of $\hw$, the set
$\tilde{\b}^0(\hw)\ni \hw'$ 
was introduced in Theorem \ref{PHIBRUHAT} 
and the corresponding {\dfont standard} decompositions of $\hw'$
obtained by deleting simple singular reflections. This
set has a natural partial ordering. Recall that
\begin{align}\label{macfilt}
&\widetilde{E}_{c}^\dag\equal
\widetilde{\Psi}_{\hw}(1)\in\oplus_{a\succeq c'}\Q_{q,t}X_a
\for c'=\hw'(\!(c)\!),
\end{align}
where $\hw'\in \tilde{\b}^0(\pi_c)$.
Without $\dag$, the $\widetilde{E}$\~polynomials are
$$
\widetilde{E}_c=\hbox{Const}\,P_r\widetilde{\Psi}_{i_l}\cdots
\widetilde{\Psi}_{i_1}(1) \for \pi_c=\pi_rs_{i_l}\cdots s_{i_1},
$$
where the reduced decomposition is used and the constant
is chosen to ensure the normalization
(\ref{ecxfilt}), i.e., the coefficient
of $X_c$ has to be $1$. 
If the decomposition of $\pi_c$ is not assumed
to be reduced (the  $\dag$\~case), then
$\widetilde{E}^\dag_c$ are defined only up to proportionality
and they can be zero. The spaces $\widetilde{V}_c$ and $\v_c$
are constructed only for the reduced decompositions of $\pi_c$ in
the following proposition and later.
\smallskip

Recall that $\widetilde{V}_c$ is the linear span of
$\widetilde{E}^\dag_{b}$ for $b=\hw'\llb 0 \rrb$
defined for all {\em standard
decompositions\,} of $\hw'\in \tilde{\b}^0(\pi_c)$,
{\em possibly non-reduced}, obtained from
a given reduced decomposition of $\pi_c$.
The space $\v_c$ is linearly generated by
$\widehat{E}_b$ for the same set $\{b\}$, however they are
defined only for {\em reduced}
decompositions of $\pi_b$; they are linearly independent
for any particular choice of these decompositions.
\smallskip

Practically, the difference between $\v_c$ and $\widetilde{V}_c$
is as follows. The squares $\ldots(\Psi_i)^2\ldots$ 
(they are $Y$\~rationals)
in the formulas for $\widetilde{E}^\dag_c$
and $\widetilde{V}_c$
are replaced by their $Y$\~expressions; such squares
are simply deleted when constructing $\v_c$.
\smallskip

Note that the linear span of
$\widetilde{\Psi}_{\hw}\widetilde{V}_{c}$
for {\em all} decompositions, possibly non-reduced,
of $\hw\in \hW^\flat$ is an $\HH^\flat$\~submodule
of $\v$.

\smallskip
As we know, the space $\v_c$ does not depend
on the choice of the reduced decomposition of
$\pi_c$. It follows from the limiting procedure
or can be readily checked using
Theorem \ref{PHIBRUHAT}\, (apply $\tau_-\si$ to $\Phi_{\hw}$
there 
and use $\widetilde{\Psi}_i=\tau_+(T_i)$ for singular
$s_i$ in $\widetilde{\Psi}_{\hw}$).

Individual $\widehat{E}_c$ and the space
$\widetilde{V}_c$ defined for
a given reduced decomposition of $\pi_c$ may change
if the homogeneous Coxeter transformations of type
(\ref{atwosing}) from Theorem \ref{PHIBRUHAT},($a$)
are applied. However formula (\ref{phitphi}) makes it
possible to control the change in this case.
\smallskip

\begin{proposition}\label{VCHAT}
The spaces $\v_c$ and $\widetilde{V}_c$
are $\Q_{q,t}[Y_b]$\~modules and
are also invariant with respect to the action of
the operators $\tau_+(T_i)$ for $s_i$ such that
$\pi_c^{-1}(\al_i)\in \tR^0$ and $l(s_i\pi_c)<l(\pi_c)$.
If $q,t$ are not roots of unity
and the whole centralizer of
$q^{-c_\#}$ in $\hW^\flat$ is generated by
simple reflections $s_i$, then
$$
\widetilde{V}_c=\v_c=\v(-c_\#)^\infty.
$$
\end{proposition}
{\em Proof}. See Theorem \ref{VINVHAT},(iii) and
Theorem \ref{STRIKOUT},($e$). The last claim follows from
the fact that $\v(-c_\#)^\infty$ is an induced module
over the {\em affine} Hecke algebra defined for $\tR^0$,
that is irreducible and $Y$\~cyclic in this case;
cf. Proposition \ref{YSTRUCTURE} below.
\sq
\medskip

\subsection{The 
\texorpdfstring{\mathversion{bold}$Y$}{\em Y}--action}
The following theorem makes the $Y$\~structure
of $\widetilde{V}_c$ or $\v_c$ as explicit
as possible (the latter is somewhat simpler to calculate
than the former). Note that
the formulas below do not give
a complete description of the $Y$\~action
because the polynomials $\widehat{E}^\dag_{c'}$ that
appear in process of calculations
can be zero or linearly dependent in $\widetilde{V}_c\subset\v_c$.

\begin{maintheorem} \label{THEORCHAIN}
Given $b\in B$, an element $c\in B$ and also its decomposition
$\pi_c=\pi_rs_{i_l}\cdots s_{i_1}$, non-necessarily reduced,
let $\{i_p\}$ be the set of 
singular indices in this decomposition
for  $p\in\{p_g>,\ldots,>p_1\}$.  Then
\begin{align}\label{ybofea}
&q^{(b,c_\#)}\,Y_b(\widetilde{E}^\dag_{c})\,=\,
\widetilde{E}^\dag_{c}
+\sum_{c'} \prod_{p'=p_j'}^{1\le j\le g'}
(t_{i_{p'}}^{1/2}-t_{i_{p'}}^{-1/2})\,
C_{c'}\widetilde{E}^\dag_{c'},
\end{align}
where $c'=\hw'(\!(0)\!)$,
$\hw'\in \tilde{\b}^0(\pi_c)$ and
$\widetilde{E}^\dag_{c'}$
are defined for the decompositions, possibly non-reduced, of
$\hw'$ obtained from the initial decomposition
of $\pi_c$ by deleting singular $s_{i_{p'}}$ for the indices
$p'$ forming a subsequence
$$\{p'_{g'}>,\ldots,> p'_{1}\}\,\subset\,\{p_g,\ldots,p_1\}.
$$
For instance, if there is only one such $p'$, i.e., $g'=1$, then
\begin{align} \label{eaviaec}
&C_{c'}=((\tbe^{p'})^\vee,b+d)
\for \tbe^{p'}\equal(\pi_r s_{i_l}\cdots 
s_{i_{p'+1}})(\al_{i_{p'}}),
\\
((\tbe^{p'})^\vee,b+d)&=-((\tal^{p'})^\vee,\pi_c^{-1}\llb b\rrb+d)
\for
\tal^{p'}=s_{i_1}\cdots s_{i_{p'-1}}(\al_{i_{p'}}).\notag
\end{align}
Generally, when the number of the indices $p'$ is $h=g'\ge 1$,
we set
$\widehat{b}=\pi_c^{-1}(\!(b)\!)$ and the formula for 
$C_{c'}$ reads as follows:
\begin{align} \label{eaviaecc}
C_{c'}=
(Y_{\hb}^{-1}&(D^Y_{\tal^{p\,'_1}}\cdots D^Y_{\tal^{p\,'_h}})
(Y_{\widehat{b}}))
\Downarrow 1^0=(D^Y_{\tal^{p\,'_1}}\cdots D^Y_{\tal^{p\,'_h}})
(Y_{\widehat{b}})
\Downarrow 1,\notag\\
\hbox{setting\ \,}
Y_{\tal}&\Downarrow 1^0=1\hbox{\ for\ }\tal\in \tR^0,\ \,
Y_{a}\Downarrow 1=1\hbox{\ for\ all\ } a\in B,\notag\\
&\hbox{where\ \ } D^Y_{\tal}\,\equal\,
(Y_{\tal}-1)^{-1}(s_{\tal}^Y-1)\for
\tal\in \tR,
\end{align}
defined in terms of the action\,
$s_{\hw}^Y(Y_a)\equal Y_{\hw(a)}$\, for $\hw\in \hW$.
\end{maintheorem}
{\em Proof}.
We use the $\tau_-\si$\~image of the relations (\ref{tixi}):
\begin{align}
&Y_b \,\tau_+(T_i)=\tau_+(T_i)\,Y_{s_i(b)} +
(t_i^{1/2}-t_i^{-1/2})\frac{s^Y_i(Y_b)-Y_b}{
Y_{\al_i}^{-1}-1}.\label{tixiy}
\end{align}
Recall that
$\tau_+(T_i)=\tau_+\tau_-^{-1}(T_i)=\tau_-\si(T_i)$
and $\tau_-\si(X_b)=Y_b^{-1}$.

Let us introduce the operators
\begin{align}\label{tytal}
&\t^Y_{\tal}\equal t_{\al}^{1/2}s_{\tal}^Y-
\frac{t_{\al}^{1/2}-t_{\al}^{-1/2}}{Y_{\tal}^{-1}-1}(s^Y_{\tal}-1),
\end{align}
satisfying the relations
\begin{align}\label{tixiyy}
&Y_b \,\t^Y_{\tal}=\t^Y_{\tal}\,Y_{s_i(b)} +
\frac{t_{\al}^{1/2}-t_{\al}^{-1/2}}{Y_{\tal}^{-1}-1}(s^Y_{\tal}-1)
(Y_b).
\end{align}
For instance, $\t^Y_i=\t^Y_{\al_i}$  
satisfy (\ref{tixiy}) for $\tau_+(T_i)$.

These relations readily give (\ref{eaviaecc}).
One needs to move $Y_b$
in $Y_b\widetilde{\Psi}_{\hw}$ through 
$\widetilde{\Psi}_{\hw}$ and calculate
the coefficients of $\widetilde{\Psi}_{\hw'}^\dag\,Y_a$ for
$\hw'\in\tilde{\b}^0(\hw)$ in the resulting decomposition.
Here using ${}^\dag$ in $\widetilde{\Psi}_{\hw'}^\dag$ indicates
that when some simple intertwiners $\widetilde{\Psi}_i$ disappear
from the product due to (\ref{tixiyy}), we leave the
product as it is (without reducing the corresponding 
decompositions).

These coefficients {\em coincide} with those obtained when $Y_b$ is
moved through $\hw^Y\r^Y_{\tal^{p_g}}\cdots \r^Y_{\tal^{p_1}}$
(instead of $\widetilde{\Psi}_{\hw}$), where
$\r^Y_{\tal}\equal s^Y_{\tal}\t^Y_{\tal}$ satisfy
\begin{align}
&Y_b \,\r^Y_{\tal}=\r^Y_{\tal}\,Y_{b} +
(t_{\al}^{1/2}-t_{\al}^{-1/2})D^Y_{\tal}(Y_b)\label{tixiyr}.
\end{align}
The product
$Y_b\hw^Y\r^Y_{\tal^{p_g}}\cdots \r^Y_{\tal^{p_1}}$
will become a linear combination of the terms 
$$
\hw^Y\r^Y_{\tal^{p_{h}'}}\cdots \r^Y_{\tal^{p_1'}}Y_a
\hbox{\ \ for\ \ }
\hw'\in \tilde{\b}^0(\hw),
$$
where we do not perform reductions
if some of these terms are linearly dependent.

Then we apply the evaluation $\Downarrow 1^0$; here
one can also take
$\Downarrow 1$, sending $Y_a\mapsto 1$ for all $a$,
because $\{\tal^{p'}\}$ are from $\tR^0$.

In the case of one $p'$, the coefficient
$C_{c'}$ from (\ref{eaviaec})
is $D_{\tal^{p'}}(Y_{\hw^{-1}\llb b\rrb})$ upon the evaluation
$\Downarrow 1$ for $\hw=\pi_c$; it equals
\begin{align*}
&-((\tal^{p'})^\vee,\hw^{-1}(\!(b)\!)+d)=
-(\hw s_{i_1}\cdots s_{i_{p'-1}}(\al_{i_{p'}})^\vee,b+d)\\
&=-(\pi_r s_{i_l}\cdots s_{i_{p'}}(\al_{i_{p'}})^\vee,b+d)=
((\tbe^{p'})^\vee,b+d).
\end{align*}
\sq
\smallskip

Note that $C_{c'}$ are integers. Their
calculation is a combinatorial problem that can be
formulated in terms of the algebra generated by
the {\dfont affine Demazure operators\,}
$D^Y_{\tal}$ defined for $\tal\in \tR^0$. They
satisfy the $r$\~matrix relations and the
quadratic ones: $D^Y_{\tal}(D^Y_{\tal}-1)=0$.
The affine Hecke algebra for the root system
$\tR^0$ is sufficient for calculating these coefficients.
The degenerate affine Hecke algebra can be used here,
which is useful for analyzing the $Y$\~cyclicity
of $\widetilde{V}_c$. Moreover, almost always a reduction
to the {\em non-affine} Hecke algebra is sufficient.  
\smallskip

The elements $\widetilde{E}^\dag_{c'}$, when considered
as vectors in $\widetilde{V}_c$, can vanish and some can become
linearly dependent. The next stage of the calculation is
when we express the $\widetilde{E}^\dag$\~polynomials
in terms of {\em reduced} $\widetilde{E}$\~polynomials.
Replacing $T_i^2$ using the quadratic relations may be necessary
and moving the $Y$\~functions $\Psi_i^2$ to the right through
$\Psi$ and $T$. Eventually, some $\widetilde{E}^\dag_{c'}\neq 0$
may vanish and some may become linearly dependent.
It is convenient to perform this calculation inside $\v_c$.

The following ``rationalization" of (\ref{ybofea})
is needed for the complete description of $\v_c$
as $Y$\~modules.

Given a rational function $P_Y$ of $Y_b$ such that $P_Y(1)$ is
well defined, setting $\widehat{P}_Y=(\hw^Y)^{-1}(P_Y)$:
\begin{align}\label{ybofeaa}
&P_Y(\widetilde{E}^\dag_{c})\,=\,
\widetilde{E}^\dag_{c}\,P_Y(1)\notag\\
+&\sum_{c'} \widetilde{E}^\dag_{c'}\,
\big(\prod_{p'=p_j'}^{1\le j\le h=g'}
(t_{i_{p'}}^{1/2}-t_{i_{p'}}^{-1/2})\big)\,
(D^Y_{\tal^{p'_1}}\cdots D^Y_{\tal^{p_{h}'}}(\widehat{P}_Y)))(1).
\end{align}

\smallskip

We can use here (\ref{lapimin}), describing the set of
$\tbe^{p}\equal(\pi_r s_{i_l}\cdots s_{i_{p+1}})(\al_{i_{p}})$
that may appear in (\ref{eaviaec}) for a given reduced
decomposition $\pi_c=\pi_r s_{i_l}\cdots s_{i_{1}}$
and singular $p$. 
Setting $\tbe=[\be,\nu_\be j]$, it is as follows:
\begin{align}\label{crossout}
&\tbe\in \la(\pi_c^{-1})\hbox{\ and\  }
q_\be^{j+(\be^\vee, c-u_c^{-1}(\rho_k))}=1, \where\notag\\
&\la(\pi_c^{-1})\ =\
\{ \tal=[\al,\nu_\al j]\in \tR_+,\  -(c,\al^\vee)>j\ge 0 \}.
\end{align}
There are two immediate applications of this description:

(1) the set $\{\tbe^p\}$ does not depend on the choice of the
reduced decomposition of $\pi_c$;

(2) the non-affine components of such $\tbe$
cannot coincide for distinct $\tbe$ for generic $q$.
\smallskip

\begin{corollary}\label{HATEECOM}
Given $\pi_c$, let us assume that the singular
roots in $\la(\pi_c)$ are pairwise orthogonal; for
instance, it holds for any $\hw$ if $\tR^0$ is a direct sum
of one dimensional root systems, respectively, $\tW^0$ is a
commutative group. Then the module $\widetilde{V}_c$
is $Y$\~cyclic generated by $\widetilde{E}_c$.
\end{corollary}
{\it Proof.} The corresponding Demazure operators for
{\em singular} roots in (\ref{ybofeaa}) are pairwise commutative
and it is possible
to find a polynomial $P^{p_j'}$ for each $p_j'$ such that
$$
D^Y_{\tal^{p_{i}'}}(\widehat{P}^{p_j'})=\de_{ij}
\hbox{\ \, for \ all\ \,} 1\le i,j\le g.
$$
One can use directly formula (\ref{eaviaec}) here. Hence, the
$Y$\~span of $\widetilde{E}_c$ contains $\widetilde{E}_{c'}^\dag$
for $g'=1$ and $c'$ as in (\ref{eaviaec}). Then we can
proceed by induction.
\sq
\smallskip

It is not known how far the $Y$\~modules $\widetilde{V}_c$
are from cyclic in general. The strongest 
{\em criterion} we give in the paper 
is Proposition \ref{VINDUCEDY} below. 
Combined with Theorem \ref{RANKTWO}, it shows that
in quite a few cases $\widetilde{V}_c$ are cyclic. 
It is not difficult to construct examples
when $\v_c$ are not $Y$\~cyclic. 
\smallskip

\begin{corollary}\label{HATEE}
(i) Given $c\in B$ and a reduced decomposition of $\pi_c,$
the polynomial $\widetilde{E}^\dag_{c'}\neq 0$
defined for the corresponding standard decomposition 
(maybe non-reduced) of \,$\hw'\in \tilde{\b}^0(\pi_c)$ 
is a $Y$-eigenvector
for $c'=\hw'\llb 0\rrb$\, if\, $\widetilde{E}^\dag_{c''}=0$
for all $\pi_{c''}\,<_0\, \pi_{c'}$.

Any $Y$\~quotient of $\v_c$ contains
at least one nonzero $Y$\~eigenvector that is the image of
$\widetilde{E}_{c'}\in \v_c$ for suitable $c'$;
one can take here $c'$ such that  $\pi_{c'}\,\le_0\, \pi_c$ is 
minimal with respect to $\,\le_0\,$ considered only 
among $\pi_{c'}$ with nonzero images of
$\widetilde{E}_{c'}$.

(ii) The kernel of
$\Psi_i$ in $\v$ or its any finite dimensional $Y$\~submodule
consists only of $Y_{\al_i}$\~eigenvectors
provided that the action of $\Psi_i$ is well defined at
such vectors and $t_i\neq \pm 1$.
Given $c\in B$ and  $0\le i\le n$, the action of $\Psi_i$ is
well defined in the space $\v_c$ unless  $s_i$ is singular
in the product $s_i\pi_c$; see (\ref{interinfty}).

Let $s_i$ be singular and $t_i\neq 1$. Then
$(\al_i^\vee,c+d)\neq 0$ and $s_i\pi_c\in \pi_B$.
Moreover, if $0\neq E_1\in \widetilde{V}_c$, then at
least one of the vectors $E_1$ and $E_2=\tau_+(T_i)E_1$
is not a $Y_{\al_i}$\~eigenvector.
\end{corollary}
\comment{
(iii) Under the same assumption for $c,i,t_i$,
let us suppose that the ideal $\i_c=$Ker$(P(Y)\mapsto
P(Y)\widetilde{E}_c)$ is $s_i^Y$\~invariant
for $P(Y)\in \Q_{q,t}[Y_a,a\in B]$.
Then the $Y$\~cyclic
module $\lan \widetilde{E}_c\ran_Y\equal$
$\Q_{q,t}[Y_a]\,\widetilde{E}_c$ is smaller than
$\lan \widetilde{E}_b\ran_Y$ for $b=s_i\llb c\rrb$;
and the $Y$\~homomorphism sending
$\widetilde{E}_b\mapsto$
$\widetilde{E}_c$ is surjective.
}
{\em Proof.}
Claim (i) follows from the fact that the $Y$\~submodule
generated by $\widetilde{E}_c$
is given in terms of $\widetilde{E}_{c'}$ for
$\pi_{c'}\,\le_0\,\pi_c$.

The first claim in (ii) results from the
formula
\begin{align*}
&(\Psi_i)^2=
\frac{(t_i^{1/2}Y_{\al_i}^{-1}-t_i^{-1/2})
(t_i^{1/2}Y_{\al_i}-t_i^{-1/2})}
{(Y_{\al_i}^{-1}-1)(Y_{\al_i}-1)}.
\end{align*}
Relation (\ref{tixiy}) can be naturally
extended to the vectors $[b,l]$
from $[B,\Z]$. Applying it to $\al_i,$
\begin{align}\label{yale}
&(Y_{\al_i}-1)E_2=
2(t_{i}^{1/2}-t_{i}^{-1/2})E_1
\end{align}
if $E$ is a $Y_{\al_i}$\~eigenvector; the
corresponding eigenvalue is $1$.
Respectively, if $\widetilde{E}$ is a
$Y$-eigenvector, then
\begin{align*}
&
(Y_{\al_i}^{-1}-1)E_1=2(t_{i}^{-1/2}-t_{i}^{1/2})E_2.
\end{align*}
It gives the rest of (ii).
\sq
\comment{
Claim (iii) is actually a generalization of (ii).
Let $E_1= \widetilde{E}_c$,
$E_2=\tau_+(\widetilde{E}_c)$; recall that
$\widetilde{E}_b=CE_2$
for $C\neq 0$. The corresponding variant of
(\ref{yale}) reads as follows:
\begin{align}\label{yalee}
&(Y_{\al_i}-1)E_2=
\Psi_i^\diamond (E_1)
+(t_{i}^{1/2}-t_{i}^{-1/2})E_1 \for\\
&\Psi_i^\diamond=
\Psi_i\cdot (Y_{\al_i}^{-1}-1)\
\hbox{\ \ from\ (\ref{Psiprime}).}\notag
\end{align}
The module
$\lan (Y_{\al_i}-1)E_2\ran_Y$ is
isomorphic to $\lan E_1\ran_Y$. Here we need to
check that
$$
P(Y)E_1^\diamond=0\Leftrightarrow
P(Y)E_1=0 \for E^\diamond=\Psi^\diamond(E),\ P(Y)\in\Q_{q,t}[Y_b].
$$
Let $P(Y)(E_1+E_1^\diamond)=0$, where
$P(Y)(E_1+E_1^\diamond)=P(Y)E_1+\Psi_i^\diamond(s_i(P(Y))E_1)$.
Setting $E'=P(Y)E_1\neq 0$, let $M>0$ be the smallest power
such that $(Y_{\al_i}-1)^M\,E'=0$, equivalently,
$(Y_{\al_i}^{-1}-1)^M\,E'=0$. Then
$\Psi_i^\diamond((Y_{\al_i}-1)^{M-1}\,E')=0$
because $\Psi_i^\diamond$ annihilates all $Y_{\al_i}$\~eigenvectors,
including $E''=(Y_{\al_i}-1)^{M-1}\,E'$.
Thus,
$$
(Y_{\al_i}-1)^{M-1}\,E'=E''+
\Psi_i^\diamond((Y_{\al_i}^{-1}-1)^{M-1})E')=E''.
$$
$(Y_{\al_i}^{-1}-1)^{M-1}\,E'=0$.
}

\begin{corollary}\label{SEMIY}
Let $t_\nu\neq 1$ for all $\nu$.
The condition $\tR^0_+ \subset R_+$, i.e., the absence
of singular $\al_i$ in all $\pi_c$,
is  necessary and sufficient
for the polynomial representation $\v$
to be $Y$\~semisimple. In this case, any reduced
chain originated at $E_0=1$ consists of
one-dimensional nonzero spaces $\widetilde{V}_c$,
even if some of the simple intertwiners $\Psi_i$ become
non-invertible in this chain; all Macdonald polynomials
$E_c$ are well defined
in such $\v$ (and coincide with  $\widetilde{E}_c$).
\end{corollary}
{\em Proof}. Use Corollary \ref{HATEE},(ii).
\medskip

\subsection{Induced representations}
The $Y$\~modules\, $\widetilde{V}_c$\, are closely
related to their counterparts for
the {\em $Y$\~induced representations
of $\HH^\flat$}. 

Let us fix a weight $\xi$. In this section $\tR^0$ is an
arbitrary {\em root subsystem} from
Section \ref{sec:RightBruhat} satisfying:
\begin{align}\label{xisingular}
&q^{(\tal,\xi+d)}=1\ \Rightarrow\
\tal\in \tR^0.
\end{align}
The most natural choice (cf. (\ref{troeq})) is
\begin{align}\label{troeqxi}
\tR^0\equal&\{\tal=[\al,\nu_\al j]\in \tR\, \mid \,
q^{\nu_\al j+(\al, \xi)}=1\}.
\end{align}

let $\i_\xi$ be the  
$\HH^\flat$\~module induced from
the corresponding representation of its subalgebra 
$\Q_{q,t}[Y_a,a\in B]$.
By definition,
$\i_\xi$ is  generated by the element $v$ with the defining
relations $Y_a(v)=q^{(a,\xi)}v$ and is naturally isomorphic to
$\h_X^\flat\,=\,\lan T_i, 1\le i\le n,\, X_b, b\in B\ran$ due
to the PBW\~theorem for the pair of subalgebras $\h_X^\flat$ and
$\Q_{q,t}[Y_b,b\in B]$.
It will be more convenient to use the identification of $\i_\xi$
with $\tau_+(\h_Y^\flat)$:
\begin{align*}
&\tau_+(\h_Y^\flat)\,=\, \tau_+\tau_-^{-1}(\h_Y^\flat)\,=
\,\tau_-\si(\h_Y^\flat) \\
&=\lan\tau_+(T_i) \,, 0\le i\le n,\
\tau_+(\pi_r) \for b_r\in \Pi^\flat\ran, \hbox{\ where\ }  \\
&\h_Y^\flat=\lan T_i,Y_b,\,b\in B \ran=\lan T_i,\, 0\le i\le n,\
\pi_r \for b_r\in \Pi^\flat\ran.
\end{align*}
Here the automorphisms $\tau_{\pm},\si$ are used;
see (\ref{taux}),(\ref{tauminax}), and
(\ref{tauintery}). The isomorphism $\i_\xi=\tau_+(\h_Y^\flat)$
(as vector spaces and as $\tau_+(\h_Y^\flat)$\~modules) is based
on the PBW\~theorem
for the pair of subalgebras $\tau_-\si(\h_Y^\flat)$ and
$\Q_{q,t}[Y_b,b\in B]=\tau_-\si(\Q_{q,t}[X_b,b\in B])$, which
is the $\tau_-\si$\~image of the standard PBW\~theorem for
$\h_Y^\flat$ and $\Q_{q,t}[X_b,b\in B]$.

We introduce the {\em $\widetilde{e}$\~elements} in $\i_\xi$, 
counterparts of $\widetilde{E}$\~polynomials, 
 as follows:
\begin{align}
&\widetilde{e}_{\hw}\equal\widetilde{\Psi}_{\hw}(v)
\hbox{\ \, for\ }
\hw\in \hW^\flat,\ \widetilde{\Psi}_{\hw}\, \hbox{\ from\ Theorem\ }
\ref{PHIBRUHAT}.
\end{align}
Recall that $\widetilde{\Psi}_{\hw}$ is constructed by
replacing $\Psi_{i_p}$ with $T_{i_p}$ for singular indices \,$p$\,
in the product for $\Psi_{\hw}$ corresponding to a
reduced decomposition of $\hw=\pi_rs_{i_l}\cdots s_{i_1}$.

Similarly, one defines the spaces $\v^\xi_{\hw}$,
$\widetilde{V}^\xi_{\hw}$. The elements $\widetilde{e}_{\hw}$
may depend on a choice of the reduced decompositions
of $\hw\in\hW^\flat$ but they are always nonzero. It follows from
the calculation of the leading coefficient of $\widetilde{e}_{\hw}$,
that is the coefficient of $\tau_+(T_{\hw})$.

The space $\v^\xi_{\hw}$ is linearly generated by
$\widetilde{e}_{\hw'}$
for $\hw'\in \b^0(\hw)$;
$\widetilde{V}^\xi_{\hw}\subset \v_{\hw}.$ It does not depend
on the choice of the reduced decomposition and can be obtained
by the same kind of limiting procedure as in the polynomial case.
However, now 
the complete right Bruhat ordering must be used instead of its
restriction to $\pi_B$.

If $\xi=-\rho_k$, then the images of $\widetilde{e}_{\hw}$ under
the $\HH^\flat$\~homomorphism $\i_\xi\to \v$ sending $v\mapsto 1$
become zero for $\hw\not\in \pi_B$, otherwise they are
nonzero and proportional to $\widetilde{E}_{b}$ for $\hw=\pi_b$.
Respectively,  $\v^\xi_{\pi_c}$ and $\widetilde{V}^\xi_{\pi_c}$
map {\em onto} $\v_{c}$ and $\widetilde{V}_{c}$.
\medskip

{\bf The case of generic weights}.
We will begin with generic $\xi$
subject to (\ref{xisingular}). Practically, it means
that the intertwines can be singular, but all are
invertible; cf. Proposition 8.11 from \cite{L}. 
Then the normalized intertwiners
$F=\tau_-\si(G)$ can be used instead of the $\Psi$.

Note that it is exactly the case where there is a complete
parallelism with Kauffman's axioms of {\em virtual links};
the {\em normalized} intertwiners and $\{T\}$ provide the
key example for Kauffman's axioms.

Thus one may switch from $\widetilde{e}_{\hw}$
to $\widetilde{e}'_{\hw}\equal F_{\hw}(v)$. Such
$\widetilde{e}'_{\hw}$
do not depend on the choice of the reduced decompositions
and are proportional to $\widetilde{e}_{\hw}$; this is due
to using the normalized intertwiners.

\begin{proposition}\label{VINDUCED}
Let $\xi$ be generic subject to (\ref{xisingular}). Then

(i) $\v^\xi_{\hw}=\widetilde{V}^\xi_{\hw}$ and therefore
the elements $\widetilde{e}_{\hw'}$ for $\hw'\in \b^0(\hw)$
form a basis in the space $\widetilde{V}^\xi_{\hw}$;

(ii) the elements $\{\widetilde{e}_{\hw'}^\dag\}$, defined for
all standard decompositions of $\hw'\in $
$\widetilde{\b}^0(\hw)$, can be linearly expressed in terms of
$\{\widetilde{e}_{\hw'}\}$;

(iii) the expressions from (ii) combined with (\ref{ybofeaa})
and the quadratic relations for $\{T\}$
give a complete description of the $Y$\~action
in $\widetilde{V}^\xi_{\hw}$. \\
\sq
\end{proposition}

The $Y$\~structure of $\widetilde{V}^\xi_{\hw}$ can be
calculated for generic $q,t$ in terms of the affine Hecke
algebras associated with the root system $\tR^0$, provided
that $\xi$ is generic from (\ref{xisingular}).

The construction requires the set of positive roots
$\tR^0_+\subset \tR^0$, the corresponding simple roots
$\{\al^0_i\}$ in $\tR^0_+$,
and the standard Bruhat ordering for $\tR^0$. Recall
that the notation for
the {\em standard Bruhat sets} for $\tR^0$ is $\b(\tu;\tR^0_+)$ where
$\tu\in \tW^0=\lan s_i^0\ran$.
Let $\widetilde{H}^0$ be the Hecke algebra defined for the system
$\tR^0$; it is generated by $\{T^0_i\}$ satisfying
the homogeneous Coxeter relations and the quadratic ones
with $t_{\al^0}=t_{\al}$
as $\al^0=[\al,\nu_{\al}j]\in \tR$. The elements
$T^0_{\tu}\in \h^0$ for $\tu\in \tW^0$ are defined naturally and
do not depend on the choice of the reduced decompositions of $\tu$:
$$
\widetilde{H}^0\, = \,
\oplus_{\tu}\Q_{q,t}T^0_{\tu} \for \tu\in \tW^0.
$$

Let $B^0\equal\{b\in B,(\xi,b)=0]$,
$\widetilde{\h}_Y^0$ be the algebra generated by
$\Q_{q,t}[Y_b,\, b\in B^0]$
and the Hecke algebra $\widetilde{H}^0$.
We impose 
\begin{align}\label{tixiyy0}
&Y_b \,T^0_{i}=T^0_{i}\,Y_{s_{\tal}(b)} +
(t_{\al}^{1/2}-t_{\al}^{-1/2})\frac
{Y_{s_{\tal}(b)}-Y_b}{Y_{\tal}^{-1}-1} \hbox{\ \ for\ }
\tal=\al^0_i;
\end{align}
to put it simply, $T_i^0$ are $\t_{\tal_i^0}$
from (\ref{tixiyy}).

The induced module
$\widetilde{\i}^0$ is defined by setting $Y_b(1)=1$ and
can be naturally identified
with $\widetilde{H}^0$ (as $\widetilde{H}^0$\~modules).

\begin{proposition}\label{YSTRUCTURE}
Provided that $q,t$ are not roots of unity and $\xi$
is \underline{generic}, the $Y$\~module
$\widetilde{V}^\xi_{\hw}$ is cyclic
for $\hw\in \hW^\flat$  if and only if the following holds
in $\widetilde{\i}^0$ for $\tu=\hw|_0\in \tW^0$ and
$\Q_{q,t}^0[Y]\equal\Q_{q,t}[Y_b,\,b\in B^0]$:
\begin{align}\label{cyclro}
&\Q_{q,t}^0[Y]\,(T^0_{\tu})=
\oplus_{\tu'\,}\Q_{q,t}\,T^0_{\tu'\,} \for 
\tu'\in \b(\tu;\tR_+^0),
\end{align}
where $\b$ is the standard Bruhat ordering in $\tW^0$.
For instance, relations (\ref{cyclro}) are satisfied
if $\tu$ is the element $\tu_0$ of maximal length in $\tW^0$.
\end{proposition}
{\em Proof.} Using
the intertwiners $F$ and the polynomials $\widetilde{e}\,'_{\hw}$,
the calculation of Theorem \ref{THEORCHAIN} can be reduced
to the case of $\widetilde{\h}_Y^0$ and $\widetilde{\i}^0$.

Let us examine the case of $\tu=\tu_0$. The reduction to
the degenerate affine Hecke algebra $\widetilde{\h}_y^0$
from \cite{Ch4}, Lemma 2.12 can be used (the analysis of
irreducibility of $\widetilde{\i}^0$). Laurent polynomials
in terms of $\{Y_b\}$ are replaced by the usual polynomials
in terms of $\{y_b\}$. The elements $T^0_{\tu}$ become $\tu$;
the Demazure operators in formulas (\ref{tixiyy0}),
(\ref{eaviaecc}) are replaced by the BGG operators
$d_{\tal}^y=y_{\tal}^{-1}(s_{\tal}-1)$.

The degeneration of $\widetilde{\i}^0$ is isomorphic
to the quotient $\widetilde{\i}_y^0$  of the standard
representation of $\widetilde{\h}_y^0$ in
the space of polynomials in terms of $y_b$ for $\,b\in B^0$
by the ideal generated by the $\tW^0$\~invariant polynomials
with zero constant term. It contains a unique $y$\~eigenvector
$d$, the {\em discriminant}, and $\tu_0(d)$ is a linear
generator of the one-dimensional space
$\widetilde{\i}_y^0/(y_b, b\in B^0)$; therefore it is
a $y$\~cyclic vector in $\widetilde{\i}_y^0$
(Nakayama's lemma). It gives the degenerate version of
(\ref{cyclro}) and results in the required claim.
\sq
\medskip

\rmk
(i) A counterexample to (\ref{cyclro}) is $\tR^0=A_3$,
$\tu=(4231)$. There are four elements in $\b(\tu)$ of length
$l(\tu)-1$, i.e., their number is greater than the rank.
This leads to a contradiction (consider the degeneration
above). One can use Proposition
\ref{BRUHPSEUDO} to construct more general counterexamples.

(ii) Note that the construction we discuss is connected
with the theory of Schubert polynomials upon the reduction
we performed when proving the proposition. Using the 
BGG\~operators for Schubert polynomials is similar 
to what we did. Generally, the combinatorics of {\em affine}
Schubert manifolds has a lot in common with that of 
nonsymmetric Macdonald polynomials.
\sq
\medskip

{\bf Arbitrary weights}.
There is a natural extension of
Proposition \ref{YSTRUCTURE} to the case of
arbitrary, {\em not only generic}, $\xi$.
We impose (\ref{xisingular}); $q,t$ will be not 
roots of unity.
The sequence $\be_g,\be_{g-1},\cdots,\be_2,\be_1$ of 
{\em simple, maybe coinciding} roots in $\tR_+^0$ constructed 
in (\ref{lambdarelprodu}),(\ref{lambdainvrel}) will
be used; it is defined in terms 
of a given reduced decomposition of $\hw\in \hW^\flat$.

We continue using the Hecke algebra $\widetilde{H}^0$.
Given a reduced decomposition of 
$\tu\in \tW^0$, we will need the elements 
$T^{0\dag}_{\tu'\,}\in \widetilde{H}^0$ for 
$\tu'\in \b(\tu;\tR_+^0)$, defined by crossing out the
corresponding $T_i^0$ in the corresponding product for
$T_{\tu}^0$. In contrast to $T_{\tu'\,}$, 
$T^{0\dag}_{\tu'\,}$ may depend on the choice
of the reduced decomposition of $\tu$ (unless 
the resulting decomposition of $\tu'$ remains reduced).
Respectively, ${\b}^\dag(\tu;\tR^0)$ is the set
of such decompositions, called {\em standard} in 
(\ref{macfilt}).

\begin{proposition}\label{VINDUCEDY}
Let all (simple) singular reflections in the decomposition
of $\hw$ belong to
a disjoint union $\{L_j,\, 1\le j\le r\}$ of consecutive
segments (connected portions) of a given 
reduced decomposition of $\hw$;\, $\cup_j L_j$ may be smaller
than $\la(\hw)$, there can be gaps between 
$L_j$, $L_{j+1}$. Following Proposition \ref{BRUHTW0}, 
we obtain the reduced decomposition
$\hw|_0=s_{\be_g}\cdots s_{\be_1}$ for
\underline{simple} $\be_j\in \tR^0$. In terms of
$L_j$, \, $\hw|_0=\tu_r\cdots\tu_1$, where 
$\tu_j=s_{\be_{b^j}}\ldots s_{\be_{a^j}}$
is the part of the decomposition of $\hw|_0$
corresponding to $L_j$ for
$a^1\le b^1<a^2\le b^2\ldots$ determining this
partition.

We impose (\ref{cyclro}) for $\tu=\hw|_0$, however, do not
assume $\xi$ to be generic.
Then $Y$\~module $\widetilde{V}^\xi_{\hw}$ is cyclic
under the assumptions\,:

(i) $q^{(\tal,\xi+d)}\neq t_{\tal}^{\pm 1}$ for
$\tal=\tal^m\in \la(\hw)$ corresponding
to the (simple) reflections\, $s_{i_m}\in \cup_j L_j$ for 
$\hw=\pi_r s_{i_l}\cdots s_{i_1}$, i.e., non-singular
intertwiners $\widetilde{\Psi}_{i_m}$ 
in $\widetilde{\Psi}_{\hw}(v)$ are all invertible 
for such $s_{i_m}$;

(ii) if $T^{0\dag}_{\tu'\,}\in \widetilde{H}^0$ are linearly 
dependent for some $\tu'\in \tU'\subset\b^\dag(\tu;\tR_+^0)$ then
the singular reflections removed to obtain
the standard decompositions of these $\tu'\in \tU'$ 
must be \underline{all} from \underline{one} $L_j$, 
i.e., can be obtained from the same $\tu_j$.
\sq
\end{proposition}

Note that (ii) holds if $\tu'\neq\tu''$ for elements 
$\tu',\tu''\in \b(\tu;\tR_+^0)$ that come from different
$L_{j}$, however (ii) is of more general nature.
For instance, if all singular reflections form a connected
segment in the reduced decomposition of $\hw$ then (i,ii)
hold. 

Assuming that there exists 
$\hv\in \hW^\flat$ sending {\em simple}
roots of $\tR^0_+$ to {\em simple} roots of $\tR_+$, the element
\begin{align}\label{hvbeal}
&\hw=\hv\tu_0=s_{\al(g)}\cdots s_{\al(1)}\,\hv \for
\hv(\{\be_j\})=\{\al(j)\},\, 1\le j\le g,
\end{align}
for the longest $\tu_0\in\tW^0$
satisfies all conditions of the proposition. The corresponding
module $\widetilde{V}^\xi_{\hw}$ is $Y$\~cyclic. Indeed,
$\hv(\be)>0$ for any $\be\in\tR_+^0$ by construction, therefore
$\la(\hv)$ contains no singular roots. Any elements
$\hv'\tu_0$ such that $l(\hv'\tu_0)=l(\hv'\hv^{-1})+l(\hv\tu_0)$
can be taken too, since
$\la(hv'\tu_0)\setminus \la(\hv\tu_0)$ cannot contain
roots from $\tR^0_+$ (they are all already in $\la(\hv\tu_0)$).
We come to the following corollary.

\begin{corollary}\label{CORYCYC}
Assuming that $q,t$ are not roots of unity and 
an element $\hv\in \hW^\flat$ exists sending simple roots
of $\tR^0_+$ to simple roots of $\tR_+$, one can find $\hw$
such that the modules $\widetilde{V}^\xi_{\hw'}$ are
$Y$\~cyclic when $l(\hw')=l(\hw'\hw^{-1})+l(\hw)$, i.e. when
$\hw'$ are sufficiently big. \sq
\end{corollary}

Note that if
$\widetilde{V}^{\xi}_{\hw}$ is cyclic for $\xi=-\rho_k$ and
$\hw=\pi_c$, then
$\widetilde{V}_c$ is cyclic in the polynomial representation;
thus we can use Corollary \ref{CORYCYC} and
Proposition \ref{VINDUCEDY} to study $\v$.
For instance, one can take $\hv\tu_0$ for $\hv$ from (\ref{hvbeal})
provided that it is in the form $\pi_c$ for certain $c$.
The latter means that $\hv=\pi_b$  and
the roots from $\tR^0_+$ have negative nonaffine
components.

\medskip
\section{The radical}
\setcounter{equation}{0}
The technique of intertwiners is expected to help in
decomposing $\v$ in terms of the irreducible constituents.
The first step in this direction is finding
{\dfont singular} $q, t_\nu$ making the 
{\em radical} of the polynomial
representation $\v$ nonzero.
In this section $q$ is generic, not a root of unity, so
the problem is in finding singular $t_\nu$ in terms of $q$.
The radical is an $\HH^\flat$\~submodule
defined for the ``evaluation
pairing" in $\v$. There are several cases
when the radical is zero but $\v$ is reducible
(this never happens in the rational case and for
simply-laced root systems);
the technique of intertwiners makes it possible to 
describe these cases.

Actually using the intertwiners alone 
is essentially sufficient to describe 
all cases when $\v$ becomes reducible 
(for instance, it is possible in the simply-laced case) 
without any reference to the radical.
However, the approach via the radical significantly
simplifies the combinatorics involved.

In the $A-D-E$\,\~\,case, the answer is simple
to formulate. It follows from
Theorem \ref{RADZERO} (see also (\ref{yofxdeg}))
and the Zigzag Lemma \ref{ZIGZAG} below. 

\begin{theorem}\label{RADADE} 
Let
$t=t_{\sht}$, $\ze_i$ be primitive $(m_i+1)$th\,
roots of unity for the classical exponents $m_i$ in
the simply-laced case;
see (\ref{tdegr}).
When  $q$ is not a root of unity,
the radical is nonzero if and only if 
$$
t=q^{-l-\frac{j}{m_i+1}}\,\ze_i^{j\,'}
\hbox{\ as\ \,} 1\le i\le n,\ 
0\le j,j\,'\le m_i,\, j+j\,'>0, \, l\in 1+\Z_+. 
$$ 
Moreover, the radical is nonzero if and only if $\v$ is 
reducible.\sq
\end{theorem}
\smallskip

\subsection{Basic properties}
Recall that the {\em evaluation pairing} from (\ref{symfg}) is
as follows:
\begin{align}
 &\{f,g\}=\{L_{\imath(f)}(g(X))\} =
 \{L_{\imath(f)}(g(X))\}(q^{-\rho_k}) \for
f,g\in \v,
\label{symfgg}
\\
&\imath(X_b)\ =\ X_{-b}\ =\ X_b^{-1},\
\imath(z)\ =\  z \for
z\in \Q_{q,t}\, .
\notag \end{align}
It induces the
$\Q_{q,t}$-linear anti-\-involution $\phi$ of
$\HH^\flat$ from (\ref{starphi}): 
\begin{align}
&\phi\equal \vep\,\star=\star\,\vep:\
X_b\mapsto Y_b^{-1},\
T_i\mapsto T_i\ (1\le i\le n).
\label{triangledown}
\end{align}

\begin{lemma}\label{RADI}
For arbitrary nonzero $q,t_{\sht},t_{\lng},$
\begin{align}
&\{f,g\}=\{g,f\} \and
 \{H(f),g\}=\{f,H^\phi\,
(g)\},\ H\in \HH^\flat.
\label{syminv}
\end{align}
The quotient $\v'$ of $\v$ by the radical
$Rad\equal$Rad$\{\, ,\,\}$
of the pairing $\{\ ,\ \}$ is an
$\HH^\flat$-module such that

(a) all $Y$-eigenspaces of $\v'$ are zero or
one-dimensional,

(b) $E(q^{-\rho_k})\neq 0$ if the image $E'$ of
$E$ in  $\v'$
is a nonzero $Y$-eigenvector.
\end{lemma}
{\em Proof.} Formulas (\ref{syminv}) are from
Theorem 2.2 of \cite{C4}. See also  \cite{C1},  Corollary 5.4.
We recall the argument from \cite{C12}.
Since  Rad$\{\, ,\,\}$  is a submodule, the
form  $\{\ ,\ \}$
is well defined and nondegenerate on $\v'.$
For any pullback $E\in \v$ of $E'\in \v',$
$E(q^{-\rho_k})=\{1,E\}=$ $\{1',E'\}.$
If $E'$ is a $Y$\~eigenvector in $\v'$ and
$E(q^{-\rho_k})$ vanishes , then
$$
\{\h_Y^\flat(1'),Q_{q,t}[Y_b](E')\}=0=
\{\v\cdot\h_Y^\flat(1'),E'\}.
$$
Therefore $ \{\v',E'\}=0,$ which is impossible.
\sq

The following lemma follows directly from the definition
of the radical.

\begin{lemma}\label{LEMKER}
The radical $Rad$ is the greatest $\HH^\flat$-submodule
in the kernel of the map
$f\mapsto \{1,f\}=f(q^{-\rho_k}).$
It is also
the intersection of the kernels of 
evaluation maps $f\mapsto f^\vee_c(e)=\{e,f\}$
for all $e\in \v(-c_\#)^\infty$:
\begin{align}\label{xicmap}
\xi_{c}:\v\ni f\mapsto f^\vee_c\in  (\v(-c_\#)^\infty)^\vee\equal
\hbox{Hom}(\v(-c_\#),\Q_{q,t}),
\end{align}
where the later linear space is supplied with the
following natural action of $\Q_{q,t}[X_b]$:
$X_b(f^\vee(e))=f^\vee(Y_b^{-1}(e))$ for $e\in  \v(-c_\#)^\infty$.
Here the maps
$\xi_c$ are $X$\~homomorphisms and their kernels 
are ideals in $\v$.
\sq
\end{lemma}
\smallskip

Using the spaces $\widetilde{V}_c$
from Section \ref{sec:tildeE}, defined for reduced decompositions
of $\pi_c$,
we can improve this statement.
Namely, $f\in Rad$ if and only if the $X$-homomorphisms
\begin{align}\label{xicmapv}
\zeta_{c}:\v\ni f\mapsto f^\vee_c\in  \widetilde{V}_c^\vee=
\hbox{Hom}(\widetilde{V}_c,\Q_{q,t})
\end{align}
are all zero for $f\mapsto f^\vee_c=\{\widetilde{E},f\}$, where
$\widetilde{E}\in \widetilde{V}_c$.
\smallskip

It is important to know how the zero-value condition 
$\zeta_c(\widetilde{E})=0$
is transformed when $c$ changes.
Let us assume that
$(\al_i^\vee,c)>0$ and that $\Psi_i^c$ is infinity.
i.e., satisfies (\ref{interinfty}).
Then the zero-value condition
for the next $\widetilde{V}_{s_i\llb c\rrb}=
\widetilde{V}_c+\tau_+(T_i)\widetilde{V}_c$
reads as follows:
$$\zeta_c(\tau_+(T_i)(\widetilde{E}))=0\hbox{\ and\ }
\zeta_c(\widetilde{E})=0.
$$
Here we use (\ref{etatxpi}):
$$ \phi(\tau_+(T_i))=(\star\cdot\tau_+\cdot\eta)(T_i)=
\tau_+(T_i) \for n\ge i\ge 0.
$$
Recall that $\tau_+(T_i)=T_i$ for $i>0$ and
$\tau_+(T_0)=X_0^{-1}T_0^{-1}$ for $X_0=qX_{\vth}^{-1}$.

The cases ($a,c$) from the definition
of $\widetilde{V}$\~spaces can be readily considered using
formulas (\ref{phiPsi}) for the action
of the $\phi$\~images of the
simple intertwiners acting on polynomials.
\smallskip

Combining Lemma \ref{LEMKER} and Lemma \ref{RADI}, we come to
the following lemma.

\begin{lemma}\label{EINRAD}
(i) A $Y$-eigenvector $E$ belongs to $Rad$ if and only
if $E(q^{-\rho_k})=0.$ It holds true if the later vanishing
condition is replaced by $E(q^{-c_\#})=0$ provided that
$\Psi_{\pi_c}$ evaluated at $q^{-\rho_k}$
is an invertible element in the subalgebra $\h_X^\flat$
of $\HH^\flat$, which is generated by $X_a (a\in B),\,$ 
and $T_i (i>0).$
For instance, $c=\pi_r(\!(0)\!)=\om_r$
can be taken instead of $c=0$ for $r\in O,\, \om_r\in B$;
similarly, 
$c=s_0(\!(0)\!)=\vth$ can be taken instead of $0$ provided that
$qt_0^{(\rho,\vth)}\neq t^{\pm 1}.$

(ii) The equality $E(q^{-\rho_k})=0$ automatically results in
the equalities
\begin{align}
E(q^{-b^\circ_\#})=0 \hbox {\ \ for\ \ }
b^\circ \in B^\star\equal \{b=b^\circ\in B\,
\mid \,  E_{b^\circ}(q^{-\rho_k})\neq 0\},
\end{align}
where $E_{b^\circ}$ is the
Macdonald polynomial for primary
$b^\circ$ defined in (\ref{primarydef}).

(iii) Expanding $E=\sum_{\mathbf{m}}
C_{\mathbf{m}} X_1^{m_1}\cdots X_n^{m_n},\ C_{\mathbf{m}}\neq 0$
for $E$ from (i) and setting $X_i=X_{b_i}$, the degree Deg$(E)$
is defined as
$$
\hbox{Deg}(E)\equal\hbox{Max}_{\mathbf{m}}\{m_1+\ldots+m_n\}-
\sum_i\hbox{Min}_{\mathbf{m}}\{m_i\}.
$$
Then Deg$(E)^n$ is no smaller than the cardinality
of the intersection
$\cap_{i=1}^n \{(c^i +B^*_\#)\cap B^*_\#\}$, where
$B^*_\#=B^*\cap B_\#$ and
the translations $c^i(E)$ of $E$ are assumed to have finitely
many common zeros; such set
$\{c^i\}$ always exists in  $B_+$, more generally,
in $u(B_+)$ for an arbitrary $u\in W$.
\end{lemma}

{\em Proof.} The first claim follows from
Lemma \ref{RADI}. Let us check that $E\in Rad$
if and only if  $\{A\,,\,E\}=0$ for {\em one}
invertible $A\in \h_X^\flat.$ Indeed,
\begin{align*}
&\{A\,\H\, Y_a(1)\,,\,\Q_{q,t}[Y](E)\}=0
\Rightarrow \{\Q_{q,t}[X]A\,\H\, Y_a(1)\,,\,E\}=0\\
&\Rightarrow \{\h_X^\flat\, Y_a(1)\,,\,E\}=0
\hbox{\ \ for\ any\ }a\in B,
\end{align*}
where the nonaffine Hecke algebra $\H\subset\HH^\flat$ is used.

Coming to (iii), an arbitrary monomial
$X_b$ in the decomposition of $E$ can be represented
as a linear combination $\sum_i\hbox{Coeff}_i\, c^i(E)$ for proper
coefficients Coeff$_i\neq 0$ and sufficiently large number of
translations $c^i$.
The translations $c^i(E)$ cannot have common zeros in $C^n$ here,
because  otherwise $X_b$ would have such a zero. Moreover,
Deg$(E)^n$
is exactly the number of common zeros, counted with multiplicities,
for $n$ generic (transversal) translations
$c^i(E)$. Here $c^i$ can be taken from $B$ or from any $u(B_+).$
\sq

As an immediate application,
we conclude that the set $B^*$ must be always
smaller than the whole $B$ if the radical is nonzero.
Setting $u=$id and using the previous lemma,
we obtain that the radical always contains {\em symmetric}
Macdonald polynomials.

\begin{lemma}\label{PINRAD}
Let us suppose that the
radical $Rad$  is nonzero. Then,
given $C_i>0\, (i>0),$
it contains at least one
Macdonald polynomial
$E_{b^\circ}$ for primary $b^\circ=b_+^\circ$
satisfying $(b^\circ,\al_i)>C_i;$ see (\ref{primarydef}).
In particular,
the corresponding symmetric Macdonald polynomial
$P_{b_-^\circ}=\P (E_{b_+^\circ})$
is well defined for $b_-^\circ=w_0(b_+^\circ)$ 
(see (\ref{symmetr}))
and belongs to the radical.
\end{lemma}
{\em Proof.} Generally,
if $E_{b_+}$ is well-defined, then so is
$P_{b_-}$, and 
\begin{align*}
&\Pi_R\,\{E_b,1\}=\{P_{b_-},1\}=
P_{b_-}(q^{-\rho_k})
\hbox{\ for the Poincar\'e\ polynomial\ } \Pi_R.
\end{align*}
Here we apply the symmetrization $\P$ to $1$,
move $\P$ to $E_b$ and use (\ref{symmetr}); 
$\P(1)=\Pi_R$ due to (\ref{tdegra}). Moreover,
if $\Pi_R\neq 0$ then 
\begin{align}\label{ebpbradical}
&E_b\in Rad\,\Rightarrow\, P_{b_-}\in Rad\,\Rightarrow\, 
\{P_{b_-},1\}\,=\,0,\\
&\{P_{b_-},1\}=0\,\Rightarrow\,
\{E_b,1\}=0\,\Rightarrow\, E_b\in Rad.\notag
\end{align}
Then we use Lemma \ref{EINRAD},(iii).\sq
\medskip

\subsection{Non-negative 
\texorpdfstring{{\mathversion{bold}$k$}}{{\em k}}}
As another application,
we come to the following theorem generalizing the
description of singular $k$ from \cite{O3};
it is equivalent
to the description below from Theorem \ref{RADZERO}
in terms of the shift-operators.

\begin{theorem}\label{THMRADI}
Assuming that $q$ is generic, 
the radical vanishes if and only if
$E_{b^o}(q^{-\rho_k})\neq 0$ for all sufficiently big
primary $b^o,$ i.e., if the product in
the right-hand side of
(\ref{ebebs}) is nonzero for all $b\in B$
with sufficiently large $(b,\al_i)$ for $i>0.$
Here $t_\nu\neq 0$ are arbitrary. 
\sq
\end{theorem}
\medskip

We are going to check the irreducibility of $\v$
in the case of non-negative $k$ in this section;
it automatically gives that $Rad=\{0\}$.
One can of course use Theorem \ref{THMRADI} to
check that  $Rad=\{0\}$ for non-negative $k$.

\begin{proposition}\label{IRRPOLYN}
Let us assume that $\Re k_{\lng}\ge 0$ and $\Re k_{\sht}\ge 0$
for $q>1.$
Then all elements $\Psi^c_{i}, c\in B$ from (\ref{Phijb}),
they belong to $\H$, are invertible,
the polynomial representation $\v$ is
$Y$-semisimple with simple $Y$\~spectrum and, moreover,
irreducible.
\end{proposition}
{\em Proof.} Setting $\al^\vee=\sum_i m_i\al_i^\vee$ for
$\al\in R,$
let
$$
\hbox{sht}(\al)=\sum_{\al_i=\sht}m_i,\
\hbox{lng}(\al)=\sum_{\al_i=\lng}m_i.
$$
Then $(\al^\vee,\rho_k)=\hbox{sht}(\al)k_{\sht}+
\hbox{lng}(\al)k_{\lng}.$
In the simply-laced case we set $\hbox{lng}(\al)=0,$ i.e.,
treat all roots as short.
Note that $\hbox{sht}(\al)>0$ for short $\al$ and
$\hbox{lng}(\al)>0$ for long $\al.$
Indeed, the intersections of $R$ with the lattices
$\oplus_{\nu_i=\nu}\Z\al_i$ for fixed $\nu=\nu_{\sht},$
$\nu=\nu_{\lng}$ contain only roots $\al$ of the same
$\nu_\al=\nu$. These root subsystems of $R$ are respectively
the sets of $\al$ such that $ \hbox{lng}(\al)=0$ and
$\hbox{sht}(\al)=0$.

Recall the conditions that are necessary and
sufficient for the intertwiner $\Psi_i^c$ to be
proportional to

\centerline{
(a) $\tau_+(T_i)+t_i^{-1/2},$\ \ \  (b) $\tau_+(T_i)+\infty,$\ \ \
(c) $\tau_+(T_i)-t_i^{1/2}$ ;}
\noindent
they are as follows:
\begin{eqnarray}\label{interplust}
&q_\al^{-k_\al+(\tal^\vee,c_-+d)-(\al^\vee, \rho_k)}=1
&\for (a): \tau_+(T_i)+t_i^{-1/2},\\
&q_\al^{(\tal^\vee,c_-+d)-(\al^\vee, \rho_k)}=1
&\for (b): \tau_+(T_i)+\infty,
\label{interinfinity}\\
&q_\al^{k_\al+(\tal^\vee,c_-+d)-(\al^\vee, \rho_k)}=1
&\for (c): \tau_+(T_i)-t_i^{1/2},
\label{interminust}
\end{eqnarray}
where we set
$$\tal=u_c(\al_i),\ \al=u_c(\al_i) \for i>0, \and
\al=u_c(-\vth) \for i=0.
$$
Here $\nu_\al=\nu_i,$
$\al<0,$ and $(\tal,c_-+d)=(\al_i,c+d)>0.$
See (\ref{aljb}).
In case ($b$), the intertwiner
$\Psi_i^c$ is not well defined. We will call 
it {\em infinity} and the corresponding $s_i$ {\em singular}
following the previous sections.

Due to the positivity assumptions,
$$
-\Re(\,k_\al+(\al^\vee, \rho_k)\,)\ge 0 \and
\Re(\,-k_\al+(\tal^\vee,c_-+d)-(\al^\vee, \rho_k)\,)>0.
$$
Thus (\ref{interplust})-(\ref{interminust}) never take place,
all intertwiners are invertible, and all $E$-polynomials
are well defined. This approach actually allows $\Re k$ to
be ``small" negative. 

Recall that the eigenvalues are
$\{q^{-c_\#}=q^{u_c^{-1}(\rho_k-c_-)}\}$. Generally,
\begin{align}\label{qucqub}
&q^{u_c^{-1}(\rho_k-c_-)}\ \neq\ q^{u_b^{-1}(\rho_k-b_-)}\for
c\neq b \in B \hbox{\ \, if}\\
&\Re(\,\rho_k-c_-\,)\,\in \{z\in \R^n,
\Re(\,(\al_i^\vee,z)\,)>0\}\ni\, \Re(\,\rho_k-b_-\,),
\notag
\end{align}
for instance, for sufficiently large $(\al_i^\vee,-c_-)$ 
for all $i>0$. Indeed, then $u_c=u_b$ and $c_-=b_-.$ 
Due to $\Re k\ge 0$, the inequalities
$\Re(\,(\al_i^\vee,z)\,)\ge 0$ are sufficient in  (\ref{qucqub})
and the latter holds for any $c,b\in B$.


Concerning the irreducibility, if $\v'\subset \v$
is an $\HH^\flat$-submodule then it contains at
least one eigenvector $v$ corresponding to $q^{-c_\#}.$
Applying the intertwiners, it must contain all
eigenvalues (including the simple ones for
sufficiently big $c$ if $\Re k$ is allowed to be small negative).
Therefore it contains at least one $E$\~polynomial. However then
the intertwines will make all of them in $\v'.$ \sq
\medskip

\rmk
The method above gives the irreducibility of
$\v$ as $\Re k>-1/h$ for the Coxeter number 
$h=1+(\rho,\vth)$ in the simply-laced case;
the inequalities are somewhat more involved for 
$B,C,F,G$. Here one can also use an analytic variant of the 
inner product (\ref{staru}) in $\v$. Let us outline an
approach via the roots of unity assuming
that $k_{\lng},k_{\sht}$ are rational numbers, which
is sufficient for the analysis of the irreducibility.

Following \cite{C4},(6.12)-(6.14),\, $\v$ can be 
supplied with a $\star$\~invariant inner product, 
where $q^N=1$ for sufficiently large $N$ (also satisfying
certain divisibility conditions). Under $k>-1/h$
and the corresponding conditions in the non-simply-laced
case, this inner product is nonzero and is
{\em positive definite} on $W$\~invariant polynomials.
Then we use that reducibility of $\v$ implies the 
reducibility of $\v^W$ under the action of the 
subalgebra of invariants of $\HH^\flat$ (for generic $q$)
and that $\v$ is in the category $\o_Y$.
\sq 

\medskip
\subsection{Using affine exponents}
Continuing to assume that $q$ is in a general
position (see the exact condition in the
theorem below), we will
examine when $Rad$ is zero via the
shift operator. As a matter of fact, the answer has 
been already obtained (Theorem \ref{THMRADI}). However, 
using the shift operator has certain technical advantages 
in the $q,t$\~case and generalizes Opdam's analysis
of the rational case 
(where it is the only known approach). Note that
Theorem \ref{THMRADI} has no {\em differential} counterparts,
rational or trigonometric.

The set of {\em singular} $k$ with $Rad\neq \{0\}$
is given by some algebraic equalities.
Since $q$ is generic (only roots of unity must
be avoided), real $q>1$ are sufficient to consider.

We will apply formulas (\ref{kconfin}),(\ref{kconfinnx}) :
\begin{align}\label{kconfinn}
&\Pi_{\tR}\,P_{b+\rho}^{q,tq}(q^{-\rho_{k+1}})
\ =\ q^{\{\cdots\}}\,\Bigl(\prod_{\al\in R_+}
1- q_\al^{k_\al+(\al^\vee,\,\rho_k-b)}\Bigl)\,
P_b^{q,t}(q^{-\rho_k}).
\end{align}
Recall the notation
$$
(tq^j)^{(\,\cdot\,,\rho^\vee)}=
\prod_\nu (t_\nu q_\nu^j)^{(\,\cdot\,,\rho_\nu^\vee)};
$$
see (\ref{yofxeval}) for the definition of $\Pi_{\tR}$.

\begin{maintheorem}\label{RADZERO}
(i) We assume that  $q\neq 0$ is not a root of unity. For
the Poincar\'e polynomial $\Pi_R$ from (\ref{tdegra}),
$Rad=$Rad$\{\, ,\,\}$ is zero if and only if 
\begin{align}\label{piinfty}
\Pi_R^{-1}\,&\prod_{l=0}^\infty
\hbox{\rm rad}(q,tq^l)\neq 0\,\hbox{\ \ for\ \ }
\hbox{\rm rad}(q,t)\equal\Pi_{\tR}
\\
=&\prod_{\al\in R_+}
\Bigl( (1- t_\al(tq)^{(\al,\rho^\vee)})
\prod_{ j=1}^{(\al^\vee,\rho)}
\frac{
(1- q_\al^{j-1}t_\al t^{(\al,\rho^\vee)})}
{(1- q_\al^{j-1}t^{(\al,\rho^\vee)})}\Bigr). \notag
\end{align}

(ii) Provided that $\rho\in B$ and 
$\Pi_R\neq 0$, 
the zeros of\, $\Pi_{\tR}$ described explicitly in
formula (\ref{yofxdeg})
in the cases $\widetilde{A},\widetilde{D},\widetilde{E},$\\
\centerline{(\ref{yofbfinal}),(\ref{yofbdenom}) 
for $\widetilde{B}$,\ \ \ \, formula
(\ref{yofxevalc}) for $\widetilde{C}$,}\\
\centerline{(\ref{g2affinen}),(\ref{g2affined}) 
for $\widetilde{G}$,
and (\ref{f4affinen}),(\ref{f4affined}) for $\widetilde{F}$}\\
constitute the set of all $q,t$
such that the $t$\~discriminant $\x^t$ belongs to $Rad$;
see formula (\ref{shiftoper}). Here $q\neq 0$ can be
arbitrary, possibly a root of unity.

(iii) If $q$ is not a root of unity ($B$ can be arbitrary), 
then the zeros of $\Pi_{\tR}/\Pi_R$ from (i) are those
listed in (ii) and their \underline{positive translations}. 
The latter are obtained from the zeros in (ii) with
$k_\nu$ replaced by \,\underline{all}\, $k'_\nu\in k_\nu+\Z_+$
disregarding the binomials that do contain the factors $q^j$ 
with $j>0$.  Moreover, if 
\begin{align}\label{qtgeneric}
&\{q^a\prod_\nu t_\nu^{b_\nu}=1\ \hbox{\ for\ } a,b\in \Z\}
\Rightarrow  \{a+\sum_\nu\nu k_\nu b_\nu=0\},\ \nu\in \nu_R,
\end{align}
then only the \,\underline{rational exponents}\,,
described in Theorem \ref{BPRIMDEG}, are 
sufficient to consider here and in (ii).
\end{maintheorem}

{\em Proof}. If the radical is nonzero
then it contains symmetric
Macdonald polynomials $P^{q,t}_b$ with
sufficiently large negative $b=b_-$
due to  Lemma \ref{PINRAD}. Recall that
$P^{q,t}_b\in Rad$ $\Leftrightarrow$ 
$P^{q,t}_b(q^{-\rho_k})=0$ provided that
$\Pi_R\neq 0$; see (\ref{ebpbradical}).
Let us assume that $\Pi_R\neq 0$ in the following
reasoning.

We apply
consecutively the shift operators changing
$$
b\mapsto b^{(m)}= b+m\rho,\ k_\nu\mapsto (k+m)_\nu= k_\nu+m,
\ t\mapsto t^{(m)}=tq^m,
$$
to make $\Re\, (k+m+1)_\nu\ge 0.$
The radical is
trivial for $k+m+1$ thanks to Proposition \ref{IRRPOLYN}
and $P\,^{q,t^{(m+1)}}_{b^{(m+1)}}(q^{-\rho_{k+m+1}})\neq 0.$
This step is similar to the approach of 
\cite{O3,DJO} in the rational 
case, although the relation between the
evaluation pairing and the shift operators
is different in the rational and 
the $q,t$\~cases. 

Due to Theorem \ref{YXPB}, see also (\ref{kconfinn}),
\begin{align*}
&\{\overline{\y}\,^{t^{(m)}}\x^{t^{(m)}}\}
P\,_{b^{(m+1)}}^{q,t^{(m+1)}}(q^{-\rho_{k+m+1}})=
\overline{g}\,^{q,t^{(m)}}(b)\,
P\,_{b^{(m)}}^{q,t^{(m)}}(q^{-\rho_{k+m}}),\\
&\overline{g}\,^{q,t^{(m)}}(b)=
\prod_{\al\in R_+}
((t_\al q_\al^m)^{-1}X_\al(q^{(b-\rho_{k})/2}) -
 X_\al(q^{(\rho_{k}-b)/2})),
\end{align*}
where $\{\overline{\y}\,^{t^{(m)}}\x^{t^{(m)}}\}$ is
$\hbox{\rm rad}(q,tq^m)$ up to
a power of $q.$ Here
$\overline{g}\,^{q,t^{(m)}}(b)$ can be supposed 
nonzero because $b$ is assumed sufficiently large negative.

We note that there are no restrictions
concerning $q,t$, when using the shift operator
(although the relation $\Pi_R\neq 0$ is necessary
for the applications to the radical).
Also the above formula holds for {\em arbitrary}
symmetric $\Q_{q,t}[Y]^W$-eigenfunctions $P\in \v$
with the $Y$\~weight coinciding with the weight of
$P_b\, (b=b_-\in B_-)$ and such that 
$P=\sum_{a\in W(b_-)}X_a$ modulo
$X_c$ for $c_-\succ b_-\,$;
see (\ref{Lf}),(\ref{order}). 

Thus all partial products
$$
\hbox{\rm rad}(q,t)\,\hbox{\rm rad}(q,tq)\,\cdots\,
 \hbox{\rm rad}(q,tq^{l}) \for l\ge m
$$
must vanish. 

Without imposing $\Pi_R\neq 0$ and using the shift-operator, 
one can use directly Theorem \ref{THMRADI}, which claims that 
the radical vanishes if and only if
$E_{b}(q^{-\rho_k})\neq 0$  for sufficiently large
$b$. Using the evaluation formula (\ref{ebebs}) 
as $\la(\pi_b)\to \{[-\al,\nu_\al j],\al\in\R_+,j>0\}$
we come exactly to (\ref{piinfty}).
\smallskip

Claim (ii). Thanks to Theorem \ref{YOFXEV},
 rad$(q,t)=\Pi_{\tR}=0$ if and only if $\{\x^t,\x^t\}=0$;
cf. \cite{O3}, \cite{DJO}. We need to check that
$\{\x^t,P\}=0$ for an arbitrary polynomial $P\in \v$.
Applying the $t$\~anti-symmetrization (here $\Pi_R=\P(1)\neq 0$ 
is needed), it suffices
to assume that $P$ is $t$\~antisymmetric. Due to the
inequalities $t_\nu\neq -1$ (which follow from $\Pi_R\neq 0$),
it must be then divisible
by $\x^t$, which gives the required. 

The other claims
from (ii,iii) follow from the definition and calculations
of the affine and rational exponents performed in the formulas
listed in (ii), (iii).
\sq
\smallskip

\rmk
Using the formula for $\overline{\y}^t(\x^t)$ with 
$k'_\nu\in k_\nu+\Z_+$ in (ii,iii) instead of $k$ 
is essentially equivalent to direct using 
the evaluation formula for the $E$-polynomials or 
$P$-polynomials. The approach without the shift-operator
is more general, since it is not necessary to assume
that $\Pi_R\neq 0$. On the other hand, 
using  $\overline{\y}^{t}(\x^t)$, which is
a regular function by construction, makes the
structure of the resulting formula more transparent
and convenient for practical finding {\em singular}\,
$k$.

The claim that $\x^t\in Rad$ (the first part of (ii))
has interesting applications including the case of 
roots of unity. Note that the radical is {\em always\,}
of finite codimension as $q$ is a root of unity. 
However, rad$(q,t)=\Pi_{\tR}$ becomes zero only for finitely
many $k$. The ``singular" $k$ form finitely many sequences 
$\{k_\nu+\Z_+\}$ such that applying the shift 
operator will eventually make $\x^t\in Rad$. 
\sq
\medskip

\section{Irreducibility of 
\texorpdfstring{{\mathversion{bold}$\v$}}{V}}
We study the irreducibility
of the polynomial representation $\v$ under the assumption
that the radical $Rad$ of the form  $\{\ ,\ \}$
is zero. Recall that $\tR^0$, $\tR^{-1}$ are
from (\ref{troeq}) and (\ref{critrootoneuminus});
$\tR^0_+$, $\tR^{-1}_+$ and 
$\tR^0_+[-]$, $\tR^{-1}_+[-]$ denote positive roots
in these sets and those with negative non-affine 
components. 

The radical
is of finite codimension when $q$ is a root of unity, so
we shall consider only generic $q.$
From now on we assume that $q=q_{\sht}$ is not a root of unity,
neither is $q_{\lng}=q^{\nu_{\lng}}$.
\smallskip

\subsection{Properties of 
\texorpdfstring{{\mathversion{bold}$\tR^0,\tR^{-1}$}}
{R0,R-minus}}
We will need to develop some tools of combinatorial 
nature.
Recall that the relation $(\al_i,c+d)>0$ is necessary and
sufficient
for $s_i\pi_c$ to be represented in the form $\pi_b$ for
$l(\pi_b)=l(\pi_c)+1$; see Proposition \ref{BSTAL}. 
$$\la(\pi_b)\setminus\la(\pi_c)=\pi_c^{-1}(\al_i)=
\tal+[0,(\al_i,c)]=[u_c(\al_i),\de_{i0}+(\al_i,c)]
$$ 
for the Kronneker delta,
where $\tal\equal u_c(\al_i), (\al_i,c)=(\tal,c_-)$. The
nonaffine component of $\tal$  must be negative since
$\tal$ belongs to $\la(\pi_b)$.

For $c_\#=\pi_c(\!(-\rho_k)\!)=c-u_c^{-1}(\rho_k),$
$$(\tal, c_- -\rho_k+d)\ =\ (\al_i,c_\#+d)\ =\
(\tal+[0,(c,\al_i)]\,,\,-\rho_k+d),$$
where we transform:\, $(\tal,c_- +d)=(u_c(\al_i),u_c(c)+d)=
(\al_i,c+d)$.

The following
lemma uses the notation from the tables from \cite{Bo}.
\smallskip
\begin{lemma}\label{ZFROMM}
For $c\in B$ and an index $0\le i\le n$ satisfying the
inequality $(\al_i,c+d)=(\al,c_-)+\de_{i0}>0$,
let
\begin{align}\label{minustc}
&t_{\al}^{-1}q^{(\tal, c_--\rho_k+d)}\equal
q_\al^{(\tal^\vee, c_--\rho_k+d)-k_\al}= 1 =
t_i^{-1}q^{(\al_i,c_\#+d)}
\end{align}
for $\tal=u_c(\al_i)=[\al,\de_{i0}],$
$\al=u_c(\al_i)$ for $i>0$ or
$\al=u_c(-\vth)$ for $i=0$. 
Equivalently, $\tal'\equal \pi_c^{-1}(\al_i)=
[\al,\de_{i0}+(\al,c_-)]
\in \tR^1$\,; see (\ref{critrootone}).
Then
$-\al$ is positive and cannot be simple.

(i) If $-\al$ is long then always $\al+\al_l=(\nu_\al/\nu_\be)\be$
for proper $-\be\in R_+$ and \underline{long} 
$\al_l\, (l>0)$, where $\be$ can be taken long unless
for $\tR=\widetilde{C}_n$ (and for long $-\al$).
In the non-simply-laced case and for short $-\al$, 
we \underline{assume} that such representation exists 
with \underline{short} $\be,\al_l$.
Let $\tal'=[\al,\nu_\al j']$, $\tbe'\equal [\be,\nu_\be j']$
for $\be=(\nu_\be/\nu_\al)(\al+\al_l)$. Then
\begin{align}\label{qbetbe}
&q_{\be}^{\ \frac{\nu_\al}{\nu_\be}(\,(\tbe')^\vee\,,\,-\rho_k+d\,)}
\ =\ 1,\ \ 
\hbox{i.e.,\ \ }  \tbe'\in \Q\lan\tR^0\ran
\end{align} 
for the $\Q$\~span $\Q\lan\tR^0\ran$ of $\tR^0$.
Here $\tbe'\in \tR^0$ unless $\nu_{\al}/\nu_{\be}>1$
(see Lemma \ref{LEMNUALBE},(i) below).
 
If $\tbe'\in \tR^0$, then $\tbe'\not\in \la(\pi_c)$.
Moreover, one can find $b\in B$ and the index
$i'\ge 0$ such that
$l(\pi_b)=l(\pi_b\pi_c^{-1})+l(\pi_c)$
and
\begin{align}\label{sumforal}
&q^{(\al_{i'}, b_\#+d)}=1, \ (\al_{i'},b+d)>0,\
(u_b(\al_{i'}),u_c(\al_i))>0, \hbox{\ i.e.,}\notag\\
&\la(\pi_b)\,\not\ni\, \tbe'=\pi_b^{-1}(\al_{i'})\,\in\, 
\la(s_{i'}\pi_b)\,=\,\{\tbe',\la(\pi_b)\},
\end{align}
where
$\be=u_b(\al_{i'})$ for $i'>0$ or  
$\be=u_b(-\vth)$ for $i'=0$. 

(ii) Provided (\ref{qtgeneric}),
the assumption from (i) holds for any long $-\al$ and
for short $-\al\in R_+$ such that 
$(\al,\rho_k)+k_{\sht}\in 1+\Z_+$
unless $\al$ is from the following list:
\begin{align*}
-\al&=\eps_{n-g},\ g=1,\ldots,n-1,\
(\al,\rho_k)+k_{\sht}=-2gk_{\lng}\in 1+\Z_+,
&(\widetilde{B}_n)\\
-\al&=\ep_{n-1}+\ep_n,\
(\al,\rho_k)+k_{\sht}=-2k_{\lng}\in 1+\Z_+,
&(\widetilde{C}_n)\\
-\al&=\al_1+\al_2,\
(\al,\rho_k)+k_{\sht}=-3k_{\lng}\in 1+\Z_+,
&(\widetilde{G}_2)\\
-\al&=\al_3+\al_2 \hbox{\  or\  } -\al=\al_3+\al_2+\al_1
\hbox{\ (  see \cite{Bo}\,),\ }
&(\widetilde{F}_4)\\
&\hbox{where\ }
(\al,\rho_k)+k_{\sht}\ =-2k_{\lng} \hbox{\ \,or\ \,} =-4k_{\lng}
\in 1+\Z_+,\\
-\al&=1221,\ \ (\al,\rho_k)+k_{\sht}=-6k_{\lng}-2k_{\sht}\in 1+\Z_+.
&(\widetilde{F}_4)
\end{align*}
In the case of long $\al$ such that 
$\nu_{\al}/\nu_{\be}>1$, (\ref{qtgeneric}) 
results in $\tbe'\in \tR^0$.

(iii)  Among the cases listed 
in (ii) subject to (\ref{qtgeneric}), the radical
can be zero only when $k_{\lng}=-j\in -1-\Z_+$, which
implies that $[-\al_p,\nu_\al j]\in \tR^0_+$ for any
long \underline{simple} root $\al_p\, (p>0)$,\, or in the 
following subcases: 
\begin{eqnarray}\label{reducases}
&\widetilde{B}_n:\  &k_{\lng}=-\frac{s}{2g} \hbox{\ \ for\ \, } 
n>g>\frac{n}{2},\, s\in \N,\, (s,2g)=1,\\
&\widetilde{G}_2:\  &k_{\lng}=-s/3, \hbox{\  where\  }
 s\in 1+\Z_+,\ (s,3)=1,\notag\\ 
&\widetilde{F}_4:\  &k_{\lng}=-s/4, \hbox{\  where\  }
 s\in 1+\Z_+,\ (s,2)=1,\notag\\ 
&\widetilde{F}_4:\  &3k_{\lng}+k_{\sht}=-s/2 \hbox{\ for\ }
 s\in 1+\Z_+,\ (s,2)=1.\notag 
\end{eqnarray}
In the case of $\widetilde{F}_4$, 
it is also possible that $\ 3k_{\lng}+k_{\sht}=-j\in -1-\Z_+$\
(then $Rad=\{0\}$ for generic $k_{\sht}$), 
but in this case the long root $[-1220,2j]$ belongs to $\tR^0_+$. 
\end{lemma}

{\em Proof.}
Generally, $t_i^{\pm 1}q^{(\al_i, c_\#+d)}=1$ is equivalent to
\begin{align}\label{qalkc}
&q_\al^{\pm k_\al+(\al^\vee, c_- -\rho_k)}=1 \for
\al=u_c(\al_i) \  (i>0),\,  \\
&q_\al^{\pm k_\al+(\al, c_- -\rho_k)+1}=1 \for
\al=u_c(-\vth) \ (i=0).\notag
\end{align}
We use the definition of $c_\#$ and
the relation $u_c(c)=c_-\in B_-.$

Recall that $\al<0$ in (\ref{qalkc}) and (\ref{minustc})
due to
\begin{align}\label{alic+d}
&(\al_i,c+d)>0\Rightarrow (u_c(\al_i),c_-+d)>0,
\end{align}
which is obvious as $i>0$. However if $i=0$ than it
may happen that $(u_c(\al_i),c_-)=0=(c,\vth)$. In the
latter case the sign of $\al$ cannot be determined
from (\ref{alic+d}) and we need to use: 
$$(c, \vth)=0\Rightarrow u_c(\vth)>0, \hbox{\ \ since\ \ }
\la(u_c)=\{\al\in R_+,(c,\al)>0\}.
$$

The root $-\al$ cannot be a simple root.
Indeed, otherwise $\al=-\al_p$ for $p>0$ and
$$
1=q_\al^{-k_\al+(\al^\vee\,,\, c_- -\rho_k)+\de_{i0}}=
q_\al^{( -\al_p^\vee\,,\, c_-)+\de_{i0}}\Rightarrow
(\al,c_-)+\de_{i0}=0,
$$
which contradicts the assumption from (i); here we
use that $q_\al$ is not a root of unity.
By the way, this argument shows that
the set from (i) is empty for $\widetilde{A_1}$.
\smallskip

Now, imposing (i), let
$\frac{\nu_\al}{\nu_\be}\be=\al+\al_l$ for certain $l>0$ such that
$\nu_{\al_l}=\nu_\al,$ i.e.,
$\al$ and $\al_l$ have to be of the same length ($\nu_\be=\nu_\al$
unless for long $\al$ in $\widetilde{C}_n$). 
This
representation guaranties that $(\be,\al)>0$ and
$s_\al(\be)=\al_l>0$ that 
directly leads to the inequality
$(u_b(\al_{i'}),u_c(\al_i))>0$ stated in (\ref{sumforal}).

Let $\tbe'=\frac{\nu_\be}{\nu_\al}(\tal'+\al_l)=[\be,\nu_\be j\,']$ 
for $\tal=[\al,\nu_\al j\,']$. Then
\begin{align*}
&q_{\be}^{\ \frac{\nu_\al}{\nu_\be}\,(\,(\tbe')^\vee\,,\,-\rho_k+d\,)}
\ =\ 
 q^{(\,\frac{\nu_\al}{\nu_\be}\,\be\,,\,-\rho_k\,)+\nu_\al j\,'}\\
&=\ q_{\al}^{\,-k_{\al}+
((\,\frac{\nu_\al}{\nu_\be}\,\be-\al_l)^\vee\,,
\,-\rho_k\,)+j\,'\,}\ =\ 
q_{\al}^{-k_{\al}+(\al^\vee,-\rho_k)+j\,'}\ =\ 1.
\end{align*}
Unless $\nu_\al/\nu_\be=2$ for $\widetilde{C}_n$, we conclude
that $\tbe'\in \tR^0$. If (\ref{qtgeneric}) is assumed, then
it holds always (we can omit the factor $\nu_\al/\nu_\be$).
Recall that $q^{(\tal,\,\cdot\,)}=q_{\al}^{(\tal^\vee,\,\cdot\,)}$
{\em by definition} and all such powers must be
expressed as products of powers of $q$ and $t_\nu$.

\comment{
$q_\be^{(\be^\vee, c_- -\rho_k)+p'}=1,$ where we
set $(\al_l^\vee,c_-)=-p,$
$$p'\equal \frac{p\nu_\al}{\nu_\be}\in \Z_+ \for i>0,\ \
p'\equal\frac{(p+1)\nu_\al}{\nu_\be}\in \N \for i=0.
$$
Here $\al$ is long if $\be$ is long, that gives
the integrality of $p'$;
recall that $\nu_{\sht}=1, \nu_{\lng}=1,2,3.$ 
\begin{align*}
&\hbox{One\ has:\ \ }
q_\al^{-k_\al+((\be-\al_l)^\vee, c_- -\rho_k)}=1
\ \Rightarrow\\
&q^{(\be, c_- -\rho_k)+    p\nu_\al}=1\ (i>0)
\hbox{\ \ or\ \ }
 q^{(\be, c_- -\rho_k)+(p+1)\nu_\al}=1\ (i=0)
\ \Rightarrow\\
& q_\be^{(\be^\vee, c_- -\rho_k)+p'}=1\ \Rightarrow \tbe'\in \tR^0,
\end{align*}
where
$\tbe'\equal[\be,\nu_\be(j+p')] \for j\equal(\be^\vee,c_-).$
}

Let us check that $\tbe'\not\in \la(\pi_c).$ It can be seen
directly from (\ref{lambpi}), but it is simpler to
use that 
$$
\frac{\nu_\al}{\nu_\be}\,\tbe'\ =\ 
\frac{\nu_\al}{\nu_\be}\,[\be,\nu_\be j']\ =\ \tal'+\al_l
$$
and apply Theorem \ref{INTRINLA},($b'$). The root
$\tbe'$ can appear in $\la(s_i\pi_c)=\{\tal',\la(\pi_c)\}$ 
only {\em after} $\tal'$, which is the last in this set,
since  $\al_l$ does not belong to $\la(\pi_c)$ (and
to any $\la(\pi_a)$).

\comment{
It is obvious when
$p'>0$ for $i>0$ or $p'>1$ for $i=0.$
See formula (\ref{lambpi}).
Let us examine the case $(\al_l,c_-)=0.$

If $i>0,$ using (\ref{lapiomu}):
\begin{align*}
&(\al_l,c_-)=0=
(u_c^{-1}(\al_l),c)\, \Rightarrow\, \al_l\not\in \la(u_c^{-1})
\Rightarrow\\
&u_c^{-1}(\al_l)>0\ \Rightarrow\
u_c^{-1}(\be)=u_c^{-1}(\al+\al_l)=\al_i+u_c^{-1}(\al_l)>0.
\end{align*}
Thus $\tbe'\not\in \la(\pi_c)$ in this case.

If $i=0,$ then we have $\nu_\be=1=\nu_l$ and
the end of the last formula reads as
$$
u_c^{-1}(\be)=u_c^{-1}(\al+\al_l)=-\vth+u_c^{-1}(\al_l)<0,
$$
where we use that $\vth$ is the maximal short root.
Thus $[\be,(\be^\vee,c_-)]$ belongs to $\la(\pi_c),$
but $\tbe'$ is not in this set. Indeed, its affine
component is  $(\be^\vee,c_-)+1.$
\smallskip
}

To construct $b,$
it suffices to use an {\em arbitrary} $a\in B$ satisfying
$\pi_a=\hw s_i\pi_c$ such that
$l(\pi_a)=l(\hw)+1+l(\pi_c)$ and
$\tbe'\in \la(\pi_a);$ we simply make $(a_-,\be)$
sufficiently large. Then for any reduced decomposition
$\pi_a=\pi_rs_{i_m}\cdots s_{i_1},$
the partial products
$s_{i_{h}}\cdots s_{i_2}s_{i_1}$ remain in $\pi_B$
for any $h\le m$. Note that these products also 
remain in $\pi_b$ 
if we multiply them by {\em any} $\pi_{r'}$ on the left.

We can assume that a reduced decomposition for $\pi_a$
extends that for $\pi_c.$
Let $s_{i'}\pi_b$ be the first partial product
``after" $\pi_c$ such that
$\tbe'\in\la(s_{i'}\pi_b).$ It gives the required $b$
since the affine root $\tbe'$
can appear in the $\la$\~sets of the partial products
only after $\la(\pi_c).$
\smallskip

Claim (ii) is checked
using the tables from \cite{Bo}.
Concerning (iii), we need to examine the
special short $\al$ listed in (ii) and
employ Theorem \ref{RADZERO}.
\sq
\medskip

\subsection{A generalization (any  
\texorpdfstring{{\mathversion{bold}$t$}}
{{\em t}})}
It is not too difficult to extend the lemma to a 
general setting when $q$ is assumed to be generic
(not a root of unity) but condition 
(\ref{qtgeneric}) is omitted in (ii),(iii). 
The main change we need is at the end of
(ii), where it is stated that
in the case of long $\al$ such that 
$\nu_{\al}/\nu_{\be}>1$ condition (\ref{qtgeneric}) 
results in $\tbe'\in \tR^0$. Now we cannot use it.
First of all, we need to rewrite the
formulas in the lemma by adding $q$.

The $\Z$\~integrality
conditions in (ii) become:
\begin{align}\label{tlngcases}
&t_{\lng}^{-g}\in q^{1+\Z_+} \for \widetilde{B}_n,\ \ 
t_{\lng}^{-1}\in q^{1+\Z_+} \for \widetilde{C}_n,\\
&t_{\lng}^{-1}\in q^{1+\Z_+} \hbox{\ in the  
case of } \widetilde{G}_2, \ \hbox{\ and for\ }
\widetilde{F}_4\,: \notag\\
&t_{\lng}^{-1}\in q^{1+\Z_+},\ \,t_{\lng}^{-2}
\in q^{1+\Z_+},\ \,
t_{\lng}^{-3}t_{\sht}^{-2}\in q^{1+\Z_+}.\notag
\end{align}
Recall that $t_{\lng}=q_{\lng}^{k_{\lng}}$
for $q_{\lng}=q^{\nu_{\lng}}$.
The condition $k_{\lng}\in -1-\Z_+$ now reads as
$t_{\lng}^{-1}\in q^{2+2\Z_+}.$

Similarly, the list from (iii) becomes:
\begin{eqnarray}\label{treducases}
&\widetilde{B}_n:\   t_{\lng}^{-g}\in 
q^{1+\Z_+}\setminus q^{2+2\Z_+},\ \ 
\widetilde{G}_2:\  
t_{\lng}^{-1}\in q^{1+\Z_+}\setminus q^{3+3\Z_+}, 
\\
&\widetilde{F}_4:\  t_{\lng}^{-2}\in 
q^{1+\Z_+}\setminus q^{2+2\Z_+} \hbox{\ \,or\ \,}  
t_{\lng}^{-3}t_{\sht}^{-2} \in
q^{1+\Z_+}\setminus q^{2+2\Z_+},\notag 
\end{eqnarray}
where $n/2<g<n$ and $t_{\lng}^{-g'}\not\in q^{1+\Z_+}$ with 
any $\Z_+\ni g'<g$  for $\widetilde{B}_n$. 
\smallskip

Let us describe the cases when $\al$ is {\em long\,}
and $\nu_\al/\nu_\be>1$ in (i); they 
are exactly those dual to the cases listed in (ii).
We continue to assume that
$q$ is not a root of unity but do not impose the
$q,t$\~generality condition 
(\ref{qtgeneric}). Then (\ref{minustc}) reads as
\begin{align*}
&q^{\,\nu_{\lng}\,(\,(\al^\vee,\,\rho_k)+k_{\lng}\,)\,}\, \in\, 
q^{-\nu_{\lng}(1+\Z_+)},\\
&1\in q^{(\al^\vee,\,\rho_k)+k_{\lng}+1+\Z_+}
\Rightarrow \tbe'\in \tR^0.
\end{align*}
Recall that the first condition becomes
$(\al^\vee,\rho_k)+k_{\lng}\in -1-\Z_+$ and always implies
$\tbe'\in \tR^0$ if (\ref{qtgeneric}) is imposed.
\medskip

\begin{lemma}\label{LEMNUALBE}
(i) The cases $\nu_\al/\nu_\be>1$ are as follows:
\begin{align*}
-\al&=2\eps_{n-g},\ 1\le g\le n-1,\
(\al^\vee,\rho_k)+k_{\lng}=-gk_{\sht},&(\widetilde{C}_n)\\
&\hbox{where\ condition\ \ (\ref{minustc})\ \ becomes\ \ } 
 t_{\sht}^{-2g}\,\in\, q^{2+2\Z_+},&\\
-\al&=\ep_{n-1}+\ep_n,\
(\al^\vee,\rho_k)+k_{\lng}=-k_{\sht}, \ t_{\sht}^{-2}\in q^{2+2\Z_+},
&(\widetilde{B}_n)\\
-\al&=3\al_1+\al_2,\
(\al^\vee,\rho_k)+k_{\lng}=-k_{\sht}, \ t_{\sht}^{-3}\in q^{3+3\Z_+},
&(\widetilde{G}_2)\\
-\al&=\al_2+2\al_3 \hbox{\ \,or\ \,} -\al=\al_2+2\al_3+2\al_4,
\hbox{\ respectively,}
&(\widetilde{F}_4)\\
\hbox{for}&\ \ (\al^\vee,\rho_k)+k_{\lng}=-gk_{\sht},\   
t_{\sht}^{-2g}\in q^{2+2\Z_+}\hbox{\  as\ \,} g=1,2,\\
-\al&=1242,\  (\al^\vee,\rho_k)+k_{\lng}=-2k_{\lng}-3k_{\sht},
\ t_{\lng}^{-2} t_{\sht}^{-6}\in q^{2+2\Z_+}.
&(\widetilde{F}_4).
\end{align*} 
Let $\tau=\tau(-\alpha)\equal t_{\sht}$ in all 
these cases but the last, where
$\tau(1242)\equal t_{\lng}\, t_{\sht}^{3}.$
We set $g=1$ unless stated otherwise 
(namely, when $\tR=\widetilde{C}_n$ and 
for $\widetilde{F}_4$ as $g=2$). In these notations,
if $\tbe'\not\in \tR^0_+$ (see above) then 
$\tau^{\,-g\,\nu}\in q^{\nu+\nu\Z_+}$ for
$\nu=\nu_{\lng}$.

(ii) Continuing (i), let $Rad=\{0\}$. Then either
$t_{\sht}\in q^{-1-\Z_+}$ or one of the following conditions
hold : 
\begin{eqnarray}\label{reducasesze}
&\widetilde{C}_n:\  &t_{\sht}^{g}\in -q^{-1-\Z_+},
\hbox{\ \ where\ \ } 
n>g>\frac{n}{2}\\
&\mbox{and}&t_{\sht}^{g'}\not\in \pm q^{-1-\Z_+}
\hbox{\ \ for\ \ } \Z_+\ni g'<g,\notag
\\ 
&\widetilde{G}_2:\  &t_{\sht}\in \rho\, q^{-1-\Z_+}
\hbox{\ \ for\ \ } \rho^2+\rho+1=0,
\notag\\ 
&\widetilde{F}_4:\  &t_{\sht}^2\in -q^{-1-\Z_+}
\hbox{\ \ \ \, for \  } -\al=0122,\ 
\notag\\ 
&\widetilde{F}_4:\  & t_{\lng}\, t_{\sht}^{3}\in 
-q^{-1-\Z_+}
\hbox{\ \ as\ } -\al=1242.\ 
\notag 
\end{eqnarray}
In the case $-\al=\al_2+2\al_3$ for $\widetilde{F}_4$ from (i), 
$t_{\lng}\, t_{\sht}^{3}\in q^{-1-\Z_+}$
 does not result in $Rad\neq\{0\}$; however, if
it holds, then $[-0121,j]\in \tR_+^0$ for proper $j\in 1+\Z_+$.
\end{lemma}
\sq
\medskip

\subsection{Zigzag paths}
Following the last section,
we add more detail concerning the combinatorial
structure of $\tR^0_+$ and $\tR^{-1}_+$, the sets 
of positive roots
in  $\tR^0$ and $\tR^{-1}$. We need to make more
transparent the procedure from Lemma \ref{ZFROMM} and 
also provide the combinatorial tools for
Theorem \ref{THMIRRV} below.

The condition that $\tal'\in \tR^{-1}_+$ corresponds to
$\tbe'\in \tR^0_+$ (in the sense of (i)) has the following
general meaning. Always, $s_{\tbe'}(\tal')$ belongs to
$\tR^1$ (it holds for any $\tw\in \tW^0$). Then 
$\tbe'=[\be',\nu_{\be'} j]$ in (i)
are precisely those satisfying:
\begin{align}\label{tbeprimeal}
&(a): \ s_{\tbe'}(\tal')=-\al_l\ (l>0),\ \, -\be'\in R_+,\, j>0,\\ 
&(b): \ \tal'=\hbox{short}\,\Rightarrow\,
\tbe'=\hbox{short}.\notag
\end{align}

To avoid misunderstanding,
we use here the notations $\be'$ and $\al'$ for the non-affine
components of $\tbe'$ and $\tal'$ instead of $\be$ and
$\al$ used in (i).
  
The set of all $\al_l$
satisfying (\ref{tbeprimeal}) for some 
$\tal'=[\al',\nu_{\al}j]\in \tR^{-1}_+[-]$,
i.e., with {\em negative\,} non-affine components $\al'$, 
is described as follows: 
$$
(a'):\ (\al_l,\tbe')=(\al_l,\be')>0, \ \,
(b'):\  \be'=\hbox{long}\,\Rightarrow\,
\al_l=\hbox{long}.
$$
Always $s_l(\tbe')=\tbe'-\al_l$ due to ($a,b$), but 
$s_l(\tbe')$ may be different from 
$\tal'=-s_{\tbe'}(\al_l)$ in the non-simply-laced case.

If ($a',b'$) hold, then we  connect $\tal'$ and $\tbe'$
by a {\em link}. 
We say that  different $\tbe',\tbe''\in \tR^0_+[-]$ are
{\em $\tR^0$\~neighbors\,} 
if they are connected with the 
same $\tal'\in \tR^{-1}_+[-]$ by {\em links\,};
similarly, $\tal',\tal''\in \tR^{-1}_+[-]$ are
called {\em $\tR^{-1}$\~neighbors\,} if they  correspond
to different $\al_{l}$ and $\al_{m}$ and the same $\tbe'$.
By a {\dfont zigzag\,}, we mean a {\em connected}
subgraph  of the graph of all roots from
$\tR^{-1},\tR^0$ taken as the vertices with the 
{\em links\,} considered as the edges; generally,
it may be a tree or a loop.

Note that if $\tal',\tal''$ are {\em neighbors\,}
in $\tR^{-1}_+[-]$, then
\begin{align}\label{allalm}
\{\,(\al_l\,,\,-\be')\,<\,0\,>\,(\al_m\,,\,-\be')\,\}, 
\hbox{\ which \ implies\ }
(\al_l,\al_m)\,=\,0;
\end{align}
recall that  $-\be'\in R_+$.
Indeed, if $\al_l$ and $\al_m$ are connected
in the Dynkin (non-affine) diagram $\Ga$ and
$\nu_l\ge \nu_m$, then
$s_m(\al_l)=\al_l+(\nu_l/\nu_m)\al_m$ belongs to $R_+$ and 
$$
s_{-\be'}(\al_l+(\nu_l/\nu_m)\al_m)=
\al_l+(\nu_l/\nu_m)\al_m+ M(-\be'),
$$
where $M=2(\nu_l/\nu_m)$ as $\nu_{\be'}=\nu_m$
and $M=1+(\nu_l/\nu_m)$ as $\nu_{\be'}\ge \nu_l$.
However, when $-\be'$ is long 
or when all three roots are of the same length,
then $M<2$, which contradicts to the above formulas
for $M$.

For instance, the orthogonality from (\ref{allalm}) gives
that the number of {\em links\,} from  
a given $\tbe'\in \tR^0_+[-]$
cannot be greater than the maximum degree of the 
vertices in $\Ga$ considered as a tree (ignoring the
multiplicities of laces), which is $2$ unless for $D,E$.
\smallskip

Actually, the positivity $-\be'>0$ was used above only in the
last inequality $M<2$. For negative $-\be'$ 
satisfying the inequalities from (\ref{allalm}),
$M$ can be $2$ for neighboring 
$\al_l,\al_m$ 
if $-\be'=-\al_l-\al_m$ and $\al_l,\al_m$
are of the same length. We obtain an inversion of (\ref{allalm}):
\begin{align}\label{sallm}
&(\al,\al_l)\,>\, 0\,<\, (\al,\al_m) \hbox{\ \,  for } \al\in R_+\,
\Rightarrow\notag\\
&(\al_l,\al_m)=0 \hbox{\ \ unless\ } \al=\al_l+\al_m
\hbox{\ and\ }\nu_l=\nu_m.
\end{align}
See (\ref{all+alm}) below. As a matter of fact,
(\ref{allalm}) and (\ref{sallm}) are
standard facts on subsystems of affine root system;
the latter can be deduced from the former by setting 
$-\be'=[-\al,\nu_\al]$. 
\medskip

{\bf Restricted zigzags}.
We will assume that $t_{\nu}\not\in q^{-\nu(1+\Z_+)}$
in the next lemma for any $\nu=\nu_{\sht},\nu_{\lng}$; 
the cases $[-1220,2j]\in \tR^0$
and $[-0121,j]\in \tR^0$ for $\widetilde{F}_4$ 
will be omitted. Let us also exclude the cases in  
(\ref{reducases}) or, more generally, 
in (\ref{treducases}) and 
(\ref{reducasesze}). We can now make ($b$) 
more restrictive:
\begin{align}\label{nube=nual}
&(b''):\ 
\nu_{\be'}=\nu_{\al'}=\nu_l \and \tal'=s_l(\tbe')=\tbe'-\al_l.
\end{align}
Respectively, we will use {\em restricted} links and zigzags. 
Note that one can find $\al_l$ to ``go" from {\em any} 
$\tal'$ to $\tbe'$ 
under ($b''$) since we excluded the cases 
where it does not hold. 
Also, 
$((\tbe')^\vee,\rho)=((\tbe'')^\vee,\rho)$ if $\tbe'$ and
$\tbe''$ are $\tR^{0}$\~{\em neighbors\,} and 
\begin{align}\label{tbereduced}
&\tbe''=\tbe'+\al_{l'\,}-\al_l,\ \be''=\be'-\al_{l'\,}+\al_{l'\,},\ 
\nu_l=\nu_{l'\,}
\end{align}
due to (\ref{allalm}),(\ref{nube=nual}),
where $\tbe'-\al_l=\tal'=\tbe''-\al_{l'\,}$.
Using (\ref{sallm}): 
\begin{align}\label{all+alm}
&(\al_l,\al_{l'\,})=0 \hbox{\ \, if\ \,}
-\al'\neq \al_l+\al_{l'\,} \hbox{\ \, for\ neighboring\ \, } 
\al_l,\al_{l'\,}.
\end{align}
We will not use this fact but it is helpful in
understanding the combinatorial structure of $\tR^{0,-1}$.

We mention that these definitions resemble
the BGG\~resolutions. Namely, 
$s_ls_m(\tbe')=\tbe'-\al_l-\al_m$
is {\em link}\~connected with $\tal'$ and $\tal''$, and
$s_ls_{l'\,}(\tal')=\tal'-\al_l-\al_{l'\,}$ is connected
with $\tbe'$ and $\tbe''$ in the corresponding cases.
This is an example of the BGG\~type squares.
\medskip

We see that the {\em zigzag-connectivity\,} is essentially 
a non-affine notion (if $q$ is not a root of unity)
with close relations to the classical theory, 
especially under (\ref{nube=nual}). 
Given $p\in 1+\Z_+$, let
$$
\tR_+^{\,0,\pm 1\,}(p) \hbox{\ \ be the set of all\ \ } 
\tga'=[\ga',p\,]\in \tR^{\,0,\pm 1\,}_+[-].
$$

\begin{lemma}\label{COMBR01}
Provided (\ref{nube=nual}), a restricted zigzag belongs
to the set $\tR^0(p)\cup \tR^{-1}(p)$ for 
proper $p\in 1+\Z_+$. If it is maximal, i.e., 
not a part of a larger zigzag, then it
contains a \, \underline{boundary}\,
root from $\tR^{0}$ that has a unique 
link to $\tR^{-1}$ by definition.
\end{lemma}
{\em Proof.}
If $\nu_{\lng}k_{\lng}=k_{\sht}$ (that includes the
simply-laced case), the non-affine components of the 
roots from $\tR_+^0(p)$ constitute a
subset $-R_{ht}$ of $R_-$ of the roots $\be$ of fixed 
height $ht\equal(-\be',\rho^\vee)$. Similarly, the set
$\tR_+^1(p)$ will be $-R_{ht+1}$. 
Note that the zigzags with $\be'$\~roots of different 
length do not intersect. 

The claim of the lemma
in this case becomes a simple statement
about {\em zigzags\,} in 
$-(R_{ht}\cup R_{ht+1})$, which is straightforward. 
As a matter of fact, this statement 
is a combinatorial version of 
(\ref{tdegr}). Let us outline the deduction
of (\ref{tdegr}) from the claim of
the lemma considered for abstract  
{\em zigzags\,} in $-(R_{ht}\cup R_{ht+1})$ for an
arbitrary $ht$. 

We take a maximal 
{\em segment of $\tR^0$\~neighbors\,} in $-R_{ht}$, a longest 
chain of consecutive {\em $\tR^0$\~neighbors\,} there.
It contains all roots in $-R_{ht}$ unless in the 
cases $D,E$, where we need to add an additional
{\em segment\,} with the corresponding zigzag that
begins at $-R_{ht}$ and ends at $-R_{ht}$ too. 
The zigzag corresponding to the maximal segment
can be with one or two endpoints from $-R_{ht}$. 
Correspondingly, $|R_{ht}|-|R_{ht+1}|$ is $0$, $1$,
or $2$ in the case of $D_{even}$. 
it gives the required structure of the right-hand side 
of (\ref{tdegr}) since the denominator is known.  
It is quite possible that this approach to proving 
(\ref{tdegr}) is known, but we found no proper 
references.
\smallskip

The sets $\tR_+^0(p)$ become the greatest 
as $\nu_{\lng}k_{\lng}=k_{\sht}$ 
(it is the so-called case of equal labels). If the latter
constraint is not imposed, then we inspect the cases of $B,C$ 
directly and use that each $\tR_+^0(p)$ can contain 
only one or two different roots for $F,G$ (if not empty).
Cf. (\ref{tdegra}).
\sq
\smallskip

The claim of the lemma is likely to hold (essentially) 
without imposing (\ref{nube=nual}). We expect that
there must be a general reason for any maximal
{\em zigzag\,} in $\tR_+^0\cup\tR_+^{-1}$ 
to contain the {\em boundary} points from $\tR_+^0$,
which are the ones with only
one link to $\tR_+^{-1}$. 
\smallskip

The following corollary
seems of general nature too, but we will state 
it assuming that the radical is 
zero and avoiding the exceptional cases from
(\ref{reducases}) or, more generally, from
(\ref{treducases}) and (\ref{reducasesze}).
It can somewhat simplify
proving Theorem \ref{THMIRRV} when being used instead of 
Lemma \ref{COMBR01}, for instance,
Lemma \ref{LEMREFLEC} and 
Corollary \ref{CORTALPRIME} can be avoided. 
\smallskip

\begin{corollary}\label{ZIGZAG}
Provided (\ref{qtgeneric}), we assume that the radical
is zero apart from the cases listed in (\ref{reducases}).
In the setting with $q$, the cases from (\ref{treducases})
and (\ref{reducasesze}) must be excluded.
We also assume that $t\neq q^{-1-\Z_+}$ in the 
simply-laced case.

Then there exists $\pi_a$ for $a\in B$
and its reduced decomposition such that $\la(\pi_a)$ contains 
the sets $\tR^0_+,\tR^{\pm 1}_+$ and satisfies the following
$\tR^{-1}\rightarrow \tR^{0}$\~\underline{alternation condition}.
Given  any  
$\tal'\in \tR^{-1}_+$, the first $\tbe'\in \tR^0_+$
after $\tal'$  in the sequence $\la(\pi_a)$ 
satisfies $(\tal',\tbe')> 0$ and no 
$\tga\,'\in \tR^{\pm 1}_+$ can be found in $\la(\pi_a)$ 
between $\tal'$ and $\tbe'$
such that $(\tga\,',\tbe)\neq 0$. 
\end{corollary}
{\em Sketch of the proof.}
Here adding the condition $t\not\in q^{-1-\Z_+}$ in
the simply-laced case is important; 
respectively, $k\not\in -1-\Z_+$ if (\ref{qtgeneric})
is assumed.
 
The simplest counterexample to the claim of
this corollary without imposing this condition
is $\widetilde{D}_4$ (this
condition is not needed for $\widetilde{A}$). 
\smallskip

\rmk
We not that if $R$ is simply-laced and
$k\not\in -1-\Z_+$, then the condition 
$Rad=\{0\}$ implies that $R^0_+$ is a set of
pairwise orthogonal roots 
unless in the following two cases:
$$
E_8 \hbox{\ \ where\ \,} k=-g/11,\ -g/13\not\in \Z \hbox{\ \,as\ \,}
g\in \N
$$
(similarly, in the setting with $q$).
Actually the alternating sequences are not needed
in the proof of Theorem \ref{THMIRRV} when $\tR^0_+$ consists
of pairwise orthogonal roots; one can use part (ii) of
the Key Lemma below.   
\sq
\smallskip

Essentially the procedure is as follows.
We take a maximal zigzag path that begins
with a proper {\em boundary}
root from $\tR^0_+$ (with a unique
link to $\tR^{-1}_+$) and transpose all consecutive pairs 
$\tbe'\rightarrow\tal'$ in this path. Then the resulting
sequence can be made a part of some $\la(\pi_a)$\~sequence.  

{\em The case of $\widetilde{E}_8$}.
The following example demonstrates the procedure.
Let $\tR=\widetilde{E}_8$ and
$k=-g/13$, where  $g\in \N$ and $(g,13)=1$. Then
$$
\tR^0=\{[-\be,g],\ ht_\be=(\be,\rho)=13\}
\cup \{[-\be,2g],\ ht_\be=26\},
$$
$\tR^{-1}=\{[-\al,g],\ ht_\al=14\}
\cup\{[-\al,2g],\ ht_\al=27\}$ and  also
$\tR^{1}=\{[-\ga,g],\ ht_\ga=12\}
\cup\{[-\ga,2g],\ ht_\al=25\}$.
\smallskip

The maximal zigzag we need is (in the notation from \cite{Bo}): 
{\small
\begin{align}\label{e8zig}
&& 12&32111^0-&-12&32211-&-12&32210^0-&-12&33210-&\\ \notag
&&   &2       &   &2     &   &2       &   &2     &\\ \notag     
-12&33210^0-&-12&33211-&-12&32211^0-&-12&32221-&-12&22221^0. &\\
\notag   
   &1       &   &1     &   &1       &   &1     &   &1        &
\end{align}
}
The $\be$\~components of
the roots from $\tR^0$ are marked by ${}^0$, the links are $-\,-$; 
here the endpoints are {\em boundary roots}
(that have unique link-connections with $\tR^{-1}$).
The remaining roots are of height $26,27$:
{\small
\begin{align*}
 24&65321^0\   =  &12&32111^0\ \ +\  &12&33210^0, &\\
   &3             &  &2              &  &1        &\\      
 24&65421\ \   =  &12&32211\ \ \ +\  &12&33210^0. &\\
   &3             &  &2              &  &1        &\\      
\end{align*}
}
The corresponding $\al,\be$\~parts of the roots from the
required 
$\la$\~sequence intersected
with $\tR^{-1,\,0}$ can be made:
{\small
\begin{align}\label{e8la}
&&\{\ 12&22221^0,& 12&32211^0,& 12&32221,& 12&33210^0,& 12&33211, &\\ 
\notag
&&      &1       &   &1       &   &1     &   &1       &   &1      &\\ 
\notag
  24&65321^0,&
  24&65421,& 12&32210^0,& 12&33210,& 12&32111^0,& 12&32211\ \}. &\\ 
\notag
    &3     &
    &3     &   &2       &   &2     &   &2       &   &2     & 
\end{align}
}

\noindent
Recall that the ordering of the roots 
in $\la$\~sets is from right to left.

Here the roots from $\tR^1_+$ can be 
inserted between the
roots from $\tR^0_+$ and the next ones from $\tR^{-1}_+$.
Recall that $|\tR^1_+|=|\tR^0_+|$ because of the assumption
$Rad=\{0\}$.  As a matter of fact, only one root from the $5$
ones of height $12$ in $\tR_+^1$ can appear before 
the last two roots in (\ref{e8la}).

\comment{
Its nonaffine component is the negative of the root
$\left(
  \begin{array}{ccccccc}
    1 & 2 & 3 & 2 & 1 & 1 & 0 \\
      &   & 2 &   &   &   &   \\
  \end{array}
\right).$
It may appear in sequence (\ref{e8la}) only after
the root with the $\be$\~part   
$\left(
  \begin{array}{ccccccc}
    1 & 2 & 3 & 2 & 2 & 1 & 0 \\
      &   & 2 &   &   &   &   \\
  \end{array}
\right).
$
}
\smallskip 

Use formula (\ref{lambpi}) or Theorem \ref{INTRINLA},(ii) to justify
that the resulting {\em sequence} can be made 
the intersection
$\la(\pi_a)\cap \{\tR^{-1}_+\cup\tR^0_+\cup\tR^1_+\}$ 
for certain $a\in B$.

Using Proposition \ref{GEOMLA}, one can construct geometrically
some $\pi_a$ satisfying the alternation condition.
Let $b_0\in \CC^a, b_1\in \CC$ for the {\em basic} affine
and nonaffine Weyl chambers. Then 
$L\equal\{tb_1+(1-t)b_0,0\le t\le 1\}\subset \CC$.
Let $[-\al,j]\in \tR^{\,-1,\,0,\,1\,}_+$. Here $j=g,2g$.
We assume that $(b_1,\al)>j$ for all such roots;
automatically, $(b_0,\al)<1\le j$.
Then the ordering of the points
\begin{align}\label{b0b1la}
&t=t(\al)=(j-(b_0,\al))/(b_1-b_0,\al), 
\hbox{\ \, where\ \,} 0<t<1,
\end{align}
gives the ordering of
the corresponding roots $\,[-\al,j]\,$ in the $\la$\~sequence
associated with the segment $L$, which is a collection of the  
positive (all possible) affine roots such that their hyperplanes
intersect $L$. 
The condition $b_1\in \CC$
is necessary and sufficient to make the total set in the form
$\la(\pi_a)$ for certain $a\in Q$. 
\smallskip

For instance, setting $\ep=0.1/g$, 
the vectors $b_0=\ep\,\om_5$ \,
and\, $b_1=\ep\,(\om_4+\om_5+\om_6+\om_7)+2\ep\,\om_8+\om_2$\,
(the notation is from \cite{Bo}) result in the following 
$\tR^{-1}\to \tR^0$\~alternating sequence:
{\small
\begin{align}\label{e8geom}
&&\{\ 12&32210^1,& 12&32111^1,& 12&22211^1,& 12&33210^0,& 12&32211^0,
&\\ \notag
    &&  &1       &   &1       &   &1       &   &1       &   &1      
&\\ \notag
 12&33211,  & 
 11&22221^1,& 12&22221^0,& 12&32221,& 24&64321^1,& 24&65321^0, 
&\\ \notag
   &1       &           
   &1       &   &1       &   &1     &   &3       &   &3       
&\\ \notag
 24&65421,  & 
 12&32110^1,& 12&32210^0,& 12&33210,& 12&32111^0,& 12&32211\ \}, 
&\\ \notag
   &3       &   
   &2       &   &2       &   &2     &   &2       &   &2     
\end{align}
\\ }
where the negatives of the nonaffine components of the roots from 
$\tR^1_+$ are marked by ${}^{1}$. 
\smallskip
 
We note a relation of the alternating sequences
to the so-called {\em non-crossing partitions\,}.
The Coxeter element and its powers are useful 
for constructing required $\pi_a$ and their
reduced decompositions,
namely, the formula $w_0=c^{h/2}$ from \cite{Bo} for 
the Coxeter element $c=s_n\cdots s_2 s_1$
and the Coxeter number $h=1+(\rho,\vth)$ if 
$R$ is not of type $A$ and
$$ w_0=s_1(s_2 s_1) \cdots 
(s_{n-1} \cdots s_1) (s_n\cdots s_1)\for A_n.
$$
\sq

\comment{
Second, the case
of {\em long} $\al$ in $\widetilde{C}_n$ when
$\tbe'\in \tR^0$ does not follow from 
(\ref{qbetbe}) must be added to (ii),(iii).
Similarly, if $(s,2)=1$ but $(s,g)\neq 1$  for  $\widetilde{B}_n$,
then it may lead to Rad=$\{0\}$ and such cases must be
added to (iii).
The other cases from (ii) that are not listed in (iii)
can be reduced to $t_{\lng}^{-1}\in q_{\lng}^{1+\Z_+}$,
that replaces $k_{\lng}\in -1-\Z_+$, or lead to Rad$\neq\{0\}$;
so here there is no change in (iii). Similarly,
if  $t_{\lng}^{-3k_{\lng}-k_{\sht}}\in q_{\lng}^{1+\Z_+}$
then $[-1220,2j]\in \tR^0$ for $\widetilde{F}_4$. 
\sq
\medskip
For instance, the
condition $k_{\lng}\in -1-\Z_+$ will be transformed
to $t_{\lng}^{-1}\in q_{\lng}^{1+\Z_+};$
the radical can be non-zero in this case now.
It occurs when the {\em pure $k_{\lng}$\~part}
of the Poincar\'e polynomials or its translations
$k_{\lng}\mapsto k_{\lng}-\Z_+$ are zero. The corresponding 
conditions to be added to (iii) are as follows:
\begin{align*}
&\widetilde{B}_n:\ \prod_{m=0}^{n-1}\prod_{j=0}^\infty
\frac{1-q^{2(m+1)k_{\lng}+2(m+1)j}}
{1-q^{2k_{\lng}+2j}}\neq 0;\\
&\widetilde{C}_n:\ 
\prod_{j=0}^\infty(1+q^{2k_{\lng}+2j})\neq 0;\ \ \widetilde{G}_2:\ 
\prod_{j=0}^\infty(1+q^{3k_{\lng}+3j})\neq 0;\\
&\widetilde{F}_4:\  
\prod_{j=0}^\infty(1+q^{2k_{\lng}+2j})\neq 0 \hbox{ \ or\ }
\prod_{j=0}^\infty(1+q^{2k_{\lng}+2j}+q^{4k_{\lng}+4j})\neq 0.
\end{align*}
}  
\smallskip

\subsection{Key Lemma}
We suppose till the end of this section that 
$t_{\lng}\neq \pm 1 \neq t_{\sht}$. Here
the condition $\neq -1$  makes the sets 
$\tR^{\pm 1}$ non-intersecting; the condition
$\neq 1$ is necessary in 
the next Key Lemma.

For $\tal'\in R_+^{-1}[-]$ from (\ref{minustc}),
we determine $\overline{c}$ from 
$\pi_{\overline{c}}=s_i\pi_c$ and suppose that 
there exists $\tbe'\in \tR^0_+$ (see there) such that
$\tbe'\not\in \la(\pi_c)$, $(\tal',\tbe')>0$ and,
moreover, $-s_{\tbe'}(\tal')>0$. We take
the pair $b,i'$ satisfying (\ref{sumforal}) assuming that
\begin{align}\label{overlinepib}
&\pi_{\overline{b}}=s_{i'}\pi_b,\ 
\la(\pi_b)\not\ni\tbe'\in\la(\pi_{\overline{b}}) \and\\
\hbox{there is no \ }&\tga\,'\in R_+^{\pm 1}[-]\cup R_+^{0}[-]
 \hbox{\ between\ }
\tal' \and \tbe' \hbox{ in } \la(\pi_{\overline{b}}).\notag
\end{align}     

Recall the notations:
$\pi_b=\hw s_i\pi_c$ as $l(\hw)+l(\pi_c)=l(\pi_b)$. We
set $\pi_{d}=\hw(i') s_i\pi_c$ for 
$\hw(i')\equal\hw^{-1}s_{i'}\hw$,
i.e., $\pi_d$ is
obtained when the inverse of the reduced decomposition
of $\hw$ is added to that of $s_{i'}\pi_b$. One can assume
that the decomposition $\hw^{-1}s_{i'}\hw$ is reduced 
(otherwise $\pi_b$ can be diminished), in particular, 
$(\tal',\tga')\neq 0$ for all $\tga'$ between $\tal'$ and 
$\tbe'$ in the sequence $\pi_b.$ However we do not suppose
that $(\hw^{-1}s_{i'}\hw) s_i\pi_c$ is reduced.
We set $\widetilde{E}^\dag_d=\Psi_{\hw(i')}\Psi_{s_i}E_c$. 
\smallskip    

\newtheorem{keylemma}[theorem]{Key Lemma}
\begin{keylemma}\label{KEYTRIPLE}
(i) Let $t_{i}\neq 1$ and dim\,$\v_c=1$. In the notation above,
one of the following holds:\\
(a):\ $-s_{\tbe'}(\tal')\in \tR^1_+[-]$,\ \,  
(b):\ $-s_{\tbe'}(\tal')=\al_l\, (l>0)$,\ \,
(c):\ $-s_{\tbe'}(\tal')<0$.\\
In either case, $E'=\Psi_i\widetilde{E}^\dag_d$ 
is a nonzero $Y$\~eigenvector. If ($b$) holds, 
then $E'$ is proportional to $\widetilde{E}_c$. Otherwise,
$E'\in \v_{\overline{d}}$ for $\pi_{\overline{d}}=s_i\pi_{d}$; 
moreover, dim\,$\v_{\overline{d}}=1$ and
$\widetilde{E}_{\overline{d}}=E_{\overline{d}}$ in case ($a$). 

(ii) Generalizing, let us allow several pairs 
$\{\tbe^{v},\tal^{v}\}, v=1,\ldots, p$ 
in (\ref{overlinepib}) between $\tal'$ and $\tbe'$
connected in the same way as 
$\tal'$ and $\tbe'$ (Lemma \ref{ZFROMM}). We also
assume that $(\tbe^v,\tbe^{v'\,})=0$ when $v\neq v'$.
The polynomial $E'$ is constructed as in (i),
but now the singular intertwiners must be dropped when we go
back from $\tbe'$.  Then $E'$ is a nonzero 
$Y$\~eigenvector of the same weight as $E_c$.
\end{keylemma} 
{\it Proof.} Here 
$\widetilde{E}^\dag_d$ is {\em not} a  $Y_{\al_i}$\~eigenvector
in $\v_d$ due to $(\tal',\tbe')\neq 0$. In the transformation
to $E\mapsto E'$, we multiply $\widetilde{E}^\dag_d$ 
by $(\tau_+(T_{i})-t_i^{1/2})$; then we use
Corollary \ref{HATEE},(ii).
The key claim is that $\Psi_i\widetilde{E}_c=$ 
const$\widetilde{E}_c$ for a nonzero constant in case ($b$);
it holds because 
$(\tau_+(T_{i})-t_i^{1/2})
\widetilde{\Psi}^\dag_{\pi_d}(\Q_{q,t}1)$ 
belongs to $\v_c$ and therefore must coincide with this space
due to dim\,$\v_c=1$.
In the case of ($a$) (we will omit ($c$)): 
$$
\widetilde{\Psi}_{\pi_{\overline{d}}}=
\widetilde{\Psi}_i\widetilde{\Psi}^\dag_{\pi_d}, \ \
\hbox{dim\,}\v_c=1\ \Rightarrow\ 
(Y_{\al_i}-t_i)
\widetilde{\Psi}_i\widetilde{\Psi}^\dag_{\pi_d}(\Q_{q,t}1)=0.
$$

Let us sketch another approach to (i)
based on Theorem \ref{RANKTWO}.
We try to collect $\al_l,\ldots,\tbe',\ldots, \tal'$ together 
in $\la(s_i\pi_d)$. Generally, it is possible only for
admissible triples; however, here $\al_l$ is simple non-affine
and therefore can be moved using the Coxeter relations  
in $\la(s_i\pi_d)$ to the first position. At a certain
moment it will become next to $\tbe'$ and we can apply
Proposition \ref{INTRINLAP},(iii). Then we can change
the order of intertwiners corresponding to
$\{\al_l,\tbe',\tal'\}$ to the opposite, 
corresponding to $\{\tal',\tbe',\al_l\}$, 
thanks to Theorem \ref{PHIBRUHAT},($c$). 

Due to the appearance of $\al_l$, the total
chain of intertwiners for $s_i\pi_d$ will result in zero
for $\{\tal',\tbe',\al_l\}$ since the corresponding
chain of partial reduced decompositions will leave the set 
$\pi_{B}$ after $\al_l$. The difference of two expressions for
$\widetilde{\Psi}^\dag_{s_i\pi_d}(1)$ corresponding to the orderings
$\{\al_l,\tbe',\tal'\}$ and  
$\{\tal',\tbe',\al_l\}$ is 
const\,$\widetilde{E}_c$ for const\,$\neq 0$
(we use that $t_i\neq 1$).

Concerning (ii), the orthogonality $(\tbe^v,\tbe^{v'})=0$
makes it possible to proceed by induction. Compare with
Corollary \ref{HATEECOM}.
\sq  

\comment{
In the simply-laced case, we can begin with any
reduced decomposition of $\pi_{\rho}$, setting 
$(-\rho)=w_0\pi_{\rho}$ and taking ``natural"
reduced decompositions for $w^0$. Namely, we
set $w_0=c^{h/2}$ for 
the Coxeter element $c=s_n\cdots s_2 s_1$
and the Coxeter number $h$ for $D,E$ and 
$$ w_0=s_1(s_2 s_1) \cdots 
(s_{n-1} \cdots s_1) (s_n\cdots s_1)\for A_n.
$$
\sq
}
\medskip

Recall that the space 
$\v(-c_\#)^\infty$
is linearly generated by the polynomials
$\widetilde{E}_b$ for $b$
such that  $q^{b_\#}=q^{c_\#}.$
They are defined for reduced chains have nonzero
leading terms and are normalized
as follows:
\begin{align}
&\widetilde{E}_b-X_b\ \in\ \oplus_{a\succ b}\Q(q,t) X_a.\
\label{macdgenn}
\end{align}
The $\widetilde{E}_b$ here is a $Y$-eigenvector
if $\widetilde{E}^\dag_{b'}=0$ for all $b'=\hw'\llb 0 \rrb$, where
$\dag$ indicate that the {\em standard decompositions\,}
of $\hw'\in \tilde{\b}_o^0(\pi_b)$ are taken, possibly non-reduced.

Also recall that any $\v(-c_\#)^\infty$ contains a unique
$Y$\~eigenvector up to proportionality if $Rad=\{0\}$;
indeed, otherwise a proper linear combination of two
different eigenvectors from this space
would belong to the radical. This eigenvector must be in
the form $E_b$ for (a unique) $\pi_b\,<_0\, \pi_c$,
i.e., for $\pi_b$ obtained from a reduced decomposition 
of $\pi_c$ by crossing out some singular simple reflections.

For {\em primary} $b=b^\circ$, the space $\v(-b_\#)^\infty$
is one-dimensional and $\widetilde{E}_{b}=E_{b}$ is its
generator. Otherwise there would exist $a\succ b$ with the
same weight, which contradicts the definition of primary elements.
Recall the definition of the primary elements:
$q^{-b^\circ_\#}=q^{-b_\#}$ and $q^{-b^\circ_\#}\neq q^{-a_\#}$
for $a\succ b^\circ$.

\smallskip
Due to the condition $Rad=0$, we have the following
implications for any $\pi_c$ (all claims hold true for primary
elements):

\noindent
(a)\ \{reduced decompositions of $\pi_c$ contain no singular
reflections\} $\Longrightarrow$

\noindent
(b)\ \{no intertwiners $\tau_+(T_j)-t^{1/2}$ appear in the
corresponding $\Psi_{\pi_c}$\} $\Longrightarrow$

\noindent
(c)\ \{ dim$\, \v_b=1$ and $\widetilde{E}_b=E_b \}$
\smallskip

More generally, if $\Psi_i$ is not infinity and
not in the form $\tau_+(T_j)-t^{1/2}$ the map
\begin{align}\label{psivavc}
&\ \v(-c_\#)^\infty\,\to\, \v(-b_\#)^\infty
\hbox {\ \ for\ reduced\ }\pi_b=s_i\pi_c
\end{align}
is {\em injective}. It suffices to check that it does not
kill $Y$-eigenvectors. This results from the following
lemma, which is a simple application of Corollary \ref{COREBEC}.

\begin{lemma}\label{LEMEPRIM}
Let  $\Psi_i^c\neq \infty.$
We do not assume that
$l(s_i\pi_c)=1+l(\pi_c).$

(i) The polynomial $E'=\Psi_{i}(E)$ belongs to 
$\Rad$ for a $Y$-eigenvector $E\not\in Rad$
of weight $-c_\#$
if and only if
$\Psi_i^c$ is proportional to
$\tau_+(T_i)-t_i^{1/2},$ equivalently,
(see (\ref{interminust}),
\begin{align*}
&q_\al^{k_\al+(\tal^\vee,c_-+d)-(\al^\vee, \rho_k)}=1
\for \tal=u_c(\al_i).
\end{align*}

(ii) Provided that $t_\nu\neq -1$ for all $\nu,$
the intertwiner
$$
\Psi_i:\ V'_c\ \to\ V'_b,
\where V'_c=\widetilde{V}_c \hbox{\ modulo\ } Rad,
$$
is injective if $\Psi_i^c$ is proportional
to $\tau_+(T_i)+t_i^{-1/2}$ for  $i\ge 0.$
\sq\end{lemma}

Note that the demonstration of (i,ii) is immediate
if  $i>0$. Say, if
$\Psi_i^c=T_i+t_i^{-1/2}$, then
$\{E',1\}=(t_i^{1/2}+t_i^{-1/2})\{E,1\}\neq 0$, so obviously
$E'\not\in Rad.$
\smallskip

Let us assume that there exists $a\in B$ such
that $\tR^{-1}\subset \la(\pi_a)$ and 
$\tR^{0,1}\cap\la(\pi_a)=\emptyset$. It is always the
case for the simply-laced root systems;
use Proposition \ref{GEOMLA} and follow (\ref{b0b1la})
for $b_0,b_1$ proportional to $\rho$.
Then Lemma \ref{LEMEPRIM} gives that an arbitrary
$\HH^\flat$\~submodule $\v\,'$ of $\v$ subject 
to $Rad=\{0\}$
contains infinitely many eigenvectors $E_b$ such that
{\em no singular intertwiners appear in the
reduced chains for $E_b$,} i.e., 
$\tR^{0,1}\cap\la(\pi_b)=\emptyset.$ We do not need
this claim too much but it simplifies the considerations
below.   
\medskip

\subsection{Main Theorem}
We are now in position to formulate and prove
the main theorem of this part of the paper;
the notations from the previous section are used.

\begin{maintheorem}\label{THMIRRV}
(i) Imposing (\ref{qtgeneric}), let 
the radical $Rad$ of
the form $\{\,,\,\}$ be zero. Then the
polynomial representation $\v$ can be reducible
only in the following cases:
\begin{align}\label{bfgcondi}
\widetilde{B}_n:\  &k_{\lng}\ =\ -\frac{s}{2g},
\hbox{\ \ \ provided\ that\ \ } (s,2g)=1,\\
&\hbox{where\ \ } n>g>n/2,\ \  g\in 1+\Z_+ \ni s, 
\notag\\
\widetilde{G}_2:\  &k_{\lng}\ =\ -s/3,
\hbox{\ as\ \ } s\in 1+\Z_+,\ (s,3)=1,
\notag\\
\widetilde{F}_4:\  &4k_{\lng}\in -1-2\Z_+\hbox{\ \ or\ \ }
6k_{\lng}+2k_{\sht}\in -1-2\Z_+.\notag
\end{align}
If $k_{\sht}$ is generic in this list, then $Rad=\{0\}$ and,
indeed, $\v$ is a reducible $Y$-semisimple $\HH^\flat$\~module; it
remains reducible for any $k_{\sht}$, however the radical may
become nonzero.

(ii) Without imposing (\ref{qtgeneric}), the 
the above
conditions from (i) become:
\begin{eqnarray}\label{conditbfg}
&\ \ \ \widetilde{B}_n:\   t_{\lng}^{-g}\in 
q^{1+\Z_+}\setminus q^{2+2\Z_+},\ \ 
\widetilde{G}_2:\  
t_{\lng}^{-1}\in q^{1+\Z_+}\setminus q^{3+3\Z_+}, 
\\
&\ \ \widetilde{F}_4:\  t_{\lng}^{-2}\in 
q^{1+\Z_+}\setminus q^{2+2\Z_+} \hbox{\ \,or\ \ }  
t_{\lng}^{-3}t_{\sht}^{-2} \in
q^{1+\Z_+}\setminus q^{2+2\Z_+},\notag 
\end{eqnarray}
where $n/2<g<n$ and $t_{\lng}^{-g'}\not\in q^{1+\Z_+}$ for 
any $g'$ such that\, $\Z_+\ni g'<g$\,  in the case of 
$\widetilde{B}_n$. 
If $t_{\sht}$ is generic in this list, then $Rad=\{0\}$ and
$\v$ is a reducible $Y$-semisimple $\HH^\flat$\~module.

(iii) Continuing (ii),
the remaining cases when $Rad$ is $\{0\}$ but 
$\v$ may be reducible are those from Lemma \ref{LEMNUALBE},(ii): 
\begin{eqnarray}\label{reducasesthm}
&\widetilde{C}_n:\  &t_{\sht}^{-g}\in -q^{\Z_+},
\hbox{\ \ where\ \ } 
n>g>\frac{n}{2}\\
&&t_{\sht}^{-g'}\not\in \pm q^{\Z_+}
\hbox{\ \ \, for\ \ \,} \Z_+\ni g'<g,
\notag \\
&\widetilde{G}_2:\  &t_{\sht}^{-1}\in \rho\, q^{\Z_+}
\hbox{\ \ for\ \ } \rho^2+\rho+1=0,
\notag\\  
&\widetilde{F}_4:\  &t_{\sht}^{-2}\in -q^{\Z_+}\hbox{\,\  or\ \,}
t_{\lng}^{-1}\, t_{\sht}^{-3}\in 
-q^{\Z_+}.
\notag
\end{eqnarray}
If $t_{\lng}$ is generic in this list, then $Rad=\{0\}$ and
$\v$ is a reducible $Y$-semisimple $\HH^\flat$\~module.
\end{maintheorem}

{\em Proof.}
We will assume that $t_{\nu}\not\in q^{-\nu(1+\Z_+)}$
for any $\nu=\nu_{\lng},\nu_{\sht}$ respectively
for $B,F,G$ (then $\nu=\nu_{\lng}$) and $C,F,G$
(then $\nu=\nu_{\sht}$). 
Technically, we will need to exclude the case 
$t\in q^{(1+\Z_+)}$ for the simply-laced $R$ too;
it will be considered separately.
 
We also need to consider separately the cases 
$[-1220,2j]\in \tR^0$ and $[-0121,j]\in \tR^0$ 
for $\widetilde{F}_4$. Without going into detail,
the irreducibility of $\v$ can be deduced 
directly from Key Lemma \ref{KEYTRIPLE} 
in these two cases.

Thus, apart from these cases and those
listed in (\ref{bfgcondi}) or, more generally, in
(\ref{conditbfg},\ref{reducasesthm}), 
let us suppose that $\v$ has a proper $\HH^\flat$\~submodule
$\v'$.
\smallskip

Using Corollary \ref{ZIGZAG}, we find the 
portion $s_i\pi_c$ of the reduced decomposition
of $\pi_a$ constructed there
such that $E_{c^\circ}\not\in \v'$ but
$E_{\overline{c}^\circ}\in \v'$
for $\pi_{\overline{c}}=s_i\pi_c$. Note that 
we examine {\em primary} $c=c^\circ$ 
($E_{\overline{c}}$ may remain outside $\v'$).
We use that $\v'$ contains all $E_b$ with
sufficiently large $\la(\pi_b)$.

Then we use Key Lemma \ref{KEYTRIPLE},(i) and find
$\pi_b=\hw s_i\pi_c$ as $l(\hw)+l(\pi_c)=l(\pi_b)$ and
$\pi_{d}=\hw(i')s_i\pi_c$ for
$\hw(i')\equal\hw^{-1}s_{i'}\hw$. The corresponding
intertwiner $\Psi_{\hw(i')}$ will be applied to
$E_{c^\circ}$. The result can be either proportional
to $E_{c^\circ}$, which contradicts to the assumption
$E_{c^\circ}\not\in \v'$ , or can be another (nonzero)
$Y$\~eigenvector for the same eigenvalue, 
which contradicts to $Rad=\{0\}$. It gives the required.
\smallskip

If claim (ii) of Key Lemma \ref{KEYTRIPLE} is used here,
then the combinatorial part (Corollary \ref{ZIGZAG})
can be reduced significantly. Namely, we do not need 
to use $\pi_a$ with the {\em strict alternation} 
$\tR^{-1}\to\tR^{0}$; weaker conditions
are sufficient. 
\smallskip

We note that Key Lemma \ref{KEYTRIPLE} can be generalized
in many ways; the following its variant is needed to manage
the case $k\in -1-Z_+$. Let $k=-g\in -1-Z_+$,
more precisely, $t=q^{-g}$; for the
sake of definiteness only the simply-laced case will be
discussed below and we will assume that $\rho\in B$.
Then $\la(\pi_{a_+})$
for the element $a=a_+=g\rho\in B_+$ contains $\tR_+^{-1}$ and
also $\la(\pi_{a_+})\cap \tR_+^{0}=\emptyset$; see
(\ref{lambpi}). If $a_-=-g\rho$\, is taken, then 
\begin{align}\label{lagrho}
\tR_+^{0}\subset \la(\pi_{a_-})\setminus \la(\pi_{a_+}),\ 
\la(\pi_{a_-})\cap \tR_+^{1}=\emptyset.
\end{align}
Using the chain of intertwiners for 
$\pi_{a_-}$ followed by 
its reverse, where the singular intertwiners
are omitted, we obtain a {\em nonzero}
eigenvector with the same eigenvalue as for $E_0=1$.
It is a variant of Key Lemma \ref{KEYTRIPLE},(ii);
the verification is similar 
to the considerations of Proposition \ref{VINDUCEDY}.
It gives the irreducibility of $\v$
for $t=q^{-g}$, which is in a sense
``the most non-semisimple" case. 
\smallskip

Let us outline a justification of the theorem based 
directly on Lemma \ref{COMBR01}, without using a ``global"
$\pi_a$ satisfying the {\em alternation property}.
\smallskip
 
{\em Step 1.}
Similar to the above consideration, one can find 
$E_{\overline{c}}$ in $\v'$  for 
{\em primary} $\overline{c}=\overline{c}^\circ$ 
such that $\pi_{\overline{c}}=s_i\pi_c$ for {\em primary}
$c=c^\circ$ and the index $i\ge 0$ satisfying 
(\ref{minustc}) from Lemma \ref{ZFROMM} and with the
corresponding $E_{c}$ {\em not} in $\v'$.
Indeed, let $\overline{c}$ be {\em primary} such that 
$\pi_{\overline{c}}=s_i\pi_c$ (the decomposition is reduced)  
for some $c$ satisfying $\widetilde{E}_c\not\in \v'$. 

We may  assume that 
$l(\pi_{\overline{c}})$ is  minimal possible
provided that $\widetilde{E}_c\not\in \v'$.
If $c$ is primary, then
Lemma \ref{LEMEPRIM} gives that the intertwiner
$\widetilde{\Psi}^c_i$ can be only of type
$(\tau_+(T_i)+t_i^{-1/2})$. 

If $c$ is not primary then
$\widetilde{E}_c$ is {\em not} a $Y$\~eigenvector. Then
$\widetilde{\Psi}^c_i$ is proportional
to $(\tau_+(T_i)-t_i^{1/2})$. Taking 
$\widetilde{\Psi}^c_i$
is not sufficient here since
$$
\widetilde{\Psi}_i\ =\ \Psi_i\ =\ 
(\tau_+(T_i)+
\frac{t_i^{1/2}-t_i^{-1/2}}{Y_{\al_i}^{-1}-1})
$$ 
involves $Y_{\al_i}$ when applied to  
$\widetilde{E}_c$. One has: 
\begin{align*}
&(\Psi_i)^2= 
\frac{(t_i^{1/2}Y_{\al_i}^{-1}-t_i^{-1/2})
(t_i^{1/2}Y_{\al_i}-t_i^{-1/2})}
{(Y_{\al_i}^{-1}-1)(Y_{\al_i}-1)} \and\\
&\widetilde{\Psi}_i^2(\widetilde{E}_c)\,=\, 
\widetilde{\Psi}_i(E_{\overline{c}})\,=\,
\hbox{const\,}E_{c^\circ} \for
\hbox{const}\neq 0.
\end{align*}
Recall that $E_{c^\circ}$ generates the space of
$Y$\~eigenvectors in $\v_c$. Therefore $E_{c^\circ}\in \v'$,
which contradicts to the minimality of $l(\pi_{\overline{c}})$. 
\smallskip


{\em Step 2.}
Now let $\tbe',i',\pi_b$ be from
Lemma \ref{ZFROMM} applied to $c=c^\circ, i$. 
We take a reduced decomposition
of $\pi_b$ extending that 
of $\pi_c$. We will also assume that 
$l(\pi_{b})$ is the minimal possible
such that $\widetilde{E}_c\not\in \v'\ni \widetilde{E}_b$.  
Then there are no 
{\em singular roots\,} $\tbe''\in \la(\pi_b)\cap \tR^0_+[-]$ 
after $\pi_c$ in $\pi_b$ of the same type as $\tbe'$;
it is straightforward to check that other singular roots 
(if any) in  $\pi_b$ can be removed by diminishing $\pi_b$.

Recall that applying $\widetilde{\Psi}$\~intertwiners 
along a given reduced 
decomposition of any element from $\pi_B$
is always nonzero. Therefore, 
$\Psi_{\hw}(E_c)=E_b$ for $\hw=\pi_b \pi_c^{-1}$,
and $b$ must be primary too due to $Rad=\{0\}$.
Moreover, intertwiners of type $(\tau_+(T)-t^{1/2})$
cannot appear in $\Psi_{\hw}$ because their ({\em nonzero})
images would belong to $Rad$ (Lemma \ref{LEMEPRIM}).
\smallskip

{\em Step 3}.
The only situation that may prevent us from using 
Key Lemma \ref{KEYTRIPLE} is the appearance 
of new $\tal''\in \la(\pi_b)\cap \tR^1_+[-]$ associated to 
the same $\tbe'$ (the 
last root in $\la(s_{i'}\pi_b)$) 
in the same way as $\tal'$. Recall that
any such $\tal''$  must appear {\em before} $\tbe'$
in any reduced decompositions of the elements from
$\pi_B$.

Due to (\ref{allalm}): 
$$
\tal'=-\tbe'(\al_l),\ \tal''=-\tbe'(\al_m),\ \where
(\al_l,\al_m)=0.
$$ 
Applying Corollary \ref{CORTALPRIME}, one
can make $\tal''$ {\em before}
$\tal'$ in $\la(\pi_b)$ using the Coxeter
transformations in this set. The procedure is
as follows. 
\smallskip

We can assume that $\hv=\,s_i\,\hw^{-1}s_{i'}\,\hw s_i\pi_c\,$
is reduced by diminishing $\pi_c$ if necessary;
here $\pi_b=\hw\pi_c$ (cf. Lemma \ref{KEYTRIPLE}). 
Since $\al_l$ is from $R_+$, the triple
$\{\,\al_l,\tbe',\tal'\,\}$  
can be made consecutive using the Coxeter transformations
in the sequence $\la(\hv)$. Using
Main Theorem \ref{RANKTWO}, the set $\la(\hv)\cap R_+'$
must contain at least two simple roots from $R_+'$
for an arbitrary subsystem $R'\supset$
$\{\,\al_l,\tbe',\tal'\,\}$ 
of type $B_3,C_3$ or $D_4$ in this theorem 
($\tbe'=\ep_1+\ep_2$ for $B,C_3$ and $\tbe'=\th^4$ for $D_4$).  

Using Coxeter transformations {\em inside} $\pi_b$,
one makes these roots {\em after}
$\tal'$, changing the position of $\tal'$ if
necessary. 
One of them (at least) must be orthogonal to $\tbe'$ and 
therefore will ``disappear" in the corresponding $\hv$, 
which contradicts to the minimality of 
$l(\pi_b)$ imposed above.

We conclude that $\{\,\al_l,\tbe',\tal'\,\}$  can be
assumed satisfying Corollary \ref{CORTALPRIME},(iii).
It gives that $\tal'\ $ can be transposed with $\tal''$
(or with all such $\tal''$ if there are several of them).
Note that since dim\,$\v_c=1$ 
and dim\,$\v_b=1$, applying Coxeter transformations 
will not change the corresponding $E$\~polynomials
(generally, it may influence the 
$\widetilde{E}$\~polynomials). 
\smallskip

If $E_c$ remains {\em not} in $\v'$, then it
concludes Step 3. Otherwise we switch to $\tal''$
and proceed by induction using the {\em links\,}
for creating a 
{\em zigzag} in $\tR^0\cup \tR^{-1}$  in the terminology of 
Lemma  \ref{COMBR01}. Eventually, we will make
this {\em zigzag} maximal. Then it will
contain a {\em boundary}
root $\tbe^*\in \tR^0_+[-]$ (with a unique link
to $\tR^{-1}_+$) and there will 
be no roots from $\tR^1_+[-]$ between
$\tbe^*$ and the previous $\tal^*\in \tR^{-1}_+[-]$.
\smallskip

{\em Step 4}.
In the absence of $\tal''\in \la(\pi_b)\cap \tR^1_+[-]$
after $\tal'$, one may apply 
Key Lemma \ref{KEYTRIPLE},(i). 
Let us recall the construction.
Setting $\pi_b=s_{i_u}\cdots s_{i_{1}}\pi_c$ for
$i=i_{1}$ and $\pi_{\overline{b}}=s_{i'}\pi_c$,
we consider the elements 
$\pi_d= s_{i_2}\cdots s_{i_u}s_{i'}\pi_b$ and
$\pi_{\overline{d}}= s_{i}\pi_d$. 
One has $\pi_{\overline{d}}=\pi_c s_{\tbe'}$ with
$s_{\tbe'}\in \tW^0$ (since $\tal'\in \tR^0$).

We do not claim
that the product for $\pi_{\overline{d}}$
is a reduced decomposition but the 
reductions can be only with the reflections 
corresponding to invertible intertwiners.
Such transformations are acceptable for this proof.

Here $\pi_{\overline{d}}$ belongs to $\pi_B$ 
if the last root in $\la( \pi_{\overline{d}})$
does. If it is the case, we obtain a $Y$\~eigenvector
$E_{\pi_{\overline{d}}}$ non-proportional to
$E_c$ and with the same $Y$\~eigenvalue; this
contradicts to $Rad=\{0\}$. Otherwise,
the last root in $\la( \pi_{\overline{d}})$
is in the form $\al_l$ for $l>0$ and
Key Lemma \ref{KEYTRIPLE} shows that 
the corresponding chain of intertwiners 
results in $E_c\in\v'$, a contradiction. 
This contradiction concludes this step and 
( a sketch of the) the proof
irreducibility of $\v$ 
without using Corollary \ref{ZIGZAG}.
Recall that $Rad=\{0\}$ and  
the exceptional cases
listed in (\ref{conditbfg}) and (\ref{reducasesthm})
were excluded.
\smallskip

We will check the claim about reducibility 
of $\v$ when the $t$\~parameters
are from these lists (for generic $t_{\sht}$
or $t_{\lng}$) in the cases of $\widetilde{B,C}_n$. 
The remaining three cases are very much similar.

Let us impose (\ref{conditbfg}) and (\ref{reducasesthm})
(see also (\ref{treducases}) 
and (\ref{reducasesze})) and verify that $\v$ is 
reducible (and semisimple) though $Rad=\{0\}$. 
This consideration extends the example of
$\widetilde{B}_n$ considered in \cite{C100}.
\smallskip

{\em Example of $\widetilde{B}_n$.}
We assume that $k_{\sht}$ is generic. 
According to (\ref{yofbfinal}), the condition
$t_{\lng}^{2m}\not\in q^{-2-2\Z_+}$ for $1<m\le n$ is
necessary and sufficient for $Rad=\{0\}$.
In this case, $\tR^0_+[-]=\emptyset=\tR^{1}_+[-]$;
therefore, $\v$ is semisimple.

The set $\tR^{-1}_+[-]$ is nonempty if and only if
$t_{\lng}^{g}=q^{-s}$ for $1\le g<n,\,s\in 1+\Z_+$.
Namely, the short root $[-\ep_{n-g},s]$ belongs to
this set: $((-\ep_{n-g})^\vee,\rho_k)+k_{\sht}=-2gk_{\lng}$
(see (\ref{yofbfinal})). Setting
$t_{\lng}=\exp(\frac{2\pi \imath j}{g})q^{-s/g}$
for $s\in 1+\Z_+,\, j\in \Z_+$, if $Rad=\{0\}$ then:
$$
(s,2)=1 \and t_{\lng}^{g'}\not\in \pm q^{-1-\Z_+} 
\for 0<g'<n,\, g'\neq g,
$$ 
in particular, $g>n/2$.

These conditions give that $Rad=\{0\}$, $\v$ is semisimple
and at least one non-invertible intertwiner exists,
which imply the reducibility of $\v$.
For instance, either the relation $(j,g)=1$ (any odd $s>0$) 
or the relation $(s,2g)=1$ (any $j\ge 0$) is sufficient as
$n>g>n/2$. Generally, $(2j,s,2g)=1$ as $n>g>n/2$
is necessary and sufficient. The simplest example is
$n=3,g=2,s=1,j=0$.
\smallskip

{\em Example of $\widetilde{C}_n$.}
We assume that $k_{\lng}$ is generic and present this
case following the previous one to emphasize that they
are dual to each other (under 
$q\leftrightarrow q'$\~duality).
According to (\ref{yofxevalc}), the condition
$t_{\sht}^{m}\not\in q^{-1-\Z_+}$ for $1<m\le n$ is
necessary and sufficient for $Rad=\{0\}$.
In this case, $\tR^0_+[-]=\emptyset=\tR^{1}_+[-]$
and $\v$ is semisimple.

The set $\tR^{-1}_+[-]$ is nonempty if and only if
$t_{\sht}^{2g}=q^{-2s}$ for $1\le g<n,\,s\in \Z_+$.
Namely, the long root $[-2\ep_{n-g},2s]$ belongs to
this set: $((-2\ep_{n-g})^\vee,\rho_k)+k_{\lng}=-gk_{\sht}$.
Setting
$t_{\sht}=\exp(\frac{2\pi \imath j}{2g})q^{-s/g}$
for $s\in \Z_+,\, j\in 1+\Z_+$, if $Rad=\{0\}$ then:
$$
(j,2)=1 \and t_{\sht}^{g'}\not\in \pm q^{-1-\Z_+} 
\for 0<g'<n,\, g'\neq g,
$$ 
in particular, $g>n/2$ must hold.

For instance, either the relation $(j,2g)=1$ (any $s\ge 0$) 
or the relation $(s,g)=1$ (any odd $j\ge 1$) is sufficient. 
Generally, the condition $(j,2s,2g)=1$ as $n>g>n/2$
is necessary and sufficient
for reducibility of $\v$ when $Rad=\{0\}$; $\v$ is
semisimple in this case (for generic $k_{\lng}$).
The simplest example is $n=3,g=2,s=0,j=1$.
\medskip

To recapitulate, let us recall that  
Main Theorem \ref{RANKTWO} and
Corollary \ref{CORTALPRIME} are not needed
if Corollary \ref{ZIGZAG} is used.
Moreover, claim (ii) of Key Lemma \ref{KEYTRIPLE} 
can simplify the combinatorial part of the proof.
\sq
\medskip

\rmk
\comment{
(i) Strictly speaking, we need to assume 
in (ii),(iii) that 
$t_{\al}\not\in \pm q_{\al}^{-1-\Z_+}$
to use the Key Lemma and related tools. However, due 
to the condition $Rad=\{0\}$ and other constraints there, 
if $t_{\sht}q^{1+l}=\pm 1$ for $l\in \Z_+$, then 
the corresponding intertwiner is calculated at a 
{\em long} root (at a short root for $t_{\lng}$).
We really need to avoid the cases when 
$t_{\sht}q^{1+l}=\pm 1$ appears when the
intertwiners are evaluated at short roots (similar
for long roots).
Thus, $t$ can be
$q$\~translations of roots of unity, provided
that they are not $q$\~translations of
the zeros of the corresponding Poincar\'e 
polynomial (the part of the product $\Pi_{\tR}$ without $q$).
The radical is nonzero at the latter. 
} 
(i) As a mater of fact, this proof can be used
to manage a more general problem of the irreducibility
of the quotient $\v/Rad$. We need a proper generalization
of Corollary \ref{ZIGZAG} and the Zigzag Lemma \ref{COMBR01}.
Note that the roots from $\tR^1_+[-]$ can now appear between 
$\tal'$ and $\tbe'$ connected by a {\em link}. 
The following observation is 
more or less sufficient to manage this problem. 
The roots from $\tR^1_+[-]$ represented in the
form $s_{\tbe'}(\al_l)$ cannot appear in these intervals;
such roots may occur only {\em after}
$\tbe'$ in {\em any} $\la(\pi_b)$. 
A natural expectation is that if  $\v/Rad$ is 
reducible then it has a $Y$\~semisimple quotient, 
although we do not have any confirmations beyond the case 
$Rad=\{0\}$.

(ii) More generally,
the same method is expected to help with checking the
irreducibility of constituents of $\v$. In the simply-laced
case, one can expect that a chain of intertwiner 
can ``enter" submodules of $\v$ only at the places where the
simple intertwiners become of type $\tau_+(T)-t^{1/2}$. 
If the latter intertwiners can be avoided,
generally, such chain can be expected to
remain in the same constituent
of $\v$. Here $\pi_a$ must be chosen similar to that 
in Corollary \ref{ZIGZAG})
or a proper variant of Lemma \ref{COMBR01} must be used.  
A natural step in this direction is a direct proof 
(without any reference to the localization functor) of 
the Kasatani conjecture; see \cite{Ka},\cite{En}.
\sq
\medskip

An important motivation of the constructive methods
we study in this paper is that, generally,
intertwiners are helpful for the theory of
{\em square integrable, tempered} and similar irreducible 
AHA\~modules. 
The analytic properties of the intertwiners acting
in a particular module are directly related 
to its analytic type. See, e.g.,\cite{O5,O11};
for instance, the so-called residual sets are directly
related to $\tR^{0,\pm 1}$ considered in our paper.
Here an explicit description in terms of intertwiners 
is needed, similar to what we do.

Paper \cite{MTa} and some other related papers 
indicate that the classification
of the square integrable AHA\~modules can be made
sufficiently explicit for the classical root systems.
The combinatorics involved and the methods from \cite{MTa}
employed are sophisticated (and the relation to \cite{KL1}
is far from simple). 
It seems important to develop the  {\em non-semisimple} 
technique of intertwiners from this paper toward
\cite{O11} and \cite{MTa}.
Hopefully, switching to the DAHA\~level can be productive
here, but there are no direct confirmations so far.
\medskip

To conclude, we will touch upon the 
$q\leftrightarrow q'$\~duality
for $\Pi_{\widehat{R}}$ in the notation from Theorem
\ref{OTHERAFFINE}. Let 
$\{t_{lng}^{l_r}t_{\sht}^{s_r},\ 1\le r\le d\}$
be the set of the $t$\~powers in the numerator
of the corresponding product formula for the 
Poincar\'e polynomial $\Pi_R$ from (\ref{tdegr}) in the case 
$t_{\lng}=t_{\sht}$ and
from (\ref{tdegra}-\ref{tdegrg}). We set $h_r=l_r+s_r$ and
denote the corresponding sum in the {\em denominator} of $\Pi_R$ 
by $n_r$. Thus $n_r$ always divides $h_r$ and
$d=n,\, n_r=1,\, h_r=m_r+1$ for all $r$ as $t_{\lng}=t_{\sht}$.
Let $\tau_r\equal(t_{lng}^{l_r}t_{\sht}^{s_r})^{1/h_r}$.

Equivalently, the $k$\~terms of the positive 
affine exponents of $\Pi_{\widehat{R}}$
from (\ref{yofxevalnew}-\ref{newf4affine})
are $\{l_r k_{lng}+s_r k_{\sht} +j,\ 1\le r\le d\}$.
Respectively, the negative affine exponents are 
divisors of the {\em non-rational} positive affine 
exponents with the corresponding ratios equal to
$h_r/n_r$. The number
$d$ coincides with the number of negative affine exponents 
such that $j=0$ from Theorem \ref{OTHERAFFINE}. 

We can now rewrite $\Pi_{\widehat{R}}$ in a
$q\leftrightarrow\exp(2\pi \imath)$\,\~symmetric form as follows
(cf. Theorem \ref{RADADE}) \,:
\begin{align}\label{qqprime}
&\Pi_{\widehat{R}}=
\prod_{r=1}^d \,\frac{\prod_{j,j\,'=0}^{h_r-1}\,
(1-q^{j/h_r}\,e^{2\pi \imath j\,'/h_r}\,\tau_r)}
{\prod_{j,j'=0}^{n_r-1}\,
(1- q^{j/n_r}\,e^{2\pi \imath j\,'/n_r}\,\tau_r)}\,,
\end{align}
which is an indication that the localization functor from
\cite{GGOR},\cite{VV1} can be extended to the general 
$q,t$\~theory.
     
\vfill
\eject

\vskip -1cm
\bibliographystyle{unsrt}

\end{document}